%% file: diss.tex
\documentstyle[german,12pt,titlepage]{article}
\begin{document}

\newtheorem{satz}{Satz}[subsection]
\newtheorem{lem}[satz]{Lemma}
\newtheorem{prop}[satz]{Proposition}
\newtheorem{kor}[satz]{Korollar}

\begin{titlepage}

\centerline{}
\bigskip
\bigskip

\LARGE
\centerline{Zur Affinit"at von Hyperfl"achenkomplementen}
\centerline{in affinen und projektiven Schemata}

\vspace{2.5cm}
\Large
\centerline{Dissertation}
\smallskip
\large
\centerline{zur Erlangung des Grades eines}
\centerline{Doktors der Naturwissenschaften}
\centerline{der Fakult"at f"ur Mathematik}
\centerline{an der Ruhr\--Universit"at Bochum}

\vspace{4.2cm}
\centerline{von}
\smallskip
\Large
\centerline{Holger Brenner}
\smallskip
\large
\centerline{aus}
\smallskip
\centerline{Wildberg}

\vspace{2cm}
\large
\centerline{Bochum 1998}

\end{titlepage}

\newpage

\pagenumbering{arabic}
\setcounter{page}{0}
\thispagestyle{empty}

\ 
\\[11.5cm]

\noindent
Diese Arbeit entstand unter Anleitung von Herrn Prof. Dr. U. Storch an der
Ruhr-Universit"at Bochum.
\par\bigskip\noindent
Dissertation eingereicht am 15. 12. 1998.
\par\bigskip\noindent
Referenten: Prof. Dr. U. Storch, Ruhr-Universit"at Bochum
\par\smallskip\noindent
\hspace{2.04cm} Prof. Dr. H. Flenner, Ruhr-Universit"at Bochum
\par\smallskip\noindent
\hspace{2.04cm} Prof. Dr. W. Bruns, Universit"at Osnabr"uck
\par\bigskip\noindent
Tag der m"undlichen Pr"ufung: 16. 2. 1999.

\newpage

\tableofcontents

\newpage
\addcontentsline{toc}{section}{Einleitung}

\input{einleit}

\newpage
\
\thispagestyle{empty}

\newpage
\input{grundla}

\newpage
\
\thispagestyle{empty}
\newpage
\input{proka}

\newpage
\input{afklaka}

\newpage
\input{suphoka}

\newpage
\input{divhoka}
\newpage
\
\thispagestyle{empty}
\newpage
\addcontentsline{toc}{section}{Probleme und offene Fragen}
\input{problem}

\newpage
\
\thispagestyle{empty}
\newpage
\addcontentsline{toc}{section}{Literaturverzeichnis}

\input{biblio}
\newpage
\
\thispagestyle{empty}
\newpage
\input{lebenslauf}

\end{document}

%% file: einleit.tex
\section*{Einleitung}
\par
\bigskip
\medskip
\noindent
\footnotesize
\hspace{7.5cm} Nicht allein in Rechnungssachen
\par\noindent
\hspace{7.5cm} Soll der Mensch sich M"uhe machen
\par\smallskip\noindent
\hspace{7.5cm} Wilhelm Busch, Max und Moritz (1865)
\normalsize
\par
\medskip
\bigskip
\noindent
Sei $Y \subseteq {\bf C}^n$ eine Hyperfl"ache, also eine abgeschlossene
algebraische Teilmenge der Kodimension eins.
Dann ist $Y$ die Nullstellenmenge eines Polynoms $F$, auf
dem offenen Komplement ${\bf C}^n-Y$ ist die rationale Funktion
$1/F$ definiert und ${\bf C}^n-Y$ ist isomorph zum Graph dieser
Funktion in $({\bf C}^n-Y) \times {\bf C}$, der abgeschlossen in
${\bf C}^n \times {\bf C}$ liegt.
Eine abgeschlossene algebraische Teilmenge eines ${\bf C}^m$ nennt man eine
affine Variet"at,
das Komplement einer Hyperfl"ache im ${\bf C}^n$ ist also affin.
\par\smallskip\noindent
Sei $X$ eine algebraische Variet"at oder allgemeiner
ein noethersches separiertes Schema und $Y$ eine Hyperfl"ache in $X$.
In dieser Arbeit untersuchen wir die Frage,
wann das offene Komplement $U=X-Y$ affin ist.
Dabei gilt unser Hauptinteresse zum einen der Situation, wo $X$ selbst eine
affine Variet"at oder ein affines Schema ist,
also Spektrum eines kommutativen Ringes,
zum anderen der Situation, wo $X$ eine projektive Variet"at ist.
Die Beschr"ankung auf Hyperfl"achen ist dabei keine
wirkliche Einschr"ankung, da eine abgeschlossene Teilmenge mit einer
Komponente der Kodimension $\geq 2$ kein affines Komplement
besitzen kann.
\par\smallskip\noindent
Es handelt sich dabei um ein Thema mit einer
Vielzahl von Aspekten, Teilproblemen und spezielleren Fragen.
Innerhalb der algebraischen Geometrie und
kommutativen Algebra bestehen
vielf"altige Bez"uge zu ger"aumigen Divisoren auf projektiven
Variet"aten, kohomologischen
Eigenschaften, endliche Erzeugtheit von globalen Schnitt\-ringen
(14. Hilbertsches Problem), Schnitttheorie, Quotientenbildung.
Auch zur komplexen Analysis ergeben sich durch den Begriff des Steinschen
Raumes Beziehungen.
\par\smallskip\noindent
Auf der einen Seite gibt es einige prominente Resultate zur Affinit"at, wie
etwa die kohomologische Charakterisierung, der Satz von Chevalley,
die Kodimensionseigenschaft oder der Satz von Nagata.
Auf der anderen Seite gibt es wenig speziellere Literatur zu diesem Thema,
erw"ahnt seien
\cite{goodman},\cite{harcodi},\cite{haramp},\cite{goodlandman},
\cite{binsto1},\cite{bingener},\cite{arezzo}.
\par\smallskip\noindent
Im Laufe der "Uberlegungen zu dieser Arbeit haben sich folgende beiden
Fragestellungen als Leitfragen herauskristallisiert, deren Behandlung
im Rahmen dieser Arbeit ein besonderes Gewicht zukommt.
\par\bigskip
\noindent
1) Inwiefern gibt es ein geometrisches Kriterium f"ur die Affinit"at
eines Hyperfl"achenkomplementes?
\par\smallskip\noindent
Ausgangspunkt f"ur diese Fragestellung ist die Beobachtung, dass man die
Nicht-Affinit"at von $X-Y \subseteq X$ h"aufig am einfachsten
dadurch erweisen kann, dass man einen affinen Morphismus
$f:X' \longrightarrow X$ angibt, unter dem das Urbild $f^{-1}(Y)$ der
Hyperfl"ache Komponenten der Kodimension $\geq 2$ besitzt.
Dies ist beispielsweise dann der Fall, wenn der Schnitt von
$Y$ mit einer abgeschlossenen
Fl"ache in $X$ isolierte Punkte besitzt.
Dieses schnittgeometrische Prinzip zum Nachweis der Nicht-Affinit"at
durchzieht die gesamte Arbeit und wird im vierten Kapitel auf seine
Reichweite und Grenzen hin untersucht.
\par\bigskip
\noindent
2) Was bedeutet die Eigenschaft, dass in einem affinen Schema s"amtliche
Hyperfl"achen affines Komplement besitzen und wie kann man die Abweichung
von dieser Eigenschaft messen?
\par\smallskip\noindent
Eine treibende Kraft innerhalb dieses Problemkreises
war dabei die Frage, ob sich diese Eigenschaft von ${\rm Spek}\, A$
auf die affine Gerade ${\rm Spek}\, A[T]$ "ubertr"agt.
Als mathematisches Vorbild diente dabei die Divisorenklassengruppe, die
die Abweichung von der Faktorialit"at misst und die sich beim "Ubergang
zum Polynomring nicht "andert.
Dazu f"uhren wir im dritten Kapitel die affine Klassengruppe ein, die eine
Restklassengruppe der Divisorenklassengruppe ist und genau dann verschwindet,
wenn jede Hyperfl"ache affines Komplement besitzt.
\par\bigskip
\noindent
Die Arbeit bem"uht sich, durch Beispiele die erzielten Resultate zu
illustrieren und auf ihre Voraussetzungen zu "uberpr"ufen.
Die Anwendbarkeit der affinen Klassengruppe zeigt sich etwa in der expliziten
Berechnung in einigen wichtigen Beispielklassen.
In eine andere Richtung gehen die offenen Fragen, die an verschiedenen Stellen
durch den Gang der Untersuchung aufgeworfen werden. Eine Problemliste
findet sich am Ende der Arbeit.
Die Verbindung zur komplexen Analysis wird in Bemerkungen verfolgt.
\par\smallskip\noindent
Wir geben eine "Ubersicht "uber Inhalt und Ergebnisse der Arbeit,
wobei die Aussagen nicht unbedingt in der vollen Sch"arfe ausgef"uhrt werden.
\par\bigskip
\noindent
Das erste Kapitel der Arbeit hat einleitenden und vorbereitenden
Charakter. Hier werden einerseits die ben"otigten Begriffe eingef"uhrt
und entwickelt, andererseits finden sich hier die wichtigen Resultate
aus der Literatur zur Affinit"at wieder.
\par\smallskip\noindent
Nach einer kurzen Beschreibung des Themas der Arbeit in einer begrifflich
unaufwendigen Sprache (1.1) und der schematheoretischen Reformulierung (1.2)
beginnen wir in 1.3 mit der Untersuchung der kanonischen Abbildung eines Schemas
auf das Spektrum des globalen Schnittringes.
Mittels dieser Abbildung werden quasiaffine Schemata und die G"ultigkeit des
globalen Hilbertschen Nullstellensatzes charakterisiert.
Diese beiden Eigenschaften zusammen charakterisieren die Affinit"at.
F"ur Variet"aten wird gezeigt, dass sich die Quasiaffinit"at bereits
aus der scheinbar schw"acheren Eigenschaft ergibt, dass sich je zwei
abgeschlossene Punkte durch eine globale Funktion trennen lassen,
als auch dadurch, dass es auf jeder Kurve globale, nicht konstante Funktionen
gibt.
\par\smallskip\noindent
In 1.4 wenden wir uns den verschiedenen Beschreibungsm"oglichkeiten
des globalen Schnittringes zu.
F"ur eine offene Teilmenge $U=D({\bf a})$ in einem affinen Schema
lautet das (bekannte) Affinit"atskriterium, dass das Erweiterungsideal
im globalen Schnittring zum Einheitsideal wird.
Es wird gezeigt, dass dies f"ur Hyperfl"achen, die lokal mengentheoretisch
durch eine Funktion beschrieben werden k"onnen, erf"ullt ist.
Damit ist in jedem eindimensionalen und in jedem lokal faktoriellen
affinen Schema die beschriebene Eigenschaft erf"ullt, dass jede
Hyperfl"ache affines Komplement besitzt.
\par\smallskip\noindent
1.5 zeigt, dass unter den nulldimensionalen Schemata sich die affinen
topologisch charakterisieren lassen. Ein kompaktes (es gen"ugt quasikompakt und
quasisepariert) nulldimensionales Schema ist affin.
\par\smallskip\noindent
In 1.6 wird der grundlegende Begriff der affinen Abbildung rekapituliert.
Wichtig ist dabei, dass unter einer affinen Abbildung das Urbild jeder
affinen Teilmenge (nicht nur aus einer "Uberdeckung)
wieder affin ist.
\par\smallskip\noindent
1.7 nimmt das Studium der kanonischen Abbildung wieder auf.
Es wird gezeigt, dass diese mit flachen Basiswechseln und insbesondere
mit Lokalisierungen vert"aglich ist.
Die Betrachtung der Lokalisierung am generischen Punkt (im integren Fall)
f"uhrt zum Begriff generisch affin, der besagt, dass die kanonische
Abbildung im generischen Punkt affin ist (und damit eine Isomorphie).
"Uberraschenderweise ist in diesem Fall bereits die Lokalisierung an jedem
Primideal der H"ohe eins affin (1.7.5). Ein verwandtes Ergebnis besagt, dass ein
separiertes, integres, noethersches Schema mit offenem generischen
Punkt affin ist (1.7.6).
\par\smallskip\noindent
In 1.8 wird gezeigt, dass die Affinit"at einer Abbildung eine offene
und (punktweise) lokale Eigenschaft ist.
Ist zu $f:X \longrightarrow Y$ und einem Punkt $y \in Y$ die lokale Faser
$X_y=X \times_Y {\rm Spek}\, {\cal O}_y$ affin, so gibt es bereits eine offene
affine Umgebung von $y$ mit affinem Urbild (1.8.1). Insbesondere kann damit die
Affinit"at einer Abbildung lokal getestet werden (1.8.2).
Dies f"uhrt bei einer offenen Teilmenge $U \subseteq {\rm Spek}\, A$
zu dem bekannten Lokalisierungsprinzip, dass $U$ genau dann affin ist,
wenn dies in jeder Lokalisierung von ${\rm Spek}\, A$ gilt.
\par\smallskip\noindent
Anhand von einigen Beispielen wird in 1.9 gezeigt, dass die echten
Fasern  $X_{\kappa (y)}=X \times_Y {\rm Spek}\, \kappa(y)$ ihr
Affinit"atsverhalten von Punkt zu Punkt "andern k"onnen. Insbesondere
kann die Affinit"at abh"angig von der Charakteristik des Restek"orpers sein.
\par\smallskip\noindent
In 1.10 wird der eindimensionale Fall abgehandelt.
Offene Teilmengen affiner noetherscher Kurven sind
affin (1.10.1). Ebenso ist ein generisch affines, separiertes Schema mit
eindimensionalem globalen Schnittring affin (1.10.4).
\par\smallskip\noindent
In 1.11 wird an den Begriff des Krullringes erinnert, den unser Beweis
der Kodimensionseigenschaft in 1.12 erfordert, die besagt,
dass das Komplement einer affinen Menge die Kodimension $\leq 1$ besitzt.
\par\smallskip\noindent
In 1.13 erw"ahnen wir die kohomologische Charakterisierung affiner Schemata
mit dem Korollar, dass man sich bei der Affinit"atsuntersuchung
auf die (reduzierten) Komponenten beschr"anken kann.
\par\smallskip\noindent
1.14 besch"aftigt sich mit der Frage, wann man bei einer surjektiven
Abbildung $f:X \longrightarrow Y$
von der Affinit"at von $X$ auf die Affinit"at von $Y$ schlie"sen kann.
Das prominente Resultat ist hierbei der Satz von Chevalley, der dies
f"ur endliche Abbildungen best"atigt.
Eine Folgerung daraus ist, dass man die Affinit"at in der
Normalisierung testen kann, so dass man sich bei der Untersuchung
h"aufig auf normale integre Schemata beschr"anken kann.
Neben endlichen Abbildungen erlauben auch Basiswechsel, die
treuflach oder ein direkter Summand sind, diesen R"uckschluss auf die
Affinit"at.
Beispiele zeigen, dass dies f"ur hom"oomorphe Abbildungen nicht sein muss.
\par\smallskip\noindent
1.15 schlie"slich beinhaltet den Satz von Nagata, dass in einem
zweidimensionalen, normalen, exzellenten Ring das Komplement jeder (reinen)
Kurve affin ist.
Insbesondere erf"ullen diese Ringe die in 2) beschriebene Eigenschaft.
Unsere Ausarbeitung des Beweises von Nagata zeigt,
dass es in einem beliebigen noetherschen Ring der Dimension zwei f"ur
eine nicht affine Teilmenge einen noetherschen Ringwechsel gibt, unter dem das
beschreibende Ideal die Kodimensionseigenschaft verletzt.
Im zweidimensionalen Fall gibt es also ein geometrisches Kriterium f"ur die
Affinit"at.
\par\bigskip
\noindent
Das zweite Kapitel besch"aftigt sich mit offenen Teilmengen auf einem
projektiven Schema und den Urbildern davon unter der Kegelabbildung.
Dabei ist der affine Kegel das Spektrum $X$ eines graduierten Ringes $A$, und
das projektive Schema $Y$ das homogene Spektrum.
Der Graduierung auf $A$ entspricht dabei eine Operation des Gruppenringes
${\rm Spek}\, A_0[T,T^{-1}]$ auf ${\rm Spek}\, A$.
Unter dieser Operation ist $V \subseteq Y$ der Quotient von $U=p^{-1}(V)$,
so dass wir in 2.1 generell Gruppenoperationen und Quotientenbildung untersuchen.
Dabei erfordert h"aufig der Nachweis, dass ein schematheoretischer Quotient
existiert, den Affinit"atsnachweis f"ur den Quotienten auf einer kleineren
invarianten Menge.
\par\smallskip\noindent
In 2.2. zeigen wir allgemein die "Aquivalenz von $D-$Graduierungen und
Operationen von ${\rm Spek}\, K[D]$, wobei $D$ eine beliebige kommutative
Gruppe ist und geben ein (trivialerweise notwendiges) Kriterium an,
unter dem das Spektrum des Ringes der nullten Stufe der Quotient
der Operation ist (2.2.2). Dieses Kriterium verallgemeinert ein Kriterium aus
\cite{SGA3}, das die Freiheit der Operation voraussetzt.
\par\smallskip\noindent
Dieses Kriterium wird in 2.3 spezialisiert auf den Fall von (nicht
unbedingt positiven) ${\bf Z}-$Graduierungen und erlaubt den Nachweis, dass
unter der Kegel\-abbildung $p$ das Urbild $p^{-1}(V)$ genau dann affin ist,
wenn $V$ affin ist (2.3.5).
Ferner beschreiben wir in diesem Abschnitt Beispiele und Ph"anomene, die bei
nicht positiven Graduierungen auftreten k"onnen.
\par\smallskip\noindent
In 2.4 werden f"ur eine offene Teilmenge $U$ in einer projektiven Variet"at
notwendige und hinreichende Kriterien angegeben, dass $U$ affin ist
(oder verwandte Eigenschaften erf"ullt).
Insbesondere wird die Beziehung zu Eigenschaften eines Divisors mit dem
Komplement von $U$ als Tr"ager untersucht.
\par\smallskip\noindent
2.5 besch"aftigt sich mit der Frage, wann die Aufblasung eines Ideals
in einem integren noetherschen Ring affin ist.
Dies ist genau dann der Fall, wenn das Erweiterungsideal in der Normalisierung
invertierbar wird, was bereits dann der Fall ist, wenn es in irgendeiner
ganzen Erweiterung invertierbar wird (2.5.2).
Die dabei entstehenden Aufblasungsringe liegen innerhalb der
Normalisierung; ist diese endlich, so kann sogar durch affine Aufblasungen
normalisiert werden (2.5.7).
\par\bigskip
\noindent
Im dritten Kapitel wird die affine Klassengruppe zu einem normalen Schema $X$
als Ma"s daf"ur eingef"uhrt, wie viele Hyperfl"achen nicht affines Komplement
besitzen.
\par\smallskip\noindent
In 3.1. wird gezeigt, dass es in einem affinen normalen Schema zu
einer Hyperfl"ache mit nicht affinem Komplement stets effektive, linear
"aquivalente Divisoren mit affinem Komplement gibt, so dass also die Affinit"at
keine Eigenschaft der Divisorenklasse ist.
Dies ist in 3.2
der Ausgangspunkt zum Begriff des koaffinen Divisors, der dadurch
definiert wird, dass jeder effektive "aquivalente Vertreter affines
Komplement besitzt.
Diese Eigenschaft kann dadurch charakterisiert werden, dass der zugeh"orige
reflexive Modul ${\cal L}_D$ ein Erzeugendensystem besitzt, dessen Nullstellen
affines Komplement haben (3.2.2). Mit Hilfe dieses Kriteriums l"asst
sich zeigen, dass die Eigenschaft koaffin unter einem flachen Ringwechsel
erhalten bleibt (3.2.3).
Dies erlaubt bereits den Schluss, dass sich die Eigenschaft eines normalen
affinen Schemas ${\rm Spek}\, A$,
dass alle Hyperfl"achenkomplemente affin sind,
auf die affine Gerade ${\rm Spek}\, A[T]$ "ubertr"agt (3.2.4).
Ist speziell $A$ zweidimensional, normal und exzellent, so folgt,
dass der Schnitt zweier Fl"achen in ${\rm Spek}\, A[T]$ keine
isolierten Punkte besitzt (3.2.5). Beispiele zeigen, dass diese Aussagen
ohne die Voraussetzung normal nicht richtig sind.
\par\smallskip\noindent
In 3.3 wird die affine Klassengruppe eingef"uhrt. Sie beruht auf dem Begriff
des affin-trivialen Divisors, der eine weitere Versch"arfung von koaffin ist.
Die affin-trivialen Divisoren bilden eine Untergruppe der
Divisorenklassengruppe, die Restklassengruppe nennen wir die affine
Klassengruppe. Sie verschwindet genau dann, wenn jede Hyperfl"ache affines
Komplement besitzt. Bei treu\-flachen Ringwechseln mit surjektiver
Abbildung der Divisorenklassengruppen "andert sich die affine Klassengruppe
nicht (3.3.3). Dies gilt insbesondere beim "Ubergang von $A$ zu
$A[T_1,...,T_n]$ (3.3.4).
\par\smallskip\noindent
3.4 beschreibt einige Folgerungen, wenn die affine Klassengruppe verschwindet.
So ist f"ur eine in der Kodimension eins definierte rationale Abbildung
in ein normales Schema mit trivialer affiner Klassengruppe
der Undefiniertheitsort im Abschluss des Urbildes von jeder Hyperfl"ache
enthalten (3.4.3).
In diesen Kontext f"allt auch der Reinheitssatz von van der Waerden (3.4.7).
Wir diskutieren, wie durch nicht affine Hyperfl"achenkomplemente
(unter gewissen Zusatzbedingungen) der Reinheitsatz verletzt wird.
\par\smallskip\noindent
3.5 kn"upft an den "Ubergang von ${\rm Spek}\, A$ zu
${\rm Spek}\, A[T]={\rm Spek}\, A \times {\bf A}^1$ an und
untersucht, wie sich die affine Klassengruppe eines Produktes
$X \times Y$ von normalen affinen Variet"aten "uber einem algebraisch
abgeschlossenen K"orper verh"alt, insbesondere im einfachsten Fall,
dass $Y=C$ eine glatte affine Kurve ist.
In dieser Situation wird zu einer rationalen Funktion
$f: X \cdots \longrightarrow C$ das Komplement des Abschlusses des Graphen
studiert (3.5.3). Darauf aufbauend kann f"ur eine normale Fl"ache gezeigt
werden, dass die affine Klassengruppe von $X \times C$ verschwindet (3.5.6).
\par\smallskip\noindent
In 3.6. wird die affine Klassengruppe zu drei prominenten Ringklassen
berechnet: f"ur Hyperbelringe, f"ur Monoidringe und f"ur
Determinantenringe. Die detaillierten Berechnungen ergeben dabei in allen
drei F"allen, dass die affine Klassengruppe gleich der Divisorenklassengruppe
modulo Torsion ist.
\par\smallskip\noindent
In 3.7 steht die Frage im Mittelpunkt, ob es zu einer affinen Teilmenge in einer
projektiven Variet"at einen ger"aumigen Divisor mit dem Komplement als Tr"ager
gibt. Dabei wird haupts"achlich "uber Ergebnisse von Goodman referiert.
\par\smallskip\noindent
3.8 beschreibt koaffine Divisoren und die affine Klassengruppe eines
positiv graduierten Ringes in Begriffen der zugeh"origen projektiven
Variet"at. Insbesondere ergibt sich, dass die affine
Klassengruppe des affinen Kegels genau dann verschwindet, wenn dies auf dem
homogenen Spektrum gilt (3.8.4). Ist die projektive Variet"at glatt, so f"uhren
numersch triviale Divisoren zu koaffinen Divisoren im Kegel (3.8.6).
\par\smallskip\noindent
3.9 untersucht entlang der Klassifikation Affinit"atseigenschaften auf
projektiven glatten Fl"achen und den affinen Kegeln dar"uber. Es werden
Bedingungen f"ur das Verschwinden der affinen Klassengruppe und Folgerungen
daraus vorgestellt. Bei geometrischen Regelfl"achen best"atigt sich der
Zusammenhang zwischen affiner Klassengruppe und numerischer Klassengruppe,
die affine Klassengruppe der Kegel ist ${\bf Z}$.
Abelsche Fl"achen mit trivialer affiner Klassengruppe werden
charakterisiert.
\par\bigskip
\noindent
Im Mittelpunkt des vierten Kapitels steht die Frage, inwiefern man die
Affinit"at einer offenen Teilmenge $U=D({\bf a}) \subseteq {\rm Spek}\, A$
dadurch sichern kann, dass das beschreibende Ideal universell die H"ohe eins
besitzt, ob man also im nicht affinen Fall einen Ringwechsel angeben kann, wo
die Kodimensionsbedingung verletzt wird.
\par\smallskip\noindent
In 4.1 wird dazu der Begriff der Superh"ohe eines Ideals besprochen, der die
maximale H"ohe bezeichnet, die das Erweiterungsideal unter einem Ringwechsel
bekommen kann. Dieser Begriff wurde in Zusammenhang mit den homologischen
Vermutungen von Hochster eingef"uhrt und wird hier in Verbindung zu
Affinit"atsfragen gebracht.
\par\smallskip\noindent
In 4.2 wird die Superh"ohe im globalen Schnittring von $U$ untersucht.
Es ergibt sich, dass ein Ideal genau dann eine affine Teilmenge beschreibt,
wenn es die Superh"ohe eins bezogen auf alle Krullbereiche besitzt (4.2.1).
Ist der globale Schnittring von endlichem Typ, so ist die Affinit"at "aquivalent
dazu, dass die endliche Superh"ohe von $\bf a$ eins ist (4.2.2).
Dies ist im Fall einer $K-$Algebra "aquivalent dazu, dass s"amtliche
(abgeschlossenen) Fl"achen in $U$ affin sind (4.2.3).
Nebenbei geben wir einen einfachen Beweis f"ur die
(zur direkten Summand Vermutung "aquivalente) Monomialvermutung im
(bekannten) zweidimensionalen Fall, der sich auf Affinit"at st"utzt.
\par\smallskip\noindent
4.3 bringt Beispiele, dass sich im Allgemeinen die Affinit"at nicht aus der
noetherschen Superh"ohe eins ergibt. Ansatzpunkt sind dabei im Fall des
komplexen Grundk"orpers Steinsche offene Mengen.
Diese erf"ullen ebenfalls die universelle Kodimensionseigenschaft,
m"ussen aber nicht affin sein. Die Superh"ohe eins l"asst sich
("ahnlich wie mit algebraischen Fl"achen) dadurch
charakterisieren, dass alle analytischen Fl"achen in $U$ Steinsch sind (4.3.5).
\par\smallskip\noindent
In 4.4 f"uhren wir die Superh"ohe f"ur ein beliebiges Schema ein, womit man
die Superh"ohe eines Ideals $\bf a$ als nur von $D({\bf a})$
abh"angig erweisen kann (4.4.1). Auf diesen Begriff st"utzt sich in 4.5
der Beweis der Aussage, dass das homogene beschreibende Ideal einer
irreduziblen Kurve auf einer glatten projektiven Fl"ache mit Selbstschnitt
null und positivem Schnitt mit den "ubrigen Kurven die Superh"ohe eins besitzt.
Dies liefert f"ur jeden algebraisch abgeschlossenen K"orper
(unabh"angig von analytischen Methoden) Beispiele von nicht affinen, offenen
Teilmengen mit noetherscher Superh"ohe eins.
\par\smallskip\noindent
4.6 behandelt zwei Beispielklassen, bei denen die Affinit"at einer offenen
Teilmenge durch das H"ohenverhalten in einer einzigen Erweiterung
getestet werden kann, und wo insbesondere die Affinit"at durch endliche
Superh"ohe eins gesichert ist, n"amlich Fl"achen mit ihrer Normalisierung
und Monoidringe mit einer Darstellung des Monoids, das die
Durchschnittseigenschaft besitzt.
\par\bigskip
\noindent
Das f"unfte Kapitel besch"aftigt sich
abschlie"send mit Eigenschaften abgeschlossener
Teilmengen h"oherer Kodimension, die analog zur Affinit"at sind, insbesondere
mit dem Verh"altnis von kohomologischer Dimension und Superh"ohe.
\par\smallskip
\noindent
5.1 erinnert an das affine Dimensionstheorem und die Verallgemeinerung
von Serre f"ur beliebige regul"are Ringe und zeigt darauf aufbauend,
dass f"ur ein Primideal in einem regul"aren Ring H"ohe und Superh"ohe
"ubereinstimmen (5.1.3).
\par\smallskip\noindent
In 5.2 wird der Begriff von Hartshorne der kohomologischen Dimension und
die wesentlichen Eigenschaften besprochen, insbesondere in Vergleich zu
Eigenschaften der Superh"ohe. In 5.3 wird basierend auf die kohomologische
Charakterisierung von Grothendieck
der punktierten lokalen Schemata gezeigt, dass die Superh"ohe
durch die kohomologische Dimension $+1$ beschr"ankt ist (5.3.1). Daneben werden
weitere H"ohenbegriffe diskutiert und mit Superh"ohe und kohomologischer
Dimension in Verbindung gebracht.
\par\smallskip\noindent
5.4 beschreibt die Absch"atzung der Superh"ohe von $\bf a$ durch die
projektive Dimension des Restklassenringes $A/{\bf a}$, die auf
dem Durchschnittstheorem von Roberts beruht (5.4.2).
In der Charakteristik $p$ beschr"ankt die projektive Dimension auch die
kohomologische Dimension $+1$ (5.4.5), ein Beispiel zeigt,
dass dies in Charakteristik null nicht richtig ist.
\par\smallskip\noindent
5.5 schlie"slich diskutiert, wie sich die verschiedenen H"ohenbegriffe unter
Ringwechseln verhalten, insbesondere unter dem Gesichtspunkt,
wann sie konstant bleiben.
\par\bigskip
\noindent
Mein Dank gilt in erster Linie
Herrn Prof. Dr. U. Storch f"ur die Anregung zu dieser Arbeit
und die intensive Betreuung.
\par\smallskip\noindent
F"ur Ideen, Kritik,
Ratschl"age, Literaturtips und Diskussionen
m"ochte ich mich bei Herrn Prof. Dr. H. Flenner, Herrn Prof. Dr.
A. Huckleberry, bei Frau Dipl.-Math. H. Wagner und den Herren
Dipl.-Math. T. Fenske, G. Rei"sner und S. Schroer bedanken.
Bei den Mitgliedern des Oberseminars Algebra und Geometrie
m"ochte ich mich f"ur die Gelegenheit zum
Vortrag und das Interesse an meiner Arbeit bedanken.
\par\smallskip\noindent
F"ur die Hilfe bei Computerproblemen bin ich Herrn Dr. H. Wiebe,
Herrn T. Fenske und Herrn T. Gorski zu Dank verpflichtet.
\par\smallskip\noindent
Dem Graduiertenkolleg Geometrie und mathematische Physik der deutschen
Forschungsgemeinschaft verdanke ich durch ein Stipendium
die M"oglichkeit, mich eineinhalb Jahre ausschlie"slich der Arbeit widmen
zu k"onnen.
\par\smallskip\noindent
Zum Schluss danke ich Sabine f"ur die sch"onen Seiten des Lebens jenseits
der Rechnungssachen.

%% file: grundla.tex
\section{Grundlagen}
\par
\bigskip
\noindent
Das erste Kapitel der Arbeit hat vorbereitenden und einleitenden Charakter.
Es soll das Thema vorgestellt und der begriffliche Rahmen abgesteckt
werden, indem sich die Behandlung des Themas abspielen soll.
Zum begrifflichen schematheoretischen
Werkzeug geh"oren neben anderem die kanonische Abbildung,
der globale Schnittring, affine Abbildungen, Basiswechsel, Krullringe,
Kohomologie.
\par\smallskip\noindent
An prominenten Resultaten begegnen hier das Lokalisierungsprinzip,
die Kodimensionseigenschaft,
kohomologische Charakterisierung, Satz von Chevalley und
Satz von Nagata. Diese grundlegenden Ergebnisse finden sich zum Teil auch
in Lehrb"uchern, d"urfen aber hier nicht fehlen. Zum "uberwiegenden
Teil ist die Darstellung dieser Resultate eigenst"andig, und zwar
in Hinblick auf Beweis, Tragweite der Aussagen und Anwendungen.
\par\smallskip\noindent
Die vorgestellten Methoden erlauben die Behandlung der Affinit"atsproblematik
in niederen Dimensionen, d.h. hier bis einschlie"slich der Dimension
zwei, so dass diese Beispielklassen in eigenen Abschnitten
schon hier besprochen werden k"onnen. 
\par
\bigskip
\noindent
\subsection{Affine Variet"aten}
\par
\bigskip
\noindent
Zu Beginn soll das Thema der Arbeit in einer vereinfachten
Situation beschrieben werden. Wir betrachten den ${\bf C}^n$ und dazu die
Polynomalgebra ${\bf C}[X_1,...,X_n]$. Eine Menge von Polynomen
$P_1,...,P_r \in {\bf C}[X_1,...,X_n]$
definiert dann in ${\bf C}^n$ die gemeinsame Nullstellenmenge
$$ V=V(P_1,...,P_n)=
\{x=(x_1,...,x_n) \in {\bf C}^n:\, P_i(x)=0 \mbox{ f"ur } i=1,...,r \}\, .$$
Ein solches Nullstellengebilde hei"st eine
affine oder affin-algebraische Variet"at.
Sei $U \subseteq V$ eine Zariski-offene Teilmenge, d.h.
$U$ ist das Komplement in $V$ eines weiteren algebraischen Nullstellengebildes
$Y=V(F_1,...,F_s)$, also $U=V-V\cap Y$.
Unsere Frage ist nun: Wann ist $U$ selbst affin, also
mittels auf $U$ definierten rationalen Funktionen
(das sind Quotienten aus polynomialen Funktionen, deren Nenner auf $U$
nullstellenfrei sind) abgeschlossen einbettbar in
einen ${\bf C}^m$?
\par\smallskip\noindent
Man kann immer von der vorgegebenen Einbettung
$i=(X_1,...,X_n): U \subseteq V \hookrightarrow {\bf C}^n$
ausgehen und versuchen, sie zu
einer Abbildung vom Typ $U \longrightarrow {\bf C}^n \times {\bf C}^k$
mit einem abgeschlossenen Bild mittels rationalen Funktionen auf $U$
zu erg"anzen.
Ist jede auf $U$ definierte rationale Funktion auf $V$ ausdehnbar,
so l"asst sich eine abgeschlossene Einbettung f"ur $U$ zu einer
abgeschlossenen Einbettung f"ur $V$ fortsetzen, was nur bei $U$ abgeschlossen
in $V$ der Fall sein kann.
\par\smallskip\noindent
Eine notwendige Bedingung f"ur die abgeschlossene Einbettbarkeit
bei $U \neq \emptyset,\, \neq V$ und $V$ zusammenh"angend ist somit,
dass es auf $U$ rationale Funktionen geben muss,
die nicht auf ganz $V$ ausdehnbar sind.
Ist $V$ normal und hat das Komplement $Y$ in $V$
eine Kodimension von zumindest
zwei, so folgt aus dem Riemannschen Hebbarkeitssatz, dass jede rationale
Funktion auf $U$ nach $V$ fortsetzbar ist. Als Komplemente offener affiner 
Teilmengen kommen daher nur Hyperfl"achen in Frage, die punktierte Ebene
${\bf C}^2-\{ P \}$ beispielsweise ist nicht affin.
\par\smallskip\noindent
Wird die Hyperfl"ache $Y$ in $V$ durch ein einziges Polynom
$F \in {\bf C}[X_1,...,X_n]$ beschrieben, so
ist auf $U$ die Funktion $1/ F$ definiert und die Gesamtabbildung
$i \times 1/F:U \longrightarrow {\bf C}^n \times {\bf C}$ ist
abgeschlossen, das Bild ist n"amlich gleich
$V(P_1,...,P_n, X_{n+1}\cdot F(X_1,...,X_n)-1)$.
Insbesondere ist in ${\bf C}^n$ selbst das Komplement jeder Hyperfl"ache
affin, da der Polynomring faktoriell ist und damit jede Hyperfl"ache durch
eine Funktion beschrieben werden kann.
\par\smallskip\noindent
Wird eine Hyperfl"ache nicht durch eine einzige Funktion beschrieben, so 
ist es im Allgemeinen schwierig zu entscheiden, ob das Komplement affin ist
oder nicht. Man muss dazu die auf $U$ und $V$ definierten rationalen Funktionen
kennen und die Affinit"at von $U$ mittels dieser Funktionen
charakterisieren.
Dies f"uhrt zu
einem intrinsischen Standpunkt, der unabh"angig von der
vorgegebenen Einbettung
$V \hookrightarrow {\bf C}^n$ ist.
\par\smallskip\noindent
Ferner ist es "au"serst hilfreich, die Situation $U \subseteq V$ unter einer
affin-algebra\-ischen Abbildung $g:V' \longrightarrow V$ zu betrachten und
die Affinit"at des Urbildes $U'=g^{-1}(U)$ zu untersuchen.
Gibt es beispielsweise eine affine Abbildung
$g:{\bf C}^2 \longrightarrow V$ mit der punktierten Ebene als Urbild
von $U$, so kann $U$ nicht affin sein.
Hierbei ist insbesondere das Schnittverhalten von $Y$ mit anderen
abgeschlossenen Untervariet"aten aufschlussreich.
\par\smallskip\noindent
Von diesem Abbildungstyp sind auch die Lokalisierungen
$V_x \longrightarrow V,\, x \in V,$
mit denen man das Affinit"atsproblem sozusagen auf den Punkt bringen kann.
Dabei ist $U$ genau dann affin, wenn dies lokal in jedem Punkt gilt.
Mit dieser Methode sieht man zum einen, dass eine offene
Teilmenge einer affinen Kurve wieder affin ist, und zum anderen,
dass in einer affinen Mannigfaltigkeit $V$
\-- also einer glatten Variet"at \-- das Komlement jeder Hyperfl"ache
affin ist. In diesem Fall sind die lokalen Ringe regul"ar und damit wie
der Polynomring faktoriell, so dass darin jede Hyperfl"ache durch eine
Funktion beschrieben wird.
\par\smallskip\noindent
Mit solchen Lokalisierungen
und allgemeineren Basiswechseln verl"asst man schnell den eingangs
beschriebenen
Rahmen von algebraischen Variet"aten und so ist es ratsam, von Anfang an in
der Kategorie der Schemata zu arbeiten.
\par
\bigskip
\noindent
\subsection{Das Thema aus schematheoretischer Sicht \--- Zur Notation}
\par
\bigskip
\noindent
Wir erinnern kurz an einige Begriffe aus der Theorie der Schemata, um
das Thema in dieser Sprache formulieren zu k"onnen, und an einige
Grundtatsachen, die wir durchgehend verwenden werden.
Ferner legen wir einige Bezeichnungen fest.
In der Terminologie orientieren wir uns an
EGA I-IV, \cite{EGAI}, \cite{EGAII}, \cite{EGAIII}, \cite{EGAIV},
Hartshorne \cite{haralg}
und f"ur die kommutative Algebra an Eisenbud \cite{eisenbud}.
\par\smallskip\noindent
Ein {\it affines Schema} ist das {\it Spektrum} eines kommutativen Ringes,
versehen mit der zubeh"origen Strukturgarbe.
Ein {\it Schema} ist ein beringter Raum, der von offenen Mengen "uberdeckt wird,
die affine Schemata sind.
\par\bigskip
\noindent
Eine offene Teilmenge eines Schemas ist wieder ein Schema. Dies gen"ugt f"ur
eine offene Menge in einem affinen Schema zu zeigen, f"ur
$U=D({\bf a}) \subseteq {\rm Spek}\, A$ ist aber
$U=\bigcup_{f \in {\bf a}}D(f)$ mit $D(f) \cong {\rm Spek}\, A_f$ affin.
Wir stellen uns die Frage:
\par\bigskip
\noindent
{\it Sei $D({\bf a})=U \subseteq X={\rm Spek}\, A$ eine
offene Teilmenge in einem affinen  Schema. Wann ist $U$ affin}?
\par\bigskip
\noindent
Wir sagen, dass eine abgeschlossene Teilmenge $Y=V({\bf a})$ eine
{\it Hyperfl"ache} ist, wenn jede Komponente von $Y$ die Kodimension eins hat.
Das ist gleichbedeutend dazu, dass alle minimalen Primoberideale zu $\bf a$
die H"ohe eins haben.
F"ur beliebige Ideale ${\bf a} \subseteq A$ halten wir uns an die Terminologie
von Nagata, dass die {\it H"ohe} des Ideals die minimale H"ohe der
minimalen Primoberideale von $\bf a$ ist, und die {\it Altitude}
die maximale H"ohe der minimalen Primoberideale.
Wir werden sp"ater allgemein sehen, dass bei $A$
noethersch nur Ideale eine affine offene Teilmenge definieren k"onnen,
wenn die Altitude davon $\leq 1$ ist.
\par\bigskip
\noindent
Eine {\it abgeschlossene Einbettung} eines Schemas ist eine topologisch
abgeschlossene Einbettung, so dass
der Garbenmorphismus der Strukturgarben surjektiv ist.
Eine abgeschlossene Einbettung in ein affines Schema ${\rm Spek}\, A$ ist
isomorph zu einer Abbildung ${\rm Spek}\, A/{\bf a} \longrightarrow {\rm Spek}\, A$.
Abgeschlossene Unterschemata eines affinen Schemas sind also affin.
\par\bigskip
\noindent 
Sind $X={\rm Spek}\, A$ und $Y={\rm Spek}\, B$ affine Schemata
"uber einem affinen Basisschema $S={\rm Spek}\, K$, so wird das Produkt
beschrieben durch $X\times_S Y= {\rm Spek}\, A\otimes_KB$ und ist insbesondere
selbst affin.
In einem beliebigen Produkt $X \times_S Y$ von Schemata bilden in kanonischer
Weise die $U \times_W V$ zu affinen, offenen Teilmengen der beteiligten
Schemata eine affine "Uberdeckung.
\par\bigskip
\noindent
Ein Schema $X$ hei"st {\it separiert}, wenn die {\it Diagonale}
$X \longrightarrow X \times _{{\rm Spek}\, {\bf Z}}X$ eine abgeschlossene
Einbettung ist. In diesem Fall ist f"ur $U,V \subseteq X$ beide affin
auch das Produkt $U \times V \subseteq X \times X$ affin und damit
auch das Urbild, und das ist $U \cap V$. Bei separierten Schemata ist also
der Durchschnitt von affinen offenen Teilmengen wieder affin.
\par\bigskip
\noindent
Sei $X$ ein Schema und $x \in X$. Mit $X_x= {\rm Spek}\, {\cal O}_{X,x}$
bezeichnen wir das {\it lokale Schema} im Punkt $x$. Dazu geh"ort eine
kanonische topologische Einbettung $X_x \hookrightarrow X$, die aber in aller
Regel nicht offen ist. Ist ein Morphismus $Y \longrightarrow X$ gegeben,
so sei mit $ Y_x = Y_{{\cal O}_x}= Y \times _X X_x$ die {\it lokale Faser}
bezeichnet.
Die echte Faser bezeichnen wir dagegen mit
$ Y_{k(x)}= Y \times_X {\rm Spek}\, k(x)$.
\par\bigskip
\noindent
Wir sagen, dass eine abgeschlossene Teilmenge $Y \subseteq {\rm Spek}\, A$
{\it geometrisch} (oder {\it mengentheoretisch})
durch eine Funktion beschrieben wird,
wenn $Y=V(f)$ ist f"ur ein $f \in A$. Ebenso sagen wir, dass zwei Ideale
geometrisch "ubereinstimmen, wenn sie dieselbe Nullstellenmenge definieren.
\par\bigskip
\noindent
Unter einer {\it Variet"at} verstehen wir ein integres, separiertes Schema von
endlichem Typ "uber einem K"orper.
\par
\bigskip
\noindent
\subsection{Die kanonische Abbildung}
\par
\bigskip
\noindent
Sei $X$ ein Schema mit globalem Schnittring $A=\Gamma(X,{\cal O}_X)$.
Dazu geh"ort dann der {\it kanonische Morphismus}
$\varphi: X \longrightarrow {\rm Spek}\, A=:X^{\rm aff}$, den man auch als
{\it Affinisierungsabbildung} bezeichnet.
Die universelle Eigenschaft dieses Morphismus ist, dass jeder andere
Morphismus in ein affines Schema dar"uber faktorisiert.
\begin{satz}
Ein Schema $X$ ist genau dann affin, wenn die Affinisierungsabbildung
eine Isomorphie ist. 
\end{satz}
{\it Beweis}. Trivial. \hfill $\Box$
\par\bigskip
\noindent
In diesem Abschnitt untersuchen wir die kanonische Abbildung und bringen
Eigenschaften des Schemas $X$ mit solchen der Abbildung in Zusammenhang.
Im Prinzip liefert jede Eigenschaft eines Morphismus eine Eigenschaft
f"ur das Schema, und jede Aussage, die es erlaubt,
durch das Zusammentreffen verschiedener Eigenschaften eines
Morphismus auf die Isomorphie zu schlie"sen, liefert bereits ein
Affinit"atskriterium. Ist umgekehrt $X$ nicht affin, so schl"agt sich das
in irgendeiner Weise als Defizit der Abbildung nieder.
\par\bigskip
\noindent
Zu einem $f \in \Gamma(X,{\cal O}_X)$ nennen wir die offene Menge
$X_f= \{ x \in X :\, f \in {\cal O}_x^\times \}$
die {\it Invertierbarkeitsmenge}
der Funktion $f$. $X_f$ ist das Urbild von $D(f)$ unter der kanonischen
Abbildung.
\begin{lem}
Sei $X$ ein quasikompaktes separiertes Schema und
$f \in A=\Gamma(X,{\cal O}_X)$. Dann gelten folgende Aussagen.
\par\smallskip\noindent
{\rm (1)}
Es ist $\Gamma(X_f,{\cal O}_X)= A_f$.
\par\smallskip\noindent
{\rm (2)} Ist $X_f$ affin, so induziert
die kanonische Abbildung $X \longrightarrow {\rm Spek}\, A$
eine Isomorphie $X_f \longrightarrow D(f)$
\end{lem}
{\it Beweis}.
(1) F"ur affines $U$ ist $X_f \cap U= D(f \vert U)$ wieder affin mit
$\Gamma (X_f \cap U,{\cal O}_X) = \Gamma(U,{\cal O}_X)_f$.
Sei $U_i,\,  i \in I,$ eine endliche affine "Uberdeckung von $X$;
wir betrachten die exakte Sequenz
$$0 \longrightarrow A \longrightarrow \prod_{i \in I}\Gamma(U_i,{\cal O}_X)
\longrightarrow \prod_{i,j \in I} \Gamma(U_i \cap U_j, {\cal O}_X) \,  ,$$
die den globalen Schnittring von $X$ festlegt. Dabei wird die letzte
Abbildung beschrieben durch
$(a_i)_{i \in I} \longmapsto (a_i-a_j)_{(i,j) \in I^2}$.
Wir machen die Nenneraufnahme zu $f$, die exakt ist. Da es sich um endliche
Produkte handelt,  
erhalten wir die exakte Sequenz
$$0 \longrightarrow A_f \longrightarrow
\prod_{i \in I}\Gamma(U_i,{\cal O}_X)_f
\longrightarrow \prod_{i,j \in I} \Gamma(U_i \cap U_j, {\cal O}_X)_f \, .$$
Wegen separiert sind
Durchschnitte affiner Mengen wieder affin, daher l"asst sich
nach der Vorbemerkung diese Sequenz
schreiben als
$$0 \longrightarrow A_f \longrightarrow
\prod_{i \in I}\Gamma(U_i \cap X_f ,{\cal O}_X)
\longrightarrow
\prod_{i,j \in I} \Gamma(U_i \cap U_j \cap X_f, {\cal O}_X) \, .$$
Somit liegt eine definierende Sequenz f"ur $\Gamma(X_f,{\cal O}_X)$ vor.
\par\smallskip\noindent
Zu (2). Wegen (1) liegt eine Isomorphie der globalen Schnittringe vor,
daher muss, da $X_f$ und $D(f)$ beide affin sind, eine Isomorphie
vorliegen. \hfill $\Box$
\begin{satz}
Sei $X$ ein Schema mit kanonischer Abbildung
$\varphi: X \longrightarrow {\rm Spek}\, A$. Dann gelten folgende Aussagen.
\par\smallskip\noindent
{\rm (1)} $\varphi$ ist dominant, d.h. $\varphi$ faktorisiert nicht durch
ein echtes abgeschlossenes Unterschema.
\par\smallskip\noindent
{\rm (2)}
Ist $U \subseteq {\rm Spek}\, A $ offen und affin und faktorisiert $\varphi$
durch $U$, so ist bereits $U={\rm Spek}\, A$.
\par\smallskip\noindent
{\rm (3)} Ist $X$ quasikompakt und separiert, so ist
$\varphi_\ast({\cal O}_X) = \tilde{A}$.
\end{satz}
{\it Beweis}.
Zu (1). Jeder Schemamorphismus in ein affines Schema,
bei dem der globale Ringhomomorphismus injektiv ist, ist dominant.
\par\smallskip\noindent
Zu (2).
Es liegt die Situation
$X \longrightarrow {\rm Spek}\, B \hookrightarrow {\rm Spek}\, A$ vor,
wobei die letzte Abbildung eine offene Einbettung ist.
Auf der Ringebene haben wir
$A \longrightarrow B \longrightarrow A=\Gamma(X,{\cal O}_X)$ vor, wobei die
erste Abbildung die Restriktionsabbildung ist, und wobei die Gesamtabbildung
die Identit"at ist.
F"ur die Spektren muss
dann in
${\rm Spek}\, A \longrightarrow {\rm Spek}\, B \longrightarrow {\rm Spek}\, A$
die hintere Abbildung surjektiv sein, also ist $U={\rm Spek}\, A$.
\par\smallskip\noindent
Zu (3).
Nach 1.3.2 ist
$\Gamma(D(h),\varphi_\ast({\cal O}_X))=\Gamma(X_h,{\cal O}_X)=A_h$ f"ur
$h \in A$.
Da $\varphi_\ast({\cal O}_X)$ eine Garbe ist, ist sie lokal auf einer Basis
festgelegt, also $= \tilde{A}$.  \hfill $\Box$
\par\bigskip
\noindent
{\bf Beispiel}
Sei $B$ ein Ring mit einem Element $f$ und einer Folge von
Elementen $s_n,\, n \in {\bf N}^\ast$ derart, dass
$f^{n-1} \cdot s_n \neq 0 $ und $ f^n \cdot s_n =0$ ist.
Sei $X$ die abz"ahlbar unendliche disjunkte Summe mit ${\rm Spek}\,B$ als
Komponenten und sei $A$ der globale Schnittring davon.
$A$ enth"alt die Elemente $f=(f,f,f,...) $ und $s=(s_1,s_2,s_3, ...) .$
Auf $U=X_f$ ist dann $s=0$, dagegen ist $s \neq 0$ in $A_f$.
Man kann also in dem Lemma auf die Voraussetzung quasikompakt nicht verzichten.
Ein Ring $B$ wie angegeben l"asst sich problemlos mit Variablen und
Relationen basteln.
$X$ ist reich beringt, da die einzelnen Komponenten durch
Funktionen, die auf dieser Komponente den Wert eins und sonst "uberall
den Wert null haben, beschrieben werden k"onnen.
\par
\bigskip
\medskip
{\it Quasiaffine Schemata}
\par
\bigskip
\noindent
Ein Schema $X$ nennen wir {\it reich beringt},
wenn die Invertierbarkeitsmengen von
globalen Funktionen die Topologie von $X$ definieren.
Ein offenes, quasikompaktes Unterschema eines affinen Schemas nennt man ein
{\it quasiaffines Schema}.
Eine invertierbare Garbe $\cal L$ auf einem separierten quasikompakten
Schema nennt man {\it ger"aumig}
(oder {\it ampel}),
wenn die Invertierbarkeitsmengen $X_s$ zu homogenen Schnitten
$s \in {\cal L}^n(X),\, n \geq 0, $ die Topologie definieren. Die beiden
letzten Definitionen halten sich an \cite{EGAII}, 4.5.2f und 5.1.1f. Es wird dabei
also jeweils die Quasikompaktheit vorausgesetzt. 
\begin{satz}
Sei $X$ ein Schema mit globalem Schnittring $A$.
Dann sind folgende Aussagen "aquivalent.
\par\smallskip\noindent
{\rm (1)} $X$ ist reich beringt.
\par\smallskip\noindent
{\rm (2)} $\varphi: X \longrightarrow {\rm Spek}\, A$ ist eine topologische
Einbettung.
\par\medskip\noindent
Ist $X$ zus"atzlich quasikompakt und separiert, so sind folgende
Aussagen "aquivalent.
\par\smallskip\noindent
{\rm (1)} $X$ ist reich beringt.
\par\smallskip\noindent
{\rm (2)} Die Strukturgarbe ${\cal O}_X$ ist ger"aumig.
\par\smallskip\noindent
{\rm (3)} Es gibt eine Familie $f_i \in A,\, i \in I,$ globaler Funktionen mit
$X_{f_i}$ affin und $X=\bigcup_{i \in I}X_{f_i}$.
\par\smallskip\noindent
{\rm (4)} Die kanonische Abbildung $X \longrightarrow {\rm Spek}\, A$ definiert
eine Isomorphie auf ein offenes Unterschema von ${\rm Spek}\, A$
\par\smallskip\noindent
{\rm (5)} $X$ ist quasiaffin.
\end{satz}
{\it Beweis}.
Der erste Teil ist klar. Zum zweiten Teil.
Die "Aquivalenz von (1) und (2) ist klar. Sei dies erf"ullt und $x \in U$ 
eine affine Umgebung. Dann gibt es eine globale Funktion $f \in A$ mit
$x \in X_f \subseteq U$. Innerhalb von $U$ ist dann $X_f=D(f)$ und damit
ist $X_f$ ebenfalls affin. Solche affinen Mengen "uberdecken dann $X$
(Das gilt schon unter reich beringt). Es gilt also (3).
Da $X$ separiert und quasikompakt ist, so ist nach Lemma 1.3.2
$X_f \longrightarrow D(f) \subseteq {\rm Spek}\, A$ eine Isomorphie und
damit identifiziert die kanonische Einbettung $X$ mit einem offenen
Unterschema von $X$. Es gilt also (4) und damit sofort (5), und aus (5) folgt
sofort wieder (1). \hfill $\Box$
\par\bigskip
\noindent
Die Quasiaffinit"at bedeutet insbesondere, dass die kanonische Abbildung
injektiv ist. Unter recht allgemeinen Bedingungen folgt umgekehrt
daraus schon die Quasiaffinit"at.
Dieses Kriterium beruht auf folgendem nicht trivialen Satz,
den wir hier zitieren und der auf der Theorie der \'{e}talen Morphismen
beruht:
{\it Ein separierter, quasiendlicher Schemamorphismus $f:X \longrightarrow Y$
ist quasiaffin}.
Dieser Satz findet sich in \cite{EGAIV}, 18.12.12.
Ein {\it quasiendlicher} Morphismus ist ein Morphismus vom endlichen Typ, dessen
Fasern nur endlich viele Punkte haben, und {\it quasiaffin} bedeutet, dass
es eine affine "Uberdeckung von $Y$ gibt, deren Urbilder alle quasiaffin sind.
Es sind dann die Urbilder aller affinen Teilmengen von $Y$ quasiaffin,
siehe \cite{EGAII}, 5.1.8.
Wir wenden den Satz auf den kanonischen Morphismus an und gelangen zu folgendem
Resultat.
\begin{satz}
Sei $X$ ein separiertes Schema von endlichem Typ "uber einem affinen
Basisschema, und $A$ sei der globale Schnittring.
Dann sind folgende Aussagen "aquivalent.
\par\smallskip\noindent
{\rm (1)} $X$ ist quasiaffin.
\par\smallskip\noindent
{\rm (2)}
Die kanonische Abbildung $X \longrightarrow {\rm Spek}\, A$ ist injektiv.
\par\smallskip\noindent
{\rm (3)} Die kanonische Abbildung ist quasiendlich.
\par\smallskip\noindent
{\rm (4)} $X$ ist quasiendlich "uber einem affinen Schema.
\par\medskip\noindent
Ist $X$ eine Variet"at "uber einem K"orper $K$,
so sind diese Bedingungen auch "aquivalent zu
\par\smallskip\noindent
{\rm (5)} Die globalen Funktionen trennen abgeschlossene Punkte.
\par\smallskip\noindent
{\rm (6)}
Zu jeder Kurve $C \subseteq X$ gibt es eine globale Funktion auf $X$,
die auf $C$ nicht konstant ist.
\end{satz}
{\it Beweis}.
Die Implikationen von (1) nach (2) und von (3) nach (4)
sind trivial.
Von (2) nach (3). Da $X$ von endlichem Typ "uber einem affinen Basisschema ist,
ist auch die kanonische Abbildung von endlichem Typ.
Sei (4) erf"ullt. Nach Voraussetzung an $X$ ist der
angegebene Morphismus auch separiert,
und der oben angef"uhrte Satz liefert (1).
\par\smallskip\noindent
Sei jetzt $X$ eine Variet"at "uber einem K"orper.
In den Formulierungen (5) und (6) fassen wir globale Funktionen
als Morphismen nach ${\rm Spek}\, K[T]={\bf A}^1_K$ auf, ein solcher
Morphismus faktorisiert nat"urlich durch die kanonische Abbildung.
\par\smallskip\noindent
Von (1) nach (5). Sind $P$ und $Q$ verschiedene abgeschlossene
Punkte, so ist $U:=X-\{ Q \}$ eine offene Umgebung von $P$ und
wegen quasiaffin gibt es ein $f \in A$ mit
$P \in X_f \subseteq U$.
Fasst man $f$ als Morphismus auf die affine Gerade auf, so wird $Q$ auf
den Nullpunkt abgebildet, $P$ aber auf einen Punkt in $D(T)$.
\par\smallskip\noindent
Von (5) nach (6). Sei $C \subseteq X$ eine Kurve, d.h. ein integres,
abgeschlossenes Unterschema der Dimension eins.
Die Kurve ist dann insbesondere wieder vom endlichen Typ "uber $K$,
und somit liegt ein Jacobsonschema vor, siehe \cite{EGAI}, § 6.4,
d.h. die abgeschlossenen Punkte
liegen dicht (und zwar in jeder offenen Menge)
und insbesondere gibt es zwei verschiedene abgeschlossene
Punkte auf $P,Q \in C$, die auch in $X$ abgeschlossene Punkte sind.
Nach Voraussetzung gibt es eine Funktion $f$
mit $f(P) \neq f(Q) \in {\rm Spek}\, K[T]={\bf A}^1_K$, und diese Funktion ist
nat"urlich nicht konstant auf $C$.
\par\smallskip\noindent
Von (6) nach (4).
Wir verwenden folgende Hilfskonstruktion, die induktiv angewendet wird.
Sei $P \in X$ ein abgeschlossener Punkt von $X$ auf einer irreduziblen
Komponente maximaler Dimension und
seien $X_1,...,X_n$ alle Komponenten von $X$ durch $P$.
All diese Komponenten haben dann eine Dimension von zumindest eins,
andernfalls ist $X$ nulldimensional.
Auf jeder Komponente gibt es dann eine Kurve $P \in C_i \subseteq X_i$ und
damit Funktionen $f_i \in A$, die auf den Kurven nicht konstant sind.
Der durch $f=(f_1,...,f_n)$ definierte Morphismus $X \longrightarrow {\bf A}^n$
ist dann insbesondere auf keiner der Komponenten konstant und damit ist
dim $f^{-1}(f(P)) < {\rm dim}\, X$. $f(P)$ ist ein abgeschlossener Punkt,
daher ist die Faser wieder eine (in aller Regel reduzible) Variet"at.
\par\smallskip\noindent
Auf dieser Faser w"ahlen wir ebenso einen Punkt und finden dann wieder
Funktionen, um die wir die Abbildung $f$ verl"angern. So sukzessive fortfahrend
findet man eine Abbildung von $X$ in einen affinen Raum ${\bf A}_K^m$,
die in einem Punkt eine endliche Faser besitzt.
\par\smallskip\noindent
Wir wollen zeigen, dass es eine quasiendliche Abbildung
$X \longrightarrow {\bf A}^k$ gibt, dazu gen"ugt es zu zeigen,
dass dies auf jeder Komponente erreicht werden kann (mit auf $X$ definierten
Funktionen).
Sei also $X$ integer und $f:X \longrightarrow {\bf A}^n$ eine Abbildung mit
$f^{-1}(f(P))$ nulldimensional.
Nach \cite{eisenbud}, Theorem 14.8, gibt es dann eine offene
Umgebung $V$ von $P$,
auf der die Abbildung nulldimensionale Fasern besitzt und
quasiendlich ist.
$X-V$ hat dann eine kleinere Dimension als $X$.
Durch Induktion folgt dann, dass es eine quasiendliche Abbildung in einen
affinen Raum gibt. \hfill $\Box$
\par\bigskip
\noindent
{\bf Bemerkung}
In der komplex-analytischen Geometrie nennt man einen komplexen
Raum holomorph separabel, wenn es analog zu der Eigenschaft
(5) zu verschiedenen Punkten eine global definierte holomorphe Funktion
gibt, die diese Punkte trennt, und holomorph ausbreitbar, wenn es
analog zu der Eigenschaft (4) zu einem Punkt holomorphe
Funktionen gibt, so dass der Punkt unter dieser Abbildung isoliert
liegt.
Gilt in holomorph separablen oder holomorph ausbreitbaren komplexen R"aumen,
dass es entsprechend (1) zu einem Punkt $P$ und einer dazu disjunkten
analytischen Teilmenge $A$ eine holomorphe Funktion gibt mit
$f(P)=1, f \vert A =0$?
Ist ein solcher Raum offen einbettbar in einen Steinschen Raum?
Wir werden gelegentlich auf solche analogen Situationen aus der komplexen
Analysis hinweisen, und orientieren uns dabei an \cite{grauert}.
\par\bigskip
\noindent
{\bf Bemerkung}
Die Eigenschaft, dass es auf $X$ keine projektiven Kurven gibt, dass
es also auf jeder Kurve auf $X$ nicht konstante Funktionen gibt,
ist nicht "aquivalent zur Eigenschaft (6) im Satz.
Ein Beispiel von dieser Art, bei dem alle globalen Funktionen
konstant sind, findet sich im 4.3.
\par\bigskip
\noindent
F"ur eine quasiaffine Variet"at ist a priori nicht klar, ob sie sich
als offene Untervariet"at einer affinen Variet"at realisieren l"asst.
Die kanonische Abbildung liefert lediglich eine offene Einbettung in ein
affines Schema, das aber nicht vom endlichen Typ sein muss, da
der globale Schnittring $A=\Gamma(U,{\cal O}_U)$ nicht von
endlichem Typ sein muss.
\begin{satz}
Sei $U$ ein Schema von endlichem Typ "uber einem noetherschen
Ring $K$. Ist $U$ quasiaffin, so ist $U$ offen einbettbar in ein affines Schema
vom endlichen Typ "uber $K$.
\end{satz}
{\it Beweis}.
Sei $U$ reich beringt und "uberdeckt von affinen Menge $U_i=U_{f_i}, i=1,...,n$
mit globalen Funktionen $f_i \in \Gamma(U,{\cal O}_U)=A$ und so,
dass die Algebren $ A_{f_i} = \Gamma(U_i,{\cal O}_U)$ vom endlichen Typ
"uber $K$ sind. Sei
$$A_i=A_{f_i} = A_i = K[a_{i,1}/ f_i,...,a_{i,m_i}/f_i]
= K[1/f_i,a_{i,1},...,a_{i, m_i}]$$
mit $a_{i,1},...,a_{i,m_i} \in A $ (dabei muss man $f_i$ durch eine geeignete
Potenz ersetzen).
Seien $a_1,...,a_m$ alle dabei beteiligten Elemente.
\par\smallskip\noindent
Wir betrachten die Abbildung
$$ F=(f_1,...,f_n,a_1,...,a_m):\ U \longrightarrow {\bf A}^{n+m}_K  \ ,$$
und wollen zeigen, dass es sich dabei um eine topologische Einbettung
handelt, deren Bild offen im Abschluss des Bildes ist.
Seien Punkte $x \neq y \in U$ gegeben. Bei $x \in U_i, y \not\in U_i $
ist $f_i(x) \neq f_i(y)=0$. Ist hingegen $x,y \in U_i={\rm Spek}\, A_i$,
so k"onnen nicht alle
Koordinatenfunktionen von $F$ darauf gleiche Werte haben, da man mit ihnen
ja ein Algebra-Erzeugendensystem von $A_i$ basteln kann. Damit ist $F$
injektiv.
\par\smallskip\noindent
Es sei $V$ der Abschluss von $F(U)$ in ${\bf A}_K^{r+m}$. Dieser ist gleich
der Vereinigung der Abschl"usse der $F(U_i)$ darin, wie eine einfache
topologische "Uberlegung zeigt. Wir haben also nur noch zu zeigen,
dass die eingeschr"ankte Abbildung
$F \vert_{U_i} \longrightarrow {\bf A}_K^{r+m}$ eine offene Einbettung
in den Abschluss hinein ist. Dies zeigt aber das folgende Lemma. \hfill $\Box$
\begin{lem}
Sei $A=K[1/f,a_1,...,a_n]$ mit $f \in A^{\times} $. Dann
ist die Abbildung
$ F=(f,a_1,...,a_n): X={\rm  Spek} \, A \longrightarrow {\bf A}^{1+n} $
eine offene Einbettung auf den Abschluss des Bildes.
\end{lem}
{\it Beweis}.
Wir betrachten
$$ X \stackrel{F}{\longrightarrow} {\bf A}^{\times} \times {\bf A}^n
\buildrel 1/ x_0 \times id \over \longrightarrow {\bf A}^{1+n} \ .$$
Dabei haben wir f"ur $F$ den Bildbereich eingeschr"ankt.
Die Gesamtabbildung ist eine abgeschlossene Einbettung, da es sich
bei $ 1/f , a_1,...,a_n$ um ein Algebra\-Erzeugendensystem von $A$
handelt, und die hintere Abbildung ist eine offene Einbettung. 
Damit ist $F(X)$ abgeschlossen in ${\bf A}^\times \times {\bf A}^n $.
Eine einfache topologische "Uberlegung zeigt nun, dass $F(X)$ im
urspr"unglichen Bildbereich offen im Abschluss $V$ ist, und zwar gilt dabei
$F(U)= V \cap D(X_0) = V \cap ({\bf A}^\times \times {\bf A}^n) \ .$
\footnote{Ist
$X \subseteq W \subseteq Y $ mit $W$ offen in $Y$ und $X$ abgeschlossen in
$W$, so ist $X$ offen im Abschluss $V$ von $X$ in $Y$, es gilt n"amlich
$X=W \cap V$, wobei $\subseteq$ trivial ist. Sei $x \in W,\, x \not\in X $.
Dann gibt es eine offene Menge $Z \subseteq W $ mit $x \in Z$ und
$Z \cap X = \emptyset$. $Z$ ist nat"urlich auch offen in $Y$, also kann
$x$ nicht zum Abschluss von $X$ geh"oren.}
\par\smallskip\noindent
Die Abbildung $B:=\Gamma(V,{\cal O}_V) \longrightarrow A$ ist injektiv,
da $V$ der schematheoretische Abschluss ist. Andererseits
ist $F(U) = D(X_0)$ (in $V$), und somit faktorisiert diese Abbildung
durch $B_{X_0} \longrightarrow A$. Diese Abbildung ist aber surjektiv, da
ja $ 1/f ,a_1,...,a_n$ ein Erzeugendensystem von $A$ bilden. Damit
ist diese letzte Abbildung ein Ringisomorphismus,
und $U \cong D(X_0)\ .$ \hfill $\Box$
\par
\bigskip
\medskip
{\it Der globale Hilbertsche Nullstellensatz}
\par
\bigskip
\noindent
Der {\it Hilbertsche Nullstellensatz} besagt f"ur
ein affines Schema, dass eine Funktion $f \in A$ zum Radikal eines Ideals
${\bf a}$ geh"ort, sobald deren
Nullstellenmenge $V(f)$ die Nullstellenmenge des Ideals
$V({\bf a})$ umfasst.
Diese Aussage ist nat"urlich f"ur jedes Schema formulierbar, wenn man f"ur
$A$ den globalen Schnittring nimmt.
\par\smallskip\noindent
Nimmt man f"ur die Funktion $f$ eine Einheit,
so ist die Nullstellenmenge $V(f)$ leer,
und der Hilbertsche Nullstellensatz besagt
dann, dass ein Ideal $\bf a$ ohne gemeinsame Nullstelle das Einheitsideal
sein muss.
\par\smallskip\noindent
Wir sagen, dass in einem Schema $X,{\cal O}_X$ mit globalem Schnittring $A$
der {\it globale Hilbertsche Nullstellensatz} gilt,
wenn globale Funktionen, die keine
gemeinsame Nullstelle in $X$ besitzen, das Einheitsideal
erzeugen.
\begin{satz}
F"ur ein Schema $X,{\cal O}_X$ mit globalem Schnittring $A$ sind
folgende Aussagen "aquivalent.
\par\smallskip\noindent
{\rm (1)} In $X,{\cal O}_X $ gilt der globale Hilbertsche Nullstellensatz.
\par\smallskip\noindent
{\rm (2)}
Aus einer "Uberdeckung $X = \bigcup_{i \in I} X_{f_i},\, f_i \in A, $ folgt
$ 1 \in ( f_i ,\, i \in  I ) $.
\par\smallskip\noindent
{\rm (3)} Unter der kanonischen Abbildung
$ \varphi : X \longrightarrow { \rm Spek}\, A $ werden alle
maximalen Ideale erreicht.
\end{satz}
{\it Beweis}.
(1) und (2) sind lediglich Umformulierungen.
Sei (2) erf"ullt und $\bf m$ ein maximales Ideal in $A$.
Dann "uberdecken nach Voraussetzung die $\{ X_f : f \in {\bf m} \}$
das Schema $X$ nicht, sei $x$ ein Punkt au"serhalb dieser Vereinigung.
Dann ist $\varphi (x)= \{ f \in A: x \not\in X_f \}={\bf m}$,
wobei die Inklusion
$\supseteq $ nach Konstruktion gilt und wegen der Maximalit"at von $\bf m$
nicht echt sein kann.
\par\smallskip\noindent
Sei (3) erf"ullt, und $X= \bigcup_{i \in I}X_{f_i}$. Angenommen, die $f_i$ 
w"urden nicht das Einheitsideal erzeugen, so w"aren sie in einem maximalen
Ideal enthalten, das dann aber nicht von einem Punkt aus $X$ getroffen werden
k"onnte. \hfill $\Box$
\par\bigskip
\noindent
{\bf Bemerkung}
In einem affinem Schema gilt der globale Hilbertsche Nullstellensatz.
Ebenso in jedem Schema, dessen globaler Schnittring ein K"orper ist, etwa
in einer projektiven integren Variet"at.
\par\bigskip
\noindent
Wenn f"ur ein Schema der Hilbertsche Nullstellensatz gilt, so folgt
daraus nicht, dass die kanonische Abbildung surjektiv ist, wie das n"achste
Beispiel zeigt.
\par\bigskip
\noindent
{\bf Beispiel}
Sei $X$ eine dreidimensionale, affine, glatte Variet"at
"uber einem algebraisch abgeschlossenen K"orper $K$,
$P$ ein abgeschlossener Punkt darauf und $C$  eine
irreduzible Kurve in $X$ durch $P$.
Wir blasen in $X$ den Punkt $P$ auf, $P$ wird also ersetzt
durch einen ${\bf P}^2$.
Die Aufblasung sei $\tilde{X}$. Wir betrachten die eigentliche
Transformierte $\tilde{C}$ von $C$, die wieder eine irreduzible Kurve ist.
Unser Beispiel ist jetzt
$$W=\tilde{X} - \tilde {C} \hookrightarrow \tilde{X}
\longrightarrow X = {\rm Spek}\, A $$
Bei diesen Abbildungen "andert sich jeweils der globale Schnittring nicht,
es handelt sich also insgesamt um eine Affinisierung.
Der Punkt $P$ wird dabei erreicht, die Faser dar"uber besteht in
$\tilde{X}$ aus einem zweidimensionalen projektiven Raum, und daraus
wird nur etwas nulldimensionales rausgenommen,
es bleiben also (viele) Punkte in der Faser "ubrig.
\par\smallskip\noindent
Dagegen wird die Kurve $C$ nicht erreicht, da au\ss erhalb von $P$
kein Punkt darauf erreicht wird.
Damit haben wir ein Beispiel, das zeigt, dass unter einer
Affinisierungsabbildung das Bild nicht abgeschlossen sein muss unter
Generalisierungen.
\par\smallskip\noindent
Lokalisiert man die Situation in $P$,
so erh"alt man ein Beispiel,
wo das maximale Ideal erreicht wird, und folglich der globale
Hilbertsche Nullstellensatz gilt, die kanonische Abbildung aber nicht
surjektiv ist (die kanonische Abbildung ist mit Lokalisierung vertr"aglich,
siehe 1.7).
\begin{satz}
F"ur ein Schema $X$ sind folgende Eigenschaften "aquivalent.
\par\smallskip\noindent
{\rm (1)} $X$ ist affin.
\par\smallskip\noindent
{\rm (2)} $X$ ist reich beringt und es gilt der globale
Hilbertsche Nullstellensatz.
\par\smallskip\noindent
{\rm (3)} Es gibt eine "Uberdeckung aus affinen Mengen der Form $X_{f_i}$ derart,
dass die $f_i$ in $\Gamma(X,{\cal O}_X)$ das Einheitsideal erzeugen.
\end{satz}
{\it Beweis}.
Von (1) nach (2).
Ist $X$ affin, so gelten diese beiden Aussagen.
Von (2) nach (3). Die Existenz einer solchen affinen "Uberdeckung folgt bereits
aus reich beringt, und der globale Hilbertsche Nullstellensatz sichert, dass
die Funktionen das Einheitsideal erzeugen.
Von (3) nach (1).
Mit dem globalen Hilbertschen Nullstellensatz ist $X$
insbesondere quasikompakt, ferner ist $X$ separiert,
da die Diagonale zur kanonischen Abbildung
zun"achst lokal und damit "uberhaupt abgeschlossen ist.
Daher ist die kanonische Abbildung
eine Isomorphie auf den beteiligten $X_f$.
Da die zugeh"origen $D(f)$ aber ${\rm Spek}\, \Gamma(X,{\cal O}_X)$
"uberdecken, liegt eine globale Isomorphie vor. \hfill $\Box$
\par\bigskip
\noindent
Im Verlauf des ersten Kapitels werden uns noch zwei weitere
Begriffe begegnen, die mittels der kanonischen Abbildung definiert
werden k"onnen. Zum einen die {\it pseudoaffinen} oder
{\it generisch affinen} Schemata, bei denen die kanonische Abbildung
generisch eine Isomorphie ist, siehe 1.7, und zum anderen
die {\it semiaffinen} Schemata,
bei denen die kanonische Abbildung eigentlich ist, siehe 1.12.
\par
\bigskip
\noindent
\subsection{Der globale Schnittring eines offenen Unterschemas}
\par
\bigskip
\noindent
Sei $A$ ein zun"achst beliebiger kommutativer Ring, $X={\rm Spek}\, A$
und $U=D({\bf a})$ eine offene Menge in $X$.
Ist $U=D(f)$, so ist der globale Schnittring dazu $A_f$, und $U$ ist isomorph
zu ${\rm Spek}\, A_f$.
Der globale Schnittring einer beliebigen offenen Teilmenge ist gegeben durch
den projektiven Limes $\lim_{D(f) \subseteq U} A_f$, wobei das Indexsystem
in nat"urlicher Weise geordnet ist. Ein Element aus $\Gamma(U,{\cal O}_X)$
ist gegeben durch eine "Uberdeckung $U= \bigcup_{i \in I} D(f_i)$ und
Elementen $(a_i)_{i \in I} \in (A_{f_i})_{i \in I}$ mit den
Vertr"aglichkeitsbedingungen $ a_i =a_j$ in $A_{f_if_j}$.
\par\smallskip\noindent
Ist $A$ integer, so gibt es eine einfachere Beschreibung des globalen
Schnittringes. In diesem Fall sind mittels den Restriktionen
und den Halmabbildungen auf den generischen Punkt
alle Schnittringe von offenen Mengen
und alle Halmringe in kanonischer Weise Unterringe des Quotientenk"orpers
$Q(A)$ (bzw. des Funktionenk"orpers, wenn es sich um ein integres Schema
handelt).
Im integren Fall werden die Schnittringe beschrieben durch
$\Gamma(U,{\cal O}_X)=  \bigcap_{x \in U} {\cal O}_x $ und f"ur eine
beliebige "Uberdeckung $U=\bigcup_{i \in I} U_i$ gilt
$\Gamma(U,{\cal O}_X)=\bigcap_{i \in I} \Gamma(U_i,{\cal O}_X)$. Insbesondere gilt
$\Gamma(D({\bf a}),{\cal O}_X)= \bigcap_{f \in {\bf a}} A_f =
\bigcap_{{\bf p} \in U}A_{\bf p}$.
Dabei sieht man sofort, dass die Durchschnittbildung der Halme eine
Garbe liefert, und man hat nur noch
$A=\bigcap_{{\bf p} \in {\rm Spek}\, A} A_{\bf p}$ zu zeigen, und dies ist ein
elementarer Schluss: ist $a/ b$ im Durchschnitt drin, so gibt es
zu jedem $\bf p$ eine
Darstellung $a/ b = {a_i/ b_i}$ mit $b_i \not\in {\bf p}$.
Die $b_i$ erzeugen dann das Einheitsideal, also $1=\sum_{i=1,..,n}c_ib_i$,
und somit gilt
$$a= \sum_{i=1,...,n} ac_ib_i = \sum_{i=1,...,n} ba_ic_i
=b(\sum_{i=1,...,n}a_ic_i) \, ,$$
also
wird $a$ von $b$ geteilt.
\par\bigskip
\noindent
F"ur $U=D({\bf a}),\, {\bf a} = (f_i,\, i \in  I)$ und eine rationale Funktion
$q \in Q(A)$ ist $q \in \Gamma(U,{\cal O}_X)$ genau dann,
wenn es f"ur jedes $i$ eine Darstellung $ q = a_i / f^{n_i}_i$ gibt,
wenn also $q \cdot f_i^{n_i} \in A$ erreicht werden kann.
Dies ist der Ausgangspunkt einer weiteren Darstellung des globalen Schnittringes
im noetherschen integren  Fall als sogenannte ${\bf a}-${\it Transformierte}
$$ T_{\bf a}(A)=\{ q \in Q(A) :
\mbox{ Es gibt ein } n \in{\bf N}\,
\mbox{ mit } \, q{\bf a}^n \subseteq A \} \, .$$
Dabei ist klar, dass rationale Funktionen der ${\bf a}-$Transformierten zum
globalen Schnitt\-ring geh"oren; umgekehrt geh"oren zun"achst nur $qf_i^{n_i}$ f"ur ein
Radikalerzeugendensystem $f_i$ zu $A$, dann findet man aber auch f"ur ein
(endliches) Erzeugendensystem $f_j,\, j \in J,$ von $\bf a$ Potenzen mit
$qf_j^{n_j} \in A$ und f"ur die Summe $n$
dieser Potenzen ist dann $q{\bf a}^n \subseteq A$.
\par\bigskip
\noindent
Im Fall eines integren Schemas $X$
hat eine rationale Funktion $q \in K(X)$ einen eindeutig
bestimmten {\it Definitionsbereich}, n"amlich
$$ {\rm Def}\, (q)=\{ x : q \in {\cal O}_x \}
=\bigcup_{f \in A:\, q = a/f^n} D(f) \, ,$$
wobei die zweite Darstellung $X={\rm Spek}\, A$ affin voraussetzt.
Dieser Definitionsbereich ist offen, und wird im affinen Fall durch das Ideal
$\{ f \in A : qf^n \in A \}$ beschrieben.
Zu einer rationalen Funktion $q \neq 0$ nennen wir
$N(q)=\{ f \in A : qf \in A \}$ das (exakte) {\it Nennerideal} und
$qN(q) = Z(q)$ das (exakte) {\it Z"ahlerideal}.
Das Nennerideal ist dabei wirklich ein Ideal, wie aus der Definition
sofort folgt, und das Z"ahlerideal ist als der Schnitt von $A$ mit dem
Bild von $N(q)$ unter dem durch $q$ definierten $A-$Modul\--
Homomorphismus $Q(A) \longrightarrow Q(A)$ ebenfalls ein Ideal.
Die zugeh"origen Radikale nennen wir dann Z"ahlerradikal und Nennerradikal.
Betrachtet man $q^{-1}$, so drehen sich die Rollen von Z"ahler- und
Nennerideal um.
\par\bigskip
\noindent
Ist $X$ integer, noethersch und separiert, so bilden die Definitionsbereiche
rationaler Funktionen eine Subbasis der Topologie von $X$,
siehe \cite{EGAI}, Cor. 8.5.8.
Die gemeinsamen Definitionsbereiche zu endlich vielen rationalen Funktionen
bilden damit eine Basis der Topologie. Diese Aussage geht wesentlich in die
Beweise der S"atze 1.7.5, 1.7.6 und 1.8.1 ein.
\par\bigskip
\noindent 
Die Vereinigung der Definitionsbereiche ${\rm Def}\, (q) \cup {\rm Def}\, (q^{-1})$
nennen wir den {\it meromorphen Definitionsbereich} von $q$,
${\rm Def}^{\rm mer}\, (q)={\rm Def}^{\rm mer}\,(q^{-1})$.
Ist $A$ eine $R-$Algebra, so legt $q$ eine meromorphe Funktion
$q: {\rm Def}^{\rm mer}\, (q) \longrightarrow {\bf P}^1_R$ fest, die sich durch
$X_0/X_1=T \longmapsto q$ und $X_1/X_0=T^{-1} \longmapsto q^{-1}$ ergibt.
Der meromorphe Definitionsbereich
wird beschrieben durch $D(N(q)+Z(q))$.
Eine wichtige Eigenschaft des meromorphen Definitionsbereiches ist, dass er
sich invariant unter der Addition mit einer ganzen Funktion verh"alt.
\begin{prop}
Sei $A$ ein Integrit"atsbereich und $q \in Q(A)$ eine rationale Funktion.
Es sei $ZN(q)$ das Radikal aus Nennern und Z"ahlern.
Dann gilt f"ur eine beliebiges $a \in A$ die Gleichheit
$ZN(q) =ZN(q+a)$.
\end{prop}
{\it Beweis}.
Es gen"ugt, die eine Inklusion zu zeigen, sei $q=  f/g  $.
Dann ist $q+a=(f+ag)/g$. Damit ist $g \in ZN(q+a)$ und $f+ag \in ZN(q+a)$.
Dann ist aber auch $f \in ZN(q+a)$. \hfill $\Box$
\par\bigskip
\noindent
Im normalen Fall gibt es weitere Darstellungen des Schnittringes und der
Definitionsbereiche von rationalen Funktionen, die wir im Abschnitt
"uber Krullringe besprechen.
\begin{satz}
Sei $D({\bf a})=U \subseteq {\rm Spek}\, A$ ein offenes Unterschema.
Ist $U$ quasikompakt, so
ist $U \longrightarrow {\rm Spek}\, \Gamma(U,{\cal O}_X)$
eine offene Einbettung mit $U \cong D({\bf a}\Gamma(U,{\cal O}_X)$.
$U$ ist genau dann affin, wenn 
${\bf a}$ in $\Gamma(U,{\cal O}_X)$ das Einheitsideal erzeugt.
\end{satz}
{\it Beweis}.
Dies folgt direkt aus 1.3.4(4) und 1.3.9.
Wir geben f"ur die zweite Aussage noch einen weiteren Beweis, der sich an den
einleitenden "Uberlegungen in 1.1 anlehnt, und dabei die
Bedeutung der in einer Darstellung der Eins $1=q_1f_1,...,q_nf_n$
auftretenden Funktionen $q_1,...,q_n$ kl"art.
Die durch diese Funktionen definierte Abbildung $q=(q_1,...,q_n)$
ist n"amlich eine abgeschlossene Einbettung
$$q: U \longrightarrow {\rm Spek}\, A \times_{{\rm Spek}\, {\bf Z}}{\bf A}_{\bf Z}^n
={\bf A}^n_{{\rm Spek}\, A}= {\rm Spek}\, A[T_1,...,T_n] \, .$$
Dabei geht $T_i$ auf $q_i$.
Es ist klar, dass eine topologische Einbettung vorliegt, da die vorgegebene
Einbettung $U \subseteq {\rm Spek}\, A$ durch $q$ faktorisiert.
Es ist ${\rm Spek}\, A[T_1,...,T_n]=
D(f_1) \cup ... \cup D(f_n) \cup D(\sum_{i=1}^n f_iT_i-1)$ und das Bild von $q$
liegt in $V(\sum_{i=1}^n f_iT_i-1)$. Die Abgeschlossenheit ist also nur noch
auf den $D(f_i)$ zu "uberpr"ufen, dort
hat aber die Abbildung die Gestalt
${\rm Spek}\, A_{f_j} \longrightarrow {\rm Spek}\, A_{f_j}[T_1,...,T_n]$
mit $q_i \in A_{f_j}$, so dass eine abgeschlossene Einbettung
vorliegt. \hfill $\Box$
\par\bigskip
\noindent
In der Argumentation haben wir folgenden Satz mitbewiesen.
\begin{satz} 
Sei $A$ ein Ring und $U \subseteq {\rm Spek}\, A$ eine offene affine
Teilmenge. Dann ist
$\Gamma(U,{\cal O}_X) $ eine $A-$Algebra von endlichem Typ. \hfill $\Box$
\end{satz}
Die folgende Aussage ist die Grundlage f"ur die Superh"ohenkriterien
in 4.2. Sie gilt im normalen Fall ohne zus"atzliche Voraussetzung
an den Schnittring, da dieser dann immer ein Krullbereich ist.
\begin{prop}
Sei $A$ ein noetherscher integrer Ring und
$U$ eine nicht affine Teilmenge mit noetherschem globalen Schnittring
$B=\Gamma(U,{\cal O}_X)$.
Dann ist die H"ohe des Erweiterungsideales ${\bf b}={\bf a}B$ $\geq 2$.
\end{prop}
{\it Beweis}.
Nach Satz 1.4.2 ist $D({\bf a}) \cong D({\bf b}) \subseteq {\rm Spek}\, B$
und das Erweiterungsideal ${\bf b}B$ ist nicht das Einheitsideal. 
Damit ist $D({\bf b})$ eine offene Teilmenge in einem noetherschen
integren Schema mit gleichem globalen Schnittring $B$.
Sei ${\bf p}$ ein Primoberideal zu
${\bf b}={\bf a}B$.
Dann ist $D({\bf b}) \subseteq D({\bf p}) \subseteq {\rm Spek}\, B$,
da $B$ integer ist muss dann auch $\Gamma(D({\bf p}),{\cal O}_B)=B$ sein.
\par\smallskip\noindent
Angenommen, $\bf p$ habe die H"ohe eins, sei ${\bf p}=(f_1,...,f_n)$.
In der Lokalisierung $A_{\bf p}$ wird ${\bf p}$ geometrisch durch eine
Funktion $0 \neq f \in {\bf p}$ beschrieben, dies f"uhrt zu Relationen
$f_i^m={a_i \over r_i} f$ mit $r_i \not\in {\bf p}$.
Durch "Ubergang zum Hauptnenner $r=r_1 \cdot ... \cdot r_n \not\in {\bf p}$
hat man auch $f_i^m={b_i \over r} f$.
Das bedeutet aber $ r/f= b_i/f_i^m$ und damit haben wir eine
auf $D({\bf p})$ definierte Funktion. Diese kann nicht zu $A$ geh"oren,
da sonst wegen $f \cdot {r \over f}=r$ auch
$r \in {\bf p}$ w"are. \hfill $\Box$
\par\bigskip
\noindent
Eine besondere Situation liegt vor, wenn in der Darstellung
$1=f_1q_1+...+f_nq_n$ mit $f_i \in {\bf a}$ und
$q_i \in \Gamma(D({\bf a}),{\cal O}_X)$ alle Summanden $f_iq_i$ zu $A$
geh"oren, also jeweils Z"ahler zu $q_i$ sind.
Dieser extreme Fall von Affinit"at liegt etwa bei $D=D(f)$ vor, dort ist
ja $1=f \cdot  1/f $.
Ist die offene Menge durch $U={\rm Def}\, (q)$ gegeben, und erzeugen Z"ahler-
und Nennerideal von $q$ zusammen das Einheitsideal, was genau dann
der Fall ist, wenn die meromorphe Funktion $q$ global definiert ist,
so liegt ebenfalls diese
Situation vor. Sei n"amlich $h+g=1$ mit $h \in N(q)={\bf a}$ und
$g \in Z(q)$ mit $g=fq, f \in N(q)$. Dann ist nat"urlich
$1=h+fq$ eine Darstellung der geforderten Art.
Der folgende Satz charakterisiert, was es hei"st, wenn Z"ahler von
rationalen Funktionen schon die Eins erzeugen.
\begin{lem}
Sei $A$ ein noetherscher Integrit"atsbereich, $\bf a$ ein Ideal und
$U=D({\bf a})$ eine
offene Menge mit $B$ als globalem Schnittring. Dann sind "aquivalent.
\par\smallskip\noindent
{\rm (1)} Es gibt eine Darstellung $1=f_1q_1+...+f_nq_n$ mit
$f_i \in {\bf a},\, q_i \in B,\,  f_iq_i \in A \, .$
\par\smallskip\noindent
{\rm (2)} In jeder Lokalisierung $A_{\bf p}$ wird ${\bf a}_{\bf p}$
als Radikal durch eine Funktion beschrieben.
\par\smallskip\noindent
{\rm (3)}
Es gibt eine offene "Uberdeckung ${\rm Spek}\, A= \bigcup_{i =1,...,n}D(g_i)$
so, dass $U \cap D(g_i)$ geometrisch durch eine Funktion $f_i \in A_{g_i}$
beschrieben wird.
\end{lem}
{\it Beweis}.
Von (1) nach (2). Die Vorraussetzungen bleiben beim "Ubergang zu $A_{\bf p}$
erhalten, sei also $A$ lokal. Da Elemente aus $A$ die Eins darstellen, muss
bereits ein Summand $f_iq_i$ eine Einheit sein, und zwar ohne Einschr"ankung
gleich eins. Dann ist einerseits $D(f) \subseteq U$ und andererseits ist
$ 1/f $ auf $U$ definiert, also gilt die Gleichheit.
\par\smallskip\noindent
Von (2) nach (3). Ein (endliches) Radikalerzeugendensystem $(f_j)$
von $({\bf a})$
aus $A$ wird im Punkt $A_{\bf p}$ dargestellt durch $f_j=q_jf $, mit
$f,q_j \in A_{\bf p}$. Diese Funktionen sind dann auch in einer Umgebung
(ihren Definitionsbereichen) von $\bf p$
definiert und erf"ullen dort den gleichen Zweck,
so dass darauf ebenfalls $V({\bf a})$ durch eine Funktion beschrieben werden
kann
\par\smallskip\noindent
Von (3) nach (1). Seien $f_i = a_i/ g_i^{n_i} $ mit $a_i \in {\bf a}$.
In ${\rm Spek}\, A_{g_i}$ ist $U=D(f_i)$,
daher ist $1/f_i$ auf $U \cap D(g_i) =U_{g_i}$ definiert.
Daher gibt es Darstellungen $1/f_i = q_i / g_i^{m_i}$ mit
$q_i \in B$.
Es ist also $ a_i/ g_i^{n_i} \cdot q_i/ g_i^{m_i}=f_i \cdot 1/f_i=1$, also
$a_iq_i = g_i^{n_i+m_i} \in A$.
Da die $g_i$ das Einheitsideal erzeugen, folgt die Behauptung. \hfill $\Box$
\par\bigskip
\noindent
In 1.8.4 werden wir allgemeiner sehen, dass generell
die Affinit"at einer offenen Teilmenge eine lokale Eigenschaft ist.
F"ur integre noethersche Kurven ist trivialerweise die Bedingung (2) erf"ullt,
so dass wir folgendes Resultat notieren k"onnen, dass im Abschnitt 1.10
"uber eindimensionale Schemata mit varierten Voraussetzung und
anderen Beweisen wieder aufgegriffen wird.
\begin{satz}
Sei $A$ ein noetherscher, integrer, eindimensionaler Ring und
$U \subseteq {\rm Spek}\, A$ eine
offene Teilmenge. Dann wird $U$ lokal durch eine Funktion beschrieben
und ist affin. \hfill $\Box$
\end{satz} 
\begin{satz}
Sei $A$ ein noetherscher, lokal faktorieller Ring. Dann besitzt
jede Hyperfl"ache affines Komplement.
Dies gilt insbesondere f"ur regul"are Ringe.
\end{satz}
{\it Beweis}.
Ein lokaler regul"arer Ring ist bekanntlich faktoriell,
siehe \cite{EGAIV}, 21.11.1 oder \cite{eisenbud}, Theorem 19.19,
lokal liegt also ein faktorieller Ring vor. In einem faktoriellen
Ring aber wird eine Hyperfl"ache durch eine
Funktion beschrieben. Ist $\bf p$ ein Primideal der H"ohe eins
und $0 \neq f \in {\bf p}$, so gibt es nach Voraussetzung eine
Primfaktorzerlegung von $f$, und dabei muss ein Primfaktor $p$ zu ${\bf p}$
geh"oren. Dann ist aber schon ${\bf p} = (p)$. \hfill $\Box$
\par\bigskip
\noindent 
{\bf Bemerkung}
Einen normalen Integrit"atsbereich $A$,
bei dem die Divisorenklassengruppe eine Torsionsgruppe
ist, nennt man {\it fastfaktoriell} (oder ${\bf Q}-${\it fak\-toriell}),
siehe \cite{storch}.
In einem fastfaktoriellen Ring wird jede Hyperfl"ache
geometrisch durch eine Funktion beschrieben.
Der Satz 1.4.7 gilt genau so, wenn man lediglich lokal fastfaktoriell
voraussetzt.
\par\bigskip
\noindent
Wir fragen uns, welche Ringe als Schnittringe in einem integren noetherschen
Ring in Frage kommen k"onnen. 
\begin{satz}
Sei $A$ ein noetherscher Integrit"atsbereich. Dann sind f"ur eine
$A-$Algebra $B$ folgende Aussagen "aquivalent.
\par\smallskip\noindent
{\rm (1)} $B$ ist der Schnittring einer offenen
affinen Teilmenge von ${\rm Spek}\, A$.
\par\smallskip\noindent
{\rm (2)} $B$ ist vom endlichen Typ, birational und flach "uber $A$.
\end{satz}
{\it Beweis}.
Von (1) nach (2). Dass $B$ den gleichen Quotientenk"orper hat ist klar,
dass es vom endlichen Typ ist wurde in 1.4.3 gezeigt.
Eine offene Einbettung ist immer flach, also ist
$ {\rm Spek}\, B=U \hookrightarrow {\rm Spek}\, A$ flach, was die Flachheit von 
$R:A \longrightarrow B$ bedeutet, da die Flachheit eine lokale Eigenschaft
ist.
Sei umgekehrt (2) erf"ullt.
Mit $A$ ist auch $B$ noethersch, und da $B$ flach und vom endlichen Typ ist,
ist nach \cite{EGAIV}, 2.4.6,
${\rm Spek}\, B \longrightarrow {\rm Spek}\, A$ offen,
sei $U$ das offene Bild davon.
Sei $D(f) \subseteq U$. Es ist dann auch $A_f \longrightarrow B_f$ flach
und nach Voraussetzung surjektiv,
also treuflach.
Eine treuflache, birationale Erweiterung $A \subseteq B \subseteq Q(A)$ ist
aber schon die Identit"at, da aus $a/b=q \in B$ mit $a,b \in A$ folgt,
dass $a$ von $b$
in $B$ geteilt wird, was dann aber schon in $A$ gilt, da eine
treuflache Abbildung rein ist, siehe \cite{scheja}, \S 88, Aufg. 29.
Damit ist ${\rm Spek}\, B \longrightarrow U$ lokal und damit auch global eine
Isomorphie. \hfill $\Box$
\par
\bigskip
\noindent
\subsection{Nulldimensionale Schemata}
\par
\bigskip
\noindent
Unter den nulldimensionalen Schemata lassen sich die affinen
topologisch charakterisieren. Ein Schema hei"st {\it quasisepariert},
wenn die Diagonale ("uber $\bf Z$) quasikompakt ist, siehe \cite{EGAI}, 6.1.3.
In diesem Fall ist der Durchschnitt zweier affiner Teilmengen wieder
quasikompakt.
\begin{satz}
Sei $X$ ein nulldimensionales Schema. Dann gelten folgende Aussagen.
\par\smallskip\noindent
{\rm (1)} Ist $X$ quasisepariert, so ist die Topologie von $X$ Hausdorffsch und $X$
ist ein reich beringtes Schema.
\par\smallskip\noindent
{\rm (2)} Es sind folgende Aussagen "aquivalent.
\par\smallskip\noindent
{\rm (i)} $X$ ist ein affines Schema.
\par\smallskip\noindent
{\rm (ii)} $X$ ist quasikompakt und quasisepariert.
\par\smallskip\noindent
{\rm (iii)} $X$ ist kompakt.
\end{satz}
{\it Beweis}.
Wir zeigen zuerst, dass ein quasisepariertes nulldimensionales Schema 
Hausdorffsch ist. Ist $X$ affin, so seien $\bf p$ und $\bf q$
verschiedene Primideale der H"ohe null. Dann gibt es $g \not\in {\bf p}$
und $g \in {\bf q}$. $(A_{\bf q})_{\rm red}$ ist dann ein K"orper
und somit ist $g^n=0$ in $A_{\bf q}$.
Dann gibt es ein $f \not\in {\bf q}$ mit $g^nf=0$, und somit sind $D(g)$ und
$D(f)$ disjunkte Umgebungen von $\bf p$ und $\bf q$.
\par\smallskip\noindent
Seien im allgemeinen Fall $x,y \in X,\, x \neq y $ vorgegeben.
Besitzen $x,y$ eine gemeinsame affine Umgebung, so finden sich darin 
trennende offene Umgebungen.
Seien jetzt $x \in U$ und $y \in V$ affine Umgebungen mit dem
Schnitt $W=U \cap V$. Wir k"onnen dann $x \not\in W$ annehmen.
Da $W$ nach Voraussetzung quasikompakt ist,
finden wir innerhalb von $U$ eine offene Umgebung
von $x$, die $W$ nicht trifft und damit "uberhaupt zu $V \ni y$ disjunkt ist.
Damit ist Hausdorffsch bewiesen, reich beringt wird sich im Verlauf 
des Beweises zum zweiten Teil ergeben.
\par\smallskip\noindent
Von (i) nach (ii) ist klar, ebenso von (ii) nach (iii) nach dem ersten Teil.
Wir zeigen zun"achst, dass aus Hausdorffsch und nulldimensional reich beringt
folgt. Sei $U \subseteq X$ eine affine, offene Teilmenge. $U$ ist wegen
(i) $\longrightarrow $ (iii) kompakt und daher auch abgeschlossen, dann ist
auch das Komplement offen und es l"asst sich durch $f \vert U =1$ und
$f=0$ au"serhalb von $U$ eine globale Funktion auf $X$ definieren mit
$U=X_f$. Damit ist $X$ reich beringt und auch Teil eins ist vollst"andig
bewiesen.
\par\smallskip\noindent
Sei jetzt (iii) erf"ullt und $X$ kompakt.
$X$ wird dann "uberdeckt durch endlich viele affine Mengen der Form $X_f$,
wobei $f$ auf $X_f$ den Wert eins und au"serhalb den Wert null habe.
Das ebenfalls offene und abgeschlossene Komplement ist dann gleich $X_{1-f}$.
Dann ist $X_g \cap X_{1-f}=X_{g(1-f)}$ affin und so kann man sukzessive
annehmen, dass $X$ "uberdeckt wird mit solchen disjunkten affinen Mengen.
Dann ist aber $X$ die endliche direkte Summe affiner Schemata und daher selbst
affin. \hfill $\Box$
\begin{kor}
Sei $X$ ein nulldimensionales, topologisch noethersches Schema.
Dann besteht $X$ aus endlich vielen Punkten, die Topologie
ist diskret und $X$ ist ein affines Schema.
\end{kor}
{\it Beweis}.
Ein noethersches Schema ist zun"achst immer quasisepariert, so dass aus der
Dimensionseigenschaft Hausdorffsch folgt.
Damit ist jede offene Menge kompakt und damit auch abgeschlossen, es liegt
also die diskrete Topologie 
vor, und daher muss $X$ endlich und affin sein. \hfill $\Box$
\par
\bigskip
\noindent
\subsection{Affine Morphismen}
\par
\bigskip
\noindent
{\bf Definition}
Ein Schemamorphismus $f:X \longrightarrow Y$ hei"st {\it affin}, wenn es eine
offene "Uberdeckung von $Y$ mit affinen Mengen $Y=\bigcup_{i \in I} Y_i$ gibt
so, dass die Urbilder $f^{-1}(Y_i)$ wieder affin sind.
\par\bigskip
\noindent
Das ist ein Standardbegriff aus der algebraischen Geometrie, vergl.
\cite{EGAI}, 9.1, und \cite{EGAII}, 1.6,
den wir hier aber wegen seiner Bedeutung
von Grund auf darstellen.
\begin{lem}
Ist $f:X \longrightarrow Y$ ein affiner Morphismus, so ist f"ur eine
offene Menge $U \subseteq Y$ auch der Morphismus 
$f^{-1}(U) \longrightarrow U $ affin.
\end{lem}
{\it Beweis}.
Sei $Y_i,\, i \in I,$ eine affine "Uberdeckung zu $Y$ mit
affinen Urbildern. Zu $y \in U$ m"ussen wir eine affine Umgebung in $Y$
finden, deren Urbild wieder affin ist. Sei $y \in Y_i$ f"ur ein $i \in I$,
also $y \in U\cap Y_i \subseteq Y_i$.
Dann k"onnen wir annehmen, dass sich alles
in $Y_i$ abspielt und also $X$ und $Y$ affin sind.
Eine offene Menge in einem affinen Schema wird aber bereits durch Basismengen
$D(f)$ "uberdeckt, und die haben unter einem Morphismus affiner Schemata
wieder diese Gestalt und sind affin. \hfill $\Box$
\begin{lem}
Sei $f:X \longrightarrow Y$ ein affiner Morphismus "uber einem
Basisschema $S$ mit Strukturmorphismen $g:X \longrightarrow S$ und
$h:Y \longrightarrow S$, und sei $p:T \longrightarrow S$
ein beliebiger Basiswechsel.
Dann ist auch der dadurch gewonnene Morphismus
$f_{(T)}: X \times_ST \longrightarrow Y \times _S T $ affin.
\end{lem}
{\it Beweis}.
Sei $ V$ eine affine Teilmenge von $S$, $W$ eine affine
Teilmenge von $T$ mit $W \subseteq p^{-1}(V)$ und $U$ eine affine Teilmenge
von $X$ mit entsprechend $U \subseteq g^{-1}(V)$.
Dann ist $U \times_V W$ eine offene affine Teilmenge von $Y_{(T)}=Y\times _S T$
und solche Mengen bilden eine affine "Uberdeckung von $Y_{(T)}$.
Nach dem vorausgehenden Lemma gibt es zu $U$ in $Y$ eine
affine "Uberdeckung $U_i,\, i \in I$, deren Urbilder alle affin sind,
und somit kann man
$Y_{(T)}$ auch mit solchen Mengen "uberdecken. Es ist
dann  $f_{(T)}^{-1}(U_i \times_V W) =f^{-1}(U_i) \times_V W$ wieder
affin. \hfill $\Box$
\begin{satz}
Ist $f: X \longrightarrow Y$ affin, so ist das Urbild jeder offenen affinen
Teilmenge wieder affin.
\end{satz}
{\it Beweis}.
Nach dem Lemma 1.6.1 k"onnen wir $Y={\rm Spek}\, A$
als affin annehmen und
haben dann die Affint"at von $X$ zu zeigen.
Sei $U$ eine affine offene Teilmenge von $Y$ mit affinem Urbild.
Dann ist auch das Urbild von Mengen der Form $D(f) \subseteq U$ affin in $X$,
und daher k"onnen wir annehmen, dass eine affine "Uberdeckung von
${\rm Spek}\, A$ mit Basismengen $D(f_i),\, i \in I,\, f_i \in A,$ vorliegt,
deren Urbilder wieder affin sind. 
Mit dem Ringhomomorphismus der globalen Schnittringe fassen wir die
$f_i \in \Gamma(X,{\cal O}_X)=:B$ auf, und die Urbilder der $D(f_i)$ sind dann
einfach die $X_{f_i}$, die nach Konstruktion affin sind.
Da die $D(f_i)$ in ${\rm Spek}\, A$ eine "Uberdeckung bilden, erzeugen die
$f_i$ in $A$ die Eins, was dann auch in $B$ gilt.
Damit sind die Bedingungen des Affinit"atskriterium 1.3.9 erf"ullt und
es folgt die Affinit"at von $X$. \hfill $\Box$
\par\bigskip
\noindent
Der folgende Satz fasst grundlegende Eigenschaften affiner Morphismen
zusammen. 
\begin{satz}
{\rm (1)} Eine abgeschlossene Einbettung ist ein affiner Morphismus.
\par\smallskip\noindent
{\rm (2)} Die Komposition zweier affiner Morphismen ist wieder affin.
\par\smallskip\noindent
{\rm (3)} Ist $f:X \longrightarrow Y$ ein affiner $S-$Morphismus, so ist auch
$f_{(T)}:X_{(T)} \longrightarrow Y_{(T)} $ affin zu jedem Basiswechsel
$T \longrightarrow S$.
\par\smallskip\noindent
{\rm (4)} Sind $f:X \longrightarrow Y $ und
$f': X' \longrightarrow Y'$ affine Morphismen "uber
$S$, so ist auch der Produktmorphismus $f \times f' :
X \times _S X' \longrightarrow Y \times _S Y'$ affin.
\par\smallskip\noindent
{\rm (5)} Sind $f:X \longrightarrow Y$ und $g:Y \longrightarrow Z$ Morphismen mit
$g$ separiert und $g \circ f$ affin, so ist bereits $f$ affin.
\par\smallskip\noindent
{\rm (6)} Ist $f$ affin, so gilt das auch f"ur $f_{\rm red}$ .
\end{satz}
{\it Beweis}.
Bekanntlich folgen aus den drei ersten Aussagen die drei letzten,
siehe \cite{EGAI}, 5.2.7. 
(1) ist die Aussage, dass ein abgeschlossenes Unterschema eines
affinen Schemas das Spektrum eines Restklassenringes und damit
affin ist, siehe etwa \cite{haralg}, Cor. II.5.10.
(2) ist klar nach 1.6.3, (3) ist der Inhalt des
Lemmas 1.6.2. \hfill $\Box$
\par\bigskip
\noindent
{\bf Bemerkung}
Oft werden wir die Aussage (3) so anwenden,
dass auf einen Morphismus
$U \longrightarrow X $ ein affiner Basiswechsel $X' \longrightarrow X$
angewendet wird, was nach (3) zu einem affinen Morphismus
$U' =U \times _X X' \longrightarrow U$ f"uhrt.
Die Basis $S$ in obiger Notation ist dabei insbesondere gleich $X$.
In dieser Situation ist dann die Affinit"at von $U$ eine hinreichende
Bedingung f"ur die Affinit"at von $U'$ und somit ist die Affinit"at von
$U'$ eine notwendige Bedingung f"ur die Affinit"at von $U$.
Dabei wird oftmals $X={\rm Spek}\, A$ ein affines Schema sein und der
Basiswechsel durch einen Ringwechsel, also einen Ringhomomorphismus
$A \longrightarrow A'$ gegeben sein. Hierbei ist an Lokalisierung,
Reduktion modulo eines Ideals, Normalisierung, Komplettierung,
Grundk"orpererweiterung etc. zu denken. Desweiteren wird $U \longrightarrow X$
einfach eine offene Einbettung sein, mit $U=D({\bf a})$. Dann ist
$U'=D({\bf a}A')$ einfach das Urbild von $U$ unter der affinen
Abbildung ${\rm Spek}\, A' \longrightarrow {\rm Spek}\, A$. 
\par\smallskip\noindent
Eine Fragestellung, die sich durch diese Arbeit durchziehen wird, ist,
inwiefern die erw"ahnte notwendige Bedingung, dass $U'$ affin ist,
bereits hinreichend f"ur die Affinit"at von $U$ ist, und zwar bezogen auf
einen bestimmten Typ von Ringwechsel oder einer Familie von Ringwechseln
$A \longrightarrow A_i,\, i \in I$.
Im "ubern"achsten Abschnitt wird die Situation, wo die $A_i$ alle
Lokalisierungen von $A$ durchl"auft, untersucht.
\begin{satz}
Sei $f:X \longrightarrow X'$ ein Schemamorphismus, der topologisch eine
Hom"oomorphie ist. Dann ist $f$ eine affine Abbildung.
\end{satz}
{\it Beweis}.
Sei $X'$ affin und sei $x \in U \subseteq X$ eine offene Umgebung.
Die gleiche (=hom"oomorphe)
Situation in $X'$ betrachtet gew"ahrleistet eine Funktion
$f \in \Gamma(X',{\cal O}_{X'})$ mit $ x' \in D(f) \subseteq U' \subseteq X'$.
Diese Inklusionen gelten dann auch f"ur das Urbild $X_f$,
und damit ist $X$ reich
beringt. Ist dabei $U$ eine affine Umgebung, so sind die darin enthaltenen
Mengen der Form $X_f$ affin.
Insbesondere finden wir eine endliche "Uberdeckung aus affinen Mengen der
Form $X_f$, wobei die Funktionen aus $\Gamma(X',{\cal O}_{X'})$ kommen.
Dort erzeugen diese Funktionen aber das Einheitsideal, was dann auch
in $\Gamma(X,{\cal O}_X)$ gilt. \hfill $\Box$
\par\bigskip
\noindent
{\bf Bemerkung}
Die Aussage des Satzes gilt umgekehrt nicht, es kann $X$ affin sein,
ohne dass das hom"oomorphe Bild $X'$ affin ist, siehe hierzu den
Abschnitt 1.14.
Die Aussage gilt auch nicht, wenn $f$ bijektiv ist, wie das folgende Beispiel
zeigt.
\par\bigskip
\noindent
{\bf Beispiel}
Betrachte eine nicht normale affine
Fl"ache $Y$, wo in der Normalisierung $\tilde{Y}$
"uber einem einzigen Punkt $P\in Y$ zwei
Punkte $Q$ und $R$ liegen und sonst au"serhalb von $P$ eine Isomorphie vorliegt.
Sei $X=\tilde{Y} - \{ Q \}$. Dann ist die Abbildung bijektiv,
das punktierte Schema $X$ ist aber nicht affin. Die Abbildung kann nach dem
vorangegangenen Satz nicht hom"oomorph sein, beispielsweise
sind Kurven durch $Q$, die $R$ nicht treffen, in $X$ abgeschlossen, ihr Bild
in $Y$ aber nicht.
\par\bigskip
\noindent
Im vorangegangenen Beispiel liegt nat"urlich nicht die kanonische Abbildung vor.
In 1.8 werden wir eine bijektive kanonische Abbildung
konstruieren, die keine Isomorphie ist,
bei der aber $X$ nicht quasikompakt ist.
Ob bei $X$ 
noethersch und separiert aus einer bijektiven kanonischen Abbildung bereits
$X$ affin folgt, ist nicht klar.
Ist $X$ separiert und von endlichem Typ "uber einem affinen Basisschema,
so folgt aus der Bijektivit"at der kanonischen Abbildung, dass $X$ affin ist,
siehe die S"atze 1.3.5, 1.3.8 und 1.3.9.
\par
\bigskip
\noindent
\subsection{Lokalisierung der kanonischen Abbildung}
\par
\bigskip
\noindent
Sei $X \longrightarrow {\rm Spek}\, A$ eine Affinisierungsabbildung und
$A \longrightarrow B$ ein Ringwechsel. Dann ist die zugeh"orige
Abbildung $X_{(B)} \longrightarrow {\rm Spek}\, B$ in aller Regel nicht mehr
die Affinisierungsabbildung, ganz einfach deshalb, weil der globale
Schnitt\-ring $\Gamma(X_{(B)},{\cal O}_{X_{(B)}} )$ nicht gleich $B$ sein muss.
Ein extremes, aber typisches Beispiel ist die
Affinisierungsabbildung der punktierten Ebene
${\bf A}_K^2 -\{ P \} \longrightarrow {\bf A}_K^2$. Durch den affinen
Basiswechsel ${\rm Spek}\, K =P \longrightarrow {\bf A}_K^2 $ entsteht 
daraus die Abbildung $\emptyset \longrightarrow P$, die nat"urlich keine
Affinisierungsabbildung ist. Unter diesem Basiswechsel wird der globale
Schnittring der Nullring, w"ahrend $K[X,Y] \otimes_{K[X,Y]}K =K$ ist.
\begin{satz}
Sei $X \longrightarrow {\rm Spek}\,A $ ein Schemamorphismus mit
$X$ quasikompakt und separiert.
$A \longrightarrow A'$ sei ein
flacher Basiswechsel und  $X' =X \times_A {\rm Spek}\, A'$.
Dann gilt f"ur den globalen Schnittring
$\Gamma(X',{\cal O}_{X'}) = \Gamma(X,{\cal O}_X) \otimes_A A'$.
\end{satz}
{\it Beweis}.
Ist $X$ selbst affin, so ist das Resultat trivial.
Sei $X= \bigcup_{i \in I}U_i$ eine endliche "Uberdeckung mit affinen
Teilmengen $U_i$. Der globale Schnittring ist dann durch die exakte Sequenz
$$0 \longrightarrow \Gamma(X,{\cal O}_X) \longrightarrow 
\prod_{i \in I} \Gamma(U_i,{\cal O}_X)  \longrightarrow
\prod_{i,j \in I} \Gamma (U_i \cap U_j,{\cal O}_X ) $$
definiert, wobei die letzte Abbildung komponentenweise die Differenzbildung ist.
Wie tensorieren diese exakte Sequenz mit der flachen $A-$Algebra $A'$ und
erhalten die exakte tensorierte Sequenz
$$0 \rightarrow \Gamma(X,{\cal O}_X)\otimes _A A' \longrightarrow 
(\prod_{i \in I} \Gamma(U_i,{\cal O}_X))\otimes_AA'  \longrightarrow
(\prod_{i,j \in I} \Gamma (U_i \cap U_j,{\cal O}_X))\otimes_AA'  .$$
Da hinten endliche Produkte stehen, k"onnen wir das Tensorieren reinziehen,
und nach der Vorbemerkung ist f"ur affines $U$ (die $U_i \cap U_j$ sind
wieder affin, da wir separiert vorausgesetzt haben)
$\Gamma(U,{\cal O}_X)\otimes_A A' =\Gamma(U',{\cal O}_{X'})$, wobei
$U'$ das Urbild von $U$ bezeichne; $U'$ ist wieder affin, da ein affiner
Basiswechsel vorliegt, daher liegt wieder eine exakte Sequenz vor, die den
Schnittring beschreibt, also ist $\Gamma(X,{\cal O}_X) \otimes_AA'$ gleich
dem Schnittring $\Gamma(X',{\cal O}_{X'})$. \hfill $\Box$
\begin{kor} 
Sei $X$ quasikompakt und separiert mit kanonischer Abbildung 
$X \longrightarrow {\rm Spek}\, A$. Vorgegeben sei ein flacher Ringwechsel
$A \longrightarrow A'$. Dann ist
$X \times_A {\rm Spek}\, A' \longrightarrow {\rm Spek}\, A'$ ebenfalls die
kanonische Abbildung.
Insbesondere ist f"ur ein Primideal ${\bf p} \in {\rm Spek}\, A$ die
Abbildung
$X \times_A {\rm Spek}\, A_{\bf p}=X_{\bf p}
\longrightarrow {\rm Spek}\, A_{\bf p}$
die kanonische Abbildung. \hfill $\Box$
\end{kor} 
Aus dem Lokalisierungsprinzip f"ur die kanonische Abbildung ergibt sich
die folgende Aussage.
\begin{satz}
Sei $X$ ein noethersches separiertes Schema mit $A=\Gamma(X,{\cal O}_X)$.
Dann werden unter der
kanonischen Abbildung $f: X \longrightarrow {\rm Spek}\, A$
alle Primideale der H"ohe null und eins von ${\rm Spek}\,A$ erreicht.
\end{satz}
{\it Beweis}.
Sei $\bf p$ ein Primideal der H"ohe null.
Wir lokalisieren daran und erhalten die entsprechende Situation wo
jetzt ${\rm Spek}\, A$ einpunktig ist.
Da eine Affinisierungsabbildung vorliegt,
kann $X$ nicht leer sein und somit gibt es Urbildpunkte.
\par\smallskip\noindent
Sei jetzt $\bf p$ ein Primideal der H"ohe eins, wir lokalisieren wieder
und haben dann die Situation einer kanonischen Abbildung eines noetherschen
separierten Schemas mit einem lokalen eindimensionalen Schnittring,
der allerdings nicht noethersch sein muss.
Sei $U=D({\bf m})$ die offene Menge aller Primideale der H"ohe null, die
alle erreicht werden, und nehmen wir an, das sei das ganze Bild der
kanonischen Abbildung. Als Bild eines noetherschen Schemas ist dann $U$
ein zumindest topologisch noethersches Schema, das nulldimensional
ist. Nach Kor. 1.5.2  ist dann $U$ insbesondere affin
und wegen Satz 1.3.3 w"are dann $U$ bereits die Affinisierung.   \hfill $\Box$
\par\bigskip
\noindent
Aufgrund von 1.3.8 kann man den globalen Hilbertschen Nullstellensatz
dadurch charakterisieren, dass unter der kanonischen Abbildung alle maximale
Ideale getroffen werden. Daher ergibt sich aus 1.7.3 sofort
das folgende Korollar.
\begin{kor}
Besitzt ein noethersches separiertes Schema einen globalen Schnittring der
Dimension $\leq 1$, so gilt in $X$ der globale Hilbertsche
Nullstellensatz. \hfill $\Box$
\end{kor}
Sei nun $X$ ein integres, noethersches, separiertes Schema. Dann ist auch
der globale Schnittring $\Gamma(X,{\cal O}_X)=A$ integer.
Der "Ubergang zum Quotientenk"orper $K=Q(A)$ entspricht der
Lokalisierung am generischen Punkt von ${\rm Spek}\, A$,
und dar"uber liegt dann als Faser der kanonischen Abbildung
das Schema $X \times_{{\rm Spek}\, A} {\rm Spek}\, Q(A)= X_{K}$. Das ist dann
ein Schema mit dem K"orper $Q(A)$ als globalem Schnittring.
Der n"achste Satz bringt verschiedene "aquivalente Formulierungen daf"ur, dass
dieses Schema affin ist.
\begin{satz}
Sei $X$ ein separiertes, integres, noethersches Schema mit globalem
Schnittring $A$, sei $Q(A)=K$ und $\eta$ der generische Punkt von $X$.
Dann sind folgende Aussagen "aquivalent.
\par\smallskip\noindent
{\rm (1)}
Es ist $X \times_{{\rm Spek}\, A} {\rm Spek}\, Q(A)$ affin, die kanonische
Abbildung ist also generisch affin.
\par\smallskip\noindent
{\rm (2)}
Die kanonische Abbildung ist im generischen Punkt von ${\rm Spek}\, A$ eine
Isomorphie.
\par\smallskip\noindent
{\rm (3)} Es ist $Q(A) = K(X)$.
\par\smallskip\noindent
{\rm (4)} Es gibt eine globale Funktion $0 \neq f \in A$ mit $X_f$ affin.
\par\smallskip\noindent
{\rm (5)} F"ur alle Primideale ${\bf p}$ von $A$ der H"ohe $\leq 1$ ist $X_{\bf p}$
affin mit $X_{\bf p} = {\rm Spek}\, A_{\bf p}$.
\end{satz}
{\it Beweis}.
Unter den ersten vier Aussagen macht man einen Ringschluss, dabei
ist nur der Schluss von (3) nach (4) nicht trivial. 
Wie in 1.4 erw"ahnt bilden in einem noetherschen, integren, separierten Schema
die gemeinsamen Definitionsbereiche zu endlich vielen rationalen
Funktionen eine Basis der Topologie, daher gibt es
$\eta \in {\rm Def}\, (q_1,...,q_n) \subseteq U$ mit $U$ affin.
Die $q_i \in K(X)=Q(A)$ haben einen gemeinsamen Nenner $0 \neq f \in A$,
und somit ist $\eta \in X_f \subseteq {\rm Def}\, (q_1,...,q_n)$
und $X_f$ ist affin.
\par\smallskip
\noindent
Wir zeigen (5) und betrachten dazu die Sache lokalisiert an einem Primideal
der H"ohe eins, dann ist $X$ ein integres, noethersches, separiertes Schema
mit einem eindimensionalen, integren, lokalen Schnittring $A$, und nach
Voraussetzung ist die kanonische Abbildung generisch isomorph, d.h.
"uber dem generischen Punkt von ${\rm Spek}\, A$
liegt nur der generische Punkt von $X$.
Daraus folgt, dass der generische Punkt von $X$ offen ist, so dass
der n"achste Satz das Resultat liefert.  \hfill $\Box$
\begin{satz}
Sei $X$ ein integres, noethersches, separiertes Schema.
Dann sind folgende Aussagen "aquivalent.
\par\smallskip\noindent
{\rm (1)}
Der generische Punkt von $X$ ist offen.
\par\smallskip\noindent
{\rm (2)} $X$ besteht nur aus endlich vielen Punkten.
\par\smallskip\noindent
{\rm (3)}
$X$ ist das Spektrum eines semilokalen Ringes der Dimension $\leq 1$.
\end{satz}
{\it Beweis}.
Von (1) nach (2).
Die Voraussetzungen  gelten dann auch auf einem affinen Teilst"uck
${\rm Spek}\, B$.
Dass der generische Punkt offen ist, hei"st, dass es eine Funktion
$0 \neq f \in B$ gibt mit $\{ \eta \} = D(f) $ (Es ist also $A_f$ bereits
der Quotientenk"orper), und $V(f)$ ist dann alles au"ser dem generischen Punkt.
Wegen noethersch besteht $V(f)$ aus endlich vielen Komponenten
$V({\bf p}_1),...,V({\bf p}_n)$ mit Primidealen der H"ohe eins.
Das sind dann alle Primideale der H"ohe eins. Das sind dann
aber "uberhaupt schon alle von null verschiedenen Primideale,
da es in einem integren noetherschen Ring der Dimension $\geq 2$ schon
unendlich viele Primideale der H"ohe eins gibt,
siehe \cite{bruske}, Folgerung 6.21.
Damit gibt es in $X$ nur enlich viele Punkte und zugleich ist gezeigt,
dass in einem solchen Schema der generische Punkt offen sein muss.
\par\smallskip\noindent
Sei jetzt (2) und (1) erf"ullt.
Wir wissen also, dass $X$ "uberdeckt wird von affinen Spektra zu
semilokalen, integren, noetherschen, eindimensionalen Ringen, und dann
auch schon von endlich vielen Spektra zu solchen lokalen Ringen, verklebt
am offenen generischen Punkt.
Seien $B_i,\, i \in I,$ die beteiligten integren, eindimensionalen, lokalen
Ringe und $A= \bigcap_{i \in I}B_i$ der globale Schnittring.
\par\smallskip\noindent
Wir nehmen nun zus"atzlich an, dass das Schema normal ist, wir werden in
1.14.2 sehen, dass aus der Affinit"at der Normalisierung
die Affinit"at selbst folgt. Zugleich "andern sich die anderen Voraussetzungen
nicht, da die Normalisierung noetherscher eindimensionaler Integrit"atsbereiche
nach dem Satz von Krull-Akizuki noethersch bleibt, siehe \cite{nagloc},
Theorem 33.2.
Die beteiligten lokalen Ringe $B_i$ sind dann diskrete Bewertungsringe.
\par\smallskip\noindent
Wir beweisen durch Induktion "uber die Anzahl der beteiligten Punkte,
dass $X$ affin ist, wobei der Fall,
dass es nur einen abgeschlossenen Punkt gibt, trivial ist.
\par\smallskip\noindent
Seien $x_1,...,x_n$ die abgeschlossenen Punkte von $X$,
sei $k \in I=\{1,...,n\}$.
Es ist dann $A_k = \bigcap_{i \neq k} B_i$ der globale Schnittring des
affinen Schemas $X-\{ x_k \}$.
Es sei $f \in A_k$ eine Funktion $\neq 0$, die an allen abgeschlossenen Punkten von
${\rm Spek}\, A_k= X-\{ x_k\}$ eine Nullstelle hat
und $g \in A_k$ eine Funktion, die an der Stelle $x_k$ eine negative
Ordnung hat. Eine solche Funktion muss es geben, da andernfalls
$A_k \subseteq B_k$ w"are, woraus $B_i =B_k$ f"ur $i \neq k$ folgen w"urde,
was bei separiert nicht sein kann, siehe hierzu \cite{EGAI}, Prop. 8.5.5.
Damit findet man dann eine Funktion $h=fg^s$, die auf $x_i,\, i \neq k,$ eine
Nullstelle hat und auf $x_k$ eine Polstelle.
Dann besitzt $1+h$ an $x_k$ ebenfalls eine Polstelle und ist in den anderen
Punkten eine Einheit. Damit ist $1/(1+h) \in A$ eine globale Funktion,
die nur an $x_k$ eine Nullstelle hat. 
\par\smallskip\noindent
Damit findet man auch umgekehrt f"ur $k \in I$ eine globale Funktion $f_k$, die
nur an der Stelle $k$ eine Einheit ist und sonst "uberall eine Nullstelle
hat. Damit ist $X_{f_k}=\{ x_k, \eta \} ={\rm Spek}\, B_k$ affin.
\par\smallskip\noindent
Wir m"ussen noch zeigen, dass diese $f_k$ in $A$ das Einheitsideal erzeugen,
dazu betrachten wir $f=f_1 +...+f_n$. 
$f$ ist in jedem diskreten Bewertungsring eine Einheit, da jeweils bis auf
einen Index alle Summanden im maximalen Ideal liegen und f"ur einen
Index eine Einheit vorliegt. Damit ist $f$ in $A$ eine Einheit.
\par\smallskip\noindent
Von (3) nach (1) ist klar. \hfill $\Box$
\par\bigskip
\noindent
In \cite{goodlandman}, Def 5.5,
wird eine Variet"at {\it pseudoaffin} genannt,
wenn sie die Bedingung (3) in 1.7.5 erf"ullt.
Mit gleichem Recht kann man solche Schemata {\it generisch affin}
nennen.
Quasiaffine Schemata sind generisch affin. Die affine Ebene,
aufgeblasen an einem Punkt, ist generisch affin (und semiaffin), aber nicht
quasiaffin.
\par\bigskip
\noindent
Im generisch affinen Fall taucht jedes Primideal der H"ohe eins von
$A=\Gamma(X,{\cal O}_X)$ genau einmal als Punkt in $X$ auf,
es kann aber nat"urlich noch weitere Punkte der H"ohe eins
in $X$ geben, die nicht auf ein Primideal unten der H"ohe eins abgebildet
werden, etwa die Aufblasung in einem Punkt.
Aber auch, wenn generisch affin vorliegt und jeder Punkt der H"ohe
eins auf einen solchen Punkt abgebildet wird, wenn sich also die Primstellen
bijektiv entsprechen, folgt daraus noch nicht quasiaffin, wie die
kleinen Kontraktionen zeigen, siehe etwa das Beispiel nach 3.7.3.
\par
\bigskip
\noindent
\subsection{Affinit"at als offene und lokale Eigenschaft}
\par
\bigskip
\noindent
Zu einem beliebigen Morphismus $\varphi: X \longrightarrow Y$
gibt es die gr"o"ste
offene Menge $V \subseteq Y$, auf der der Morphismus affin ist.
Sie ist einfach die Vereinigung aller affinen Teilmengen von $Y$, deren
Urbilder affin sind. Ist $Y ={\rm Spek}\, A$ affin, so wird diese offene
Menge durch das Ideal
$${\bf c}=\{ f \in A: X_f=\varphi^{-1}(D(f))\, {\rm affin} \}$$
beschrieben.
F"ur einen Punkt $ y \in V=D({\bf c})$ ist dann nat"urlich
die lokale Faser $X_y=X \times_Y {\rm Spek}\,{\cal O}_y$ affin.
In Satz 1.7.5 wurde gezeigt, dass es im generisch affinen Fall schon eine
nicht leere affine Teilmenge gibt vom Typ $X_f$, dass also bereits in einer
offenen Umgebung des generischen Punktes die Abbildung affin ist.
Das Hauptergebnis dieses Abschnitts verallgemeinert diese Aussage.
\begin{satz}
Sei $X$ ein noethersches separiertes Schema
und $\varphi:X \longrightarrow Y$ ein Schemamorphismus.
Zu $y \in Y$ sei die lokale Faser
$X_y=X\times_Y {\rm Spek}\, {\cal O}_{Y,y}$
affin. Dann gibt es bereits eine offene affine
Umgebung $y \in V$ mit $\varphi^{-1}(V)$ affin.
\par\smallskip\noindent
Die Menge der Punkte von $Y$ mit affiner lokaler Faser ist gleich der
Vereinigung aller affinen Teilmengen mit affinem Urbild und insbesondere
offen.
\end{satz}
{\it Beweis}.
Die Menge der Punkte $y \in Y$ mit affiner lokaler Faser ist nat"urlich
abgeschlossen unter Generalisierungen, doch folgt daraus noch nicht die
Offenheit.
Ohne Einschr"ankung kann man $Y={\rm Spek}\, A$ als affin annehmen.
\par\smallskip\noindent
Sei zun"achst $X$ integer und $\varphi: X \longrightarrow Y$ die kanonische 
Abbildung, mit $y \in Y$ und $X_y$ affin. Dann ist $X_y \longrightarrow Y_y$
eine Isomorphie und $y$ hat einen eindeutig bestimmten Urbildpunkt $x \in X$
und die zugeh"orige Halmabbildung
$A_{\bf p}={\cal O}_{Y,y} \longrightarrow {\cal O}_{X,x}$ ist eine Bijektion.
Insbesondere ist die Abbildung generisch affin und der Quotientenk"orper
von $A$ ist gleich dem Funktionenk"orper von $X$.
Sei $x \in U \subseteq X$ eine affine Umgebung von $x$, dann gibt es
nach 1.4
rationale Funktionen mit $x \in {\rm Def}\, (q_1,...,q_n) \subseteq U$,
also $q_1,...,q_n \in {\cal O}_x = A_{\bf p}$. Damit gibt es einen
Hauptnenner $f \in A,\, f \not\in {\bf p}$ der $q_i$.
Es ist dann $ 1/f  \in {\cal O}_x$ und also
$x \in X_f \subseteq {\rm Def}\, (q_1,...,q_n) \subseteq U$, also ist
$X_f$ affin.
Damit haben wir in $D(f)$ eine offene Umgebung gefunden, worauf die
Abbildung affin ist.
\par\smallskip\noindent
Sei weiterhin $X$ integer, $B=\Gamma(X,{\cal O}_X)$,
aber $Y$ nicht mehr die Affinisierung von $X$.
Dann gibt es eine Faktorisierung durch die kanonische Abbildung
$X \longrightarrow {\rm Spek}\, B \longrightarrow {\rm Spek}\, A$.
Sei ${\bf p} \in {\rm Spek}\, A$ mit $X_{\bf p}$ affin.
Sei $ \psi : A \longrightarrow B$ der Ringhomomorphismus und
$ S= \psi (A-{\bf p})$ das multiplikative System in $B$.
Es ist
$X_S=X \times_{{\rm Spek}\, B} {\rm Spek}\, B_S = X \times_{{\rm Spek}\, A} A_{\bf p}$
affin und damit ist f"ur alle
${\bf q} \in {\rm Spek}\, B, \, {\bf q} \cap S = \emptyset$ auch $X_{\bf q}$ affin.
Dann liegt ${\rm Spek}\, B_S$ innerhalb der offenen Menge $U=D({\bf c})$,
auf der die kanonische Abbildung affin ist. 
${\rm Spek}\, B_S \subseteq D({\bf c})$ bedeutet
$S \cap {\bf c} \neq \emptyset $,
sei $g' \in S \cap {\bf c}$. Dann kommt $g'$ von einem Element
$g \in A-{\bf p}$
und dann ist wegen $(\psi^*)^{-1}D(g) =D(g') \subseteq U$ in $D(g)$ eine
affine Umgebung von $\bf p$ in ${\rm Spek}\, A$ gefunden, deren Urbild in $X$
affin ist.
\par\smallskip\noindent
F"ur den nicht integren Fall m"ussen wir auf die Aussage,
dass man die Affinit"at auf den Komponenten testen kann, vorgreifen,
die wir im Abschnitt 1.13 kohomologisch beweisen.
Seien $X_i,\, i=1,...,n$,
die integren Komponenten von $X$ und $\varphi_i:X_i \longrightarrow Y$ die
eingeschr"ankten Abbildungen. Mit $X_y$ sind
nat"urlich auch die $(X_i)_y$ affin, man findet also eine offene, affine
Umgebung $U$ von $y$ mit $\varphi_i^{-1}(U)$ ist affin in $X_i$ f"ur
$i=1,...,n$.
Dann ist aber nach 1.13.2 schon $\varphi^{-1}(U) \subseteq X$ affin. \hfill $\Box$
\par\bigskip
\noindent
Die Offenheit der Affinit"at liefert sofort ein lokales Affinit"atskriterium.
\begin{satz}
Sei $\varphi: X \longrightarrow Y$ ein Schemamorphismus,
$X$ sei noethersch und separiert. Dann sind folgende Aussagen "aquivalent.
\par\smallskip
\noindent
{\rm (1)} $\varphi$ ist ein affiner Morphismus.
\par\smallskip\noindent
{\rm (2)} F"ur jedes $y \in Y$ ist $X_y$ affin.
\end{satz}
{\it Beweis}.
Die Aussage folgt aus dem Satz. Wir geben noch einen anderen Beweis.
Sei (2) erf"ullt,
wir k"onnen annehmen, dass $Y={\rm Spek}\, A$ affin ist und haben dann
die Affinit"at von $X$ zu zeigen.
\par\smallskip\noindent
Es gibt eine Faktorisierung
$X \longrightarrow {\rm Spek}\, B
\longrightarrow {\rm Spek}\, A $
durch die Affinisierung $Z={\rm Spek}\, B$, sei $z \in {\rm Spek}\, B$ 
mit Bildpunkt $y$. Dazu geh"ort dann der Morphismus der
lokalen Schemata $Z_z \longrightarrow Y_y $, der affin ist.
$X \times_Z Z_z \longrightarrow X \times _Y Y_y $ ist dann ein affiner
Morphismus,
und somit ist $X_z$ affin f"ur jedes $z \in Z$. Wir k"onnen also annehmen,
dass die Eigenschaft (2) bez"uglich der Affinisierungsabbildung
erf"ullt ist.
Da die Affinisierungsabbildung im noetherschen separierten
Fall mit Lokalisierung
vertr"aglich ist, sind die Abbildungen
$X_z \longrightarrow {\rm Spek}\,B_z $ ebenfalls die
Affinisierungsabbildungen und da $X_z$ nach Voraussetzung affin ist, liegt
sogar punktweise lokal in $Z$ eine Isomorphie vor.
Dann sichert das folgende Lemma, dass bereits global eine
Isomorphie vorliegt. \hfill $\Box$
\begin{lem}
Sei $\varphi :X \longrightarrow Y$ ein Schemamorphismus, so dass f"ur alle
$y \in Y $ die induzierte Abbildung $\varphi_y: X_y \longrightarrow Y_y$ eine
Isomorphie ist. $X$ sei noethersch. Dann ist $\varphi$ bereits eine
Isomorphie.
\end{lem}
{\it Beweis}.
Ist $\varphi$ sogar lokal isomorph in dem Sinn, dass es eine "Uberdeckung
von $Y$ mit offenen Mengen $U_i$ gibt mit
$\varphi^{-1}(U_i) \longrightarrow U_i$
isomorph, so folgt daraus generell die Isomorphie, auch ohne die
noethersche Voraussetzung. Hier muss man aber sorgf"altiger argumentieren,
wir zeigen zun"achst, dass eine Hom"oomorphie vorliegt.
\par\smallskip\noindent
Die Abbildung ist nach den Voraussetzungen insbesondere bijektiv,
wir wollen $\varphi$ als abgeschlossen nachweisen.
Wegen $X$ noethersch sind die abgeschlossenen Teilmengen von $X$ endliche
Vereinigung irreduzibler abgeschlossener Mengen, es gen"ugt daher, nur die
Bilder irreduzibler abgeschlossener Mengen $Z$ mit generischem Punkt $\eta $
zu untersuchen. Zum Abschluss des Bildes geh"ort der generische Punkt
$\varphi({\eta})$. Sei $y \in \overline{\varphi(Z)}$ im Abschluss des Bildes.
Dann ist $\varphi({\eta}) \in Y_y= {\rm Spek}\, {\cal O}_{Y,y}$ und die gleiche
Situation liegt in dem eindeutig bestimmten Urbildpunkt $y'$ von $y$ 
vor, dort besitzt insbesondere $\varphi({\eta})$ einen eindeutig bestimmten
Urbildpunkt, der dann aber schon mit $\eta $ "ubereinstimmen muss.
Also ist auch $y' \in \bar{\eta} $ und daher $y' \in Z$ und
somit geh"ort $y$ selbst zum Bild von $Z$.
\par\smallskip\noindent
Ist nun $U \subseteq Y$ affin, so ist wegen hom"oomorph auch das Urbild affin,
und f"ur einen Morphismus zwischen Spektra gen"ugt die punktweise Isomorphie,
und das eingangs erw"ahnte lokal-offene Kriterium sichert "uberhaupt
die Isomorphie.  \hfill $\Box$
\par\bigskip
\noindent
{\bf Bemerkung}
Die punktweise Isomorphie ist keine offene Eigenschaft, die
typischen Beispiele hierf"ur sind Lokalisierungen
${\rm Spek}\, A_{\bf p} \longrightarrow {\rm Spek}\, A$. Das lokale Schema
ist ja nur in Ausnahmef"allen offen in ${\rm Spek}\, A$. Ist die Abbildung vom
endlichen Typ, und liegt in einem Punkt lokal eine Isomorphie vor,
so besitzt der Punkt auch eine offene, isomorphe Umgebung.
\par\bigskip
\noindent 
{\bf Beispiel}
Wir wollen anhand eines Beispiels zeigen, dass im Allgemeinen,
ohne die Voraussetzung noethersch, die kanonische Abbildung lokal
punktweise eine Isomorphie, also insbesondere punktweise affin und
bijektiv sein kann, ohne dass das Schema selbst affin ist,
ja ohne das irgendeine Urbildmenge einer affinen, nicht leeren Teilmenge
affin ist. Das zeigt auch,
dass ohne die Voraussetzung vom endlichen Typ die Aussage 1.3.5 falsch ist.
\par\smallskip\noindent
Sei hierzu $ X:= {\rm Spek}\, {\bf Z} $ als Menge, aber
versehen mit folgender Topologie: neben der leeren Menge sei jede Teilmenge
von $X$ offen, die das Nullideal enth"alt.
Als Garbe definieren wir
$$ \Gamma(U, {\cal O}_X) := \bigcap_{(p) \in U}{\bf Z}_{(p)}={\bf Z}_{F(U)}
\mbox{ mit }
F(U) = \{ f \in {\bf Z} : U \subseteq D(f) \} \, .$$
Die Halme sind dabei die des richtigen Spektrums von ${\bf Z}$.
F"ur eine Primzahl $p$ ist
$ U= \{ ( 0 ) , ( p ) \} $ eine offene Umgebung mit
dem Schnitt\-ring
$\Gamma(U,{\cal O}_X) = {\bf Z}_{ ({\bf Z}-{\bf Z}p) } $,
und damit
ist $U$ isomorph zum Spektrum dieses diskreten Bewertungsringes, und damit ist
$X$ ein Schema.
\par\smallskip\noindent
Der globale Schnittring dieses Schemas ist ${\bf Z}$, die kanonische Abbildung
$X,{\cal O}_X \longrightarrow {\rm Spek}\, {\bf Z} $
ist bijektiv, aber keine Einbettung, da etwa das Nullideal im Spektrum von
$\bf Z$ nicht offen ist. $X$ ist nicht quasikompakt und nicht noethersch.
Die globalen Funktionen trennen die Punkte von $X$, $X$ ist aber nicht
reich beringt.
\par\bigskip
\noindent
Einen Spezialfall des obigen lokalen Kriteriums f"ur affine Abbildungen
er\-w"ahn\-en wir hier explizit.
\begin{kor}
Sei $U \subseteq X$ eine offene, noethersche, separierte
Teilmenge eines Schemas. Dann sind folgende Aussagen "aquivalent.
\par\smallskip\noindent
{\rm (1)} $U \hookrightarrow X$ ist eine affine Abbildung.
\par\smallskip\noindent
{\rm (2)} F"ur alle $x \in X$ sind in den lokalen Schemata $X_x$
die $U_x$ affin.
\par\smallskip\noindent
Insbesondere ist ein offenes Unterschema $U \subseteq {\rm Spek}\, A$ in einem
noetherschen Ring $A$ genau
dann affin, wenn dies lokal f"ur $U_{\bf p} \subseteq {\rm Spek}\, A_{\bf p}$
gilt. Dabei kann man sich auf Punkte $\bf p$ au"serhalb von $U$ beschr"anken.
\end{kor}
{\it Beweis}.
Nur der Zusatz ist neu, die Affinit"at von $U_{\bf p}$ ist aber f"ur
${\bf p} \in U$ wegen $U_{\bf p}= {\rm Spek}\, A_{\bf p}$ trivial. \hfill $\Box$
\par\bigskip
\noindent
Offene Unterschemata, die die zweite Bedingung des Satzes erf"ullen, wollen
wir {\it lokal affin} nennen.
\par\bigskip
\noindent 
{\bf Beispiel}.
Sei $X$ ein Schema und $\cal L$ ein Geradenb"undel auf $X$ und
$s \in {\cal L}(X)$ ein Schnitt mit der zugeh"origen Invertierbarkeitsmenge
$X_s$. Dann ist $X_s$ lokal affin. Lokal ist ja $\cal L$ trivial und
daher entspricht $s$ einem Schnitt im trivialen B"undel, also einer
Funktion $f \in {\cal O}_X(U)$, wobei $U$ eine affine Menge sei,
und darauf ist $X_s$ gleich der Basismenge $D(f)$ und daher affin.
Nat"urlich muss $X_s$ selbst nicht affin sein.
\par\smallskip\noindent
Ein effektiver Divisor $D$ in einem lokal faktoriellen,
noetherschen Schema entspricht einem Schnitt
$s $ im Geradenb"undel ${\cal L}_D$. Dabei ist dann
$X_s$ das offene Komplement des Tr"agers des Divisors. Insbesondere ist also
in dieser Situation das Komplement eines Divisors lokal affin,
und die
Affinit"at eines Hyperfl"achenkomplementes ist ein globales Problem.
Diese Voraussetzungen werden
in der Regel vorliegen, wenn wir nach Bedingungen f"ur $U$ affin suchen in der
Situation, wo $U$ eine offene Teilmenge
in einem projektiven Schema ist.
\par\bigskip
\noindent
Sei $D({\bf a})=U \subseteq X={\rm Spek}\, A$
ein offenes Unterschema eines affinen noetherschen Schemas.
Diese offene Einbettung faktorisiert durch
$U \longrightarrow U^{\rm aff}={\rm Spek}\, \Gamma(U,{\cal O}_X)
\longrightarrow {\rm Spek}\, A$, wobei vorne die kanonische, offene Einbettung
steht und hinten der zur Restriktionsabbildung geh"orende Morphismus
affiner Schemata. Wie in Satz 1.4.2 gezeigt, ist $U$ genau dann affin, wenn
in $B= \Gamma(U,{\cal O}_X)$ das Erweiterungsideal ${\bf a}B$ das
Einheitsideal ist. Dies ist genau dann der Fall, wenn das zur"uckgeschnittene
Ideal ${\bf b}:=R^{-1}({\bf a}B)$ das Einheitsideal ist. Auch bei $U$ nicht
affin ist dieses Ideal von Interesse und liefert eine weitere
Beschreibungsm"oglichkeit des Ortes, wo sich $U$ nicht affin verh"alt.
Es sei wieder ${\bf c}=\{ f \in A :\, U_f= U \cap D(f) \, {\rm affin} \}$.
Mit diesen Bezeichnungen gilt der folgende Satz.
\begin{satz}
Der nicht affine Ort wird beschrieben durch
$${\rm Naf}\,(U,X) := \{ x \in X : U_x {\rm \, nicht \, affin \, } \}
= V({\bf c}) =V({\bf b}) = \overline{R^*(V({\bf a}B))} \, .$$
Die Abbildung $R^*: {\rm Spek}\, B \longrightarrow {\rm Spek}\, A$ ist auf
$D({\bf a}B)$ eine offene Einbettung und das Bild von
$V({\bf a}B)$ liegt in ${\rm Naf}\,(U,X)$.
\par\smallskip\noindent
Es ist $U \cap {\rm Naf}\, (U,X) = \emptyset$.
Ist $X$ integer und $U={\rm Def}\, (q)$, so ist
${\rm Naf}\,(U,X) \subseteq V(N(q)+Z(q))$. 
\end{satz}
{\it Beweis}.
Die erste Gleichung folgt direkt aus Satz 1.8.1.
Die zweite Gleichung ergibt sich aus den "Aquivalenzen
$U_{\bf p} $ affin gdw. $({\bf a}B)_{\bf p}=(1)$ gdw.
$R(A-{\bf p}) \cap {\bf a}B \neq \emptyset $ gdw.
$(A-{\bf p}) \cap R^{-1}({\bf a}B) $ gdw.
${\bf b} \not\subseteq {\bf p}$.
Die letzte Gleichung gilt immer, die Zus"atze sind klar.
\par\bigskip
\noindent 
{\bf Bemerkung}
 In \cite{EGAIV}, 21.12.5, wird anders als hier
unter ${\rm Daf} (U/X)$ das Bild von $V({\bf b}B)$ verstanden.
Die hier vertretene Bezeichnungsweise findet sich auch
in \cite{binsto1}, \S 1.
\par
\bigskip
\noindent
\subsection{Affinit"at unter Deformationen}
\par
\bigskip
\noindent
Nachdem wir in den beiden vorausgegangenen Abschnitten zu einem
Schemamorphismus $f: X \longrightarrow Y$ die lokalen Fasern $X_y,\, y \in Y,$
in Hinblick auf die Affinit"at untersucht haben und insbesondere die
Affinit"at als eine offene Eigenschaft nachgewiesen haben,
wenden wir uns jetzt den echten Fasern $X_{k(y)},\, y \in Y,$ zu. Diese Schemata bilden
dann eine durch $Y$ parametrisierte Familie von Schemata. Dieser
Sprachgebrauch ist insbesondere in der Situation "ublich, wo $Y$ eine glatte,
integre Kurve ist, $X$ integer und die Abbildung dominant. Dann ist
die Abbildung flach, siehe \cite{haralg}, III.9.7.
\par\smallskip\noindent
Ist $f:X \longrightarrow Y$ affin, so sind nat"urlich alle Fasern affin.
F"ur den generischen Punkt von $Y$ bei $Y$ integer ist die Faser gleich
der lokalen Faser.
Ist diese Faser affin, so ist schon auf einer offenen Umgebung die
Abbildung affin und damit sind darauf auch die Fasern affin.
Ansonsten kann sich aber das Affinit"atsverhalten von Faser zu Faser "andern,
wie folgende Beispiele zeigen.
\par\bigskip
\noindent
{\bf Beispiel}
Schon ganz einfache Beispiele zeigen, dass man die Affinit"at einer
Abbildung nicht durch die Affinit"at aller echten Fasern erhalten kann,
dazu betrachtet man etwa die offene Einbettung eines
quasiaffinen, nicht affinen Schemas in ein affines Schema.
Die Fasern sind entweder leer oder einpunktig.
\par\bigskip
\noindent
{\bf Beispiel}
Wir betrachten die Projektion einer punktierten Ebene auf eine Gerade,
${\rm Spek}\, K[T,S] - \{ (T,S) \} \longrightarrow {\rm Spek}\, K[T]$.
Die punktierte Ebene ist nicht affin,
jede Faser ist affin, die Abbildung ist auf $D(T)$ affin.
Die lokale Faser "uber dem Nullpunkt ist
${\rm Spek}\, K[T]_{(T)}[S] - \{ (T,S) \}$ und damit nicht affin.
\par\bigskip
\noindent
{\bf Beispiel}
Ein Beispiel einer birationalen
Affinisierungsabbildung, die surjektiv ist und wo die Fasern alle affin sind,
sieht so aus. Man startet mit der affinen Ebene und bl"ast in einem Punkt auf,
und nimmt dann aus der exzeptionellen Faser wieder einen Punkt heraus.
Dieses Schema ist nicht affin, "uber dem aufgeblasenen Punkt liegt eine
affine Gerade als Faser. Eine Affinisierungsabbildung liegt vor, da sich bei
den Operationen der globale Schnittring nicht "andert.
\par\bigskip
\noindent
{\bf Beispiel}
Sei $X:={\bf A}_K^3 - \{ 0 \} $ der punktierte affine Raum und
$p: X \longrightarrow {\bf A}^1_K={\rm Spek}\, K[T]$
die (glatte) Projektion auf eine affine Gerade.
F"ur  $T \neq 0$ ist die Faser eine affine Ebene, f"ur $T=0$ eine punktierte
Ebene, also nicht affin. Betrachtet man die Situation lokalisiert am
Nullpunkt $(T)$, so erh"alt man ein Beispiel einer Abbildung auf einen
diskreten Bewertungsring, wo die generische Faser affin ist und die
spezielle Faser nicht. 
\par\bigskip
\noindent
{\bf Beispiel}  Sei $X= {\bf A}^3- C $, wobei $C$ eine in einer affinen
Ebene liegende Hyperbel ist, also in Gleichungen $C=V(x,yz-1)$.
Wir betrachten die durch $y$ definierte Abbildung auf die affine Gerade.
F"ur $y = 0$ ist die Faser eine affine Ebene, f"ur $y \neq 0$ sind
die Fasern punktierte Ebenen. Lokalisiert man am Nullpunkt, so erh"alt man
eine Abbildung auf einen diskreten Bewertungsring, wo die spezielle Faser
affin ist, die generische Faser nicht.
\par\bigskip
\noindent
{\bf Beispiel}
Ein "ahnliches Beispiel ist die Projektion von
$X={\bf A}^1_K \times_K {\bf P}^1_K - \{ P \}$ auf ${\bf A}^1_K$,
die die kanonische Abbildung ist.
Fast alle Fasern sind projektive Geraden, die Ausnahmefaser ist
eine affine Gerade.
\par\bigskip
\noindent
{\bf Beispiel}
Sei $X$ ein separiertes noethersches Schema "uber ${\rm Spek}\, {\bf Z}$.
Zu jedem Primideal, also zu $(0)$ und zu jeder Primzahl $(p)$, bzw. zu jedem
zugeh"origen K"orper $K= {\bf Q},\, {\bf Z}/(p)$
gibt es dann wieder das Schema $X_K$.
Ist $X_{\bf Q}$ affin, so gilt das wieder
in einer offenen Menge, das hei"st f"ur fast alle Primzahlen ist dann
auch $X_{{\bf Z}_{(p)}}$ und insbesondere $X_{\kappa(p)}$ affin.
Das muss aber nicht f"ur alle Primzahlen
gelten, in diesem Sinn ist die Affinit"at abh"angig von der Charakteristik,
wie folgendes Beispiel zeigt.
\par\smallskip\noindent
Sei $B={\bf Z}[S,T,U]/(ST-nU) $, mit einer nat"urlichen
Zahl $n \neq 0$, und sei $X= {\rm Spek}\, B - V(S,U)$.
Die Fasern "uber einem K"orper ${\bf Z} \longrightarrow K$ werden durch die
entsprechenden Gleichungen beschrieben.
Wird $n$ eine Einheit in $K$, so ist
$B \otimes _{\bf Z}K= K[S,T] $, da man die Relation nach $U$ aufl"osen kann.
Es ist dann $V(S,U)=V(S)$ und dessen Komplement ist affin.
Ist dagegen die Charakteristik von $K$ ein Teiler von $n$,
so ist die Faser ${\rm Spek}\, K[S,T,U]/(ST) -\, V(S,U) $.
Das ist die affine Gerade "uber dem Achsenkreuz, wo man eine Gerade auf
einer Komponente herausnimmt. Das ist nicht affin, weil dies auf der anderen
Komponente eine punktierte Ebene ist.
\par
\bigskip
\noindent
\subsection{Eindimensionale Schemata}
\par
\bigskip
\noindent
Schon in 1.4.6 haben wir gesehen, dass jede offene Teilmenge
eines noetherschen, integren, eindimensionalen, affinen Schemas wieder affin
ist. Wir zeigen hier verschiedene weitere
Affinit"atsaussagen f"ur den eindimensionalen
Fall, die einfach aus den "Uberlegungen der Abschnitte 1.7 und 1.8 folgen.
\begin{satz}
Sei $X={\rm Spek}\, A$ ein affines noethersches Schema der Dimension $\leq 1$ und
$U \subseteq X$ eine offene Teilmenge. Dann ist $U$ affin.
\end{satz}
{\it Beweis}.
Wir k"onnen aufgrund des Lokalisierungsprinzips 1.8.4 annehmen,
dass $A$ ein lokaler
noetherscher Ring der Dimension $\leq 1$ ist.
Ist $A$ nulldimensional, so ist die Aussage klar, sei $A$ eindimensional.
Eine offene Teilmenge $\neq X$ enth"alt den abgeschlossenen Punkt nicht
und ist daher ein nulldimensionales noethersches Schema und
nach Kor. 1.5.2 affin. \hfill $\Box$
\begin{kor}
Sei $A$ ein noetherscher Ring und $U \subseteq X={\rm Spek}\, A$
eine offene Teilmenge.
Dann besitzt ${\rm Naf}\,(U,X)$ eine Kodimension $\geq 2$.
\end{kor}
{\it Beweis}.
Dies folgt direkt aus 1.10.1 und 1.8.5. \hfill $\Box$
\begin{satz}
Sei $X$ ein eindimensionales, noethersches, quasiaffines Schema mit
einem eindimensionalen Schnittring $\Gamma(X,{\cal O}_X)=A$.
Dann ist $X$ affin.
\end{satz}
{\it Beweis}
Wir haben den globalen Hilbertschen Nullstellensatz nachzuweisen,
aufgrund von 1.7.3 ist aber die kanonische Abbildung
sogar surjektiv. \hfill $\Box$
\begin{satz}
Sei $X$ ein eindimensionales, noethersches, integres, separiertes und
generisch affines Schema mit
einem eindimensionalen globalen Schnittring $\Gamma(X,{\cal O}_X)=A$.
Dann ist $X$ affin.
\end{satz}
{\it Beweis}.
Die Voraussetzung generisch affin besagt, dass die kanonische
Abbildung $X \longrightarrow {\rm Spek}\, A$ im generischen
Punkt eine Isomorphie ist.
Nach Satz 1.7.4 ist dann die kanonische Abbildung in jedem Punkt der H"ohe eins
von ${\rm Spek}\, A$ eine Isomorphie und damit ist
nach 1.8.2 $X$ affin.   \hfill $\Box$
\par\bigskip
\noindent
{\bf Beispiel}
Verdoppelt man auf einer affinen Gerade einen Punkt,
d.h verklebt man zwei affine Geraden durch Identifizieren von
${\bf A}^1_K-\{ (0) \}$, so erh"alt man ein integres, eindimensionales,
noethersches Schema, das nicht separiert und nicht affin ist.
\par\bigskip
\noindent
{\bf Beispiel}
Eindimensionale, noethersche, quasiaffine Schemata, die nicht affin sind,
m"ussen nach Satz 1.10.2 einen zumindest zweidimensionalen
Schnittring besitzen.
Solche erh"alt man, wenn man mit einem lokalen, noetherschen,
integren, zweidimensionalen Ring $A$ startet und das
Schema $X={\rm Spek}\, A -\{ {\bf m} \}$
betrachtet. $X$ ist eindimensional, aber nicht affin (siehe 1.12.1).
Ist $A$ normal, so ist $\Gamma(X,{\cal O})=A$.
Insbesondere gibt es auch regul"are Schemata von diesem Typ.
\par\bigskip
\noindent
{\bf Beispiel}
Ein eindimensionales Schema mit nulldimensionalen Schnittring kann
nat"urlich nicht affin sein.
Die typischen Beispiele hierf"ur sind projektive Kurven "uber einem
algebraisch abgeschlossenen K"orper.
Durch "Ubergang zur Normalisierung kann man mit dem Satz von Riemann-Roch
zeigen, dass eine echte offene Teilmenge einer projektiven Kurve
affin ist, siehe den Abschnitt 2.4.
F"ur $U=C-\{P\}$ etwa sichert der Satz von Riemann-Roch die Existenz einer
nicht konstanten Funktion auf $U$, die genau in $P$ einen Pol hat.
Fasst man $q$ als meromorphe Funktion nach ${\bf P}^1_K$ auf,
so ist $U=q^{-1}({\bf A}^1_K)$ und die Affinit"at von $q$ folgt etwa aus
1.8.2. Aus der Separiertheit folgt die Affinit"at dann f"ur jede
Teilmenge.
\par\bigskip
\noindent
{\bf Beispiel}
Es gibt regul"are, quasiaffine, noethersche Kurven, aus der man einen
abgeschlossenen Punkt herausnehmen kann, so dass der Restraum nicht affin ist.
Hierzu betrachtet man einen zweidimensionalen, noetherschen, integren, lokalen
Ring $A$, der regul"ar in der Kodimension eins ist, aber nicht normal, und
zwar so, dass in der Normalisierung mehrere Punkte "uber dem abgeschlossenen
Punkt liegen.
In ${\rm Spek}\, A$ gibt es dann Kurven $C=V({\bf p})$, deren Komplement
nicht affin ist, siehe Abschnitt 4.5.
$X=D({\bf m})$ ist eindimensional und regul"ar, darin ist $\bf p$ ein
abgeschlossener Punkt, und dessen Komplement ist $D({\bf p})$, also nicht
affin.
\par\bigskip
\noindent
{\bf Bemerkung}
Es ist offen, ob es eindimensionale, noethersche, integre, separierte
Schemta mit eindimensionalem globalen Schnittring gibt, die nicht affin sind.
Dabei kann man annehmen, dass das Schema normal und der globale
Schnittring ein diskreter Bewertungsring ist.
Nach Satz 1.10.3 darf die generische Faser dann nicht affin sein, sie
ist dann ein eindimensionales Schema mit dem Quotientenk"orper des
Bewertungsringes als globalem Schnittring.
\par
\bigskip
\noindent
\subsection{Krullringe}
\par
\bigskip
\noindent
Die Klasse der noetherschen Ringe ist in verschiedener Hinsicht zu restriktiv,
die Normalisierung eines integren noetherschen Ringes und der globale
Schnittring einer offenen Menge eines noetherschen normalen Schemas
m"ussen nicht noethersch sein,
haben aber mit normalen noetherschen Ringen viele Eigenschaften
gemeinsam, die im Begriff des Krullringes zusammengefasst werden.
Wir erw"ahnen hier nur einige Grundtatsachen und orientieren uns
an \cite{nagloc}, V.33 und \cite{fossum}, Ch.1.
Im n"achsten Abschnitt beweisen wir via
Krullbereichen die Kodimensionseigenschaft affiner Teilmengen
f"ur noethersche separierte Schemata.
\par\bigskip
\noindent
{\bf Definition}
Ein {\it Krullring} $R$ ist ein Integrit"atsbereich mit folgenden drei
Eigenschaften:
\par\smallskip
\noindent
{\rm (1)} Die Lokalisierungen an Primidealen der H"ohe eins sind diskrete
Bewegungsringe.
\par\smallskip
\noindent
{\rm (2)} $R$ ist der Durchschnitt aller Lokalisierungen der H"ohe eins.
\par\smallskip
\noindent
{\rm (3)} Jedes von Null verschiedene Element liegt nur in endlich vielen
Primidealen der H"ohe eins.
\par\bigskip
\noindent
{\bf Bemerkung}
Ein noetherscher normaler Integrit"atsbereich ist ein Krullbereich,
und allgemeiner ist
die Normalisierung eines noetherschen Integrit"atsbereiches ein
Krullbereich, siehe \cite{nagloc}, 33.4 und  33.10.
\par\smallskip\noindent
"Aquivalent dazu ist folgende Charakterisierung: $R$ ist der Durchschnitt
von einer Familie $F$ von diskreten Bewertungsringen in einem
K"orper $K$ derart, dass f"ur $0 \neq a \in R$ gilt: $r$ ist in
fast allen Bewertungsringen der Familie eine Einheit,
siehe \cite{fossum} I, Prop. 1.9.
\begin{prop}
Sei $A$ ein Krullbereich und $S \subseteq A$ ein multiplikatives System.
Dann ist $A_S$ ebenfalls ein Krullbereich, und zwar ist
$A_S =\bigcap_{ht\, {\bf p}=1 :\, {\bf p} \cap S = \emptyset} A_{\bf p}$.
\par\smallskip\noindent
Insbesondere sind Lokalisierungen krullscher Ringe wieder
krullsch.
\end{prop}
{\it Beweis}.
Siehe hierzu \cite{nagloc}, Theorem 33.5. \hfill $\Box$
\begin{satz}
Sei $A$ ein Krullbereich und
$D({\bf a})=U \subseteq {\rm Spek}\, A =X$ eine offene Menge.
Dann ist der globale Schnittring $B=\Gamma(U,{\cal O}_X)$ von $U$
wieder ein Krullbereich.
\par\smallskip\noindent
Ist $U$ nicht affin, so besitzt das Erweiterungsideal ${\bf a}B$ eine
H"ohe $\geq 2$.
\end{satz}
{\it Beweis}.
Es ist
$B = \bigcap_{ {\bf p} \in U: \, {\rm ht}\, {\bf p} = 1} A_{\bf p}$.
Dies gilt zun"achst f"ur $U=D(f)$ und folgt daraus sofort allgemein.
Damit ist $B$ der Durchschnitt von diskreten Bewertungsringen. Sei
$0 \neq r =a/b \in B$ mit $a,b \in A$. Dann sind $a$ und $b$ in
fast allen Lokalisierungen von $A$ der H"ohe eins invertierbar und das gilt
dann auch f"ur $r$ bez"uglich der kleineren Familie von diskreten
Bewertungsringen. Die obige zweite Charakterisierung zeigt, dass $B$
ein Krullbereich ist.
\par\smallskip\noindent
Ist $U$ nicht affin, so ist ${\bf a}B$ nicht das Einheitsideal, ${\bf a}B$
beschreibt damit eine echte Teilmenge von ${\rm Spek}\, B$, dessen globaler
Schnittring $B$ ist. Es ist zu zeigen, dass eine solche Teilmenge alle
Primideale der H"ohe eins enthalten muss.
Sei angenommen, $\bf p$ geh"ore nicht dazu, sei $p$ eine Ortsuniformisierende
von $B_{\bf p}$, $q=p^{-1}$ und ${\bf q}_1,...,{\bf q}_n$ die weiteren
Polstellen von $q$. Es gibt ein $f \in B$ mit $f \in {\bf q}_i,\, i=1,...,n,$
und $f \not\in {\bf p}$. Dann ist f"ur eine geeignete Potenz
$f^nq$ in allen diskreten Bewertungsringen au"ser in $B_{\bf p}$ definiert,
also eine nicht ausdehnbare Funktion auf der offenen Menge.  \hfill $\Box$
\begin{satz}
Sei $A$ ein Krullbereich und $U \subseteq X={\rm Spek}\, A$ eine affine
offene Teilmenge. Dann besitzt das Komplement von $U$ die reine
Kodimension eins, als affine Teilmengen kommen nur
Komplemente von effektiven Divisoren in Frage.
\end{satz}
{\it Beweis}.
Sei $U=D({\bf a})$ und ${\bf a} \subseteq {\bf p}$ minimal,
wobei $\bf p$ eine H"ohe $\geq 2$ habe.
Wir lokalisieren an ${\bf p}$.
$A_{\bf p}$ ist ein lokaler Krullbereich der Dimension $\geq 2$,
$U_{\bf p}=D({\bf p}A_{\bf p})$ enth"alt alle Primideale der H"ohe eins
und damit ist
$\Gamma(U_{\bf p},{\cal O})=A_{\bf p}$.
Damit ist $U_{\bf p} \hookrightarrow {\rm Spek}\, A_{\bf p}$ die
Affinisierungsabbildung von $U_{\bf p}$, $U_{\bf p}$ und $U$ sind
somit nicht affin. \hfill $\Box$
\par
\bigskip
\noindent
\subsection{Die Kodimensionseigenschaft affiner Teilmengen}
\par
\bigskip
\noindent
Mit der Aussage f"ur Krullringe k"onnen wir nun das entsprechende Resultat f"ur
noethersche Ringe beweisen.
\begin{satz}
Sei $U$ eine offene Teilmenge eines noetherschen Schemas $X$.
Ist die Inklusion $U \hookrightarrow X$ ein affiner Morphismus, so besitzt
jede Komponente von $Y=X-U$ die Kodimension $\leq 1$.
\par\smallskip\noindent
Ist $X$ separiert, so gilt das insbesondere, wenn $U$ selbst affin ist.
\end{satz}
{\it Beweis}.
Der Zusatz ist klar, da im separierten Fall der Durchschnitt
zweier affiner Mengen wieder affin ist.
Sei also $U \hookrightarrow X$ affin und
$Z$ eine Komponente von $X-U=Y$ der Kodimension $r$,
wir haben $r \leq 1$ zu zeigen. Wir lokalisieren am generischen
Punkt zu $Z$ und erhalten das Spektrum eines lokalen Ringes $R$
der Dimension $r$, in
dem alle anderen Komponenten von $Y$ herausgeschnitten wurden und
folglich die offene Menge gleich $D({\bf m})$ ist, wobei $\bf m$ das
maximale Ideal von $R$ ist.
$D({\bf m})$ ist nach Voraussetzung und 1.6.2 affin.
\par\smallskip\noindent
Die H"ohe von $\bf m$ werde oberhalb des Primideals
$\bf p$ der H"ohe null angenommen, wir betrachten den
affinen Basiswechsel $V({\bf p}) \longrightarrow X$ auf die zugeh"orige
Komponente. Unter diesem affinen Morphismus bleiben
die Voraussetzungen und die H"ohe des maximalen Ideals erhalten, und
folglich k"onnen wir annehmen, dass ein lokaler, integrer,
noetherscher Ring $A$ vorliegt, indem $D({\bf m})$ affin ist.
\par\smallskip\noindent
Wir gehen zur Normalisierung $A \longrightarrow A'$ "uber.
Bei einer ganzen Erweiterung gilt das going up, das hei"st
insbesondere, dass "uber jeder Primidealkette
eine Primidealkette liegt. Insbesondere gibt es eine solche "uber einer
Kette, die die H"ohe von ${\bf m}$ realisiert, es sei $\bf n$ ein
so gewonnenes Primoberideal zu $\bf m$ mit (zumindest) der gleichen H"ohe
wie $\bf m$. 
Dann ist ${\bf m}A' \subseteq {\bf n}$, und es kann da auch kein weiteres
Primideal dazwischen liegen, da bei einer ganzen Erweiterung
die Fasern nulldimensional sind,
siehe \cite{nagloc}, 10.7-10.10, oder \cite{eisenbud}, Cor. 4.18.
$\bf n$ ist also eine Komponente des Erweiterungsideals ${\bf m}A'$,
dessen H"ohe $\geq r$ ist. Da die Normalisierung ein Krullbereich ist,
folgt $r \leq 1$ aus 1.11.3.  \hfill $\Box$
\par\bigskip
\noindent
{\bf Bemerkung}
Ohne die Voraussetzung der Separiertheit ist die
Aussage nicht richtig. Hierzu betrachte man eine affine
Ebene, auf der man einen Punkt verdoppelt,
vergleiche \cite{haralg}, II, example 2.3.6, so ist dieses Schema ohne
einen dieser Punkte, also nach Herausnahme
eines Punktes der Kodimension zwei, eine affine Ebene.
\par\smallskip
\noindent
Die Aussage gilt auch nicht ohne die noethersche Voraussetzung.
Dazu betrachte man einen Ring, in dem der Krull'sche Hauptidealsatz
verletzt ist, wo also eine Funktion $f \in A$ ein minimales
Primoberideal der H"ohe $\geq 2$ hat. Nach Lokalisierung gilt dann sogar,
dass geometrisch $f$ den abgeschlossenen Punkt beschreibt, also
$D(f)=D({\bf m})$. $D(f)$ ist aber nat"urlich immer affin.
Beispiel f"ur diese Situation finden sich in \cite{scheja}, \S59, Aufg.34,35.
\par\bigskip
\noindent
Das folgende Korollar wird in 3.2 verwendet.
\begin{kor}
Sei $X$ ein noethersches separiertes Schema und $W \subseteq X$ eine
offene Teilmenge, die alle Punkte der Kodimension eins umfasst, $Z=X-W$.
F"ur eine abgeschlossene Menge $Y$ gilt dann: ist $W-Y \cap W$ affin,
so ist $Z \subseteq Y$.
\end{kor} 
{\it Beweis}.
Die Menge $W- Y \cap W$ ist eine affine, offene Teilmenge von $X$, deren
Komplement gleich $Y \cup Z$ ist.
Da $Z$ eine Kodimension $\geq 2$ hat, muss $Z$ zu $Y$ geh"oren. \hfill $\Box$
\par\bigskip
\noindent
Die im Satz angef"uhrte Kodimensionseigenschaft besitzen neben den affinen
Schemata auch die umfassendere Klasse der semiaffinen Schemata,
siehe \cite{goodlandman} f"ur eine genaue Untersuchung dieser Schemata.
\par\bigskip
\noindent
{\bf Definition}
Ein Schema hei"st {\it semiaffin}, wenn es einen eigentlichen Morphismus
$X \longrightarrow {\rm Spek}\, A$ in ein affines Schema gibt.
\par\bigskip
\noindent
{\bf Bemerkung}
Ist $X \longrightarrow {\rm Spek}\, A$ ein eigentlicher Morphismus, so
faktorisiert dieser durch
$X \longrightarrow {\rm Spek}\, \Gamma(X,{\cal O}_X)
\longrightarrow {\rm Spek}\, A$.
Das zeigt, das zu einem semiaffinen Schema bereits der
Affinisierungsmorphismus eigentlich ist.
Affine Schemata und eigentliche Schemata sind nat"urlich semiaffin.
Ebenso sind Aufblasungen von affinen Schemata semiaffin.
Man kann semiaffine Variet"aten als analog zu holomorph konvexen R"aumen
ansehen,
ein komplexer Raum ist n"amlich genau dann holomorph konvex,
wenn es eine eigentliche Abbildung auf einen Steinschen Raum gibt,
siehe \cite{gunning}, VII. Theorem 9.
Die Analogie hinkt nat"urlich insofern, dass die Analogie zwischen affin und
steinsch hinkt. Ist das Komplement einer analytischen Teilmenge $Y$ in einem
komplexen Raum holomorph konvex, so muss $Y$ eine Hyperfl"ache sein,
siehe \cite{grauert}, V.§1 Satz 4. Das gilt ebenso f"ur semiaffine Teilmengen.
\begin{satz}
Sei $X$ separiert und $U \hookrightarrow X$ eine offene
semiaffine Teilmenge. Dann ist diese Einbettung eine affine Einbettung.
Ist $X$ zus"atzlich noethersch, so besitzt jede Komponente von $Y=X-U$
eine Kodimension $\leq 1$.
\end{satz}
{\it Beweis}.
Sei $p: U \longrightarrow {\rm Spek}\, A$ der eigentliche
Strukturmorphismus,
und sei $i:U \hookrightarrow X$ die offene Einbettung.
Wir betrachten den Produktmorphismus
$f=i \times p: U \longrightarrow
X \times _{{\rm Spek}\, {\bf Z}}{\rm Spek}\, A$.
Da $X$ ("uber $\bf Z$) separiert ist, gilt dies auch f"ur die Projektion
$X \times _{\bf Z} {\rm Spek}\, A \longrightarrow {\rm Spek}\, A$, daher ist
mit der Gesamtabbildung auch $f$ eigentlich. Betrachtet man die
Projektion auf $X$, so sieht man, dass auch $f$ eine topologische
Einbettung sein muss.
$f$ ist also eine abgeschlossene Einbettung und somit eine affine Abbildung.
Da auch die Projektion
$X \times_{\rm Spek {\bf Z}} {\rm Spek}\, A \longrightarrow X$ affin ist,
folgt die Affinit"at der Einbettung
$U \hookrightarrow X$.  \hfill $\Box$
\par\bigskip
\noindent
Wir werden uns im Folgenden bei der Untersuchung der Affinit"at offener
Teilmengen auf Komplemente von Hyperfl"achen beschr"anken k"onnen.
Es gibt nat"urlich Hpyerfl"achen, deren Komplemente nicht affin sind.
Die Kodimensionseigenschaft liefert ein wichtiges notwendiges Kriterium
f"ur die Affinit"at eines offenen Unterschemas
$U=D({\bf a}) \subseteq{\rm Spek}\, A$:
\par\smallskip\noindent
F"ur alle Ringwechsel $A \longrightarrow A'$ mit einem noetherschen Ring
(oder einem Krullbereich) $A'$ ist das Erweiterungsideal ${\bf a}'={\bf a}A'$
das Einheitsideal oder jede Komponente davon besitzt die H"ohe $\leq 1$.
\par\smallskip\noindent
Aufgrund der Affinit"at eines Ringwechsels muss ja die Kodimensionsbedingung
in dem beschriebenen Sinn universell gelten.
Mit diesem Kriterium er"ubrigen sich h"aufig Berechnungen des Schnittringes,
um nachzuweisen, dass eine offene Menge nicht affin ist, stattdessen
gibt man einen Ringwechsel an, wo das Ideal die Kodimension zwei bekommt.
\par\smallskip\noindent
Wir geben daf"ur ein Beispiel, das in verschiedener Hinsicht typisch ist und
sp"ater unter verschiedenen Aspekten wieder auftaucht, als
Hyperbelring, als Monoidring, als Minorenring.
\par\bigskip
\noindent
{\bf Beispiel}
Sei $K$ ein K"orper und $A=K[X,Y,U,V]/(UX-VY)$.
$A$ ist ein normaler dreidimensionaler Integrit"atsbereich.
Das Ideal $( X , Y ) $ ist prim mit dem Resklassenring
$K[U,V]$, so dass es die H"ohe eins besitzt.
Wir betrachten die Reduktion modulo $( U,V ) $, also die
Abbildung $A \longrightarrow A'=K[X,Y]$ mit $U,V \longmapsto 0$.
Das Erweiterungsideal ist wieder $( X , Y) $, und das hat die
H"ohe zwei (es ist ein Punkt auf einer Ebene). Daher ist $D(X,Y)$ nicht
affin.
\par\bigskip
\noindent
Im vorangegangen Beispiel konnte man mit einem abgeschlossenen Unterschema
(kleinerer Dimension), also mit einem Restklassenring,
die Nichtaffinit"at erweisen. Zweidimensionale abgeschlossene Unterschemata,
also Fl"a\-chen, sind immer die ersten Kandidaten, um f"ur eine Hyperfl"ache
zu zeigen, dass das Komplement nicht affin sein kann, n"amlich dann,
wenn sich Hyperfl"ache und Fl"ache nur in einem Punkt schneiden.
Das offene Komplement enth"alt dann eine punktierte Fl"ache als abgeschlossenes
Unterschema.
\par\bigskip
\noindent
Nat"urlich ist die Frage zu stellen, inwiefern man stets zu einer
offenen Menge durch einen geeigneten Ringwechsel, bei dem die H"ohe eins
verloren geht, zeigen kann, dass sie nicht affin ist.
Diese Fragestellung ist der Inhalt von Kapitel 4, f"ur den zweidimensionalen
Fall siehe 1.15.
\par
\bigskip
\noindent
\subsection{Kohomologische Charakterisierung affiner Schemata}
\par
\bigskip
\noindent
Eine wichtige und fundamentale Charakterisierung affiner Schemata 
beruht auf dem Verschwinden aller h"oheren Kohomologiegruppen zu
quasikoh"arenten Garben. Sie erlaubt es, affine Schemata als die
kohomologiefreien Bausteine beliebiger Schemata anzusehen, mittels denen man 
\v{c}echkohomologisch die Kohomologien 
von separierten Schemata berechnen kann.
\par\smallskip\noindent
In diesem Abschnitt besprechen wir die kohomologische Charakterisierung,
die auf Serre zur"uckgeht, so wie einige typische Folgerungen f"ur die
Affinit"at, wie das
Verhalten bei Reduktion und beim "Ubergang zu den Komponenten.
Dies und der im n"achsten Abschnitt zu besprechende Satz von Chevalley erlauben
es, bei der Untersuchung der Affinit"at von noetherschen Schemata
sich auf normale integre Schemata zu beschr"anken. Weiterhin erlaubt
die kohomologische Charakterisierung, sich im Fall, dass ein Schema vom
endlichen Typ "uber einem K"orper vorliegt, zum algebraischen Abschluss
"uberzugehen.
\begin{satz}
Sei $X$ ein noethersches Schema. Dann sind folgende Aussagen "aquivalent.
\par\smallskip\noindent
{\rm (1)} $X$ ist affin.
\par\smallskip\noindent
{\rm (2)} Es ist $H^r(X,{\cal F})=0$ f"ur $r \geq 1$ und alle
quasikoh"arenten Garben $\cal F$.
\par\smallskip\noindent
{\rm (3)} Der Funktor $\Gamma(X,\, .\, )$ ist exakt auf den quasikoh"arenten Garben auf
$X$.
\par\smallskip\noindent
{\rm (4)} Es ist $H^1(X,{\cal F})=0$ f"ur alle quasikoh"arenten Garben $\cal F$.
\par\smallskip\noindent
{\rm (5)} Es ist $H^1(X,{\cal I})=0$ f"ur alle quasikoh"arenten Idealgarben
$\cal I$.
\end{satz}
Zum {\it Beweis} siehe \cite{haralg}, III.3.7
oder \cite{EGAII}, Th\'{e}or\`{e}me 5.2 und \cite{EGAIII}, 1.
f"ur aus (1) folgt (2).  \hfill $\Box$
\begin{satz}
{\rm (1)} Sei $X$ ein noethersches Schema, ${\cal I}$ und ${\cal J}$
Idealgarben so, dass die durch ${\cal I}$ und $\cal J$ definierten
abgeschlossenen Unterschemata affin sind.
Dann ist auch das durch ${\cal I}{\cal J}$ definierte Unterschema affin.
\par\smallskip\noindent
{\rm (2)}
Ein noethersches Schema ist genau dann affin, wenn das reduzierte
Schema affin ist.
\par\smallskip\noindent
{\rm (3)} Ein noethersches Schema ist genau dann affin, wenn alle Komponenten
affin sind.
\end{satz}
{\it Beweis}.
Zu (1).
Die Affinit"at des durch $\cal I$ definierten abgeschlossenen Unterschemas
bedeutet, dass $H^1(X,{\cal G}/{\cal I}{\cal G})=0$ gilt
f"ur jede quasikoh"arente Garbe $\cal G$ auf $X$.
Ohne Einschr"ankung sei ${\cal I}{\cal J}=0$, es ist also die Affinit"at
von $X$ zu zeigen.
F"ur eine quasikoh"arente Garbe $\cal F$ auf $X$ betrachten wir die Sequenz
$$ 0 \longrightarrow {\cal I}{\cal F} \longrightarrow
{\cal F} \longrightarrow {\cal F}/{\cal I}{\cal F} \longrightarrow 0 $$
mit der zugeh"origen langen exakten Kohomologiesequenz
$$  H^1(X,{\cal I}{\cal F}) \longrightarrow H^1(X,{\cal F}) \longrightarrow
H^1(X,{\cal F}/{\cal I}{\cal F}) \, .$$
Dabei verschwindet die hintere Kohomologiegruppe, da das durch $\cal I$
definierte Schema affin ist.
Es verschwindet aber auch die vordere Kohomologiegruppe, da
${\cal I}{\cal F}$ durch ${\cal J}$ annulliert wird
und daher bereits ein ${\cal O}_X/{\cal J}-$Modul ist und die Kohomologie
auf dem durch ${\cal J}$ definierten Unterschema berechnet werden kann.
Damit ist $H^1(X,{\cal F})=0$.
\par\smallskip\noindent
Zu (2).
Sei $\cal I$ die Idealgarbe der nilpotenten
Elemente auf $X$ und ${\cal I}^d=0$.
Aus Teil (1) folgt durch Induktion, dass das durch ${\cal I}^n$ definierte
Unterschema affin ist f"ur jedes $n$,
also ist auch $X$ selbst affin.
\par\smallskip\noindent
Zu (3). Die Komponenten seien beschrieben durch Idealgarben
${\cal I}_1,...,{\cal I}_k$.
Durch Induktion folgt aus Teil (1), dass die durch
${\cal I}_1 ... {\cal I}_k$ definierten Unterschemata affin sind. \hfill $\Box$
\begin{kor}
Sei $U$ ein noethersches Schema "uber einem {\rm (}nicht
notwendigerweise noetherschen{\rm )} affinen Basisschema $S$. Sind dann
f"ur alle irreduziblen Komponenten $S_i$ von $S$ die
Schemata $U_i =S_i \times _S U$ affin, so ist bereits $U$ affin.
\end{kor}
{\it Beweis}.
Die $U_i$ sind abgeschlossene Unterschemata und "uberdecken $U$.
Da $U$ nur endlich viele Komponenten hat, liegen die Bildpunkte der
generischen Punkte
in endlich vielen $S_i$ und deren Urbilder bilden bereits eine
abgeschlossene "Uberdeckung von $U$.
Aus dem Satz folgt dann die Aussage. \hfill $\Box$
\par\bigskip
\noindent
Man kann den Satz auch so sehen: die Familie der Basiswechsel, die den
"Ubergang zu den irreduziblen Komponenten beschreibt, erlaubt den
R"uckschluss auf die Affinit"at.
So wie das Lokalisierungsprinzip im affinen Basisfall es erlaubt,
im wesentlichen nur lokale Ringe zu betrachten, erlaubt es die Aussage,
sich auf integre Ringe zu beschr"anken.
\par
\bigskip
\noindent
\subsection{Der Satz von Chevalley}
\par
\bigskip
\noindent
Sei $X \longrightarrow Y$ ein Schemamorphismus.
Wir fragen, unter welchen Bedingungen
aus $X$ affin folgt, dass $Y$ ebenfalls affin ist. Die Fragestellung macht
nat"urlich wenig Sinn, wenn die Abbildung nicht surjektiv ist, und somit
ist zu fragen, wann Bilder affiner Schemata wieder affin sind.
Dabei interessieren wir uns insbesondere f"ur die Situation, wo $Y$ "uber
einem Basissschema $S$ gegeben ist, und $X$ durch einem Basiswechsel aus
$Y$ entsteht. Dabei ist vor allem wieder die Situation
$Y \subseteq S={\rm Spek}\, A$
offen und einem Ringwechsel $A \longrightarrow A'$ mit surjektiver
Spektrumsabbildung von Interesse. Wir fragen uns im Einklang mit dem
nach Satz 1.6.4 beschriebenen Konzept, inwiefern man auf die Affinit"at durch
die Affinit"at nach einem einzigen Basiswechsel zur"uckschlie"sen kann.
Der erste Satz in diese Richtung ist der folgende Satz von Chevalley.
\par\bigskip
\noindent
\begin{satz}
Sei $X$ ein affines Schema, $Y$ ein noethersches Schema und
$f:X \longrightarrow Y$ eine endliche surjektive Abbildung. Dann ist auch
$Y$ affin.
\end{satz}
{\it Beweis}.
Wir geben im wesentlichen den Beweis aus \cite{EGAII}, 6.7.1.
Man kann nach dem Satz 1.13.2 annehmen, dass $Y$ integer ist, da beim
"Ubergang zu einer reduzierten Komponente die Voraussetzungen erhalten bleiben.
Ferner muss dann zumindest eine Komponente von $X$, die man wieder als
reduziert annehmen kann, $Y$ dominieren (und damit surjektiv auf $Y$ sein),
und somit k"onnen wir auch $X$ als integer annehmen.
\par\smallskip\noindent
F"ur jede echte abgeschlossene Teilmenge $Z$ von $Y$ erhalten sich unter
dem Basiswechsel $Z \longrightarrow Y$ die Voraussetzungen und 
somit k"onnen wir noethersche Induktion anwenden und annehmen, dass jede
solche Teilmenge affin ist (Ist $Y$ nulldimensional, so ist es nat"urlich
affin).
Kohomologisch interpretiert hei"st das, dass alle koh"arenten Garben,
deren Tr"ager eine echte abgeschlossene Teilmenge ist, nullkohomolog sind,
da die Kohomologie
auf diesen kleineren Schemata berechnet werden k"onnen, die nach Annahme
affin sind.
\par\smallskip\noindent
Seien also $X$ und $Y$ integer, es sei ${\cal A}={\cal O}_Y$ und
${\cal B}:=f_*({\cal O}_X)$ die endliche ${\cal A}-$Algebra dar"uber.
Auf einer hinreichend kleinen affinen offenen Menge $U \subseteq Y$ ist dann
${\cal B}$ frei vom endlichen Rang $m$, und somit gibt es auf
$U$ einen ${\cal A}-$Modul Isomorphismus ${\cal A}^m \longrightarrow {\cal B}$.
Dieser Morphismus ist definiert durch $m$ Schnitte
$s_i \in {\cal B}(U)= \Gamma(f^{-1}(U),{\cal O}_X)$. Da $X$ affin ist, gibt es
ein $0 \neq g \in \Gamma(X,{\cal O}_X)$ mit $D(g) \subseteq U$.
Dann gibt es nach \cite{haralg}, Theorem II.5.3, eine Potenz $g^n$,
so dass $g^ns_i$ von globalen
Schnitten $t_i$ herkommen.
Diese globalen Schnitte definieren nun einen globalen Garbenmorphismus
${\cal A}^m \longrightarrow {\cal B}$, der im generischen Punkt von
$Y$ eine Isomorphie ist.
\par\smallskip
\noindent
Es wird auf $Y$ weiterargumentiert, ausgehend von dem generisch isomorphen
Garbenmorphismus ${\cal A}^m \longrightarrow {\cal B} \longrightarrow {\cal T}$,
wobei $\cal T$ als Kokern der Abbildung eine Torsionsgarbe ist, d.h. im
generischen Punkt und damit auch auf einer offenen Umgebung verschwindet.
Wendet man auf diese Sequenz den linksexakten Funktor
${\cal H}om_{\cal A}(\, . \, ,{\cal F})$ an,
so erh"alt man die exakte Sequenz 
$$0 \longrightarrow {\cal H}om_{\cal A}({\cal T},{\cal F}) \longrightarrow
{\cal G}={\cal H}om_{\cal A}({\cal B},{\cal F}) \longrightarrow
{\cal H}om_{\cal A}({\cal A}^m,{\cal F})= {\cal F}^m \, .$$
Sei nun $\cal F$ torsionsfrei, was insbesondere f"ur Idealgarben erf"ullt ist,
so ist dabei ${\cal H}om_{\cal A}({\cal T},{\cal F}) = 0$ und wir haben eine
kurze exakte Sequenz
$$0 \longrightarrow {\cal G} \longrightarrow  {\cal F}^m 
\longrightarrow {\cal K}=Coker \longrightarrow 0    \, .$$
Als ${\cal B}-$Modul ist dabei $\cal G$ nullkohomolog, und da der
Garbenmorphismus ${\cal G} \longrightarrow {\cal F}^m$ generisch eine
Isomorphie ist, ist  $\cal K$ eine (koh"arente) Garbe mit einem Tr"ager
$Z \subset Y$, deren Kohomologie auf $Z$ berechnet werden kann und die
nach Voraussetzung null ist. Mit der Kohomologiesequenz folgt daraus aber
$H^1(X,{\cal F}^m)=0$ und damit die Aussage $H^1(X,{\cal I})=0$
f"ur jede Idealgarbe $\cal I$.  \hfill $\Box$
\begin{kor}
Sei $X$ ein noethersches integres Schema mit affiner Normalisierung $\tilde{X}$.
Dann ist bereits $X$ selbst affin. Insbesondere
gilt f"ur eine offene Menge $U \subseteq {\rm Spek}\, A$ in einem noetherschen
Integrit"atsbereich $A$, dass $U$ genau dann affin ist, wenn das Urbild in der
Normalisierung ${\rm Spek}\, \tilde{A}$ affin ist.
\end{kor}
{\it Beweis}.
Ist die Normalisierung $\tilde{X}$ endlich "uber $X$
(dazu gen"ugt vom endlichen Typ), so folgt die
Aussage direkt aus dem Satz von Chevalley.
Im allgemeinen sei ${\cal A}$ die Strukturgarbe auf $X$ und $\cal B$ die (quasikoh"arente)
Normalisierungsalgebra darauf, deren (Algebra-)Spektrum im Sinne von
\cite{EGAI}, \S 9.1, oder \cite{haralg}, Exc. II.5.17, die Normalisierung
$\tilde{X}$ definiert.
Sei $f \in {\cal B}(X)= \Gamma(\tilde{X},{\cal O}_{\tilde{X}})$ eine globale
Funktion. Dies definiert eine Algebra ${\cal A}[f] \subseteq {\cal B}$, die
als Untergarbe von $\cal B$ lokal auf einem affinen St"uck durch die
Bedingung ${\cal A}[f]={\cal A}(U)[f]$ festgelegt ist. 
Zu dieser Algebra gibt es dann wieder ein Spektrum, das vom endlichen Typ
und wegen der Ganzheit schon endlich "uber $X$ ist, und auf dem $f$ global
definiert ist.
\par\smallskip\noindent
Sei $X=\bigcup_{i=1,...,n}U_i$ eine affine "Uberdeckung von $X$ und
$\tilde{U}_i$ die affinen Urbilder der $U_i$ in der Normalisierung.
Dann gibt es globale Funktionen
$f_{i,j} \in \Gamma(\tilde{X},{\cal O}_{\tilde{X}})$ 
mit $\tilde{U}_i = \bigcup_{j \in J_i}{\tilde X}_{f_{i,j}}$. 
Durch sukzessive Hinzunahme der beteiligten, endlich vielen
Funktion $f_{i,j}$ in der
eingangs beschriebenen Weise erhalten wir ein Schema $X'$, auf dem diese
Funktionen global definiert sind. Wegen der Surjektivit"at von
$\tilde{X} \longrightarrow X'$ gilt auch $X'_{f_{i,j}} \subseteq U_i'$, wobei
$U'_i$ das affine Urbild von $U_i$ ist. Damit ist $X'$ reich beringt.
Sei nun $X'_{f_k},\, k=1,...,n,$ eine affine "Uberdeckung von $X'$ mit globalen Funktionen
$f_k \in \Gamma(X',{\cal O}_{X'})$. Die entsprechende Situation liegt dann auch
in der Normalisierung vor und wegen deren Affinit"at gibt es dann
Funktionen $a_k \in \Gamma(\tilde{X},{\cal O}_{\tilde{X}})$ mit
$\sum_{k=1,...,n}a_kf_k=1$, und durch Hinzunahme der $a_k$ kann man diese
Relation auch in einer endlichen Erweiterung $X''$ von $X'$ erhalten.
Damit ist $X''$ affin und wegen dem Satz von Chevalley auch $X$.  \hfill $\Box$
\par\bigskip
\noindent
{\bf Bemerkung}
Die Aussage im Satz gilt nicht, wenn man affin durch quasiaffin
ersetzt. Es gibt Beispiele von algebraischen Schemata "uber einem K"orper,
deren Normalisierung quasiaffin ist, sie selbst aber nicht,
siehe \cite{EGAII}, 6.6.13.
\par\bigskip
\noindent
Wir geben einige Beispiele von surjektiven Schemamorphismen (mit weiteren
Eigenschaften)
$f:X \longrightarrow Y$ mit $X$ affin, aber $Y$ nicht.
\par\bigskip
\noindent
{\bf Beispiel}
Eine hom"oomorphe Abbildung $f: X \longrightarrow Y$ mit
$X$ affin, aber $Y$ nicht.
Sei $K$ ein algebraisch abgeschlossener K"orper und
$$ U= {\rm Spek}\, K[X,Y] -\{ (X,Y) \} \longrightarrow {\bf P}_K^1$$
die Kegelabbildung auf die projektive Gerade.
Sei $F \subseteq K[X,Y]$ das multiplikative System, das von allen Primelementen
mit konstantem Term ungleich null und allen Primelementen
vom Grad (bzgl. einer Variablen) $\geq 2$ erzeugt wird. Sei $B=K[X,Y]_F$. In ${\rm Spek}\, B$
"uberleben also nur die Primelemente vom Typ $aX+bY$.
Alle maximalen Ideale $(X-c,Y-d)$ von $K[X,Y]$ enthalten dann
Einheiten, etwa $(X-c)^2+(Y-d)^3$,
und kommen in ${\rm Spek}\, B$ nicht mehr vor, ${\rm Spek}\, B$ ist
also ein eindimensionaler Integrit"atsbereich.
Wir betrachten die Abbildung
${\rm Spek}\, B \hookrightarrow U \longrightarrow {\bf P}_K^1$. Die Faser
"uber einem abgeschlossenen Punkt der projektiven Gerade in $U$ ist eine
punktierte Gerade samt generischem Punkt, und in ${\rm Spek}\, B$
bleibt davon genau
der generische Punkt "ubrig. Umgekehrt ist au"ser dem Nullideal jedes
Primideal in ${\rm Spek}\, B$ ein generischer Punkt einer solchen Gerade, und daher
ist die Abbildung bijektiv. Da beide Schemata eindimensional und noethersch
sind, liegt bereits eine Hom"oomorphie vor; ${\rm Spek}\, B$
ist affin, ${\bf P}^1_K$
nicht.
\par\bigskip
\noindent
{\bf Beispiel}
Ein Ringwechsel $A \longrightarrow A'$ zwischen endlich erzeugten $K-$Algebren
kann eine bijektive Spektrumsabbildung besitzen, ohne dass zu einer offenen
Teilmenge $U \subseteq {\rm Spek}\, A$ aus der Affinit"at des Urbildes $U'$
die Affinit"at von $U$ folgt.
\par\smallskip\noindent
Wir starten mit der Normalisierungsabbildung zu einer Knotenkurve
${\bf A} \longrightarrow  C $, wobei genau zwei Punkte
identifiziert werden sollen.
Punktiert man die affine Gerade an einem der beiden Punkte, so erh"alt
man eine bijektive Abbildung ${\bf A}^\times \longrightarrow C$,
die auch eine Hom"oomorphie ist.
Wir betrachten den affinen Zylinder "uber dieser Abbildung, also
die bijektive (aber nicht hom"oomorphe) Gesamtabbildung
$$ {\rm Spek}\, A' = {\bf A}_K \times {\bf A}_K^\times
\hookrightarrow {\bf A}_K \times {\bf A}_K
\longrightarrow {\bf A}_K \times C ={\rm Spek}\, A \, .$$
Die beiden parallelen Geraden, die unter der hinteren Abbildung
(der Normalisierung) identifiziert
werden, seien mit $G_0$ und $G_1$ bezeichnet, die erste Abbildung sei die
Herausnahme von $G_0$.
\par\smallskip\noindent
Sei $H$ eine Hyperbel auf der affinen Ebene
${\bf A}_K^2={\bf A}_K \times {\bf A}_K$,
die zu $G_0$ disjunkt ist und $G_1$ in einem Punkt schneidet.
Das Bild von $H$ ist eine abgeschlossene Kurve $D$ auf ${\bf A} \times C$,
und deren Urbild besteht aus $H$ und einem isolierten Punkt auf $G_0$.
Mit dem Urbild in der Normalisierung ist das Komplemet von $D$ nicht affin.
Dagegen hat $H$ auf der geschlitzen Ebene affines Komplement, da dort dieser
st"orende Punkt wieder fehlt.
\par\bigskip
\noindent
{\bf Bemerkung}
Ist wie im Beispiel die Abbildung $f: X \longrightarrow Y$ surjektiv,
birational, von endlichem Typ und die Basis $Y$ zus"atzlich normal
(und separiert), so handelt es sich bereits um eine Isomorphie,
insbesondere folgt aus $X$ affin, dass auch $Y$ affin sein muss. 
Ist n"amlich $V$ eine offene affine Teilmenge von $Y$, so ist 
das Urbild $f^{-1}(V)$ ebenfalls affin, und aufgrund von \cite{nagloc}, 33.1,
liegt lokal und damit global eine Isomorphie vor.
\par\smallskip\noindent
Es ist zu fragen, ob sich bei einer Hom"oomorphie $f:X \longrightarrow Y$
zwischen Variet"aten die Affinit"at von $X$ auf $Y$ "ubertr"agt.
Die beiden Beispiele zeigen, dass die Bijektivit"at allein nicht ausreicht
und dass man auf die Voraussetzung von endlichem Typ nicht verzichten kann.
\par\bigskip
\noindent
{\bf Beispiel}
Einfach ergeben sich Beispiele f"ur
surjektive, abgeschlossene, offene, treuflache,
quasiendliche Abbildungen $f:X \longrightarrow Y$ zwischen glatten algebraischen
Kurven mit $X$ affin und $Y$ projektiv.
Dazu startet man mit einer surjektiven Abbildung vom Grad $\geq 2$
einer glatten projektiven Kurve auf eine andere.
Eine solche Abbildung erf"ullt alle angef"uhrten Eigenschaften und ist endlich.
Nimmt man nun aus einer nicht einelementigen Faser einen Punkt heraus,
so wird die obere Kurve affin, die Endlichkeit (und die universelle
Abgeschlossenheit) geht verloren, die anderen Eigenschaften bleiben erhalten.
\par\bigskip
\noindent
{\bf Beispiel}
Jedes quasikompakte Schema $Y$ wird treuflach
"uberdeckt durch ein affines Schema.
Dazu startet man mit einer endlichen affinen "Uberdeckung $Y_i$ von $Y$
und bildet die disjunkte Vereinigung $X=\biguplus_{i=1,..,n} Y_i$.
$X$ ist nat"urlich affin und die Abbildung $X \longrightarrow Y$ wird auf
den disjunkten Teilst"ucken durch die nat"urlichen Einbettungen gegeben.
Die Abbildung ist \'{e}tale und surjektiv. Ist $Y$ separiert, so ist die Abbildung 
auch affin.
Die Abbildung besitzt lokal einen Schnitt und bei $Y$ separiert
liegt dann lokal ein direkter Summand vor, d.h. es
gibt eine offene affine "Uberdeckung $Y_i$ von $Y$ mit affinen Urbildern
und so, dass der globale Schnittring von $p^{-1}(Y_i)$ den von $Y_i$
als direkten Summanden enth"alt.
Ein nicht affines Schema kann also treuflach durch ein affines "uberdeckt
werden. Dies kann nicht sein, wenn sich alles "uber einem Basisschema abspielt
und ein treuflacher Basiswechsel vorgenommen wird.
\begin{satz}
Sei $X$ ein separiertes Schema von endlichem Typ
"uber einem affinen noetherschen Schema ${\rm Spek}\, A$.
$A \longrightarrow A'$ sei ein treuflacher Ringhomomorphismus und
$A'$ ebenfalls noethersch.
Ist dann $X'=X \times_{{\rm Spek}\, A} {\rm Spek}\, A'$ affin, so auch $X$.
\end{satz}
{\it Beweis}.
Bei einem flachen Basiswechsel gilt f"ur eine quasikoh"arente Garbe $\cal F$
auf $X$ und die zur"uckgenommene Garbe ${\cal F}'$ auf $X'$
die Beziehung
$H^{i}(X,{\cal F})\otimes_A A' =H^{i}(X',{\cal F}')$.  
Dies beruht darauf, dass sich der \v{C}ech-Komplex zur Berechnung der
Kohomologie
von $\cal F$ liftet und zur Berechnung auf $X'$ verwendet werden kann,
siehe \cite{haralg}, Theorem III.9.3.
Ist nun $X$ affin, so annulliert $\otimes_A A'$ die Kohomologie
$H^1(X,{\cal F})$, und bei $A \longrightarrow A'$ treuflach folgt daraus
$H^1(X,{\cal F})=0$.  \hfill $\Box$
\par\bigskip
\noindent
{\bf Bemerkung}
Die Aussage ist nat"urlich wieder insbesondere f"ur ein offenes Unterschema
$U=D({\bf a}) \subseteq {\rm Spek}\, A$ mit $A$ noethersch von Interesse.
Als treu\-flache Basiswechsel $A \longrightarrow A'$ spielen f"ur uns
folgende Abbildungen eine Rolle.
\par\medskip\noindent
(1) Ist $\bf m$ ein Ideal in einem noetherschen Ring $A$, so ist die
${\bf m}-$adische Komplettierung treuflach,
siehe \cite{eisenbud}, Theorem 7.2.
Dies ist insbesondere der Fall, wenn $A$ lokal ist mit maximalem
Ideal $\bf m$.
Davon werden wir im n"achsten Abschnitt bei der Behandlung 
zweidimensionaler affiner Schemata Gebrauch machen.
\par\smallskip\noindent
(2) Sei $X$ ein separiertes, algebraisches
$K-$Schema, wobei $K$ ein K"orper sei. Dann ist jeder Basiswechsel
$K \subseteq L$ treuflach, und die Affinit"at von $X$ folgt aus
der von $X_{(L)}$.
Insbesondere kann man zu einem algebraischen Abschluss von $K$ "ubergehen
und damit kann man h"aufig gleich als Grundk"orper einen algebraisch
abgeschlossenen K"orper annehmen.
\par\smallskip\noindent
Bei einem solchen Grundk"orperwechsel k"onnen nat"urlich
auch viele fundamentale Eigenschaften wie integer, normal oder faktoriell
verloren gehen.
In 3.3 wird ein Beispiel angegeben mit ${\rm Spek}\, A$ faktoriell "uber $K$,
aber so, dass es in ${\rm Spek}\, A\otimes_K L$
Hyperfl"achen mit nicht affinem Komplement gibt. 
\par\smallskip\noindent
(3) Sei $A$ ein noetherscher Ring. Mit dem multiplikativen System
$$S=\bigcap_{{\bf p} \in {\rm Spek}\, A} A[X]-{\bf p}A[X]$$
bastelt man sich die sogenannte Kroneckererweiterung
$A(X)=A[X]_S$. Ist $A$ ein K"orper, so ist $A(X)$ der Funktionenk"orper in
einer Variablen.
Die Abbildung $A \longrightarrow A(X)$ ist flach, und in $A(X)$ ist
${\bf p}A(X)$ ein Primideal, so dass die Spektrumsabbildung surjektiv ist
und damit treuflach.
Ist $A$ lokal mit maximalem Ideal $\bf m$,
so ist $A(X)$ ein lokaler Ring unver"andeter H"ohe mit
dem maximalen Ideal ${\bf m}A(X)$.
\par\smallskip\noindent
Ist $A$ lokal, integer, normal und zumindest zweidimensional,
so gibt es in $A(X)$ unendlich viele Elemente $aX+b$ mit $a,b \in {\bf m}$
und $a,b$ teilerfremd; diese Elemente sind dann prim.
Dieser "Ubergang und diese Zusatz"uberlegung
werden im Beweis des Satzes von Nagata im n"achsten
Abschnitt eine Rolle spielen.
\begin{prop}
Sei $A$ ein direkter Summand einer $A-$Algebra $B=A \oplus V$,
wobei $V$ ein $A-$Modul sei.
Ist f"ur $U \subseteq {\rm Spek}\, A$ das Urbild
$U' \subseteq {\rm Spek}\, B$ affin, so ist auch $U$ affin.
\end{prop} 
{\it Beweis}.
Auf einer offenen Menge der Form $D(f),\, f \in A, $
ist $B_f=A_f \oplus V_f$,  und da diese Zerlegungen
mit den Restriktionsabbildungen vertr"aglich sind,
gilt
$\Gamma(U',\tilde{B})=\Gamma(U,\tilde{ A}) \oplus \Gamma(U,\tilde{V})$
f"ur jede offene Menge $U \subseteq {\rm Spek}\, A$.
Ist nun $U'$ affin, so erzeugt in dessen Schnittring das definierende
Ideal $\bf a$ das Einheitsideal, dass muss dann aber schon im Summanden
$\Gamma(U, \tilde{A})$ gelten, also ist $U$ affin. \hfill $\Box$
\par\bigskip
\noindent
{\bf Beispiel}
Besitzt eine Abbildung $f:X \longrightarrow Y$
einen globalen affinen Schnitt, so "ubertr"agt
sich nat"urlich die Affinit"at von $X$ auf $Y$.
Ist etwa $E$ ein Vektorraumb"undel "uber $Y$, also ein
Schemamorphismus $p:E \longrightarrow Y$, der lokal vom Typ
${\bf A}^k\times V \longrightarrow V$ ist mit linearen "Ubergangsabbildungen,
vergl. \cite{haralg}, Exc. II.5.18,
so gibt es immer den Nullschnitt.
Ist der Totalraum $E$ affin, so auch $Y$.
\par\bigskip
\noindent
{\bf Bemerkung}
Im n"achsten Kapitel werden wir sehen, dass bei einer graduierten Algebra
die Kegelabbildung auch den R"uckschluss auf die Affinit"at zul"asst.
In 5.5 wird f"ur die hier betrachteten Morphismen diskutiert, ob sich
die kohomologische Dimension und die Superh"ohe unter ihnen "andern kann.
\par
\bigskip
\noindent
\subsection{Zweidimensionale Ringe}
\par
\bigskip
\noindent
In diesem Abschnitt besprechen wir den Satz von Nagata, dass in einem
zweidimensionalen, normalen,
exzellenten, integren, affinen Schema das Komplement
jeder Hyperfl"ache, also jeder Kurve, affin ist.
Nagata bewies diesen Satz (f"ur einen affinen, normalen, zweidimensionalen
Ring "uber einem
Dedekindbereich) in Zusammenhang mit dem 14. Problem von Hilbert und
der damit zusammenh"angenden Frage, ob der globale Schnittring einer
offenen Menge endlich erzeugt ist.
Zum 14. Hilbertschen Problem vergleiche \cite{hilbertprob},
\cite{zariskiprob},
\cite{naghil}, \cite{rees}.
Im vierten Kapitel werden wir offenen Mengen in affinen Schemata
begegnen, deren Schnittring nicht endlich erzeugt ist.
\par\bigskip
\noindent
Im Beweis des Satzes folgen wir im wesentlichen Nagata, wobei wir aber zuerst
eine Superh"oheneigenschaft herausarbeiten, da dies Thema der beiden letzten
Kapitel der Arbeit ist. Die noethersche Superh"ohe
eines Ideals ist die maximale H"ohe des Erweiterungsideals in einem
noetherschen Ring.
\begin{satz}
Sei $A$ ein zweidimensionaler noetherscher Ring.
Ist die offene Teilmenge $U=D({\bf a}) \subseteq {\rm Spek}\, A=X$
nicht affin, so gibt es eine Abbildung $A \longrightarrow R$ in einen
noetherschen lokalen Ring der Dimension zwei mit $V({\bf a}R)=V({\bf m}_R)$.
\par\smallskip\noindent
Dabei kann erreicht werden, dass die Abbildung
faktorisiert
durch eine ganze Abbildung, eine Lokalisierung,
eine treuflache Abbildung, die im wesentlichen von endlichem Typ ist,
und schlie"slich eine injektive, lokale, birationale
Abbildung, unter der die Erweiterung der Restek"orper endlich ist.
\par\smallskip\noindent
Insbesondere ist $D({\bf a})$ genau dann affin, wenn die noethersche
Superh"ohe von $\bf a$ $\leq 1$ ist.
\end{satz}
{\it Beweis}.
Der Zusatz ist klar, da im Beweis der Hauptaussage im nicht affinen Fall
ein noetherscher Ring
(mit weiteren Eigenschaften) konstruiert wird, bei dem das Erweiterungsideal
die H"ohe zwei bekommt.
Die Umkehrung ist dabei durch die Kodimensionseigenschaft 1.12.1 gegeben.
Die beschriebene Faktorisierung ergibt sich im Verlauf des
Beweises, sie ist in Hinblick auf den n"achsten Satz von Bedeutung.
\par\smallskip\noindent
Sei also $U=D({\bf a})$ nicht affin.
Dann gibt es auch eine Komponente von $X$, auf der $U$ nicht affin ist,
und nach Reduktion gibt es dann einen integren Restklassenring
von $A$, auf dem das Urbild von $U$ nicht affin ist.
Nach 1.14.2 ist dann auch das Urbild in der Normalisierung nicht affin.
Nach \cite{nagloc}, Theorem 33.12, ist die Normalisierung eines
noetherschen zweidimensionalen Ringes wieder noethersch.
Wir k"onnen also annehmen, dass $A$ zweidimensional, normal, integer ist.
Nach dem Lokalit"atsprinzip k"onnen wir weiter annehmen, dass $A$ lokal
ist mit maximalen Ideal $\bf m$.
Durch den treuflachen "Ubergang zu $A(X)$ k"onnen wir schlie"slich
annehmen, dass es in
$A$ unendlich viele Primelemente
gibt, vergleiche den Punkt (3) der letzten
Bemerkung des letzten Abschnittes.
\par\smallskip\noindent
Sei $B=\Gamma(U,{\cal O}_X)$. Da $U$ nicht affin ist, haben nach
1.11.2 die Komponenten des Erweiterungsideals ${\bf b}={\bf a}B$
eine H"ohe von zumindest zwei.
Sei ${\bf m}'$ ein minimales Primoberideal zu $\bf b$.
Unter der Abbildung ${\rm Spek}\, B \longrightarrow {\rm Spek}\, A$ wird
$D({\bf b})$ isomorph auf $D({\bf a})=U$ abgebildet.
${\bf m}'$ wird dabei, wie ganz $V({\bf b})$, auf den nicht affinen Ort
abgebildet, also auf das maximale Ideal $\bf m$ von $A$.
Sei $R=B_{{\bf m}'}$. Dann ist $D({\bf m}')=D({\bf b}R)$ eindimensional
und damit ist $R$ zweidimensional.
\par\smallskip\noindent
Wir zeigen, dass $R$ noethersch ist. Nach \cite{EGAI}, Ch. 0, Prop. 6.4.7,
gen"ugt es hierzu zu zeigen, dass jedes Primideal endlich erzeugt ist.
Sei ${\bf q}' \in D({\bf m}')$ und ${\bf q}' \cap A={\bf q}$.
Dann ist $B_{{\bf q}'}=A_{\bf q}$ und damit
$Q(A/{\bf q})=\kappa({\bf q})=\kappa({{\bf q}'})=Q(B/{\bf q}')$.
Damit ist $A/{\bf q} \hookrightarrow B/{\bf q}' \hookrightarrow Q(A/{\bf q})$
und nach dem Satz von Krull-Akizuki, siehe \cite{nagloc}, 33.2 oder
\cite{bourbaki}, 7 \S2 Prop.5,
ist $B/{\bf q}'$ ein noetherscher
Ring und  $B/{\bf m}'\, (=R/{\bf m}')$
ist ein $A/{\bf m}$ Modul endlicher L"ange, in den
Restek"orpern liegt also eine endliche K"orpererweiterung vor.
\par\smallskip\noindent
Nach Konstruktion gibt es ein Primhauptideal $(p) \in U$.
$p$ ist dann auch in $B$ ein Primelement. Im Krullbereich $B$
ist ja ein Element genau dann
ein Primelement, wenn der zugeh"orige Hauptdivisor ein Basisdivisor ist, und
diese Eigenschaft wird von $A$ her "ubertragen.
Da $B/p$ noethersch ist, ist darin ${\bf m}'$ endlich erzeugt, und damit
ist auch ${\bf m}'$ in $B$ und in $R$ endlich erzeugt.
\par\smallskip\noindent
Sei weiter ${\bf q}'$ ein Primideal der H"ohe eins von $R$ mit dem
Urbildideal $\bf q$. Das maximale, endlich erzeugte Ideal von $R$ sei auch
mit ${\bf m}'$ bezeichnet.
Es ist ${\bf q}R \subseteq {\bf q}'$ und das ist das einzige Primoberideal.
Stimmen diese beiden Ideale "uberein, so ist ${\bf q}'$ endlich erzeugt.
Ansonsten ist
$${\bf p}=({\bf q}R : {\bf q}')
=\{  r \in R:\, r{\bf q}' \subseteq {\bf q}R \} $$
prim"ar zu ${\bf m}'$, da auf $D({\bf m}')$ die beiden Ideale "ubereinstimmen.
Wir haben
${\bf q}' \cdot {\bf p} \subseteq {\bf q}R \subseteq {\bf q}' \cap {\bf p}$.
\par\smallskip\noindent
Da ${\bf m}'$ endlich erzeugt ist, liegt eine Potenz
$({\bf m}')^n$ innerhalb von $\bf p$.
Dann ist $R/({\bf m}')^n$ artinsch (und noethersch), $\bf p$ ist modulo
$({\bf m}')^n$ endlich erzeugt, und damit ist auch $\bf p$ endlich erzeugt
(als Mittelglied einer kurzen exakten Sequenz).
\par\smallskip\noindent
Da $R/{\bf q}'$ noethersch ist, ist mit ${\bf p}/{\bf q}'{\bf p}$
auch
$({\bf q}' \cap {\bf p})/ {\bf q}'{\bf p} \subseteq {\bf p}/{\bf q}'{\bf p}$
ein endlich erzeugter $R/{\bf q}'-$Modul.
Also ist ${\bf q}' \cap {\bf p}$ modulo ${\bf q}'{\bf p}$ endlich erzeugt
und damit auch modulo ${\bf q}R$.
Damit ist bereits ${\bf q}' \cap {\bf p}$ selbst endlich erzeugt.
\par\smallskip\noindent
Da $R/{\bf q}'$ und $R/{\bf p}$ noethersch sind, ist auch
$R/({\bf q}' \cap {\bf p})$ noethersch,
siehe \cite{EGAI}, Ch. 0, Prop. 6.4.2.
Damit ist ${\bf q}'$ endlich erzeugt modulo ${\bf q}' \cap {\bf p}$
und damit ist ${\bf q}'$ selbst endlich erzeugt.
Also ist $R$ noethersch. \hfill $\Box$
\par\bigskip
\noindent
Vor dem eigentlichen Satz von Nagata erinnern wir an den Begriff
{\it exzellent}, vergl. \cite{EGAIV}, §7.
F"ur den Satz muss man lediglich wissen, dass komplette lokale Ringe
exzellent sind und dass die Komplettierung eines lokalen, normalen,
exzellenten Ringes integer ist, zu dieser Aussage siehe
\cite{EGAIV}, § 7.6.2. Insofern gilt der Satz auch unter etwas schw"acheren
Bedingungen, die wir hier aber nicht diskutieren.
\begin{satz}
Sei $A$ ein zweidimensionaler, noetherscher, integrer Ring.
$A$ sei normal und exzellent, oder lokal und komplett.
Dann ist in $X={\rm Spek}\, A$ das Komplement jeder Hyperfl"ache
{\rm (}=Kurve{\rm )}
affin.
\end{satz}
{\it Beweis}.
Ein lokaler kompletter Ring ist exzellent, daher ist die Normalisierung
nach \cite{EGAI}, 7.6.2, wieder komplett und
exzellent. Ferner beh"alt das Erweiterungsideal die H"ohe eins, und
nach 1.14.2 kann man die Affinit"at in der Normalisierung
nachweisen.
Wir m"ussen also nur den ersten Teil zeigen, wobei wir
$A$ als lokal mit maximalen Ideal $\bf m$ annehmen k"onnen.
\par\smallskip\noindent
Sei $U=D({\bf a})$ nicht affin. Wir zeigen, dass dann nicht jede Komponente
von $\bf a$ die H"ohe eins haben kann, was gleichbedeutend ist
zu $V({\bf a})=V({\bf m})$.
Nach dem Satz gibt es eine Abbildung $A \longrightarrow A' \longrightarrow R$,
mit $A \longrightarrow A'$ treuflach (und im wesentlichen von endlichem Typ),
$A'$ ein ebenfalls zweidimensionaler, normaler, lokaler, exzellenter
und integrer Ring,
und mit $R$ lokal und noethersch der Dimension zwei, mit
$V({\bf a}R)=V({\bf m})$, und so, dass zu
$A' \longrightarrow R$ die Erweiterung der Restek"orper endlich ist.
Da $A \longrightarrow A'$ treu\-flach ist, gen"ugt es zu zeigen, dass
$V({\bf a}A')$ nur aus dem abgeschlossenen Punkt besteht. Wir k"onnen $A$
durch $A'$ ersetzen, da $A'$ wieder exzellent ist.
\par\smallskip\noindent
Wir betrachten die Komplettierung
$\hat {{\bf A}} \longrightarrow \hat{R}$, bei der
ebenfalls die Erweiterung der Restek"orper endlich ist.
Damit ist auch $\hat{R}/{\bf m}_{\hat{R}}^n$ ein endlicher
$\hat{A}/{\bf m}_{\hat{A}}-$Modul
f"ur alle $n \in {\bf N}$ und wegen
$({\bf m}_{\hat{R}})^n \subseteq {\bf a}\hat{R}$ f"ur ein $n$
ist auch $\hat{R}/{\bf a}\hat{R}$ ein endlicher
$\hat{A}/{\bf m}_{\hat{A}}-$Modul.
Nach \cite{nagloc}, Theorem 30.6, ist dann schon
$\hat{R}$ endlich "uber $\hat{A}$.
\par\smallskip\noindent
Das Bild von ${\rm Spek}\, \hat{R}$ ist dann eine abgeschlossene Teilmenge
von ${\rm Spek}\, \hat{A}$, die zweidimensional sein muss,
da sie andernfalls bereits
endlich w"are, unter ${\rm Spek}\, R \longrightarrow {\rm Spek}\,A$ aber
alle Primideale aus $D({\bf a})$ erreicht werden.
Da $A$ exzellent und normal ist,
ist $\hat{A}$ irreduzibel. Also ist die Spektrumsabbildung surjektiv.
Da das Urbild von $V({\bf a}\hat{A})$ nur aus dem abgeschlossenen Punkt
besteht, muss auch $V({\bf a}\hat{A})$ und damit auch
$V({\bf a}A)$ einpunktig sein. \hfill $\Box$
\begin{kor}
Sei $A$ ein integrer normaler exzellenter Ring, $Y$ eine Hyperfl"ache
in $X={\rm Spek}\, A$ und $U=X-Y$.
Dann besitzt ${\rm Naf}\,(U,X)$ eine Kodimension $\geq 3$.
\end{kor}
{\it Beweis}.
Nach dem Satz ist $U$ in allen Punkten der H"ohe $\geq 2$ affin,
so dass die Aussage aus 1.8.5. folgt.    \hfill $\Box$

%% file: proka.tex
\section{Affine Teilmengen in projektiven Schemata und im affinen Kegel}
\par
\bigskip
\noindent
In diesem Kapitel untersuchen wir Affinit"atseigenschaften von offenen Mengen
in einer projektiven Variet"at und damit zusammenh"angend im Spektrum eines
graduierten Ringes.
Es ergibt sich, dass eine offene Teilmenge
der Variet"at genau dann affin ist, wenn das f"ur ihr Urbild unter der
Kegelabbildung gilt. Dies wiederum kann im lokalen Ring der Kegelspitze
getestet werden. Um die Korrespondenz zwischen graduierten Ringen
und projektiven Variet"aten lokaler fassen zu k"onnen,
empfiehlt es sich, auch nicht notwendig positive Graduierungen zuzulassen.
Da Graduierungen und Gruppenoperationen eng zusammenh"angen, werden
auch Affinit"atseigenschaften von Quotienten unter
Gruppenoperationen hier betrachtet.
\par\smallskip\noindent
Zum Schluss wird gefragt, wann die Aufblasung eines affinen Schemas nach einem
Ideal affin bleibt, und gezeigt, dass man mit affinen Aufblasungen
normalisieren kann.
\par
\bigskip
\noindent
\subsection{Gruppenoperationen und Quotientenbildung}
\par
\bigskip
\noindent
Sei $G$ ein Gruppenschema "uber $S$, das auf einem Schema $X$ "uber $S$
operiert, f"ur diese Begriffe siehe etwa \cite{SGA3}.
Es ist im Allgemeinen nicht m"oglich,
einen Quotienten in der Kategorie der Schemata f"ur die Operation anzugeben,
und es ist keineswegs so, dass ein Quotient eines affinen Schemas wieder
affin sein muss. Es gibt viele verschiedene Begriffe von Quotienten zu
einer Gruppenoperation (und zu allgemeineren Situationen),
wir untersuchen hier nur den kategoriellen Quotienten in der Kategorie aller
beringten R"aume.
\par\smallskip\noindent
In der Kategorie der beringten R"aume
ist der Quotient unproblematisch als Kokernobjekt zu bilden, man identifiziert
Punkte $x,x'$ aus $X$, wenn sie unter den Abbildungen
$p_2, m: G \times_S X \longrightarrow X$ von einem Urbildpunkt erreicht
werden, versieht den Restklassenraum $Y$ mit der 
Quotiententopologie und nimmt den Equalizer von
$p_2^*, m^*:{\cal O}_X \longrightarrow {\cal O}_{G \times X}$
als Strukturgarbe auf den invarianten Mengen,
also ${\cal O}_Y(V)=\{ f \in {\cal O}_X(p_2^{-1}(V)):\, p_2^*(f)=m^*(f) \}$,
siehe hierzu \cite{SGA3}, Exp. 5.1, sowie generell Exp. 5,
wo der Quotient in der Sprache der Gruppoide untersucht wird.
\par\smallskip\noindent
Ist $p:X \longrightarrow Y$ die Quotientenabbildung, so ist f"ur eine offene
Menge $V \subseteq Y$ das Urbild $U=p^{-1}(V)$ invariant, d.h. man kann die
Gruppenoperation auf $G \times U \longrightarrow U$ einschr"anken. Dabei ist
dann $V$ der Quotient von $U$ unter der Operation. Der Quotient in der
Kategorie der beringten R"aume ist also ein durch lokale Daten gegebenes
Objekt.
\par\smallskip\noindent
In dieser Arbeit betrachten wir ausschlie"slich den Quotienten in der
Kategorie aller beringten R"aume.
Wenn wir bei einer Operation eines Gruppenschemas
auf einem Schema sagen, dass der Quotient existiert, so ist
damit gemeint, dass dieser Quotient ein Schema ist.
\par\bigskip
\noindent
Wir erw"ahnen hier einige Situationen von Operationen
von affinen Gruppen auf Schemata derart, dass der Quotient
existiert, und insbesondere solche, wo der Quotient eines affinen Schemas
wieder affin ist.
Aussagen von diesem Typ sind wiederum die
Grundlage von allgemeineren Aussagen zur Existenz von Quotienten.
Besitzt eine solche Operation n"amlich eine "Uberdeckung aus
affinen invarianten Teilmengen, so ist der Quotient "uberdeckt von
den einzelnen affinen Quotienten und daher ein Schema.
\begin{satz}
Sei $G$ eine endliche Gruppe "uber einem K"orper $K$, die auf einem
affinen $K-$Schema $X={\rm Spek}\, A$ operiert. Dann ist das Spektrum des
Invariantenringes ${\rm Spek}\, A^G$ der Quotient.
\end{satz}
Zu einem Beweis siehe etwa \cite{mumab}, III.12, Theorem 1,
oder \cite{SGA3}, Exp. 5.4, Theorem 4.1. \hfill $\Box$
\par\bigskip
\noindent
Operiert $G$ auf einem separierten $K-$Schema, so findet man eine
affine "Uberdeckung aus invarianten Teilmengen, und daher existiert
dann ebenfalls der Quotient.
\begin{satz}
Sei $X$ eine affine algebraische Gruppe "uber einem K"orper $K$
und $G$ eine normale algebraische Untergruppe von $X$.
Dann existiert der Quotient und ist ebenfalls eine affine algebraische
Gruppe.
\end{satz}
Zum Beweis siehe \cite{demazure}, III.3.5.6 oder
\cite{springer}, 5.2.5. \hfill $\Box$
\par\bigskip
\noindent
{\bf Definition}
Ein {\it Hauptfaserb"undel} zum Gruppenschema $G$ "uber $S$
ist ein $S-$ Morphismus $p:X \longrightarrow Y$ und eine Gruppenoperation
von $G$ auf $X$ derart,
dass es eine "Uberdeckung $V_i,\, i\in I,$ von $Y$ gibt so,
dass $p^{-1}(V_i) \cong G \times V_i$ ist und bei dieser Identifizierung
$G$ durch Gruppenmultiplikation auf der ersten Komponente wirkt.
\begin{satz}
Sei $p:X \longrightarrow Y$ ein Hauptfaserb"undel zum Gruppenschema $G$.
Dann ist $Y$ der Quotient der Gruppenoperation.
\end{satz}
{\it Beweis}.
Da der Quotient durch lokale Daten gegeben ist, kann man
annehmen, dass die Operation
trivial ist, sei also $X=G \times Y$.
Die Operation von $G$ auf $G \times Y$ entspricht der Gruppenmultiplikation
der Gruppe $G \times Y=G_Y$ mit dem neuen Basisschema $Y$, und man hat zu
zeigen, dass bei einem Gruppenschema $G$ "uber $Y$ die Basis der Quotient ist.
Da $p:G \longrightarrow Y$ das neutrale Element als Schnitt besitzt, ist die
Abbildung surjektiv, und $Y$ tr"agt die Quotiententopologie.
\par\smallskip\noindent
Sind $P$ und $Q \in G$ Punkte, die "uber einem Punkt $R \in Y$ liegen,
so findet man wegen $\kappa(R) \longrightarrow \kappa(P), \kappa(Q)$ einen
gemeinsamen Oberk"orper $L$ und damit einen Punkt aus $G \times G$,
der unter den beiden Projektionen $P$ und $Q$ trifft.
Andererseits ist
$(m,p_2):G \times G \longrightarrow G \times G$ eine Isomorphie,
daher gibt es auch 
einen Punkt aus $G \times G$, der unter der Multiplikation auf $P$ und unter
der Projektion auf $Q$ geht.
\par\smallskip\noindent
Wir haben noch zu zeigen, dass die Strukturgarbe auf $Y$ die Eigenschaft
der Strukturgarbe des Quotienten besitzt,
d.h. es ist zu zeigen, dass die Untergarbe
$p^*{\cal O}_Y \hookrightarrow {\cal O}_G$ der Equalizer ist.
Eine Untergarbe liegt vor, da das neutrale Element $e:Y \longrightarrow G$
zu einem Linksinversen $e^*$ zu $p^*$ f"uhrt.
Die Abbildung $G \longrightarrow G \times G \longrightarrow G$
mit ${\rm id}_G \times e$ und $m$ ist insgesamt die Identit"at.
Ist $m^*(f)=p_2^*(f)$ f"ur einen Schnitt $f$ "uber einer offenen Menge
$U=p^{-1}(V),\, V \subseteq Y,$ von $G$,
so ist insgesamt
$f=(id \times e^*)m^*(f)=(id \times e^*)p_2^*(f)= e^*(f)$,
es ist also $f$
ein Element aus $p^*{\cal O}_Y$. \hfill $\Box$
\par\bigskip
\noindent
Die Hauptfaserb"undel zur multiplikativen Gruppe ${\bf A}_K^\times$ werden wir
ab dem n"achsten Abschnitt studieren.
Dort wird allgemeiner f"ur diagonalisierbare Gruppen gezeigt, dass
ein affiner Totalraum zu einem affinen Quotienten f"uhrt, falls der Quotient
"uberhaupt existiert.
F"ur beliebige Hauptfaserb"undel muss das nicht gelten, wie das folgende
Beispiel zeigt.
\par\bigskip
\noindent
{\bf Beispiel}
Wir zeigen, dass bei einem Hauptfaserb"undel zu einem affinen Gruppenschema
mit affinem Totalraum der Quotient nicht affin sein
muss.\footnote{Dieses Beispiel zeigt zugleich, dass die
in \cite{haramp}, III. problem 3.14. ge"au"serte Vermutung, dass bei
einem lokal
trivialen Faserb"undel die kohomologische Dimension des Totalraumes gleich
der Summe der kohomologischen Dimensionen der Faser und des Quotienten ist,
nicht richtig ist.}
\par\smallskip\noindent
Sei $K$ ein kommutativer Ring und
$R= K[a,b,c,d] / (ad-cb-1),\, X={\rm Spek}\, R$.
Wir versehen $X$ mit der Gruppenstruktur ("uber ${\rm Spek }\, K$), die der
Gruppe $SL_2(K)$ entspricht, also der Matrizenmultiplikation
auf den Matrizen
$$\left(\matrix{ a & b \cr  c & d \cr}\right)$$
mit Determinante $ad - cb =1$.
Auf $SL_2(K)$ operiert
von links die Untergruppe G der oberen Dreiecksmatrizen
$$ \left(\matrix{ x & y \cr  0  & z \cr} \right) $$
mit $xz=1$ durch
$$\left(\matrix{ x &y \cr 0&z \cr}\right)  \cdot 
\left(\matrix{ a&b\cr c&d \cr}\right) =
\left(\matrix{ xa+yc&xb+yd\cr zc&zd \cr}\right) \, .$$
Die Abbildung
$$ p: X \longrightarrow {\bf P}^1_K \mbox{ mit }
\left(\matrix{ a & b \cr  c & d \cr}\right) \longmapsto  (c:d) $$
ist ein Hauptfaserb"undel.
Zu $D_+(d)= {\bf A}^1_K \subseteq {\bf P}^1_K$ ist
$$G \times {\bf A}^1_K  \longrightarrow D(d) \subseteq X \mbox{ mit }
\left(\matrix{ x & y \cr  0  & x^{-1} \cr} \right), t \, \longmapsto
\left(\matrix{x+yt & y \cr x^{-1}t & x^{-1} \cr}\right) $$
eine Trivialisierung der Operation, wie man schnell ausrechnet,
und entsprechend gibt es eine
Trivialisierung auf $D_+(c)$.
Die operierende Gruppe und der Totalraum ist affin,
der Quotient ${\bf P}^1_K$ nicht.
\par
\bigskip
\noindent
\subsection{Graduierte Algebren und diagonalisierbare Gruppen}
\par
\bigskip
\noindent
Sei $D$ eine (additiv geschriebene) kommutative Gruppe und $K$ ein
kommutativer Ring. Der {\it Gruppenring} $K[D]$ ist der
freie $K-$Modul mit den Elementen aus $D$ als Basis, versehen mit
der Algebra-Struktur, die die Verkn"upfung auf $D$ fortsetzt.
F"ur ein Element 
$d \in D$ schreiben wir als Element im Gruppenring $T^d$. Es gilt dann
$T^d \cdot T^e = T^{d+e}$.
\par\smallskip\noindent
$G={\rm Spek}\, K[D]$ ist ein Gruppenschema mit der Komultiplikation
$$ c: K[D] \longrightarrow K[D] \otimes_K K[D] = K[D \times D],\,
T^d \longmapsto T^d \otimes T^d \, ,$$
der Koinversenbildung 
$$ inv^*: T^d \longmapsto T^{-d}$$
und dem neutralen Element 
$$ e^*: T^d \longmapsto 1 \, .$$
Man verifiziert unmittelbar, dass durch diese Abbildungen ein affines
Gruppenschema ${\rm Spek}\, K[D]$
"uber ${\rm Spek}\, K$ gegeben ist.
Gruppenschemata von diesem Typ werden {\it diagonalisierbar} genannt, diese
Bezeichnungsweise findet sich in \cite{wat}, 4.6, begr"undet.
Zu diagonalisierbaren
Gruppen siehe ansonsten \cite{SGA3}, Exp.1.4.4,
\par\smallskip\noindent
Sei nun $A= \sum_{d \in D}A_d$ eine durch $D$
graduierte Algebra und $K \subseteq A_0$ sei ein Unterring (etwa $A_0$
selbst).
Alle homogenen Stufen $A_d$ sind dann insbesondere $K-$Moduln.
Der Zusammenhang zwischen einer Operation des Gruppenschemas
${\rm Spek}\, K[D]$ auf ${\rm Spek}\, A$ und einer $D-$Graduierung auf $A$ wird durch folgenden
Satz beschrieben, vergl. auch \cite{SGA3}, Exp.I, 4.7.
\begin{satz}
Sei $K$ ein kommutativer Ring, $A$ eine $K-$Algebra, $X={\rm Spek}\, A$,
$D$ eine kommutative Gruppe und $G={\rm Spek}\, K[D]$ das Gruppenschema.
Dann gelten folgende Aussagen.
\par\smallskip\noindent
{\rm (1)}
Folgende Strukturen sind "aquivalent.
\par\smallskip\noindent
{\rm (i)} $A$ ist ein $D-$graduierter Ring mit $K \subseteq A_0 $
\par\smallskip\noindent
{\rm (ii)} Es liegt eine Operation von ${\rm Spek}\, K[D]$
auf ${\rm Spek}\, A$ vor.
\par\smallskip\noindent
{\rm (2)}
Eine offene Teilmenge $D({\bf a})$ ist genau dann invariant unter der
Operation, wenn es ein homogenes Ideal $\bf b$ gibt mit
$D({\bf a})=D({\bf b})$.
\par\smallskip\noindent
{\rm (3)}
Die Operation ist genau dann frei, wenn in $A$ jede Stufe das Einheitsideal
erzeugt.
\end{satz}
{\it Beweis}.
Zu (1).
Von (i) nach (ii). Die Graduierung liefert in kanonischer Weise
eine Komultiplikation $c: A \longrightarrow A[D] = K[D] \otimes_K A $, wobei
ein homogenes Element $a \in A$ vom Grad $d \in D$ auf $a T^d $ abgebildet
wird. $c$ ist damit homogen vom Grad null, wenn man $A[D]$ mit der
$D-$Graduierung versieht und die Graduierung auf $A$ au"ser Acht l"asst.
Wir haben zu zeigen, dass diese Komultiplikation eine Gruppenoperation
von $G$ auf $X={\rm Spek}\, A$ ist.
\par\smallskip
\noindent
Die Gesamtabbildung
$$A \stackrel{c}{\longrightarrow} A[D] \cong A \otimes_K K[D]
\stackrel{{\rm id} \otimes e^*}{\longrightarrow} A$$
ist offenbar die Identit"at, und dieser Abbildung enspricht auf den Spektra
die Abbildung
$$X \cong {\rm Spek}\, K  \times X
\stackrel{e \times id}{\longrightarrow} G \times X
\stackrel{m}{\longrightarrow} X \, .$$
Damit ist das neutrale Element mit der Operation vertr"aglich.
\par\smallskip\noindent
Die Abbildungen ${\rm id}_{K[D]} \otimes c$ und
$c_G \otimes {\rm id}_A$ von $A[D]=K[D]\otimes_KA$ nach
$K[D] \otimes_K K[D] \otimes_K A$ sind auf dem Bild von
$c:A \longrightarrow A[D]$ identisch. Das liefert
direkt die Vertauschbarkeit von Operation und Gruppenmultiplikation.
\par\smallskip\noindent
Von (ii) nach (i).
Sei nun umgekehrt eine Operation von ${\rm Spek}\, K[D]$
auf $X$ gegeben, es liegt als eine Gruppenkomultiplikation
$c: A \longrightarrow A[D]$ vor, und wir haben zu zeigen, dass dies zu einer
Graduierung auf $A$ f"uhrt, so dass $K$ im Unterring der nullten Stufe liegt.
Wir definieren die d-te Stufe von A 
als das Urbild der d-ten Stufe von A[D] unter $c$.
Da $c$ ein $K-$Algebra-Morphismus ist, geht ein Element aus $K$ auf
ein Element vom Grad null, also ist $K \subseteq A_0$. Es gilt ferner
$A_d \cdot A_e \subseteq A_{d+e}$.
Da $c$ ein Rechtsinverses besitzt, ist $c$ injektiv, und damit ist die Summe
dieser Stufen direkt, wir haben
noch zu zeigen, dass A die Summe der $A_d $ ist. Sei hierzu
$a \in A$ gegeben mit $ c(a) = a_1 T^{d_1} + ... + a_r T^{d_r} $.
In A gilt dann wegen dem ersten Operationsaxiom $a = a_1 +...+a_r$,
und es ist zu zeigen, dass die $a_i$ homogen sind vom Grad $d_i$.
Hierzu zieht man in $A[D][D] $ (wobei wir T und S schreiben) die Gleichheit
\begin{eqnarray*}
c(a_1)T^{d_1} +...+c(a_r)T^{d_r} & =& ({\rm id}_{K[D]} \otimes c) (c(a)) \\
&= &(c_G \otimes {\rm id}_A)(c(a))\\
&= &a_1 T^{d_1} S^{d_1} + ...+ a_r T^{d_r}S^{d_r}
\end{eqnarray*}
heran, aus der $c(a_i)=a_iS^{d_i}$ hervorgeht.
\par\smallskip\noindent
Offenbar liefert die Hintereinanderausf"uhrung der Konstruktionen wieder die
Ausgangsobjekte.
\par\smallskip\noindent
Zu (2).
Die Invarianz einer Teilmenge $U=D({\bf a})$ bedeutet, dass unter der
Gruppenoperation $m:G \times X\longrightarrow X$ gilt
$m^{-1}(U)= G \times U$.
Dies bedeutet auf der Ringebene, dass die Erweiterungsideale von $\bf a$
in $A[D]$ 
unter der Komultiplikation und unter der Injektion das gleiche Radikal
definieren.
\par\smallskip\noindent
F"ur ein homogenes Element $h$ vom Grad $d$ ist $i(h)=h$ und $c(h)=hT^d$.
Da $T^d$ eine Einheit ist, gilt f"ur ein homogenes Ideal $\bf a$ die
Gleichheit $i({\bf a})=c({\bf a})$.
Wird ein Ideal $\bf a$ als Radikal durch ein homogenes Ideal beschrieben,
so gilt diese Gleichheit immernoch im Radikalsinn.
\par\smallskip\noindent
Umgekehrt beschreibe $\bf a$ eine invariante Teilmenge, also
$c({\bf a})=i({\bf a})$ (als Radikal).
Sei ${\bf b} \subseteq {\bf a}$ das durch die homogenen Elemente von $\bf a$
erzeugte Ideal.
Sei $f \in {\bf a}$. Dann ist
$c(f^k) \in i({\bf a})$ f"ur ein $k \in {\bf N}$, und mit
der homogenen Zerlegung
$f^k=f_{d_1}+...+f_{d_r} \in {\bf a}$ ist also
$f_{d_1}T^{d_1}+....+f_{d_r}T^{d_r} \in i({\bf a})$.
Damit ist $f_{d_j} \in {\bf a}$ f"ur $j=1,...,r$
und damit ist $f^k \in {\bf b}$.
\par\smallskip\noindent
Zu (3).
Dass die Operation frei ist bedeutet, dass der Morphismus
$m \times p_2: G \times X \longrightarrow X \times X$
eine abgeschlossene Einbettung ist.
\footnote{Mit dieser Terminologie halten wir uns an
\cite{GIT}, Chapter 0,3, Def. 0.8, w"ahrend in \cite{SGA3} unter
frei verstanden wird, dass diese Abbildung ein Monomorphismus ist.}
In unserem Kontext bedeutet das die Surjektivit"at der Abbildung
$A \otimes_K A \longrightarrow A[D],\, b \otimes a_d \longmapsto ba_dT^d$,
die wiederum gesichert ist,
wenn alle Unbestimmten $T^d,\, d\in D,$ erreicht werden.
Das Bild der homogenen Abbildung ist aber gerade
gleich $(A_d ) T^d$.      \hfill   $\Box$     
\par\bigskip
\noindent
{\bf Beispiel}
Bei einer ${\bf Z}-$Graduierung ist das Radikal eines homogenen Ideals
wieder homogen, so dass in diesem Fall eine invariante Menge durch ein
homogenes Radikal beschrieben wird.
Im allgemeinen ist aber das Radikal einer invarianten Teilmenge nicht selbst
homogen.
Als Beispiel sei $K$ ein K"orper der Charakteristik zwei, $D={\bf Z}/2$
und $A=K[D]$. Dann ist $(1+T)$ nilpotent, die homogenen Bestandteile aber
nicht. 
\par\bigskip
\noindent
Wir wollen wissen, wann zu der durch eine Graduierung gegebenen Operation 
der Quotient existiert und wann er affin ist.
In \cite{SGA3}, Tome 2, Exp. VIII.5,
Th\'{e}or\`{e}me 5.1 wird gezeigt, dass bei einer freien Operation
der Quotient durch ${\rm Spek}\, A_0$ gegeben ist.
\footnote{Frei im Sinne von Monomorphismus, in dem dortigen
Beweis wird aber zugleich
gezeigt, dass Freiheit im Sinn von abgeschlossener Einbettung vorliegt.}
Aber schon bei Graduierungen, die zu den gewichteten projektiven R"aumen
f"uhren, liegt keine freie Operation vor, weshalb diese Bedingung zu
restriktiv ist.
\par\smallskip\noindent
Ist $X \longrightarrow Y$ ein Quotient, so ist jede invariante
offene Teilmenge $U \subseteq X$ Urbild einer offenen Menge $V \subseteq Y$.
Damit also das Spektrum des Invariantenringes $A^G$ der Quotient ist, muss
jede offene invariante Teilmenge durch invariante Funktionen beschrieben
werden k"onnen.  
Im graduierten Fall werden die invarianten Teilmengen durch homogene
Ideale beschrieben, und diese notwendige Bedingung besagt dann, dass
jedes homogene Ideal als Radikal durch Elemente der nullten Stufe
beschrieben  werden kann.
\par\smallskip\noindent
Ist die Graduierung sogar frei, so erzeugt nach 1.2.1
jede Stufe $A_d$ das Einheitsideal, sagen wir
$1=\sum_{i=1,...,n}a_ib_i$ mit $a_i \in A_d,\, b_i \in A_{-d}$.
Dann ist f"ur homogenes $f \in A_d$ aber $f= \sum_{i=1,...,n} a_i (fb_i)$ mit
$fb_i \in (f) \cap A_0 $, also $f \in ((f) \cap A_0)$, und diese notwendige
Bedingung ist erf"ullt. Sie ist ebenso erf"ullt, wenn nicht jede Stufe,
sondern nur zu jedem Grad $d$ eine Stufe $A_{dn},\, n \geq 1,$ die Eins erzeugt.
Dies wiederum ist im Spezialfall $D={\bf Z}$ bereits dann erf"ullt, wenn
irgendeine Stufe $A_d,\, d \neq 0 ,$ das Einheitsideal erzeugt, und insbesondere,
wenn es eine homogene Einheit vom Grad $\neq 0$ gibt.
Der folgende Satz zeigt, dass diese notwendige Bedingung bereits
hinreichend ist f"ur die Existenz des Quotienten,
und daher in all den soeben beschriebenen Situationen
das Spektrum des Invariantenringes $A_0$ den Quotienten liefert.
\begin{satz}
Sei $A=\sum_{d \in D}A_d$ ein $D-$graduierter kommutativer
Ring mit der Eigenschaft, dass jedes homogene Radikal als Radikal
von der nullten Stufe erzeugt wird.
Dann ist ${\rm Spek}\, A_0$ der Quotient der
zugeh"origen Gruppenoperation in der Kategorie aller lokal beringten R"aume.
\end{satz}
{\it Beweis}.
Wir haben zu zeigen, dass ${\rm Spek}\, A_0$ die Eigenschaften des
Quotienten erf"ullt.
Die Abbildung ${\rm Spek}\, A \longrightarrow {\rm Spek}\, A_0$
ist surjektiv, da $A_0$ ein direkter Summand von $A$ ist.
Wir haben f"ur zwei Primideale
${\bf p},{\bf q} \in {\rm Spek}\, A=X$ mit ${\bf p}_0={\bf q}_0$ die
Existenz eines Punktes aus ${\rm Spek}\, A[D]$ nachzuweisen,
der unter $m^*$ auf den einen und unter $p_2^*$ auf den anderen abgebildet
wird.
\par\smallskip\noindent
Nach Voraussetzung enthalten dann $\bf p$ und $\bf q$ "uberhaupt die gleichen
homogenen Elemente.
Wir machen die homogenen Elemente aus $\bf p$ zu null und die au"serhalb
von $\bf p$ zu Einheiten, und erhalten einen graduierten Ring, in dem
alle homogenen Elemente bis auf die Null Einheiten sind. Es ist dann
insbesondere $A_0=K$ ein K"orper.
Wir m"ussen zeigen, dass die Fasern zur Projektion "uber dem einen Punkt
und die Faser zur Gruppenoperation zum anderen Punkt
nichtleeren Durchschnitt haben, oder anders gesagt,
dass die Operation auf einer Projektionsfaser surjektiv ist.
Ist $L=k({\bf p})$ mit $\varphi:A \longrightarrow L$,
so ist die Faser "uber der Projektion
(die als Ringabbildung die Einbettung ist) einfach gleich
dem Spektrum zu $A[D] \otimes_A L=L[D]$,
und wir haben die Surjektivit"at von
$${\rm Spek}\, L[D] \hookrightarrow {\rm Spek}\, A[D]
\stackrel{m}{ \longrightarrow} {\rm Spek}\, A$$
zu zeigen, die auf der Ringebene durch
$$ A \stackrel{m^*}{\longrightarrow} A[D]
\longrightarrow L[D] \mbox{ mit }
s_d \longmapsto s_dX^d \longmapsto \varphi(s_d)X^d $$
gegeben ist.
Wir behaupten, dass diese Abbildung injektiv ist und ein direkter
Summand vorliegt, woraus die Surjektivit"at der Spektrumsabbildung folgt.
Es liegt ein homogener Ringhomomorphismus vom Grad null vor, sei $s_d$ ein
homogenes Element, dass auf die Null abgebildet wird. Dann ist
$\varphi(s_d) =0$ und $s_d$ kann keine Einheit in $A$ sein, also muss es
aufgrund der Voraussetzung schon gleich null sein.
Also ist $A \subseteq L[D]$ ein $K-$Untervektorraum. 
Sei $V$ der $K-$Untervektorraum von $L[D]$, der von allen homogenen
Elementen erzeugt wird, die nicht zu $A$ geh"oren. Dann ist
$L[D]= A + V$ eine $K-$Vektorraumzerlegung, und $V$ ist ein $A-$Modul.
Dazu hat man $A \cdot V \subseteq V$ zu zeigen, und kann sich auf
homogene Elemente beschr"anken. Sei also $0 \neq a_d \in A_d$ homogen
und $\lambda X^e \in V,\, \lambda \in L \, .$ Dann kann auch
$(a_d) \lambda X^e = \varphi(a_d)X^d \cdot \lambda X^e $ nicht zu $A$
geh"oren, da $\varphi(a_d)X^d$ eine Einheit in $A$ ist.
\par\smallskip\noindent
Zur Quotiententopologie.
Sei $T \subseteq {\rm Spek}\, A_0$ eine Teilmenge, deren Urbild offen
ist, also $p^{-1}(T)=D({\bf a})$. Da dieses Urbild invariant ist,
ist $\bf a$ homogen und damit ist nach Voraussetzung
$p^{-1}(T)=D({\bf a}_0)$. Damit ist
$T =p(D({\bf a}_0))=p(p^{-1}(D({\bf a}_0)))=D({\bf a}_0)$.
\par\smallskip\noindent
F"ur $f \in A_0$ ist $(A_0)_f = (A_f)_0$, und das sind die
invarianten Funktionen auf
$D(f) \subseteq {\rm Spek}\, A$, also die Funktionen, die
unter der Komultiplikation und unter der Standardeinbettung nach $A_f[D]$
gleiches Bild haben.     \hfill   $\Box$
\par
\bigskip
\noindent
\subsection{Das projektive Spektrum ${\bf Z}-$graduierter Algebren}
\par
\bigskip
\noindent
Ist $A$ ein $D-$graduierter Ring, so gilt die Voraussetzung des Satzes 2.2.2
im Allgemeinen nicht, und die Gruppenoperation hat in der Regel keinen
Quotienten. Es ist auch schwierig, eine invariante offene Teilmenge
$U \subseteq {\rm Spek}\, A$ anzugeben, auf der die Operation hinreichend gut
ist, um ein Schema als Quotienten zu erhalten. 
Man kann nat"urlich einfach den Quotienten in der Kategorie
der beringten R"aume betrachten und dort die Vereinigung aller offenen affinen
Schemata nehmen, was ein Schema ergibt. Da ist es jedoch schwierig, diese
Urbildmenge explizit zu beschreiben.
\par\smallskip\noindent
Einfacher ist die Situation bei $D={\bf Z}$.
Auch hier, man denke etwa an den positiv graduierten Fall, gibt es
f"ur die Gesamtoperation keinen Quotienten,
daf"ur aber gibt es eine gro"se invariante Teilmenge, den punktierten Kegel,
auf der die Operation einen Quotienten besitzt.
\begin{satz}
Sei $A$ ein ${\bf Z}-$graduierter Ring, in dem die homogenen Elemente
vom Grad ungleich null das Einheitsideal erzeugen.
Dann ist ${\rm Spek}\, A_0$ der affine Quotient.
\end{satz}
{\it Beweis}.
Sei $1 =\sum a_i h_i $ mit $h_i$ homogen vom Grad $d_i \neq 0$.
Wir k"onnen sofort annehmen, dass alle $a_i$ ebenfalls homogen sind.
Da $1$ zur nullten Stufe geh"ort, erzeugen bereits alle Summanden vom
Grad $0$ die Eins, somit k"onnen wir $a_i$ vom Grad $-d_i$ annehmen.
Durch eventuelle Vertauschung von $a_i$ und $h_i$ kann man alle
$d_i$ als positiv (oder alle als negativ) annehmen.
Sei jetzt $h \in A_d$ homogen von negativen Grad.
Dann gibt es positive nat"urliche Zahlen $m,m_i$ mit $-dm =d_im_i$.
Es erzeugen auch die $h_i^{m_i}$ das Einheitsideal,
sagen wir $\sum_{i=1,...,n} h_i^{m_i}b_i=1$,
$h^m h_i^{m_i}$ ist vom Grad null, und es ist 
$h^m= \sum_{i=1,...,n}h^mh_i^{m_i}b_i $. Also werden die homogenen
Radikale durch die neutrale Stufe gegeben und daher folgt aus Satz 2.2.2
die Behauptung.     \hfill   $\Box$
\par\bigskip
\noindent
Die Voraussetzung des Satzes ist insbesondere erf"ullt, wenn es in irgendeiner
Stufe von Grad $\neq 0$ eine Einheit gibt.
Diese Eigenschaft besitzt insbesondere zu einem ${\bf Z}-$graduierten Ring
$A$ jede Nenneraufnahme $A_f$ zu einem homogenen Element $f$ vom
Grad $d \neq 0$.
\begin{kor}
Sei $A$ ein ${\bf Z}-$graduierter Ring und $H$ sei die Menge aller
homogenen Elemente von nicht neutralem Grad.
Dann besitzt die Operation auf $U=D(H)$ ein Schema als Quotient,
das von den affinen Mengen ${\rm Spek}\, (A_h)_0=:D_+(h)$ "uberdeckt wird.
\end{kor}
{\it Beweis}.
Auf einer affinen Teilmenge ${\rm Spek}\, A_h=D(h) \subseteq U$ mit $h \in H$
homogen ist $h$ eine Einheit und daher ist die Voraussetzung des
Satzes 2.3.1 erf"ullt, also existiert der Quotient
und ist gegeben durch ${\rm Spek}\, (A_h)_0$.
Der Quotient von $U$ in der Kategorie der beringten R"aume wird also
durch diese affinen Schemata "uberdeckt und ist daher ein Schema. \hfill $\Box$
\par\bigskip
\noindent
{\bf Definition}
Ist $A$ ein ${\bf Z}-$graduierter Ring und $H$ die Menge aller
homogenen Elemente von nicht neutralem Grad, so
nennt man den Quotienten das {\it projektive Spektrum} von $A$,
${\rm Proj}\, A$.
$({\rm Spek}\, A)^\times=D(H)$ nennt man den {\it punktierten Kegel} und
$p:D(H) \longrightarrow {\rm Proj}\, A$ die Kegelabbildung.
\par\bigskip
\noindent
{\bf Bemerkung}
Im positiven Fall ist das nichts anderes
als das "uber die homogenen Primideale ${\bf p}$ mit
$A_+ \not\subseteq {\bf p}$ definierte ${\rm Proj}\, A$.
Auch im nicht positiv graduierten Fall sind die Punkte des
Quotienten die homogenen Primideale, die $H$ nicht umfassen.
\begin{satz}
Sei $A$ ein ${\bf Z}-$graduierter Ring.
Die Abbildung $p : ({\rm Spek}\, A )^\times \longrightarrow P={\rm Proj}\,A $
ist affin, und eine offene Teilmenge $ V \subseteq P $ genau dann
affin, wenn dies f"ur das Urbild $ p^{-1}(V) $ gilt.
\end{satz}
{\it Beweis}.
Die Affinit"at ergibt sich aus $p^{-1}(D_+(h) ) = D(h) $ affin und da die
$D_+(h)$ zu $h$ homogen vom Grad $\neq 0$ eine affine "Uberdeckung von
${\rm Proj}\, A$ bilden.
\par\smallskip\noindent
Sei umgekehrt $U=p^{-1}(V) $ affin. Dann ist $V$ der Quotient von $U$ unter
der Gruppenoperation. 
Da $U$ affin ist, entspricht dieser Operation die
Graduierung auf $\Gamma(U,{\cal O}_X)$, und die Restriktionsabbildung
$R:A \longrightarrow \Gamma(U,{\cal O}_X)$ ist homogen vom Grad null. 
Es ist $U=D(h_1,...,h_n)$ mit homogenen Elementen vom Grad $\neq 0$, und
diese erzeugen wegen der Affinit"at in $\Gamma(U,{\cal O}_X)$ das
Einheitsideal, es liegt also die Situation des Satzes 2.3.1 vor.
Also ist $V$ affin und isomorph
zu ${\rm Spek}\, \Gamma(U,{\cal O}_X)_0$.    \hfill   $\Box$
\par\bigskip
\noindent
{\bf Bemerkung}. Zu einem graduierten Ring kann das Proj affin sein, ohne dass
im graduierten Ring die homogenen Elemente von nichtneutralem Grad das
Einheitsideal erzeugen. Sie erzeugen dann aber auf dem punktierten Kegel
das Einheitsideal. Ist $A$ positiv graduiert mit affinem projektiven
Spektrum, so liegt auf dem punktierten Kegel eine nicht positive 
Graduierung vor. Ein Beispiel hierf"ur ist $A=K[T]$ mit der
Standardgraduierung. 
\par\bigskip
\noindent
{\bf Bemerkung}. Ein anderer Beweis f"ur die Aussage, dass aus der
Affinit"at von $p^{-1}(V)$ die von $U$ folgt, beruht f"ur den Fall, dass
die positiv graduierte Algebra von der ersten Stufe endlich erzeugt wird,
darauf, dass alle quasikoh"arenten Garben auf einer offenen
Menge $V \subseteq P$ von einer quasikoh"arenten Garbe auf $P$ herkommen
und diese vom Typ $\tilde{M}$ mit einem graduierten $A-$Modul $M$ sind.
Verschwindet dann die Kohomologie auf $M\vert p^{-1}(V)$, so auch die
von $\tilde{M}$ auf $V$, vergleiche \cite{bingener}, 1.8.
Ebenso sieht man, dass die kohomologische Dimension von $V$ und $p^{-1}(V)$
"ubereinstimmt.
\begin{lem}
Sei $A$ ein noetherscher $D-$graduierter Ring und sei
${\bf a}$ ein homogenes Ideal. Dann ist der affine Ort von $U=D({\bf a})$
ebenfalls homogen.
\par\smallskip\noindent
Ist $D={\bf Z}$ und ${\bf a}$ in allen Lokalisierungen zu
homogenen Primidealen affin, so ist $U$ affin. 
\end{lem}
{\it Beweis}.
Nach 1.8.5 wird der affine Ort durch das Ideal
$R^{-1}({\bf a}\Gamma(D({\bf a}),{\cal O}_X))$ beschrieben.
Der Schnittring
von $D({\bf a})$ ist graduiert, die Restriktionsabbildung ist
homogen, also ist das Erweiterungsideal homogen.
Damit ist auch das Urbildideal homogen.
\par\smallskip\noindent
Bei $D={\bf Z}$ gibt es "uber jedem homogenen Ideal $\neq A$ homogene
Primideale. Beschreibt das homogene Ideal den nicht leeren nichtaffinen
Ort, so gibt es also darin homogenen Primideale, an deren Lokalisierung
$U$ nicht affin ist. \hfill    $ \Box$
\begin{kor}
Ist $A$ ein positiv graduierter Ring "uber einem K"orper
$K=A_0$, so sind f"ur ein homogenes Ideal ${\bf a} \subseteq A_+$
folgende Aussagen "aquivalent.
\par\smallskip\noindent
{\rm (1)} Es ist $D_+({\bf a}) \subseteq {\rm Proj}\, A$ affin.
\par\smallskip\noindent
{\rm (2)} Es ist $D({\bf a}) \subseteq {\rm Spek}\, A$ affin.
\par\smallskip\noindent
{\rm (3)} Es ist $D({\bf a})$ affin in der Spitze $A_{(A_+)}$.
\end{kor}
{\it Beweis}. Die "Aquivalenz von (1) und (2) folgt aus Satz 2.3.3,
Da nach Voraussetzung $D({\bf a})$ im punktierten Kegel $D(A_+)$ liegt.
Aus (2) folgt sofort (3),
sei umgekehrt (3) erf"ullt.
Das homogene Primideal $A_+$ umfasst jedes andere homogene
Primideal, daher ist die offene Menge in jeder Lokalisierung zu
einem homogenen Primideal affin und nach dem Lemma 2.3.4 ist dann $D({\bf a})$
selbst affin. \hfill   $\Box$
\par\bigskip
\noindent
{\bf Beispiel}
Die Aussage des Satzes gilt nicht unter schw"acheren Bedingungen.
Ist etwa $A_0$ kein K"orper, so kann man die Affinit"at nicht in der
Spitze testen.
Sei hierzu $A=K[u,v,X,Y]/(uX-vY)$ mit $A_0=K[u,v]$ und $X$ und $Y$ vom Grad
eins. Dann ist $D(A_+)=D(X,Y)$ und ${\rm Proj}\, A$ nicht affin, in der
Spitze $A_{(X,Y)}$ ist aber $D(X,Y)$ affin.
\par\bigskip
\noindent
Wir wollen die Geometrie der projektiven Spektra zu nicht positiv graduierten
Algebren genauer kennenlernen.
Zu dem graduierten Ring $A$ seien $ A^+$  und $ A^- $ die
Unterringe von $A$, die aus den nichtnegativen bzw. den nicht positiven
Stufen bestehen. Das sind dann einseitig graduierte Ringe. Es seien $H_+$ alle
homogenen Elemente mit positivem Grad und $H_-$ alle homogenen Elemente
mit negativem
Grad. Dann ist ${\rm Proj}\, A= D_+(H_+) \cup D_+(H_-)$ und es ist
$D_+(H_+) = {\rm Proj}\, A^{+}$ und $D_+(H_-) ={\rm Proj}\, A^{-}$.
Dies ergibt sich einfach daraus, dass f"ur ein homogenes Element $h$ von
sagen wir positivem Grad $d$ die Ringabbildung
$ ((A^+)_h)_0 \longrightarrow (A_h)_0 $ eine Isomorphie ist.
\par\smallskip\noindent
$H$ erzeugt genau dann das Einheitsideal in $A$, wenn $H_+$ das 
Einheitsideal in $A$ erzeugt (Daraus folgt nat"urlich nicht, dass $H_+$
in $A^+$ das Einheitsideal erzeugt).
In diesem Fall ist ${\rm Proj}\, A={\rm Proj}\, A_+ = {\rm Spek}\, A_0$.
Diese Gleichheit gilt bereits dann,
wenn $H_- \subseteq (H_+)A$ gilt.
Ist dies nicht der Fall, so ist ${\rm Proj}\, A$ nicht separiert.
\begin{satz}
Sei $A$ ein graduierter Ring, $A^+$ und $A^-$ seien
von endlichem Typ "uber $A_0$ mit den projektiven Spektren $U_+$ und $U_-$. 
Ist $P$ separiert und zusammenh"angend, so ist $U_+ ={\rm Proj}\, A$ oder
$U_-={\rm Proj}\, A$.
\end{satz}
{\it Beweis}.
Als projektive Spektren zu positiv graduierten Algebren von endlichem Typ
sind $U_+$ und $U_-$ projektive, insbesondere eigentliche Schemata
"uber ${\rm Spek}\, A_0$.
Ist $P$ separiert, so sind diese beiden Schemata in $P$
auch abgeschlossen.
Damit sind die beiden Mengen leer oder gleich $P$,
und eine muss gleich $P$ sein.     \hfill   $\Box$
\par\bigskip
\noindent
{\bf Beispiel}
Sei $A=K[X_0,...,X_m,Y_0,...,Y_n]$, wobei $X_i$ den Grad 1 und $Y_j$ den
Grad $-1$ habe.
Als Operation von ${\bf A}_K^\times$ auf ${\bf A}_K^{m+n+2}$ entspricht dem
die Abbildung
$$ z,(x_0,...,x_m,y_0,...,y_n)
\longmapsto (zx_0,...,zx_m,z^{-1}y_0,...,z^{-1}y_n) \, .$$
Es ist $A_0=K[X_iY_j: \, i=1,..,m,\, j=1,...,n ]$, das ist der Minorenring
der Untermatritzen vom Rang 2 zu einer $(m+1)\times (n+1)-$Matrix.
Es ist
$A^+=A_0[X_0,...,X_m]$ und $A^-=A_0[Y_0,...,Y_n]$  (jeweils mit Relationen).
${\rm Proj}\, A^+$ ist isomorph zur Aufblasung des Minorenringes $A_0$ an einem
Spaltenideal, etwa an $(X_0Y_0,...,X_mY_0)$.
Macht man in $A^+$ ein Element $X_iY_j$ zu einer Einheit, so ist
$X_r=(Y_jX_r/Y_jX_i) \cdot X_i$ und damit ist $A^+_{X_iY_j}=(A_0)_{X_iY_j}[X_i]$.
Davon ist das Proj einfach $D(X_iY_j) \subseteq {\rm Spek}\, A_0$.
Insgesamt ist oberhalb von $D(X_iY_j: \, i=1,..,m,\, j=1,...,n )$ das 
projektive Spektrum von $A$ gleich dem von $A^+$ (und gleich dem von $A^-$)
und gleich ${\rm Spek}\, A_0$. An dieser offenen Menge werden die beiden Mengen
verklebt.
\par\smallskip\noindent
"Uber dem Nullpunkt liegt dagegen zum einen ein ${\bf P}^m_K$ und zum andern
ein ${\bf P}^n_K$, da bei $X_iY_j=0$ f"ur alle $i,j$ der graduierte Ring
$A_+$ zu einem Polynomring "uber $K$ wird (ohne weitere Relationen).
Diese Schemata sind nicht separiert.
\par\smallskip\noindent
Ist beispielsweise $m=n=0$,
so ist das projektive Spektrum die Verklebung zweier affiner Geraden
an der punktierten Gerade, und damit
ergibt sich die im Nullpunkt verdoppelte, nicht separierte Gerade.
Die Bahnen der Operation sind dabei Hyperbeln und die
beiden Achsen, und die Hyperbeln konvergieren gegen beide Achsen und haben
keinen eindeutigen Limespunkt.
\par\bigskip
\noindent
{\bf Beispiel}
Sei $A=K[X_1,...,X_m,Y_1,...,Y_n]/( X_iY_j,\, i=1,..,m,\, j=1,...,n )$
mit ${\rm grad}\, X_i=1 $ und ${\rm grad}\, Y_j=-1$
Da ist ${\rm Proj}\, A$ einfach die disjunkte Vereinigung
aus ${\bf P}_K^m$ und
${\bf P}_K^n$.
\par\bigskip
\noindent
{\bf Beispiel}
Sei $A=K[X_1,...,X_n,Y_1,...,Y_n]/(\sum_{i=1,...,n}X_iY_i-1)$
mit ${\rm grad}\, X_i=d_i > 0,\, {\rm grad}\, Y_i=-d_i$.
Hier erzeugen die homogenen Elemente vom
Grad ungleich null das Einheitsideal, ohne dass es (bei $n \geq 2$)
homogene Einheiten gibt.
Ist $K$ faktoriell, so ist auch $A$ faktoriell,
da es vom Typ $R[U,V]/UV-p$ mit $R$ faktoriell und
$p \in R$ prim ist.
Insbesondere ist in ${\rm Spek}\, A$ das Komplement jeder Hyperfl"ache affin,
und das gilt
dann auch in ${\rm Proj}\, A={\rm Spek}\, A_0$.
$A_0$ ist nicht immer einfach zu bestimmen, und muss nicht faktoriell sein.
Ist etwa $n=2$ und $d_1=d_2=1$, so ist
\begin{eqnarray*}
A_0 & = &K[X_1Y_1,X_1Y_2,X_2Y_1]/( (X_1Y_1)(1-X_1Y_1) - (X_1Y_2)(X_2Y_1)) \\
 & \cong & K[r,s,t]/r(1-r)-st \, .
\end{eqnarray*}
\par\bigskip
\noindent
{\bf Beispiel}
Mit beliebigen ${\bf Z}-$Graduierungen kann man "ubersichtlich die
ge\-twisteten B"undel auf einer projektiven Variet"at $V$ beschreiben, die durch
eine positiv graduierte $K-$Al\-ge\-bra $A$ gegeben ist, die von der ersten Stufe
endlich erzeugt werde, also $A=K[x_0,...,x_n]$ und $V={\rm Proj}\, A$.
Dann wird f"ur $d \in {\bf Z}$ das Geradenb"undel zur
invertierbaren Garbe ${\cal L}={\cal O}(d)= \widetilde{A(d)}$
beschrieben durch
$$ {\rm Proj}\, A[T] \supseteq D_+(x_0,....,x_n)
\longrightarrow {\rm Proj}\, A=V \, ,$$
wobei die Unbestimmte $T$ den Grad $d$ bekommt.
\par\smallskip\noindent
Die angegebene Abbildung hat lokal die Gestalt
$${\rm Proj}\, A[T]_{x_i}={\rm Spek}\, (A[T]_{x_i})_0
\longrightarrow D_+(x_i)={\rm Spek}\, (A_{x_i})_0  \,  .$$
Andererseits wird das Geradenb"undel, dessen Schnitte die invertierbare
Garbe $\cal L$ ergeben soll, erhalten durch ${\rm Spek}\, S({\cal L}^{-1})$,
wobei $S$ f"ur die symmetrische Algebra steht,
siehe \cite{haralg}, Exc. II.5.18.
Es ist $ S({\cal O}(-d)) = \sum_{k \in {\bf N}} {\cal O}(-kd)$
und damit ist
$$ \Gamma(D_+(x_i), S({\cal O}(-d))) =
\Gamma(D_+(x_i), \sum_{k \in {\bf N}} {\cal O}(-kd) )
= \sum_{k \in {\bf N}} (A_{x_i})_{-kd}
= (A_{x_i}[T])_0 $$
mit ${\rm grad}\, T=d$ wie behauptet.
\par\smallskip\noindent
Damit sieht man etwa, dass f"ur den ${\bf P}_K^n$ das ger"aumige B"undel
${\cal O}(1)$ der
Projektion eines ${\bf P}_K^{n+1}$ von einem Punkt aus entspricht.
Bei ${\rm grad}(T)=-1$ hingegen hat man mit $Y_i=TX_i$ (vom Grad null) die
homogene Ringabbildung
$K[Y_0,...,Y_n,X_0,...,X_n]/(Y_iX_j - Y_jX_i) \longrightarrow
K[X_0,...,X_n,T]$, die einen Morphismus
$$ D_+(X_0,...,X_n) \longrightarrow
{\rm Proj}\, K[Y_0,...,Y_n,X_0,...,X_n]/(Y_iX_j - Y_jX_i) $$ stiftet, den
man schnell als Isomorphie nachweist (Betrachte die Nenneraufnahmen zu $X_i$).
${\cal O}(-1)$ entspricht also als Geradenb"undel die Aufblasung des affinen
Kegels an der Spitze.
\par
\bigskip
\noindent
\subsection{Affine Teilmengen in projektiven Variet"aten}
\par
\bigskip
\noindent
Sei $X$ eine projektive Variet"at "uber einem K"orper $K$
und $U \subseteq X$ eine offene Teilmenge.
Wir besprechen hier einige notwendige und hinreichende Bedingungen f"ur die
Affinit"at von $U$, die in vielen Beispielen zu entscheiden erlauben,
ob Affinit"at vorliegt oder nicht. Diese Kriterien werden im dritten Kapitel
(insbesondere in 3.7) weiter diskutiert und auf ihre Grenzen hin untersucht.
\par\bigskip
\noindent
Eine Vielzahl von affinen Teilmengen wird durch ger"aumige Divisoren
gestiftet.
\begin{satz}
Sei $X$ eine projektive Variet"at "uber einem noetherschen Ring $K$
und $\cal L$ ein ger"aumiges Geradenb"undel. Dann ist zu jedem Schnitt
$0 \neq s \in {\cal L}(X)$ das Komplement der Nullstellenmenge $V(s)$
affin.
\end{satz}
{\it Beweis}. Zum Schnitt $s$ besitzt der Schnitt $s^n \in {\cal L}(X)$
die gleiche Nullstellenmenge, von daher k"onnen wir zu einer beliebigen
positiven Potenz "ubergehen und annehmen, dass ${\cal L}$ sehr ger"aumig ist,
siehe \cite{haralg}, Theorem II.7.6.
Dann besitzen die globalen Schnitte ein
endliches (\cite{haralg}, Theorem III.5.2),
basispunktfreies $K-$Erzeugendensystem, das eine abgeschlossene
Einbettung in einen ${\bf P}_K^n$ definiert, und dann ist $D(s)$
Urbild einer affinen Teilmenge unter einer affinen
Abbildung, also affin. \hfill $\Box$
\par\bigskip
\noindent
{\bf Bemerkung}
Die Aussage gilt nicht ohne die Voraussetzung, dass $X$ eigentlich ist.
In einem
quasiaffinen, nicht affinen Schema sind alle invertierbaren Garben ger"aumig,
aber keineswegs besitzt da zu einem Schnitt die Nullstelle affines
Komplement.
\par\bigskip
\noindent
Eine echte offene Teilmenge $U$ einer projektiven Kurve ist immer affin,
und ein echt effektiver Divisor darauf ist immer ger"aumig.
Die Hauptschwierigkeit bei dieser Aussage liegt darin, nicht konstante
Funktionen auf $U$ zu finden, was der Satz von Riemann-Roch sichert.
Auf diese eingebettete Situation kann man den Fall eines separierten
Schemas ohne projektive Komponenten zur"uckf"uhren.
\begin{satz}
Sei $C$ ein eindimensionales separiertes Schema von endlichem Typ "uber
einem K"orper $K$, wobei keine eindimensionale Komponente
eigentlich sei. Dann ist $C$ affin.
\end{satz}
{\it Beweis}.
Zun"achts bemerken wir, dass wenn $C' \longrightarrow C$ surjektiv ist
und $C' \longrightarrow {\rm Spek}\, K$ universell abgeschlossen, dies
auch f"ur $C$ "uber $K$ gilt. F"ur ein beliebiges $K-$Schema $X$ ist dann
bei
$$C' \times_K X \longrightarrow C \times_K X \longrightarrow X$$
die erste Abbildung surjektiv und die Gesamtabbildung abgeschlossen.
Dann ist aber auch die zweite Abbildung abgeschlossen.
\par\smallskip\noindent
Sei $C$ wie im Satz beschrieben.
Aufgrund von Satz 1.13.2 k"onnen wir $C$ als irreduzibel annehmen.
Sei $K \subseteq L$ der algebraische Abschluss von $K$,
siehe \cite{scheja}, Satz 56.22. 
Dann kann keine eindimensionale Komponente von $C_{(L)}$ eigentlich sein.
Sei n"amlich $C'$ eine eindimensionale Komponente von $C_{(L)}$
und angenommen, dass
$C' \longrightarrow {\rm Spek}\, L$ eigentlich ist.
Da $K \subseteq L$ ganz ist, ist
${\rm Spek}\,L \longrightarrow {\rm Spek}\, K$
universell abgeschlossen und damit ist $C'$ universell abgeschlossen "uber $K$.
Da $C_{(L)} \longrightarrow C$ ganz ist,
sind zum einen die Fasern nulldimensional, weshalb $C'$ dominant ist,
zum anderen ist die Abbildung abgeschlossen, so dass
$C' \longrightarrow C$ surjektiv ist.
Damit ist nach der Vor"uberlegung auch $C$ schon universell abgeschlossen "uber
$K$ und damit eigentlich, was aber nicht der Fall ist.
\par\smallskip\noindent
Wir k"onnen also aufgrund von Satz 1.14.3 den treuflachen Basiswechsel
zum algebraischen Abschluss von $K$ vornehmen,
da sich dabei die Voraussetzungen nicht "andern.
Sei also $K$ algebraisch abgeschlossen.
Aufgrund von Satz 1.13.2 k"onnen wir annehmen, dass
$C$ eine separierte, integre, nicht eigentliche
Kurve "uber $K$ algebraisch abgeschlossen ist.
Geht man zur Normalisierung "uber, so bleibt nach der Vor"uberlegung
die Kurve nicht eigentlich, nach Satz 1.14.2 k"onnen wir also
$C$ noch als normal und damit regul"ar annehmen.
\par\smallskip\noindent
Die abgeschlossenen Punkte von $C$ sind dann (wegen separiert) in
eindeutiger Weise als Bewertungen des zugeh"origen Funktionenk"orpers
aufzufassen. Ist $R=K[q_1,...,q_m]$ der globale Schnittring einer
affinen Teilmenge von $C$, so liegt $R$ in fast allen Bewertungsringen
des Funktionenk"orpers drin, siehe \cite{haralg}, Lemma I.6.5,
so dass in $C$ bis auf endlich viele alle Bewertungen auftauchen.
Die Menge aller Bewertungen bildet nach \cite{haralg}, Theorem I.6.9,
eine projektive glatte Kurve.
\par\smallskip\noindent
Damit ist also $C$ eine echte offene Teilmenge einer glatten projektiven
Kurve $D$, das Komplement bestehe aus den Punkten $P_1,...,P_n$,
und es gen"ugt zu zeigen, dass das Komplement eines einzigen Punktes
affin ist.
Nach dem Satz von Riemann-Roch ist
$ l(kP) \geq =k +1-g $ und f"ur $k$ hinreichend gro"s ist
$l(kP) >0$. Das hei"st, dass es auf $D-\{ P \}$ eine nicht konstante
Funktionen $q$ gibt. Diese hat in $P$ einen Pol,
fasst man $q$ als meromorphe Funktion nach ${\bf P}_K^1$ auf, so ist
$P$ das Urbild von $\infty$ und $D - \{P \}$ das Urbild von ${\bf A}_K^1$.
Die Abbildung ist eigentlich und quasiendlich,
also nach \cite{EGAIII}, 4.4.2 oder \cite{haralg}, Exc. III.11.2
endlich und affin, weshalb das Komplement von $P$ affin ist
(Die Affinit"at folgt auch aus 1.8.2.)
\par\smallskip\noindent
Ebenfalls mit Riemann-Roch kann man zeigen, dass ein Divisor mit positivem
Grad ger"aumig ist und damit auf die Affinit"at schlie"sen,
siehe \cite{haralg}, Cor. IV.3.2.
Diese Aussage steckt auch im allgemeinen Schnittkriterium f"ur ger"aumig
drin. Die Schnittbedingung
$H^{{\rm dim} \, Y}.Y >0$ f"ur jedes integre Unterschema $Y$ reduziert sich
zu ${\rm deg}\, H= H.D >0$ und
$h^0({\cal O}_D/{\cal I}_Q)=H^0.Q >0$, was erf"ullt ist. \hfill $\Box$
\par\bigskip
\noindent
Im Fall einer projektiven glatten Fl"ache $S$ besagt das
Schnittkriterium, dass ein Divisor $D$ ger"aumig ist, wenn $D^2 >0$ und
$D.C >0$ ist f"ur jede irreduzible Kurve $C$.
Mit diesem Kriterium kann man oft zeigen, dass es
zu einer offenen Menge $U\subseteq S$ einen
ger"aumigen Divisor mit dem Komplement
als T"rager gibt und dass $U$ affin ist.
Typisch daf"ur ist die n"achste Aussage, die wir in 3.5 benutzen
werden.
\begin{satz}
Sei $S$ eine glatte projektive Fl"ache "uber einem algebraisch
abgeschlossenen K"orper. Es sei $p:S \longrightarrow C$ ein projektiver,
nicht konstanter Morphismus auf eine glatte irreduzible Kurve $C$ derart,
dass alle Fasern irreduzibel sind.
Es sei $F$ eine Faser "uber einem Punkt und $D$ sei eine irreduzible
Kurve auf $S$, auf der der Morphismus nicht konstant ist {\rm (}die also $C$
dominiert{\rm )}. Dann ist $S-(F \cup D)$ affin.
\end{satz}
{\it Beweis}.
Wir zeigen, dass ein geeigneter effektiver Divisor mit den beiden Kurven
als Tr"ager ger"aumig ist, wozu wir das Schnittkriterium anwenden.
Da der Divisor effektiv ist, m"ussen wir zeigen, dass der Schnitt mit jeder
irreduziblen Kurve positiv ist.
Sei $D^2=n$. Es ist $D.F > 0$, w"ahle $k \in {\bf N}$ so, dass
$(D+kF).D$ positiv ist. Sei $H=D+kF$. Sei $E$ eine beliebige, irreduzible,
effektive
Kurve auf $S$. Liegt $E$ in einer Faser, so ist sie bereits die ganze
Faser und hat mit $D$ einen positiven Schnitt, also auch mit $H$.
Ist $E$ nicht in einer Faser, so dominiert sie die Bildkurve und hat mit
der Faser $F$ einen positiven Schnitt. Bei $E \neq D$ folgt daraus sofort,
dass $H.E > 0$ ist und bei $E=D$ gilt das nach Konstruktion.  \hfill   $\Box$
\par\bigskip
\noindent
{\bf Bemerkung}
Ohne die Voraussetzung der irreduziblen Fasern gilt die
Aussage nicht.
Ist $F \cup D$ wie oben, und bl"ast man einen Punkt aus $S-(F \cup D)$ auf,
so geht die Affinit"at verloren, und eine Faser wird reduzibel.
\par\bigskip
\noindent
{\bf Beispiel}
Die einfachsten F"alle, wo die beschriebene Situation
vorliegt, sind Produkte von zwei glatten projektiven Kurven $X=D \times C$.
Der Satz besagt dann, dass der Restraum affin ist,
wenn man eine Faser herausnimmt und eine dazu
querliegende Kurve.
Ist etwa $D=C$, und nimmt man eine Faser und die Diagonale heraus, so
ist der Restraum affin. Diese Aussage ist f"ur uns in Abschnitt 3.5 wichtig,
siehe insbesondere Satz 3.5.3.
\par\smallskip\noindent
Das Komplement der Diagonalen allein ist nur im rationalen Fall affin,
aber nicht, wenn die Kurve vom Geschlecht $\geq 1$ ist.
\par\bigskip
\noindent
Die Affinit"at einer offenen Menge hat f"ur einen Divisor mit dem
Komplement von $U$ als Tr"ager Schnitteigenschaften zur Folge,
die dem Schnittkriterium f"ur die Ger"aumigkeit nahestehen,
aber schw"acher sind. Die Frage, wann es im Komplement einer affinen Teilmenge
wirklich einen ger"aumigen Divisor gibt, wird in Abschnitt 3.7. untersucht.
\begin{satz}
Sei $U \subseteq X$ eine offene Teilmenge in einer eigentlichen
Variet"at "uber einem algebraisch abgeschlossenen K"orper,
$Y=X-U$ sei das abgeschlossene Komplement.
Ist dann $U$ affin, so ist $Y$ eine Hyperfl"ache, die jede
Kurve schneidet.
Ist $X$ glatt, so hat $Y$ mit jeder Kurve $C$,
die nicht in $Y$ liegt, positiven Schnitt.
\end{satz}
{\it Beweis}.
Eine Kurve $C$ auf einem eigentlichen Schema ist als abgeschlossene Teilmenge
wieder eigentlich und damit projektiv.
Da affine Schemata keine projektiven Kurven 
enthalten, ist $C \not\subseteq U$, also $C \cap Y \neq \emptyset$.
Ist $C \not\subseteq Y$, so ist damit
der schnitttheoretische Schnitt von $Y$ und
$C$ positiv, also $Y.C >0$. \hfill $\Box$
\par\bigskip
\noindent
{\bf Bemerkung}
Sei $X={\rm Proj}\, A$ gegeben durch eine positiv graduierte $K-$Al\-gebra,
die in der nullten Stufe konstant ist und von der
ersten Stufe erzeugt wird mit dem punktierten Kegel $X'={\rm Spek}\, A- \{A_+\}$
und den entsprechenden Mengen $U'$ und $Y'$. Eine Kurve $C$ von $X$ ist stets
projektiv und das Urbild davon ist eine Fl"ache. Genau dann liegt $C$
innerhalb von $U$, wenn $C'$ und $Y'$ disjunkt im punktierten Kegel
sind, was genau dann der Fall ist, wenn sie sich im Kegel nur in der
Kegelspitze treffen. 
Dem Nachweis, dass $U$ nicht affin ist aufgrund einer projektiven Kurve
$C \subseteq U$, entspricht also im Kegel der Nachweis, dass $U'$
nicht affin ist aufgrund einer punktierten homogenen Fl"ache als
abgeschlossenem Unterschema. F"ur die Hyperfl"ache $Y'$ bedeutet dies, dass 
sie mit einer Fl"ache lediglich einen Punktschnitt hat, und deshalb die
universelle Kodimensionseigenschaft verletzt ist.
\par\bigskip
\noindent
Die Eigenschaft einer affinen Teilmenge $U$, dass keine projektiven Kurven
in ihr enthalten sind, gilt auch f"ur quasiaffine Schemata.
Ist aber zus"atzlich $U$ semiaffin, siehe den Abschnitt 1.12,
so folgt daraus bereits affin.
\begin{prop}
Sei $X$ eine integre, normale, eigentliche Variet"at "uber
einem algebraisch abgeschlossenen K"orper $K$.
Dann gelten folgende Aussagen.
\par\smallskip\noindent
{\rm (1)} Ist $X \supseteq U$ semiaffin und gibt es
in $U$ keine projektiven Kurven, so ist $U$ affin.
\par\smallskip\noindent
{\rm (2)} Ist $D$ ein {\rm (}potentiell{\rm )} basispunktfreies
Geradenb"undel auf
$X$, so ist das Komplement der Nullstelle eines Schnittes daraus semiaffin.
\par\smallskip\noindent
{\rm (3)} Ist $X$ eine glatte Fl"ache, und ist $U={\rm Def}(q)$
der Definitionsbereich
einer rationalen Funktion, so ist $U$ semiaffin.
\end{prop}
{\it Beweis}. 
Zu (1). Ist $U$ eine semiaffine Variet"at,
so ist auch $U \longrightarrow V$ eigentlich
mit einer affinen Variet"at, siehe \cite{goodlandman}
Da $U$ keine projektiven Kurven enth"alt,
sind alle Fasern "uber den abgeschlossenen Punkten endlich,
der Morphismus ist quasiendlich und eigentlich,
also endlich, und $U$ ist affin.
\par\smallskip\noindent
Zu (2). Durch "Ubergang zu einer Potenz nehmen wir das Geradenb"undel als
basispunktfrei an und betrachten die durch ein "uberdeckendes System
$(s,s_1,...,s_n)$ von Schnitten
global definierte Abbildung in den ${\bf P}^n_K$.
Diese ist eigentlich und daher ist auch die
Einschr"ankung $D(s) \longrightarrow D_+(x)$ eigentlich, was $D(s)$ als
semiaffin nachweist (dies gilt f"ur jeden Ring $K$).
\par\smallskip\noindent
Zu (3). Das Komplement von $U$ bestehe aus den Kurven
$C_1,...,C_k$ und der Polstellendivisor zu $q$ sei
$D=n_1C_1+...+n_kC_k$ mit $n_k >0$. Der Nullstellendivisor ist dazu "aquivalent,
hat aber keine gemeinsamen Komponenten.
Damit ist das Geradenb"undel zu $D$ komponentenfrei und hat daher allenfalls
abgeschlossene Punkte als Basispunkte. Nach dem Satz von Zariski,
siehe \cite{zariskiriemann}, Theorem 6.1, ist dann aber schon ein Vielfaches
"uberhaupt basispunktfrei und die Aussage folgt aus Teil (2). \hfill $\Box$
\par\bigskip
\noindent
Eine offene Menge $U$ in einer eigentlichen Variet"at kann man manchmal
dadurch als nicht affin nachweisen, indem man zeigt, dass es auf $U$
zu wenig globale Funktionen gibt.
Begriffe, die die Reichhaltigkeit der Schnittalgebra
$\Gamma (U,{\cal O}_X)$ zum Ausdruck bringen, und die bei
$U \subseteq {\rm Spek}\,A$ automatisch erf"ullt sind,
sind quasiaffin und generisch affin, siehe den Abschnitt 1.7.
Ist das Komplement von $U$ ein Divisor, so dr"uckt die Iitaka-Dimension
des Divisors ebenfalls etwas "uber die Reichhaltigkeit aus,
vergleiche \cite{iitaka}, § 10.1, insbesondere Prop. 10.1 und
The\-o\-rem 10.2.
\par\bigskip
\noindent
{\bf Definition}
Sei $D$ ein (Weil-) Divisor auf einer eigentlichen normalen Variet"at $X$
"uber einem algebraisch abgeschlossenen K"orper $K$.
Dann bezeichnet man das Maximum
${\rm max}_{m \in {\bf N}}\, {\rm dim}\, ({\rm Bild}(\Phi_{mD}))$
als die $D-${\it Dimension} oder die {\it Iitaka-Dimension}
des  Divisors $D$.
\par\bigskip
\noindent
{\bf Bemerkung}
Dem Divisor $mD$ entspricht auf einer offenen Teilmenge
von $X$ eine invertierbare Garbe ${\cal L}_{mD}$.
Ein Erzeugendensystem von ${\cal L}_{mD}$ legt damit auf einer
offenen Teilmenge eine Abbildung $\Phi_{mD}$ in einen projektiven Raum
fest, siehe \cite{iitaka}, 2.9.d oder \cite{haralg}, II.7.
Nach \cite{iitaka}, The\-o\-rem 10.2, gilt f"ur einen ef\-fek\-tiven Divi\-sor
mit der Iitaka-Dimen\-sion $k$ die Absch"atzung (f"ur $m$ hinreichend gro"s)
$$ \alpha m^k \leq h^0(mD) \leq \beta m^k \mbox{ mit } \alpha, \beta >0 \, ,$$
d.h. die Iitaka-Dimension beschreibt das asymptotische Verhalten der globalen
Schnitte.
Die Kodairadimension einer projektiven glatten
Variet"at ist die Iitaka-Dimension des
kanonischen Divisors $K$. 
\begin{prop}
Sei $X$ eine eigentlich normale Variet"at der Dimension $d$
"uber einem algebraisch abgeschlossenen K"orper $K$.
Sei $U \subseteq X$ das Komplement einer Hyperfl"ache
$X-U=Y=Y_1 \cup...\cup Y_k$
und $D=n_1Y_1+...+n_kY_k$ ein Divisor mit $n_i >0$.
Dann ist $U$ genau dann generisch affin, wenn die Iitaka-Dimension von $D$ maximal ist,
also gleich der Dimension von $X$.
\end{prop}
{\it Beweis}.
Sei $U$ generisch affin und $0 \neq f \in \Gamma(U,{\cal O}_X)=R$
mit $U_f$ affin.
Es ist $\Gamma(U_f,{\cal O}_X)=R_f$.
$R_f$ ist eine endlich erzeugte $K-$Algebra, sagen wir
$R_f=K[a_1/f,...,a_r/f]$ mit $a_j \in R$. Setze $a_0=f$.
Die $a_i$ haben in den Punkten der H"ohe eins aus $U$ nichtnegative Ordnung
und k"onnen allenfalls an den $Y_i$ Polstellen haben.
F"ur $m$ hinreichend gro"s ist dann
$$a_j \in {\cal L}_{mD}(X) =\{ q \in K(X):\, (q) +mD \geq 0 \} \, .$$
f"ur alle $j=0,...,r$.
Sei $\Phi=\Phi_{mD}$ und $(a_0,...,a_r,....,a_s)$ ein Erzeugendensystem von
${\cal L}_{mD}(X)$.
Es ist $\Phi$ auf $X_f=X_{a_0}$ definiert.
Wegen ${\cal L}_{mD} \vert_U={\cal O}_X \vert _U$ ist $U_f \subseteq X_f$.
Auf $U_f$ ist $\Phi: U_f \longrightarrow D_{x_0} \subseteq {\bf P}^s_K$
gegeben durch  $a_0/f,....,a_s/f$. Da dies ein $K-$Algebra Erzeugendensystem
von $R_f$ enth"alt, ist $\Phi$ eingeschr"ankt auf $U_f$ eine
abgeschlossene Einbettung und damit ist die Dimension des Bildes so gro"s
wie die von $X$.
\par\smallskip\noindent
Sei umgekehrt ${\cal L}={\cal L}_{mD}$ so, dass die durch die Schnitte
$1=s_0,s_1,...,s_r$ definierte Abbildung $\Phi$ ein Bild der Dimension
$d$ hat.
Da $X$ normal ist, kann man $s_i$ auffassen als Funktionen, die auf $U$
definiert sind, und $\Phi$ hat auf $U$ die Gestalt
$$ U \stackrel{(f_1,...,f_r)}{\longrightarrow} Z=V({\bf p})
\subseteq {\bf A}_K^r=D(x_0) \subseteq {\bf P}^r_K \, .$$
Dabei ist $Z$ der (integre) Abschluss des Bildes von $\Phi \vert_U$ und
damit $d-$dimen\-sional.
Es ist also $U \longrightarrow Z$ dominant von gleicher Dimension, und
damit liegt "uber dem generischen Punkt von $Z$ nur der generische Punkt
von $U$. Damit ist die Abbildung generisch affin und nach Satz 1.8.1 gibt
es eine offene Umgebung
$\emptyset \neq D(h) \subseteq Z,\, h \in K[T_1,...,T_r]$ mit affinem Urbild,
d.h. es ist $U_{h(f_1,...,f_r)}$ affin und $U$ generisch affin. \hfill $\Box$
\par\bigskip
\noindent
{\bf Bemerkung}
Ist $X$ eine eigentliche glatte Fl"ache "uber einem algebraisch abgeschlossenen
K"orper und $D$ ein effektiver Divisor mit positivem
Selbstschnitt $D^2 >0$, so ist $X-D$ generisch affin.
Da $D$ effektiv ist, ist $h^0(K-mD)$ beschr"ankt und nach Riemann-Roch,
siehe \cite{haralg}, Theorem V.1.6, ist
$$h^0(mD) \geq {1 \over 2} m^2 D^2-{1 \over 2} m D.K +c \,\,
(c \mbox{ Konstante }) \, .$$
Da $h^0(mD)$ durch ein quadratisches Polynom nach unten abgesch"atzt wird,
muss die Iitaka-Dimension 2 sein.
\par\smallskip\noindent
Die Umkehrung davon muss nicht gelten. Ist etwa $X$ die in einem Punkt $Q$
aufgeblasene projektive Ebene, und $C={\bf P}^1 \subseteq X, \, Q \not\in C$,
so ist $X-(E \cup C)$ quasiaffin und insbesondere generisch affin,
aber $(E+C)^2=E^2+C^2=0$.
\par\bigskip
\noindent
{\bf Beispiel}
Wir wollen an unserem Standardbeispiel $A=K[x,y,u,v]/(ux-vy)$ "uber einem
K"orper $K$ illustrieren,
dass sich je nach Graduierung die Nicht-Affinit"at der
offenen Teilmenge $D(x,y)=U \subseteq {\rm Spek}\, A$
f"ur $V=D_+(x,y) \subseteq {\rm Proj}\, A$
in ganz verschiedener Weise auswirkt. Mal enth"alt $V$ projektive Kurven,
mal gibt es auf $V$ zu wenig Funktionen (nicht generisch affin), mal ist
$V$ nicht separiert, mal nicht semiaffin.
Die diskutierten ${\bf Z}-$Graduierungen werden alle durch die Grade von
$x,y,u,v$ festgelegt.
\par\medskip\noindent
Sei
${\rm grad}\, (u,x,y,v)=1$.
Dann ist ${\rm Proj}\, A={\bf P}^1_K \times {\bf P}^1_K$, die Abbildung
wird durch die beiden meromorphen Funktionen $q=u/ y=v/ x$
und $p=u/v=y/x$ gegeben.
Es ist $V_+(x,y)=\{ \infty \} \times {\bf P}_K^1 $ und daher
$D_+(x,y)={\bf A}_K^1 \times {\bf P}_K^1 $. Die Nichtaffinit"at zeigt sich z.B.
darin, dass es projektive Kurven gibt, oder darin, dass der Schnittring
eindimensional ist.
\par\medskip
\noindent
Sei ${\rm grad}\, (x,y)=1,\, {\rm grad}\, (u,v)=2$.
Das ist eine Hyperfl"ache in einem gewichteten projektiven Raum. 
Um $D_+(x,y)$ zu beschreiben, betrachten wir
$$ D(x,y) \subseteq {\rm Spek}\, \Gamma(D(x,y),{\cal O}_X)=
{\rm Spek}\, K[x,y,q]\, ,$$
wobei $q=u/y=v/x$  vom Grad eins ist.
Auf $D(x,y)$ gibt es nur die konstanten Funktionen als homogene Elemente
von Grad null, und damit besteht der globale Schnittring
von $D_+(x,y)$ nur aus den Konstanten.
$D(x,y)$ ist ein standardgraduierter affiner Raum ohne eine Gerade,
daher ist $D_+(x,y) = {\bf P}_K^2 - \{ P\} $.
Dies ist nicht semiaffin.
\par\medskip
\noindent
Sei ${\rm grad}\, (u,v)= 0,\, {\rm grad}\, (x,y)=1$.
Da ist ${\rm Proj}\, A$ die affine Ebene im Nullpunkt aufgeblasen, und dies ist
auch gleich $D_+(x,y)$. Der globale Schnittring ist zweidimensional,
es handelt sich um ein semiaffines, generisch affines Schema, das aber eine
projektive Kurve enth"alt und daher nicht affin sein kann.
\par\medskip
\noindent
Sei ${\rm grad}\, (v,x)=1,\, {\rm grad}\, (u,y)=-1 $.
Es ist $D_+(x,y) \subseteq {\rm Proj}\, K[q,x,y]$ mit
${\rm grad}\, q=0,\, {\rm grad}\, x=1,\, {\rm grad}\, y=-1 $.
Das ist dann die am Nullschnitt verdoppelte
affine Gerade "uber ${\rm Spek}\, K[q]$,
das ist eine affine Ebene mit einer Geraden verdoppelt, und damit nicht
separiert und nicht affin.
\par\medskip
\noindent
Sei ${\rm grad}\, (x,y)=1,\, {\rm grad}\, (u,v)=-1$.
$D(x,y)$ ist da der positive Teil $U_+$ von ${\rm Proj}\, A$,
also gleich ${\rm Proj}\, K[ux,vx,uy,vy][X,Y]$
(mit den zugeh"origen Relationen).
Dabei ist der Ring der nullten Stufe isomorph zu
$A_0=K[r,s,t,w]/(rs-tw)$ und ${\rm Proj} \, A$ ist isomorph
zur Aufblasung von $A_0$ am Ideal $(r,t)$.
Das ist somit ein separiertes Schema, das "uber dem Nullpunkt von $A_0$
eine projektive Kurve enth"alt und daher nicht affin ist.
\par
\bigskip
\noindent
\subsection{Affine Aufblasungen und Normalisierung}
\par
\bigskip
\noindent
Sei $X$ ein Schema, ${\cal I}$ eine
Idealgarbe und $\tilde{X}={\rm Proj} \sum_{n \in {\bf N}}{\cal I}^n$
die dadurch definierte
{\it Aufblasung l"angs} oder {\it an} $I$,
siehe \cite{haralg}, II.7. Wir stellen uns die
Frage, wann die Abbildung $\tilde{X} \longrightarrow X$ affin ist?
\par\smallskip\noindent
Dabei nehmen wir an, dass $X={\rm Spek}\, A$ affin ist, also
$\tilde{X}= {\rm Proj}\, \sum_{n \in {\bf N}} I^n$ mit einem Ideal
$I \subseteq A$.
Aufgrund der universellen Eigenschaft der Aufblasung ist das Erweiterungsideal
{\it invertierbar}, d.h. es ist endlich erzeugt und wird lokal
durch einen Nichtnullteiler beschrieben, vergleiche \cite{eisenbud}, 11.3.
Ist $\tilde{X}$ affin, so ist das Komplement dieses
invertierbaren Erweiterungsideals affin, und da die Aufblasung auf $D(I)$
isomorph ist, muss bereits $D(I) \subseteq X$ affin sein.
Insbesondere ist eine affine Aufblasung nur bei Idealen der
H"ohe eins m"oglich.
\begin{lem}
Sei $A$ ein Krullbereich. Dann ist die Aufblasung eines endlich
erzeugten Ideals nur dann affin, wenn das Ideal invertierbar oder das
Nullideal ist, wenn also die Aufblasung die Identit"at ist oder die
leere Abbildung.
\end{lem}
{\it Beweis}.
Das Ideal $I\neq0$ besitze eine affine Aufblasung.
$I$ ist in jedem Punkt der H"ohe eins invertierbar, und da es endlich
erzeugt ist, bereits in einer offenen Menge $U$, die all diese Punkte
der H"ohe eins umfasst. Damit ist auf $U$ die
Aufblasungsabbildung eine Isomorphie, und wir haben die Inklusionen
$$ A \longrightarrow \Gamma(\tilde{X},{\cal O}_{\tilde{X}}) 
\longrightarrow \Gamma(\tilde{U},{\cal O}_{\tilde{X}})
\longrightarrow \Gamma(U,{\cal O}_X) =A \, .$$
Damit ist der
globale Schnitt\-ring der Aufblasung gleich $A$. Da $\tilde{X}$ affin ist,
liegt eine Isomorphie vor.     \hfill   $\Box$
\par\bigskip
\noindent
{\bf Bemerkung}
Ist $A$ noethersch, so ist "uberhaupt bei jeder Aufblasung der globale
Schnittring endlich "uber $A$.
Da die Abbildung projektiv ist, ist nach
\cite{haralg}, Theorem 5.2, $B=\Gamma(\tilde{X},{\cal O}_{\tilde{X}})$
ein endlicher erzeugter $A-$Modul.
\par\bigskip
\noindent
Ein {\it gebrochenes Ideal} ist ein $A-$Untermodul des Quotientenk"orpers
$Q(A)$, siehe \cite{eisenbud}, 11.3. Ein gebrochenes Ideal $I$ ist genau dann
invertierbar, wenn es ein anderes gebrochenes Ideal $I'$ gibt mit
$II'=A$. Dabei ist dann $I'=I^{-1}=\{t \in Q(A): tI \subseteq A \}$, siehe
\cite{eisenbud}, Theorem 11.6.
\begin{satz}
Sei $A$ ein integrer noetherscher Ring
und $0 \neq I \subseteq A $
ein Ideal. Dann sind folgende Aussagen "aquivalent.
\par\smallskip
\noindent
{\rm (1)} Die Aufblasung  $ \tilde{X}$ von $X={\rm Spek }\, A $
l"angs $I$ ist affin.
\par\smallskip\noindent
{\rm (2)} Das Erweiterungsideal $J=I A_{\rm nor} $ in der
Normalisierung ist invertierbar.
\par\smallskip\noindent
{\rm (3)} Es gibt eine ganze integre  Erweiterung $A \subseteq B $ so,
dass das Erweiterungsideal $IB$ invertierbar ist.
\par\smallskip\noindent
{\rm (4)} Es gibt eine endliche integre Erweiterung $A \subseteq R$ so,
dass das Erweiterungsideal $IR$ invertierbar ist.
\end{satz}
{\it Beweis}.
$(1) \Longrightarrow (2)$.
Ist $A \longrightarrow A'$ ein Ringhomomorphismus, so ist die Abbildung
zwischen der zu $I$ und der zum Erweiterungsideal $I'=IA'$
geh"orenden Aufblasung affin.
Ist also die Aufblasung von $I$ affin, so gilt das auch f"ur jede
Ringerweiterung, insbesondere f"ur die Normalisierung.
Die Normalisierung eines noetherschen Bereiches ist aber
nach \cite{nagloc}, Theorem 33.10, ein Krullbereich
und daher ist nach
dem Lemma die Aufblasung eine Isomorphie und das Erweiterungsideal ist
in der Normalisierung invertierbar.
\par\smallskip\noindent
$(2) \Longrightarrow (3) $ ist trivial (in (3) wird $B$ nicht als
noethersch vorausgesetzt).
\par\smallskip\noindent
$(3) \Longrightarrow (4) $.
Sei $a_j \in I,\, j=1,...,n $, ein endliches Erzeugendensystem von $I$.
Nach Voraussetzung ist das Erweiterungsideal $J=IB$ invertierbar,
dann gibt es ein gebrochenes (endlich erzeugtes) $B-$Ideal
$G=(b_1,...,b_r)$ mit $J \cdot G = B$ und
$\sum_{j,k} c_{j,k} a_jb_k=1$ mit $c_{j,k} \mbox{ und } a_jb_k \in B$.
Sei $R:=A[c_{j,k}, a_jb_k] \subseteq B$.
$R$ ist dann ganz und endlich erzeugt, also endlich "uber $A$.
Man sieht sofort, dass auch in $R$ das Erweiterungsideal zu $I$
invertierbar ist mit dem Inversen $H=R(b_1,...,b_n)$.
\par\smallskip
\noindent
$(4) \Longrightarrow (1)$
Man braucht lediglich, dass zu einer endlichen Ringerweiterung auch die
zugeh"orige Projabbildung endlich ist, siehe das folgende Lemma.
Da $IR$ invertierbar ist, liegt die Situation
${\rm Spek}\, R \longrightarrow \tilde{X} \longrightarrow {\rm Spek}\, A$
vor, und die erste Abbildung ist diese endliche Projabbildung. Der Satz von
Chevalley 1.14.1 liefert dann die Affinit"at von $\tilde{X}$.  \hfill  $\Box $
\par\bigskip
\noindent
Wir haben soeben folgendes Lemma verwendet.
\begin{lem}
Sei $A \subseteq R $ eine endliche Erweiterung von
Integrit"atsbereichen und $I$ ein Ideal in $A$. Dann ist auch die
Abbildung zwischen den korrespondierenden Aufblasungen endlich und surjektiv.
\end{lem}
{\it Beweis}.
Sei $S = A+I+I^2+...$ und $a \in I = ( a_j,\, j \in J )$.
Es ist
$\Gamma (D_+(a),{\cal O}_{\tilde X})= \{ f/ a^r :\, f \in I^r \}
=A[a_j/ a ,\, j \in J] \, $.
Auf diesem affinen St"uck ist die Abbildung einfach gleich
$A[ a_j/ a] \subseteq C[a_j/ a]$.
Das sind endliche, injektive Erweiterungen, und die $D_+(a)$ "uberdecken das
hintere Aufblasungsschema.     \hfill   $\Box$
\begin{kor}
Sei $A$ ein integrer, normaler, noetherscher Ring und $I$ ein Ideal.
Ist dann $IB$ in einer ganzen integren Erweiterung $A \subseteq B$
invertierbar, so schon in $A$. \hfill $\Box$
\end{kor}
{\bf Bemerkung}
Sei $I$ ein divisorielles Ideal in $A$ und $A \longrightarrow B$
ein Ringhomomorphismus $A,B$ normal und noethersch.
Dann hat man zu unterscheiden zwischen dem Erweiterungsideal $IB$
und dem zur"uckgenommenen Divisor, siehe \cite{bourbaki},\S1.10.
Dieser Divisor kann invertierbar werden, ohne dass $I$ invertierbar ist,
allerdings nur dann, wenn im Erweiterungsideal
eingebettete Komponenten auftreten.
\begin{kor}
Sei $A$ ein integrer noetherscher Ring der
Dimension eins. Dann ist jede Aufblasung von $X={\rm Spek}\, A$ affin.
\end{kor}
{\it Beweis}.
In der Normalisierung der Kurve ist jedes von null verschiedene
Ideal invertierbar.
\par\smallskip\noindent
Zweiter Beweis. Man kann sich auf den lokalen Fall beschr"anken.
Die Abbildung $\tilde{X} \longrightarrow X$ ist auf dem offenen,
generischen Punkt eine Isomorphie, und insbesondere ist $\tilde{X}$ ein
noethersches, integres, separiertes Schema mit offenem, generischen Punkt.
Damit ist $\tilde{X}$ nach Satz 1.7.6 affin.  \hfill   $\Box$
\par\bigskip
\noindent   
Der Satz 2.5.2 hat insbesondere gezeigt, dass die Normalisierung
eines inte\-gren noetherschen Ringes durch jede affine Aufblasung faktorisiert
und dass somit jeder globale Schnittring einer affinen Aufblasung zwischen dem
Ausgangsring und der Normalisierung liegt. Wir fragen uns nun,
inwiefern man die Normalisierung durch affine Aufblasungen erhalten kann.
\begin{lem}
Sei $A$ ein integrer noetherscher Ring, $I \subseteq A$
ein Ideal und $\tilde{X} $ die Aufblasung von $X={\rm Spek}\, A $ zum
Ideal $I$. Es sei ein Zwischenring $B$ mit
$A \subseteq B \subseteq \Gamma(\tilde{X},{\cal O}_{\tilde{X}})$ gegeben.
Dann ist die kanonische Abbildung
$\tilde{X} \longrightarrow {\rm Spek}\, B$ auch die Aufblasung
des Erweiterungsideals $J=IB$.
\end{lem}
{\it Beweis}.
Die Aufblasung von $Y={\rm Spek}\, B$ nach $J$ sei mit $\tilde{Y}$
bezeichnet und dessen globaler Schnittring mit $C$.
Sei $\Psi: \tilde{Y} \longrightarrow \tilde{X} $ die
Abbildung zwischen den Aufblasungen.
Dann ist die Aufblasung $\tilde{Y} \longrightarrow Y$ gleich
$f \circ \Psi$, wobei $f :\tilde{X} \longrightarrow Y$ durch
die Ringinklusion gegeben ist.
Dies kann man auf der Ringebene zeigen, und ist klar, da sich alles im
Quotientenk"orper abspielt.
\par\smallskip\noindent
Da das geliftete Ideal von $J$ in $\tilde{X}$ invertierbar ist, gibt
es aufgrund der universellen Eigenschaft von $\tilde{Y}$
eine Abbildung $\Phi : \tilde{X} \longrightarrow \tilde{Y}$.
Dabei sind $\Psi \circ \Phi$ und $ \Phi \circ \Psi $
die Identit"at, da sie auf den separierten integren
Schemata $\tilde{X}$ und $ \tilde{Y} $
auf einer offenen, nicht leeren Teilmenge mit der Identit"at "ubereinstimmen.
Also ist $\tilde{X}=\tilde{Y}$. \hfill   $\Box$
\begin{satz}
Sei $A$ ein integrer noetherscher Ring. Dann gelten folgende Aussagen.
\par\smallskip\noindent
{\rm (1)}
Ist $q \in A_{\rm nor} $ und $q= x/ y$ eine Darstellung
mit Elementen aus $A$, so ist das Ideal $I=( x,y ) $ in der
Normalisierung ein Hauptideal mit erzeugendem Element $y$, die Aufblasung von
$A$ l"angs $I$ ist affin mit globalem Schnittring $A[q]$.
\par\smallskip\noindent
{\rm (2)}
Jede endliche birationale Erweiterung von $A$ ist der globale Schnittring
einer Folge von affinen Aufblasungen.
\par\smallskip\noindent
{\rm (3)}
Genau dann gibt es eine {\rm (}endliche{\rm )} Folge von affinen
Aufblasungen, deren letzte die Normalisierung ergibt, wenn die
Normalisierung endlich "uber dem Grundring ist.
Insbesondere k"onnen
integre exzellente Ringe durch affine Aufblasungen normalisiert werden.
\end{satz}
{\it Beweis}.
Zum ersten Teil ist nur noch die Aussage "uber den Schnittring zu best"atigen.
Man zieht eine Ganzheitsgleichung f"ur $q$ heran,
$$q^n +a_{n-1}q^{n-1} +...+a_1q+a_0 =0 \, .$$
Multiplikation mit $y^n$ ergibt
$$ x^n+a_{n-1}x^{n-1}y+...+a_1xy^{n-1}+a_0y^n =0 \,  ,$$
woraus sich wiederum
$q= x/ y = (-a_{n-1}x^{n-1}-...-a_0y^{n-1}) /x^{n-1} $
ergibt.
Daraus liest man ab, dass $q$ "uberall eine Darstellung mit Z"ahler und
Nenner von gleichem Grad hat, daher ist $q$ im
globalen Schnittring der Aufblasung definiert. Es liegt dann die Situation
$\tilde{X} \longrightarrow {\rm Spek}\, A[q] \longrightarrow {\rm Spek}\, A $
vor, das Ideal ist aber schon in der Mitte invertierbar und nach dem 
obigen Lemma ist dann die erste Abbildung eine Isomorphie.
\par\smallskip\noindent
Zu (2). Sei $B=A[q_1,...,q_n]$ eine endliche birationale Erweiterung
von $A$, $q_i \in Q(A)$.
Da $B$ endlich ist, erf"ullen die $q_i$ eine Ganzheitsgleichung,
geh"oren also zur Normalisierung.
In Teil (1) wurde gezeigt, wie man $A[q_1]$ als affine Aufblasung
erhalten kann, durch Iteration dieses Verfahrens gelangt man zu $B$.
(3) ergibt sich direkt aus (2).      \hfill   $\Box$
\par\bigskip
\noindent
{\bf Beispiel}
Sei $A$ integer, exzellent und eindimensional. Dann ist die Normalisierung
dar"uber endlich und man kann nach dem Satz die Kurve mit affinen
Aufblasungen normalisieren und damit desingularisieren.
Die Desingularisierung von ebenen algebraischen Kurven wurde erstmals
von Max Noether, \cite{noether}, mittels elementaren Variablentransformationen
beschrieben, die der Hinzunahme von Funktionen in der Aufblasung
entsprechen.
\par\smallskip\noindent
Im Kurvenfall gen"ugt es sogar, nur an maximalen Idealen aufzublasen,
um zur Normalisierung zu gelangen.
Dabei bl"ast man einfach nicht-invertierbare, maximale Ideale auf.
Da die Normalisierung ein noetherscher $A-$Modul
ist, wird dieser Prozess station"ar, und man gelangt so zu einem
eindimensionalen noetherschen Integrit"atsbereich, bei dem alle maximalen
Ideale invertierbar sind. Dann sind aber alle Lokalisierungen
regul"ar und es liegt bereits die Normalisierung vor.
\par\bigskip
\noindent
{\bf Beispiel}
Die {\it Seminormalisierung} eines noetherschen integren Ringes besteht aus
allen Funktionen der Normalisierung, die verschiedene Punkte zu einem
gemeinsamen Basispunkt
nicht trennen, und ist die maximale ganze Erweiterung derart, dass die
Abbildung der Spektren eine Hom"oomorphie ist.
Schon im eindimensionalen Fall l"asst sich die Seminormalisierung 
nicht durch Aufblasen an maximalen
Idealen allein erreichen.
\par\smallskip\noindent
Betrachte hierzu den eindimensionalen Ring
$$ A=K[x,y]/(x^5-(y^2-x^2)^2)=K[t^2-1)^2,(t^2-1)^2 t ]
\subseteq K[t]=A_{\rm nor} \, .$$
Wegen $t=y / x$ f"uhrt die
Aufblasung am maximalen Ideal ${\bf m}=( x,y ) $
bereits zur Normalisierung, die nicht hom"oomorph ist, und an den anderen
abgeschlossenen Punkten ist die Kurve nicht singul"ar, an maximalen Idealen
aufgeblasen tut sich also nichts.
$A$ selbst ist aber nicht seminormal, da die Funktion $ (t^2-1) $ die
Urbildpunkte in der Normalisierung nicht trennt, aber nicht zu $A$ geh"ort.
Diese Funktion ist wegen $t^2-1={y^2-x^2 \over x^2}$ im
globalen Schnittring der Aufblasung zum Ideal $( x^2,y^2 )$
enthalten.
\par\bigskip
\noindent
{\bf Beispiel}
Sei $K[M]$ ein Monoidring zum numerischen Monoid
$M \subseteq \Gamma M={\bf Z}^n$. Die Normalisierung ist
gegeben durch $K[M']$, wobei $M'$ die
Normalisierung des Monoids ist, also
$$ M' =  \{ m \in \Gamma M :\,
\mbox{ Es gibt ein } r \in {\bf N} \, \mbox{ mit } \, rm \in M \} \,  .$$
F"ur $m \in M'$ ist dann $X^m= X^s/ X^t$ mit $s,t \in M$ und $rm \in M$,
und diese ganze Funktion kommt durch das Aufblasen des Ideals
$( X^s, X^t ) $ dazu.
Insbesondere gelangt man zur Normalisierung mittels
Aufblasen von monomialen Idealen.
\par\bigskip
\noindent
{\bf Beispiel}
Sei $A=R[X]/(X^2-p^k) $ eine quadratische Erweiterung eines faktoriellen
Ringes $R$ und $p \in R$ ein Primelement. $k$ sei ungerade $ \geq 3$.
Dann ist $(X /p)^2 = p^{k-2} $ eine Ganzheitsgleichung f"ur $X/p$ und die
Aufblasung an $( X, p ) $ hat als globalen Schnittring mit
$X'=X / p$ die Gestalt
$A[X']=R[X]/(X^2-p^k) [X'] = R[X']/ (X'^2-p^{k-2})$.
So f"ahrt man fort, bis man 
$k=1$ erreicht hat.
In vielen F"allen ist dies dann normal, etwa bei $R=K[T_1,...,T_n]$,
siehe \cite{haralg}, Exc. II.6.5
\par\bigskip
\noindent
Wenn man wei"s, dass die Normalisierung endlich ist, so
gelangt man mit affinen Aufblasungen schlie"slich zur Normalisierung,
unabh"angig davon, welche Ideale man aufbl"ast, solange man immer noch
Ideale findet, die nicht triviale, aber affine Aufblasungen besitzen. Man
muss hier also nicht in irgendeiner Weise geschickt vorgehen. Dies liegt
daran, dass die Aufblasungsringe eine aufsteigende Kette von $A-$Untermoduln
im noetherschen $A-$Modul $A_{\rm nor}$ liefert, die station"ar wird und
somit letztendlich die Normalisierung ist. Dagegen gibt es aber wohl keine
systematische Methode im Gegensatz zum eindimensionalen Fall, affin
aufblasbare Ideale zu finden. Insbesondere kann man nicht mit Aufblasungen
von Idealen mit irreduziblem Tr"ager allein zur Normalisierung gelangen,
wie das Beispiel weiter unten zeigt.
\par\bigskip
\noindent
Es sei daran erinnert, dass die Serre'sche Bedingung $S_2$ f"ur einen
integren Ring bedeutet, dass
f"ur ein Ideal $\bf a$ der H"ohe zwei die Restriktion
$A \longrightarrow \Gamma(D({\bf a}),{\cal O}_X)$
bijektiv ist. 
\begin{lem}
Sei $A$ ein integrer noetherscher Ring, der die Bedingung $S_2$ erf"ullt.
$D({\bf a}) $ sei normal und ${\bf a} +{\bf b} $ habe die H"ohe von
zumindest zwei. Ist dann die Aufblasung l"angs $\bf b$ affin, so ist sie
trivial.
\end{lem}
{\it Beweis}.
Die Aufblasung ist auf dem normalen Ort und auch au"serhalb der
Nullstellenmenge
des Aufblasungsideals isomorph. Wegen der H"oheneigenschaft, wegen der
Serre'schen Bedingung und wegen
$$ A \longrightarrow \Gamma(\tilde{X}, {\cal O}_{\tilde{X}})
\stackrel{R}{\longrightarrow} \Gamma(D({\bf a}+{\bf b}),{\cal O}_X)=A $$
k"onnen keine Funktionen im
globalen Schnittring der Aufblasung hinzukommen, woraus wegen der
Affinit"at die globale Isomorphie folgt.     \hfill   $\Box$
\par\bigskip
\noindent
{\bf Beispiel}
Wir betrachten den Monoidring
$$A=K[Y,W,Z,V,U,X]/( UX-VY,\, Y^2Z-W^2) = K[Y,TY,T^2,V,U,VY/U] \, .$$
Die Gleichheit ergibt sich daraus, dass beide Ringe vierdimensional sind und
die Restklassendarstellung schon integer ist, was man so sieht: zun"achst
ist $B:=K[Y,X,U,V]/(UX-VY) [Z] $ integer und es ist $A=B[W]/(b-W^2) $ mit
$b=Y^2Z$. Eine solche quadratische Erweiterung ist aber integer,
wegen der Inklusion $A \subseteq Q(B)[W]/(b-W^2) $ und weil $b$ in
$Q(B)$ keine Wurzel hat, da auch $Z$ keine Wurzel hat.
\par\smallskip\noindent
$A$ ist als vollst"andiger Durchschnitt Cohen-Macauley und
insbesondere $S_2$.
Die Normalisierung wird erreicht durch affines Aufblasen am Ideal
$ {\bf a} = ( Y , TY ) =( Y,W ) $ und die
Normalisierung ist $A[T]=K[Y,T,V,U,X]/(UX-VY)$.
Insbesondere ist $D({\bf a})$ normal.
\par\smallskip\noindent
Wir zeigen, dass man mit Aufblasungen von Idealen ${\bf b}$, die
irreduziblen Tr"ager haben, nicht "uber
den Ring selbst hinauskommt, und schon gar nicht zur Normalisierung.
Sei also ${\bf b }$ ein in der Normalisierung
invertierbares Ideal mit dem primen Radikal ${\bf p}$ in $A$.
$\bf p$ hat dann die H"ohe eins. Wenn ${\bf p} \in D({\bf a})$ ist,
also ${\bf a} \not\subseteq {\bf p}$, so hat ${\bf a}+ {\bf b} $ eine
H"ohe $\geq 2$, und die Aufblasung ist nach dem Lemma die Identit"at.
Ist hingegen ${\bf a} \subseteq {\bf p}$, so gibt es hierf"ur
nur zwei M"oglichkeiten, n"amlich
${\bf p} = ( Y, W ,U )  $
oder $\ = ( Y,W,X ) $. Diese werden aber in der
Normalisierung zu $( Y,U ) $ und zu $( Y, X ) $,
und diese beschreiben nicht-affine offene Mengen. 
Daher kann eine solche Aufblasung nicht affin sein.

%% file: afklaka.tex
\section{Die affine Klassengruppe}
\par
\bigskip
\noindent
In diesem Kapitel zeigen wir zun"achst, dass in einem normalen Schema
die Affinit"at des Komplementes eines effektiven Divisors
keine Eigenschaft der Divisorklasse ist.
Dies motiviert den Begriff der koaffinen Divisorklasse, bei der jeder
effektive Vertreter affines Komplement besitzt.
Dieser Begriff erlaubt es zu zeigen, dass sich die Eigenschaft eines affinen
normalen Schemas, dass jede Hyperfl"ache affines Komplement besitzt,
auf die affine Gerade "ubertr"agt.
Eine weitere Versch"arfung dieser Eigenschaft f"uhrt zu affin-trivialen
Divisoren, die eine Untergruppe der Divisoren(-klassen) bilden,
und die bei Affinit"atsbetrachtungen in gewissem Sinn
au"ser Acht gelassen werden k"onnen.
\par\smallskip\noindent
Die Restklassengruppe nennen wir die affine Klassengruppe, und diese
Gruppe beschreibt das Affinit"atsverhalten von Hyperfl"achen
in einem Schema. Sie ist in einem affinen Schema
genau dann null ist, wenn jede Hyperfl"ache affines Komplement besitzt.
Die affine Klassengruppe "andert sich beim "Ubergang zu einem affinen
Zylinder nicht.
F"ur Hyperbeln, Monoidringe, Determinantenringe und die homogenen
Koordinatenringe von geometrischen Regelfl"achen
wird die affine Klassengruppe explizit berechnet.
Ferner werden Zusammenh"ange zwischen ger"aumigen und koaffinen Divisoren
und das Verh"altnis von affiner Klassengruppe zur numerischen Klassengruppe
auf glatten projektiven Variet"aten bzw. den affinen Kegeln dar"uber
untersucht.
\par
\bigskip
\noindent
\subsection{Affinit"atseigenschaften von linear "aquivalenten Divisoren}
\par
\bigskip
\noindent
Sei $A$ ein Krullbereich. Es seien
${\bf p}_i,\, i \in I,$ die Menge der Primideale
der H"ohe eins, $B_i=A_{{\bf p}_i} ,\, i \in I,$
die zugeh"origen diskreten Bewertungsringe und $\nu_i,\, i \in I,$ die
zugeh"origen Bewertungen.
Eine rationale Funktion $0 \neq q \in Q(A)$ besitzt dann zu jedem $i \in I$
die Ordnung $\nu_i(q) \in {\bf Z} $, die man im positiven Fall die
Nullstellenordnung und im negativen Fall die Polstellenordnung nennt.
In fast allen diskreten Bewertungsringen ist $q$ eine Einheit,
also $\nu_i(q)=0$.
Die formale Summe $\sum_{i \in I} \nu_i(q) \in {\bf Z}^{(I)}$
hei"st der Hauptdivisor von $q$.
Er ist zerlegt in einen effektiven und einen negativen Summanden, den man den
Nullstellen- bzw. den Polstellendivisor nennt.
Die Restklassengruppe der Divisoren modulo
den Hauptdivisoren nennt man die {\it Divisorenklassengruppe} ${\rm DKG}\, A$,
sie ist entsprechend f"ur jedes integre, normale, noethersche, separierte Schema
definiert, (siehe \cite{haralg}, II.6).
Divisoren mit gleicher Klasse nennt man linear "aquivalent.
Insbesondere sind zu einer
rationalen Funktion Pol- und Nullstellendivisor zueinander "aquivalent.
\par\bigskip
\noindent
In einem Krullbereich gibt es zu jeder Vorgabe von ganzen Zahlen
an endlich vielen Primidealen der H"ohe eins
eine Funktion, deren Hauptdivisor diese Ordnungen annimmt und au"serhalb
dieser Stellen nicht negative Ordnung hat, siehe etwa \cite{fossum} , I.5.8.
Insbesondere gibt es zu jedem effektiven Divisor einen "aquivalenten
effektiven Divisor ohne gemeinsame Komponente.
\begin{prop}
Sei $D$ ein effektiver Divisor in einem Krullbereich $A$. Dann gibt es
zu $D$ einen "aquivalenten effektiven Divisor $D'$ ohne gemeinsame
Komponente mit $D$, dessen Komplement affin ist.
\end{prop}
{\it Beweis}.
Sei $E$ ein zu $D$ "aquivalenter, effektiver Divisor ohne gemeinsame
Komponente. Aufgrund des Satzes "uber die Primvermeidung,
siehe \cite{scheja}, Lemma 58.2,
gibt es eine offene Umgebung $D(g)$ der Primdivisoren von $E$, in der die 
Primdivisoren von $D$ nicht liegen.
$X-{\rm supp}\, E$ ist in einer Menge der Kodimension eins affin,
d.h. zu jedem Primdivisor gibt es eine Umgebung $D(f)$ so,
dass $({\rm supp}\, E) \cap D(f)$ in $D(f)$ affines Komplement besitzt.
Dann gibt es auch ein $D(f)$ mit dieser Eigenschaft, das die Primdivisoren
von $D$ enth"alt.
Der Divisor $D'=E\cdot(fg)$ ist zu $D$ "aquivalent und
besitzt keine Komponente von $D$.
Es besitzt $({\rm supp}\, E) \cap D(fg)$ affines Komplement in $D(fg)$ und daher
besitzt ${\rm supp}\, (E \cdot (fg))= {\rm supp}\, E + V(fg) $ affines
Komplement in ${\rm Spek}\, A$ . \hfill $\Box$
\par\bigskip
\noindent
Ab dem n"achsten Abschnitt wird die Eigenschaft eines Divisors untersucht, dass
alle effektiven Vertreter affines Komplement besitzen. Hier untersuchen wir
noch die Affinit"at des Definitionsbereiches einer rationalen Funktion $q$
und die Affinit"at der durch $q$ definierten
meromorphen Funktion nach ${\bf P}_K^1$.
\begin{prop}
Sei $A$ ein Krullbereich und $q$ eine rationale Funktion mit dem Hauptdivisor
$(q)=m_1{\bf p}_1+...+m_r{\bf p}_r+n_1{\bf q}_1+...+n_s{\bf q}_s$ mit
$m_i >0$ und $n_j<0$..
Dann gelten folgende Aussagen.
\par\smallskip\noindent
{\rm (1)}
${\rm Def}\, (q)=
D({\bf q}_1) \cap ... \cap D({\bf q}_s) $ und
${\rm Def}\, (q^{-1})=D({\bf p}_1) \cap ... \cap D({\bf p}_r)$.
\par\smallskip\noindent
{\rm (2)}
${\rm Def}^{\rm mer}\, (q)= D({\bf b})$ mit
${\bf b}:={\bf p}_1 \cap ... \cap {\bf p}_r+{\bf q}_1 \cap ... \cap {\bf q}_s $.
Insbesondere ist $q$ als meromorphe Funktion in der Kodimension $1$ definiert.
\par\smallskip\noindent
{\rm (3)}
F"ur $a \in A$ ist
${\rm Def}^{\rm mer}\, (q) = {\rm Def}^{\rm mer}\, (q+a)$
\par\smallskip\noindent
{\rm (4)}
Sei $V(q-a)$ die Nullstelle von $q-a$ in ${\rm Def}\, (q)$.
Dann liegt $V({\bf b})$ im Abschluss von $V(q-a)$.
\end{prop}
{\it Beweis}.
(1) Sei $x \in {\rm Def}\, (q)$.
D.h. $q \in {\cal O}_x =\bigcap_{{\bf p}_i \subseteq {\bf p}_x} B_i $
und $q$ geh"ort zu allen Bewertungsringen zu Hyperfl"achen,
auf denen $x$ liegt.
Wegen $q \not\in B_{{\bf q}_j}$ ist $x \not\in V({\bf q}_j)$
f"ur $j=1,...,s $.
Ist umgekehrt $x \not\in {\rm Def}\, (q)$, so ist
$x \not\in D({\bf q}_1) \cap ... \cap D({\bf q}_s)$, da auf dieser
Menge $q$ definiert ist wegen
$q \in \bigcap_{{\bf p}_i \neq {\bf q}_j} B_i=
\Gamma(D({\bf q}_1) \cap ... \cap D({\bf q}_s) ,{\cal O}_X)$.
\par\smallskip\noindent
(2) folgt sofort aus (1), (3) ist eine Wiederholung von 1.4.1
\par\smallskip\noindent
Wegen (3) und wegen ${\rm Def}\, (q) ={\rm Def}\, (q+a)$ kann man in (4) $a=0$
annehmen. Es ist dann $V(q) \subseteq {\rm Def}\, (q)$ gegeben durch den
Abschluss der Punkte ${\bf p}_1,...,{\bf p}_r$ und zum Abschluss
in ${\rm Spek}\, A$ geh"ort
$V({\bf p}_1 \cap ... \cap {\bf p}_r) \supseteq V({\bf b}) \, .$ \hfill $\Box$
\par\bigskip
\noindent
Jedes Hyperfl"achenkomplement in einem affinen normalen Schema ist
Definitionsbereich einer rationalen Funktion,
da es rationale Funktionen gibt, die genau an den vorgegebenen Stellen
einen Pol haben.
Ist $A$ eine $R-$Algebra (etwa $R=A$), so ist genau dann die
meromorphe Abbildung
$$\varphi=\varphi_q: {\rm Def}^{\rm mer}\, (q) \longrightarrow {\bf P}_R^1 $$
ein affiner Morphismus,
wenn sowohl ${\rm Def}\, (q)$ als auch ${\rm Def}\, (q^{-1})$ affin ist.
Ist diese Abbildung durch $T_0/T_1 \longmapsto q$ festgelegt,
so sind n"amlich ${\rm Def}\, (q)= \varphi^{-1}(D_+(T_1))$
und ${\rm Def}\, (q^{-1})= \varphi^{-1}(D_+(T_0))$ Urbilder einer affinen
"Uberdeckung von ${\bf P}^1_R$.
Allgemeiner ist f"ur $a \in R$
\begin{eqnarray*}
\varphi ^{-1}(D_+(T_0-aT_1)) &=& 
\varphi ^{-1}(D_+(T_0)\cap D_+(T_0-aT_1)) \\
&\ & \cup\, \,  \varphi ^{-1}(D_+(T_1)\cap D_+(T_0-aT_1)) \\
&=& {\rm Def}\, (q^{-1})\cap D(1-aq^{-1})
\, \cup \, {\rm Def}\, (q) \cap D(q-a) \, .
\end{eqnarray*}
Das ist das Komplement in ${\rm Def}^{\rm mer}\, (q)$ der oben
definierten Menge $V(q-a)$.
Wegen (4) ist dieses Komplement auch gleich dem Komplement
von $\overline{V(q-a)}$
in ${\rm Spek}\, A$.
Sind also ${\rm Def}\, (q)$ und ${\rm Def}\, (q^{-1})$ affin,
so ist f"ur jedes $a \in A$ das Komplement von $\overline{V(q-a)}$
affin.
\par\smallskip\noindent
$q$ ist als meromorphe Funktion auf ganz $X$ definiert genau dann,
wenn Z"ahler- und Nennerideal das Einheitsideal erzeugen,
in diesem Fall ist die Abbildung
${\rm Spek}\, A \longrightarrow {\bf P}_R^1$ affin, vergl. vor 1.4.5.
Ist $A$ ein eindimensionaler Krullbereich, so ist jede Funktion als
meromorphe Funktion global definiert, was erneut zeigt,
dass jede offene Teilmenge davon wieder affin ist.
\par\smallskip\noindent
Ist hingegen der meromorphe Undefiniertheitsort
${\rm Pol}\, (q) \cap {\rm Pol}\, (q^{-1}) \neq \emptyset$, so hat er
eine Kodimension $\geq 2$. Ist diese Kodimension sogar $\geq 3$,
so kann die Abbildung nicht affin sein,
da dann etwa ${\rm Pol}\, (q)$ in ${\rm Pol}\, (q^{-1})$ eine Kodimension von
zumindest zwei hat, was ${\rm Def}\, (q)$ als nicht affin erweist.
\par\bigskip
\noindent
{\bf Beispiel}
Zu einer rationalen Funktion kann ${\rm Def}\, (q)$ affin sein, aber
${\rm Def}\, (q^{-1})$ nicht.
Betrachte etwa auf dem ${\bf P}^2$ eine (nicht-konstante) meromorphe Funktion
$q$, dann schneiden sich Pol- und Nullstellendivisor in einzelnen Punkten,
wo die meromorphe Funktion nicht definiert ist.
Bl"ast man nun einen Punkt auf, der auf dem Polstellendivisor liegt,
aber nicht auf dem Nullstellendivisor, so besteht auf der Aufblasung $X$
der Polstellendivisor von $q$ aus der
totalen Transformierten des Polstellendivisors unten und enth"alt die
exzeptionelle Faser. Das Komplement, also ${\rm Def}\, (q)$, ist dann
unver"andert affin,
im Komplement des Nullstellendivisors liegt aber die projektive exzeptionelle
Gerade, daher ist ${\rm Def}\, (q^{-1})$ nicht affin.
\par\bigskip
\noindent
{\bf Bemerkung}
Ist $A$ eine $R-$Algebra und
$q:{\rm Spek}\, A \supseteq {\rm Def}^{\rm mer}\, (q)
\longrightarrow {\bf P}_R^1$
ein affiner Morphismus, so ist auch die Abbildung
$q \times id : {\rm Def}^{\rm mer}\, (q) \times_R {\bf P}_R^1
\longrightarrow {\bf P}^1_R \times_R {\bf P}^1_R$ affin.
Da das Komplement der Diagonalen von ${\bf P}^1_R \times_R {\bf P}^1_R$
affin ist, ist
wegen $(q \times id)^{-1}(\triangle)= \Gamma_q	$
dann auch das Komplement des Graphen affin.
Diese Situation wird in 3.5 allgemeiner verfolgt.
\par
\bigskip
\noindent
\subsection{Koaffine Divisorenklassen}
\par
\bigskip
\noindent
{\bf Definition}
Sei $X$ ein noethersches, separiertes, integres, normales Schema.
Wir nennen einen Weil-Divisor $D$ bzw. seine Klasse {\it koaffin},
wenn jeder effektive Vertreter der Klasse ein affines Komplement besitzt.
\par\bigskip
\noindent
Es wird also verlangt, dass f"ur jede Funktion
$0 \neq q \in K=K(X)$ mit $D'=D+(q) \geq 0$ gilt,
dass der Tr"ager von $D'$, der die
Vereinigung der irreduziblen Hyperfl"achen ist, an denen $D'$ einen
positiven Wert hat, ein affines Komplement besitzt.
\begin{prop}
Sei $X$ ein noethersches, separiertes, integres, normales Schema.
Dann gelten folgende Aussagen.
\par\smallskip\noindent
{\rm (1)}
Genau dann ist der triviale Divisor koaffin, wenn $X$ affin ist.
\par\smallskip\noindent
{\rm (2)}
Ist ein positives Vielfaches von $D$ koaffin, so bereits $D$ selbst.
\par\smallskip\noindent
{\rm (3)}
In einem affinen Schema ist die Eigenschaft, ein koaffiner Divisor zu
sein, lokal.
Insbesondere sind Divisoren, von denen ein Vielfaches
ein Cartierdivisor ist, koaffin.
\par\smallskip\noindent
{\rm (4)}
Umfasst $W \subseteq X$ die Kodimension eins, so folgt aus der
Koaffinit"at von $D$ auf $W$ die von $D$ auf $X$.
\end{prop}
{\it Beweis}.
Zu (1). Ist der triviale Divisor koaffin, so ist insbesondere auch
der Einsschnitt, also der triviale Divisor, der leeren Tr"ager hat, koaffin,
also muss $X$ selbst affin sein.
Ein Hauptdivisor wird gegeben durch ein $(q)$ zu einer
rationalen Funktion. Wenn er effektiv ist, so ist $q$ global definiert, also 
$q \in \Gamma(X,{\cal O}_X)=A$, da $X$ normal ist.
Der Tr"ager des Divisors ist dann einfach die
Nullstelle dieser globalen Funktion.
Ist nun $X$ affin, so besitzen diese affines Komplement.
\par\smallskip\noindent
Zu (2). Sei $D+(q)$ effektiv. Dann ist auch $n(D+(q))=nD+(q^n)$ effektiv mit
gleichem Tr"ager, woraus bei $nD$ koaffin folgt, dass auch $D$ koaffin ist.
\par\smallskip\noindent
Zu (3). Sei $X={\rm Spek}\, A$.
Ist $D$ nicht koaffin, so gibt es $E=D+(q)$ effektiv mit nicht affinem
Komplement. Dann gibt es auch einen Punkt $x \in X$, wo der Tr"ager
von $E_x$ nicht affines Komplement besitzt,
und dieser wird beschrieben durch
$E_x=(D+(q))_x=D_x+(q)_x$,
also ist auch $D_x$ nicht koaffin.
\par\smallskip\noindent
Sei umgekehrt die Divisorklasse $D_x$ in ${\cal O}_x$ nicht koaffin, und es sei
$D_x +(q)_x$ effektiv in ${\rm Spek}\, {\cal O}_x$ mit nicht affinem
Komplement.
Wir ersetzen $D$ durch $E=D+(q)$ und haben also einen Divisor,
der in $x$ effektiv und dessen Komplement nichtaffin ist.
Sei $\bf p$ ein Primdivisor von $E$, der nicht durch $x$ verl"auft.
Dann gibt es eine globale Funktion $f \in A$, die in $x$ eine Einheit ist
und in $\bf p$ eine Nullstelle hat (${\bf p} \not \subseteq {\bf p}_x$).
Mit solchen Funktionen gelangt man zu einem zu $E$ "aquivalenten,
effektiven Divisor, dessen Lokalisierung an $x$ sich nicht "andert, und
dessen Komplement in $X$ daher ebenfalls nicht affin sein kann.
\par\smallskip\noindent
F"ur den Zusatz k"onnen wir $A$ als lokal annehmen,
dann ist also ein Vielfaches von $D$ bereits ein Hauptdivisor, der koaffin ist,
und dies gilt nach (2) dann auch f"ur $D$ selbst.    
\par\smallskip\noindent
Zu (4). Sei $D$ effektiv. Es ist stets
$X-{\rm supp}_X\, D \supseteq W-{\rm supp}_W \,  D$.
Bei $U=W-{\rm supp}_W \, D$ affin
gilt aber auch die Umkehrung, da
$W-{\rm supp}_W \,  D=(X-{\rm supp}_X\,  D) \cap W$ ist
und $U=V \cap W$ nach 1.12.2 nur dann affin sein kann, wenn
$V$ affin ist. Damit ist also $X-{\rm supp}_X\, D$ affin. \hfill $ \Box$
\par\bigskip
\noindent
{\bf Bemerkung}
Besitzt ein Divisor "uberhaupt keinen effektiven Vertreter, so
gilt er nach unserer Definition als koaffin.
Im quasiaffinen hat jeder Divisor effektive Verteter, eine Polstelle des
Divisors kann man ja durch Multiplikation mit einer global definierten
Funktion, die darin eine Nullstelle hat, beheben.
\par\bigskip
\noindent
{\bf Beispiel}
In einem eigentlichen Schema "uber einem K"orper $K$
ist der zu einer ger"aumigen invertierbaren Garbe
geh"orende Weildivisor koaffin, siehe Satz 2.4.1.
Dies gilt dann zun"achst auch f"ur alle positiven Potenzen, aber auch
f"ur alle negativen, da diese keine Schnitte besitzen, es sei denn $\cal L$
ist trivial. In jedem Fall erf"ullt ein ger"aumige Garbe die Eigenschaft,
dass alle Vielfachen davon koaffin oder trivial sind.
\par\bigskip
\noindent
Sei $D$ ein Weildivisor des normalen Schemas $X$. 
Ein Punkt $x \in X$ hei"st {\it Basispunkt} der Divisorklasse, wenn f"ur alle
effektiven Vertreter $D'$ von $D$ gilt $x \in {\rm supp}\, D'$. Im folgenden
Satz taucht die Voraussetzung auf, dass ein Divisor keinen Punkt der H"ohe eins
als Basispunkt haben darf, was im Fall eines affinen Schemas $X$ stets erf"ullt
ist.
\par\smallskip\noindent
Eine wichtige Rolle spielt die zum Weildivisor geh"orende
Untergarbe der konstanten Garbe des Funktionenk"orpers,
die durch die Zuordnung
$$ U \longmapsto {\cal L}_D(U)
=\{ q \in K : D+(q) \geq 0 \, \mbox{ auf }\, U \} $$
gegeben ist.
Die Garbeneigenschaft ist dabei trivial, ebenso, dass es sich dabei um einen
${\cal O}_X-$Modul handelt, der allerdings nur im lokal faktoriellen
Fall invertierbar ist.
Ist $x$ ein Punkt der H"ohe eins, so wird die Garbe dort durch eine
einzige Funktion erzeugt und ist daher dort invertierbar. Ist $x \in U=X-D'$
kein Basispunkt, ($D'$ effektiver Vertreter), so wird auf $U$ die Garbe
durch eine Funktion beschrieben, und auf $U$ ist die Garbe invertierbar.
Au"serhalb des Basisortes k"onnen wir also eine Divisorklasse mit einer
invertierbaren Garbe identifizieren. Im Allgemeinen ist die auf ganz $X$
definierte zugeh"orige Garbe nur ein reflexiver Modul.
Die Elemente aus ${\cal L}_D(U)=\{ q \in K : D+(q) \geq 0  \}$ liefern alle
effektiven, zu $D$ "aquivalenten Divisoren.
Eine invertierbare Garbe $\cal L$ auf einem Schema $X$ nennen
wir koaffin, wenn jede Nullstelle $V(s)$ zu einem globalen
Schnitt $s \in \Gamma(X,{\cal L})$
affines Komplement besitzt
\begin{satz}
Sei $X$ ein noethersches, separiertes, integres, normales Schema mit globalem
Schnittring $A$ und $D$ ein Weildivisor. 
Wir setzen voraus, dass $D$ keine fixen Komponenten hat, also in der
Kodimension eins keine Basispunkte besitzt. Das Komplement des Basisortes
sei mit $W$ bezeichnet, und ${\cal L}={\cal L}_D \vert_W$
sei die zugeh"orige invertierbare Garbe auf $W$.
Dann sind folgende Aussagen "aquivalent.
\par\smallskip\noindent
{\rm (1)}
$D$ ist koaffin auf $X$.
\par\smallskip\noindent
{\rm (2)}
$D$ ist aufgefasst auf $W$ koaffin.
\par\smallskip\noindent
{\rm (3)}
Zu jedem Schnitt $s$ auf
$W$ von $\cal L$ besitzt
die Nullstellenmenge $V(s)$ in $W$ affines Komplement.
\par\smallskip\noindent
{\rm (4)}
Zu jedem auf $W$ basispunktfreien linearen System
$(s_0,...,s_n)$
ist die Abbildung $W \longrightarrow {\bf P}_A^n$ ein affiner Morphismus.
\par\smallskip\noindent
{\rm (5)}
Der $A-$Modul ${\cal L}_D(X)$ besitzt ein
$A$-Erzeugendensystem aus Schnitten mit affinem
Komplement.
\end{satz}
{\it Beweis}.
Da $W$ die Kodimension eins umfasst, entsprechen sich die Divisoren auf
$X$ und auf $W$, und ebenso entsprechen sich "aquivalent und effektiv.
Sei $D'$ ein effektiver Vertreter von $D$. Dann ist nat"urlich
$$X-{\rm supp}_X\,   D' \supseteq (X-{\rm supp}_X\,   D')\cap W
= W-{\rm supp}_W \,   D' \, .$$
Da der Basisort $X-W$ ganz in jedem effektiven $D'$ enthalten ist,
stimmen die beiden Mengen "uberein, was (1) $ \Longrightarrow $ (2) liefert.
Die Umkehrung wurde schon in 3.2.1 gezeigt.
\par\smallskip\noindent
Die "Aquivalenz von (2) und (3).
Die Nullstellenmenge von $s \in {\cal L}_D(W)$ ist einfach der Tr"ager des
zugeh"origen Divisors.
\par\smallskip\noindent
Von (3) nach (4).
Es ist also $\cal L$ eine basispunktfreie Garbe und
$s_0,...,s_n \in {\cal L}(U)$
eine basispunktfreie Familie von Schnitten. Dadurch ist eine
Abbildung $W \longrightarrow {\bf P}^n_A$ gegeben, und dabei ist
$f^{-1}(D_+(x_i) = W_{s_i}=W-V(s_i)$ affin nach Voraussetzung,
also ist der Morphismus affin.
\par\smallskip\noindent
Von (4) nach (5) ist klar, da es Erzeugendensysteme gibt und diese
basispunktfrei auf $W$ sind.
\par\smallskip\noindent
Von (5) nach (3).
Sei $0 \neq s \in {\cal L}(W)$ ein beliebiger effektiver Schnitt.
Nach Voraussetzung ist dann $s=a_0s_0+...+a_ks_k$ mit $s_i$ aus dem
Erzeugendensystem. Diese Schnitte $(s_0,...,s_k)$ verl"angert man
zu einer $W$ "uberdeckenden Familie aus dem Erzeugendensystem
und erh"alt eine affine Abbildung von $W$ nach ${\bf P}_A^r$, die
$W_s= \varphi^{-1}(D_+(a_0x_0+...+a_kx_k))$
als affin erweist. \hfill $ \Box $
\par\bigskip
\noindent
{\bf Bemerkung} Ist $X$ affin, so gilt die im Satz gemachte Voraussetzung
"uber den Basisort f"ur jeden Divisor, da es zu jedem Divisor einen
"aquivalenten, effektiven Divisor ohne gemeinsame Komponente gibt.
Das gilt ebenso, wenn $X$ quasiaffin ist.
\par\smallskip\noindent
Die Bedingungen an den Divisor $D$ in 3.2.2
sind in der Regel nicht "aquivalent zu 
\par\smallskip\noindent
(6) Es gibt einen affinen Morphismus $f: W \longrightarrow {\bf P}_A^n$
mit ${\cal L}_D=f^\ast {\cal O}(1)$.
\par\smallskip\noindent
Es kann in der zur"uckgenommenen Garbe Schnitte geben, die durch die
Abbildung nicht erfasst werden.
Sei hierzu $X$ ein quasiaffines regul"ares Schema und $D$ eine
Hyperfl"ache ohne affines Komplement.
$D$ ist dann Nullstelle eines Schnittes in der zugeh"origen
Garbe $\cal L$, die wegen regul"ar "uberall invertierbar ist.
Eine invertierbare Garbe auf einem quasiaffinen Schema ist aber
immer ger"aumig und besitzt eine "uberdeckende Familie von affinen
Teilmengen $X_s,\, s \in {\cal L}(X)$, und mit diesen Schnitten bastelt
man sich einen affinen Morphismus nach ${\bf P}_A^n$.
Ein Beispiel hierzu ist ${\cal L}={\cal O}_W$ mit dem Einsschnitt.
\par\smallskip\noindent
Betrachte die durch $ux-vy$ gegebene Variet"at.
Die rationale Funktion $q=v(v+y)/x(x+u)=u(v+y) / y(x+u)$
hat den Pol $( x,y ) +(x+u)$ und die Nullstelle
$( v,u ) +(v+y)$, die beide affines Komplement haben.
Ihr Schnitt ist eine Kurve $C$.
Dann liegt in
$q: {\rm Spek}\, K[u,x,y,v]/(ux-vy) -C \longrightarrow {\bf P}^1_K$
eine affine Abbildung vor.
Der Divisor
besitzt aber auch den Vertreter $V(x,y)$, dessen Komplement nicht affin ist.
\begin{satz}
Sei $p:X'={\rm Spek}\, A' \longrightarrow X={\rm Spek}\, A$ die
Spektrumsabbildung zu einer flachen Ringerweiterung
normaler, noetherscher, integrer Ringe $A \longrightarrow A'$.
Dann ist mit einem koaffinen Divisor $D$ auf ${\rm Spek}\, A$
auch der zur"uckgenommene Divisor $D'=p^*D$ koaffin.
Ist die Abbildung treuflach, so gilt auch die Umkehrung.
\end{satz}
{\it Beweis}.
Zur Definition von $p^*D$ im flachen Fall siehe \cite{bourbaki}, 7.\S 1.10,
mit $X^{(1)}$ sei die Menge der Primideale der H"ohe $1$ bezeichnet.
Sei $W \subseteq X$ der basispunktfreie Ort zu $D$ mit
$X^{(1)} \subseteq W$ und ${\cal L}={\cal L}_D$ invertierbar auf $W$, sei
$W'=p^{-1}(W)$.
Wegen dem going down wird ein Primideal der H"ohe eins
von $X'$ auf den generischen Punkt oder auf ein Primideal der H"ohe eins von 
$X$ abgebildet, also ist $X'^{(1)} \subseteq W'$. Wegen Punkt (4) aus
der Proposition 3.2.1 gen"ugt es, die Koaffinit"at von ${\cal L}'$
auf $W'$ zu zeigen.
\par\smallskip\noindent
Wegen der Flachheit des Ringwechsels gilt f"ur $\cal L$ die Gleichung
${\cal L}'(W')={\cal L}(W) \otimes_AA'$,
siehe den Beweis von 1.14.3.
Insbesondere wird ein $A-$Erzeugen\-densystem $s_i,\, i \in I,$ von
${\cal L }(W) (={\cal L}(X))$ zu einem
$A'=\Gamma(W',{\cal O}_{X'})-$Er\-zeu\-gen\-densystem
$s_i'=s_i \otimes 1,\, i \in I,$ von ${\cal L}'(W')$.
Mit $s_i$ sind auch die $s_i'$ Schnitte mit affinem Komplement,
da die Nullstelle von $s_i'$
das Urbild der Nullstelle von $s_i$ ist, was man sieht, wenn man
die Sache lokal betrachtet.
Damit ist das Kriterium (5) des Satzes 3.2.2
erf"ullt und $D'$ ist koaffin.
\par\smallskip\noindent
Sei die Abbildung nun treuflach und $p^*D$ koaffin in ${\rm Spek}\, A'$, und
$E$ sei ein effektiver Vertreter von $D$. Dann ist $p^*E$ ein effektiver
Vertreter von $p^*D$ und dessen Tr"ager ist das Urbild des Tr"agers von $E$
(wegen flach). Dann ist aber mit $X'-{\rm supp}\, p^*E$
aufgrund von 1.14.3 auch $X-{\rm supp}\, E$ affin.   \hfill $ \Box $
\par\bigskip
\noindent
{\bf Bemerkung}
Die beteiligten Schemata m"ussen affin sein. F"ur die
Kegelabbildung zu einer projektiven Variet"at braucht man weitere
Bedingungen, damit der zur"uckgenommene Divisor koaffin ist, siehe
den Abschnitt 3.8.
\par\smallskip\noindent
Ohne die Flachheit gilt die Aussage vermutlich nicht, da dann der
zur"uckgenommene Modul Schnitte haben kann, die nicht von unten her
kontrollierbar sind.
\begin{satz}
Sei $A$ ein noetherscher, normaler, integrer Ring und $D$ eine koaffine
Divisorklasse von $X={\rm Spek}\, A$. Dann ist auch $D$ als Divisorklasse
aufgefasst in der affinen Gerade
$X \times {\bf A} ={\rm Spek}\, A[T]$ koaffin. Sind in $A$ alle
Hyperfl"achen koaffin, so auch in $A[T]$.
\end{satz}
{\it Beweis}.
Die erste Aussage folgt direkt aus 3.2.3, zur zweiten muss man lediglich
beachten, dass jede Divisorklasse von $A[T]$ "aquivalent zu
einer von unten ist.   \hfill $ \Box$
\par\bigskip
\noindent
{\bf Bemerkung}
Die Aussage gilt entsprechend auch f"ur einen affinen Raum
${\rm Spek}\, A[T_1,...,T_0]$ "uber ${\rm Spek}\, A$.
\par\smallskip\noindent
Sind in $A$ alle Hyperfl"achen
koaffin, so "ubertr"agt sich das nicht auf ${\bf P}^n_A$, da darin
Hyperfl"achen, die von unten kommen, nicht affines Komplement haben.
Allerdings haben Hyperfl"achen, die die Basis ${\rm Spek}\, A$ dominieren,
affines Komplement: dies kann man n"amlich im punktierten Kegel dar"uber
testen und das Komplement dort ist gleich dem Komplement im vollen
Kegel, da die Spitze zum Urbild von dominanten Hyperfl"achen  geh"ort.
\par\bigskip
\noindent
{\bf Beispiel}
Ohne die Voraussetzung der Normalit"at gilt die Aussage von 3.2.4 nicht.
Sei hierzu $C={\rm Spek}\, A$ eine affine Kurve "uber einem K"orper $K$
mit einem Doppelpunkt, sagen wir die Kurve, die entsteht, wenn
man auf der affinen Gerade zwei Punkte identifiziert.
Als eindimensionales noethersches Schema haben da alle Hyperfl"achen
affines Komplement. Betrachte $S=C \times {\bf A}_K^1={\rm Spek}\, A[T]$.
Diese Fl"a\-che entsteht aus ihrer Normalisierung ${\bf A}_K^2$ durch
Identifizieren zweier paralleler Geraden. Eine dazu schiefe (nicht parallele,
nicht senkrechte) Gerade in ${\bf A}^2_K$ hat dann auf der Fl"a\-che eine
Bildkurve, deren Komplement nicht affin sein kann, da das Urbild isolierte
Punkte besitzt.
\begin{kor}
Sei $A$ ein zweidimensionaler, integrer, normaler, exzellenter Ring.
Dann besitzt in ${\rm Spek}\, A[T]$ der Schnitt zweier irreduzibler Fl"a\-chen
keine isolierten Punkte.
\end{kor}
{\it Beweis}. Nach dem Satz von Nagata 1.15.2 besitzt in ${\rm Spek}\, A$
jede Kurve affines Komplement. Dann gilt dies nach dem Satz auch in der
affinen Gerade dar"uber, dort besitzt also jede Fl"a\-che affines Komplement.
Der Schnitt mit einer anderen Fl"a\-che muss daher leer sein oder aber die
Kodimension eins haben. \hfill $\Box$
\par\bigskip
\noindent
{\bf Beispiel}
Ist $A$ zweidimensional, aber nicht normal, so gibt es in der affinen Gerade
dar"uber Fl"a\-chen, die sich lediglich in einem Punkt schneiden.
Sei hierzu $C$ wieder die nodale Kurve von oben und betrachte jetzt
$C \times {\bf A}_K ^2=S \times {\bf A}^1_K$.
Diese dreidimensionale Variet"at entsteht aus
${\bf A}^3_K$ durch Identifizieren zweier paralleler Ebenen.
Man findet nun disjunkte Fl"a\-chen im ${\bf A}_K^3$ derart, dass die
jeweiligen Schnitte mit den Ebenen nur einzelne Punkte gemeinsam haben.
Um konkret zu sein, sei in ${\bf A}^3_K$ die eine Fl"a\-che $S_1$ durch
$z-y=0$ und die andere $S_2$ durch $xz-xy+1=0$ gegeben.
Die Fl"a\-chen sind disjunkt,
da die beiden Funktionen die Eins darstellen.
Unsere Variet"at entstehe daraus durch Identifizieren der Ebenen $z=0$
und $z=1$. Es ist einerseits
$S_1 \cap \{ z =0 \} =V(y,z)$ und $S_1 \cap \{ z =1 \}= V(y-1,z-1)$
und andererseits
$S_2 \cap \{ z=0 \} =V(-xy+1,z)$ und $S_2 \cap \{ z=1 \}=V(x-xy+1,z-1)$.
Das eine sind Geraden, das andere Hyperbeln, weshalb nur Punktschnitte
vorliegen k"onnen. Identifiziert man die beiden Ebenen, so haben die
Bildfl"achen von $S_1$ und $S_2$ nur einzelne Punkte gemeinsam.
\par
\bigskip
\noindent
\subsection{Affin-triviale Divisoren und affine Klassengruppe}
\par
\bigskip
\noindent
Wie Beispiele zeigen, siehe Abschnitt 3.6, muss ein Vielfaches
eines koaffinen Divisors nicht koaffin sein.
Unser Ziel ist es, diejenigen Divisoren in einem noetherschen, integren,
separierten, normalen Schema zu charakterisieren, die
Koaffinit"atseigenschaften von Divisoren (bei Addition) nicht
ver"andern und zwar so, dass sie eine Untergruppe der Divisorenklassengruppe
bilden. Hierzu sind verschiedene Ans"atze denkbar,
wir beschr"anken uns jedoch auf die Durchf"uhrung folgender Variante,
bei der mit folgender Eigenschaft einer Divisorenklasse $E$ gearbeitet wird.
\par\smallskip\noindent
(K) F"ur alle $n \in {\bf Z}$ ist $nE$ koaffin oder trivial
(d.h. ein Hauptdivisor).
\par\smallskip\noindent
Ist $X$ affin, so ist nat"urlich der triviale Divisor koaffin, so dass sich
die hintere Alternative er"ubrigt.
\par\bigskip
\noindent
{\bf Definition}
Sei $X$ ein noethersches, separiertes, integres, normales Schema
und $D$ ein Divisor.
Wir sagen, dass $D$ {\it affin-trivial} ist,
wenn er die Eigenschaft (K) erh"alt,
wenn also folgendes gilt:
F"ur jede Divisorklasse $E$, die (K) erf"ullt, erf"ullt auch
$D+E$ diese Eigenschaft.
\par\smallskip\noindent
Eine Divisorklasse nennen wir affin-trivial, wenn einer und damit jeder
Vertreter affin-trivial ist.
\begin{satz}
Sei $G \subseteq {\rm Div}\, X$
die Menge aller affin-trivialen Divisoren. Dann gelten folgende Aussagen.
\par\smallskip\noindent
{\rm (1)} $G$ ist eine Gruppe, die alle Hauptdivisoren umfasst.
\par\smallskip\noindent
{\rm (2)} Ein affin-trivialer Divisor besitzt die Eigenschaft {\rm (K)}
und ist insbesondere selbst koaffin oder trivial.
\par\smallskip\noindent
{\rm (3)} Ist $X$ affin und $kD$ affin-trivial, $k \neq 0$, so auch $D$.
\par\smallskip\noindent
{\rm (4)} Ist $X$ affin und ist $D$ in jedem Punkt $x \in X$ affin-trivial,
so ist $D$ affin-trivial.
Speziell ist jeder Cartierdivisor affin-trivial.
\end{satz}
{\it Beweis}.
Zu (1). F"ur einen Hauptdivisor ist die Bedingung erf"ullt,
da die Eigenschaft eine solche der Klasse ist.
Sind $D$ und $D'$ affin-trivial, und erf"ullt $E$ die Eigenschaft (K), so auch
$D'+E$ und damit auch $D+D'+E$.
Sei $D$ affin-trivial, wir m"ussen zeigen, dass auch $-D$ affin-trivial ist.
Mit $E$ erf"ullt auch $-E$ die Eigenschaft (K) und nach
Voraussetzung $D-E$, also auch $-D+E$.
\par\smallskip\noindent
Zu (2). Sei $D$ affin-trivial, dann ist nach Definition mit $E=0$ auch
$nD+nE=nD$ koaffin oder aber trivial.
\par\smallskip\noindent
Zu (3). Sei $kD$ affin-trivial und $E$ habe die Eigenschaft (K).
Wir k"onnen $k >0$ annehmen. Mit $E$ hat auch
$kE$ die Eigenschaft (K) und damit auch
$k(D+E)=kD+kE$.
Zu jedem $n \in {\bf Z}$ ist dann also $kn(D+E)$ koaffin, und dies
"ubertr"agt sich dann zur"uck auf $n(D+E)$ nach 3.2.1.
\par\smallskip\noindent
Zu (4). 
In einem affinen Schema $X$ ist ein Divisor genau dann koaffin, wenn er
lokal in jedem Punkt koaffin ist, und er erf"ullt genau dann die Eigenschaft
(K), wenn er dies lokal in jedem Punkt erf"ullt.
Ist nun $D_x$ f"ur alle $x \in X$ affin-trivial und
$E$ ein Divisor mit der Eigenschaft (K),
so besitzt $(D+E)_x=D_x+E_x$ die Eigenschaft (K) und damit
auch $D+E$.
\par\smallskip\noindent
Ein Cartierdivisor ist lokal ein Hauptdivisor und daher aufgrund
dieses Lokalit"atsprinzip affin-trivial. \hfill $\Box$
\par\bigskip
\noindent
{\bf Definition}
Sei $X$ ein noethersches, separiertes, integres, normales Schema.
Dann nennen
wir die Restklassengruppe der Divisoren nach der
Untergruppe der affin-trivialen
Divisoren die {\it affine Klassengruppe} von $X$ und bezeichnen sie mit
$$ {\rm AKG}\, X \, .$$
\par\bigskip
\noindent
{\bf Bemerkung}
Wegen der Aussage (1) des Satzes 3.3.1 ist die affine Klassengruppe
${\rm AGK}\, X= {\rm Div}\, X /G$
eine Restklassengruppe der Divisorenklassengruppe
${\rm DKG}\, X={\rm Div}\, X /{\rm PDiv}\, X$.
Im affinen Fall ist sie sogar eine
Restklassengruppe der Weildivisoren modulo der Cartierdivisoren
${\rm DKG}\,X/{\rm Pic}\, X$.
Wegen der Aussage (3) ist die affine Klassengruppe eines affinen
Schemas torsionsfrei, was im allgemeinen nicht sein muss,
siehe das Beispiel nach 3.8.2.
Man kann vermuten, dass sie im Fall eines lokalen Ringes endlich erzeugt ist.
Zu einem Divisor nennen wir das Bild in der affinen Klassengruppe einfach die
{\it Affinit"atsklasse} des Divisors. "Ahnliche Sprechweisen verstehen sich von
selbst. Beispielsweise ist die Eigenschaft (K)
jetzt eine Eigenschaft der Affinit"atsklasse.
\begin{satz}
Sei $X$ ein noethersches, separiertes, integres, normales Schema.
Das Verschwinden der affinen Klassengruppe kann man wie folgt
charakterisieren.
\par\smallskip\noindent
{\rm (1)}
Die affine Klassengruppe von $X$ verschwindet genau dann, wenn f"ur jeden
effektiven Divisor auf $X$ gilt, dass er trivial ist oder aber
das Komplement des Tr"agers affin ist.
\par\smallskip\noindent
{\rm (2)}
Ist $X$ affin, so besitzt jede Hyperfl"ache genau dann
affines Komplement, wenn die affine Klassengruppe von $X$ verschwindet.
\par\smallskip\noindent
{\rm (3)}
Ist $X$ das punktierte Schema zu einem lokalen Ring
oder ein eigentliches Schema "uber einem K"orper, so ist genau dann
${\rm AKG}\, X=0$, wenn darin das Komplement einer
jeden echten {\rm (}d.h. nichtleeren{\rm )} Hyperfl"ache affin ist.
\end{satz}
{\it Beweis}. 
Ist in $X$ jeder effektive Divisor trivial oder koaffin, so ist f"ur jeden
Divisor die
Bedingung f"ur die affine Trivialit"at automatisch erf"ullt und die Gruppe
verschwindet.
Ist umgekehrt jeder Divisor affin-trivial, so ist jeder trivial oder
koaffin. 
\par\smallskip\noindent
Zu (3) ist lediglich zu zeigen, dass der Tr"ager eines effektiven
Hauptdivisors leer ist oder affines Komplement besitzt.
Im eigentlichen Fall ist nur der Nulldivisor ein effektiver Hauptdivisor.
F"ur das punktierte Schema eines lokalen Ringes $A$ ist die Aussage bei
${\rm dim}\, A \leq 1$ klar, bei ${\rm dim}\, A \geq 2$ geh"ort
der abgeschlossene Punkt zum Tr"ager eines vom Nulldivisor verschiedenen
effektiven Hauptdivisors, so dass das Komplement im punktierten
Schema gleich dem affinen Komplement im affinen Schema ist.  \hfill $\Box$
\par\bigskip
\noindent
{\bf Bemerkung}
Mit der neuen Terminologie k"onnen wir sagen, dass f"ur einem Ring,
der lokal fastfaktoriell ist, und insbesondere f"ur einen regul"aren Ring
die affine Klassengruppe verschwindet.
Den Satz von Nagata 1.15.2 kann man so ausdr"ucken, dass in einem normalen,
exzellenten, zweidimensionalen Ring die affine Klassengruppe trivial ist.
\par\bigskip
\noindent
Es ist unklar, in welcher Allgemeinheit ein Ringhomomorphismus zu
einem Morphismus der affinen Klassengruppe f"uhrt.
\par\smallskip\noindent
Schon bei einer Nenneraufnahme $A \longrightarrow A_f$
ist nicht klar, ob ein affin-trivialer Divisor $E$ von $A$ auf
$U=D(f)$ affin-trivial bleibt.
$E$ besitzt zwar die Eigenschaft (K), es ist aber nicht klar, ob er einen
Divisor auf $U$, der diese Eigenschaft hat, in einen solchen Divisor
"uberf"uhrt.
\par\smallskip\noindent
Es besteht auch keine Aussicht, analog zur Divisorenklassengruppe eine
exakte Sequenz aufzustellen, die die affine Klassengruppe von $X$ mit
der von $U$ in Verbindung bringt.
\par\smallskip\noindent
Betrachte hierzu etwa eine isolierte dreidimensionale Singularit"at $A$ mit
nicht trivialer affiner Klassengruppe und so, dass es ein Primelement $p$
gibt.
Die Divisorenklassengruppe von $A$ und von $A_p$ stimmen "uberein.
$U=D(p)$ ist regul"ar und hat triviale affine Klassengruppe.
Das Primideal $(p)$ ist in der Divisorenklassengruppe
und erst recht in der affinen Klassengruppe von $A$ trivial, der nicht-triviale
Kern von ${\rm AKG}\, A \longrightarrow {\rm AKG}\, A_p=0$
wird also nicht durch $(p)$ beschrieben.
Ein Beispiel hierf"ur ist $A=K[x,y,u,v]/(ux-vy)$ mit $p=u+x$.
\begin{satz}
Sei $A \longrightarrow A' $ eine treuflache Erweiterung
noetherscher normaler Integrit"atsbereiche $A$ und $A'$.
Ist dann der Morphismus der Divisorenklassengruppen
$p^*:{\rm DKG}\, A \longrightarrow {\rm DKG}\, A'$ surjektiv,
so ist mit $D$ auch $D'=p^*D$
affin-trivial und die affinen Klassengruppen stimmen "uberein.
\end{satz}
{\it Beweis}. Sei $D$ aus ${\rm Spek}\, A$ affin-trivial, und $E'$ erf"ulle die
Eigenschaft (K).
Nach Voraussetzung ist $E'=p^*E$ und nach 3.2.3 muss auch $E$
die Eigenschaft (K) erf"ullen.
Dann ist auch $n(D+E)$ koaffin f"ur alle $n \in {\bf Z}$ 
und dasgleiche gilt f"ur $n(D'+E')$. Also ist $D'$ affin-trivial und es gibt
einen Morphismus ${\rm AKG}\, A \longrightarrow {\rm AKG}\, A'$. Dieser ist
surjektiv, da dies f"ur die Divisorenklassen gilt, aber auch injektiv,
da ein Divisor $D$ schon selbst affin-trivial sein muss, wenn dies f"ur $D'$
unter einer treuflachen Erweiterung gilt.  \hfill $\Box $
\par\bigskip
\noindent
{\bf Bemerkung}
Ohne die im Satz gemachte Voraussetzung, dass alle
Divisorklassen von unten kommen,
ist es nicht klar, ob affin-triviale
Divisoren affin-trivial bleiben und ob es etwa zu einem lokalen normalen
Ring $A$ eine Abbildung von ${\rm AKG}\, A$ nach ${\rm AKG}\, A^{\rm komp}$
gibt.
\begin{kor}
Sei $A$ ein normaler noetherscher Integrit"atsbereich und $A[T]$ der
Polynomring dar"uber.
Dann stimmen die affinen Klassengruppen von $A$ und von $A[T]$
"uberein.
\end{kor}
{\it Beweis}.
Die Erweiterung ist treuflach und die
Divisorenklassengruppen stimmen "uberein.    \hfill $ \Box$
\begin{satz}
Sei $A$ ein normaler noetherscher Integrit"atsbereich und $E$ ein Vektorb"undel
vom Rang $r$ darauf, also ein Schema $E \longrightarrow {\rm Spek}\, A$,
das lokal vom Typ ${\rm Spek}\, A[T_1,...,T_r]$ ist
{\rm (}mit linearen "Ubergangsabbildungen{\rm )}\footnote{Solche
Vektorb"undel korrespondieren zu freien Garben vom Rang $r$,
siehe \cite{haralg}, Exc. II.5.18.}.
Dann stimmen die affinen Klassengruppen von $A$ und $E$ "uberein.
\end{satz}
{\it Beweis}.
Die Abbildung ist treuflach, wir haben zu zeigen, dass jeder Divisor von unten
kommt. Sei $D$ ein irreduzibler Divisor, der ${\rm Spek}\, A$ dominiert.
Auf $\emptyset \neq U \subseteq {\rm Spek}\, A$ sei das B"undel trivial.
Dann gibt es eine rationale Funktion mit $D+(q)=D'$ auf $U \times {\bf A}^r$,
wobei $D'$ ein Divisor von $U \subseteq {\rm Spek}\, A$ ist.
Wir k"onnen damit annehmen, dass $D$ keine in $U \times {\bf A}^r$
liegenden Primdivisoren besitzt.
Dann ist $D$ ein Divisor, der ${\rm Spek}\, A$ nicht dominiert,
solche Divisoren kommen aber direkt von unten.  \hfill $\Box$
\par\bigskip
\noindent
Ein etwas anderer und allgemeinerer Zugang zu den beiden letzten Aussagen
wird durch folgendes Lemma geliefert. Die Voraussetzung, dass Urbilder
irreduzibler Hyperfl"achen wieder irreduzibel sind, ist f"ur jedes
Vektorb"undel, aber auch f"ur Abbildungen vom Produkttyp
$X \times Y \longrightarrow X$ mit irreduziblem $Y$ ("uber $k$ algebraisch
abgeschlossen) erf"ullt. F"ur ein Vektorraumb"undel ist auch die
zweite Voraussetzung, dass die generische Faser faktoriell ist, erf"ullt,
dagegen nicht f"ur allgemeine (affine, glatte)
Kurven $X \times C \longrightarrow X$, da
$C_{K(X)}$ beim Geschlecht $\geq 1$ nicht faktoriell ist.
F"ur die Komplettierung eines lokalen Ringes sind in der Regel beide
Bedingungen nicht erf"ullt.
\begin{lem}
Sei $A \longrightarrow B$
eine treuflache Erweiterung normaler noetherscher Integrit"atsbereiche
derart, dass die Urbilder irreduzibler Hyperfl"achen wieder irreduzibel
sind.
Wird dann die Hyperfl"ache $Y \subseteq {\rm Spek}\, B$ in der generischen
Faser ${\rm Spek}\, B \otimes_AQ(A)$ durch eine Funktion beschrieben,
so ist der Divisor linear "aquivalent zu einem Divisor von unten.
Ist ${\rm AKG}\, A=0$, so ist $Y$ koaffin.
\par\smallskip\noindent
Ist die generische Faser zus"atzlich faktoriell, so ist die
Abbildung der Divisorenklassengruppen surjektiv, und die affinen Klassengruppen
stimmen "uber\-ein.
\end{lem}
{\it Beweis}.
Sei $q \in B_{A^*}=B \otimes_A Q(A)$ die beschreibende Funktion von $Y$.
Der Hauptdivisor von $q$ ist $(q)=Y+n_1Y_1+...n_rY_r$, wobei alle $Y_i$
die generische Faser nicht treffen.
$Y$ ist also linear "aquivalent zu
einem Divisor, dessen Komponenten die generische Faser nicht treffen.
Das Bild von $Y_i$ liegt dann dicht in einer Hyperfl"ache
$H_i$ von ${\rm Spek}\, A$ (da die Abbildung flach ist, kann die H"ohe
nicht steigen)
und damit ist $Y_i \subseteq \varphi^{-1}(H_i)$ und wegen der Voraussetzung gilt
die Gleichheit. Damit ist $Y$ linear "aquivalent zu
einem Divisor, der von unten kommt, und wenn dieser koaffin ist,
ist nach 3.2.3 $Y$ koaffin.
Wenn die generische Faser faktoriell ist, so ist jeder Divisor linear
"aquivalent zu einem von unten.  \hfill $\Box$
\par\bigskip
\noindent
Wir wenden uns der affinen Klassengruppe des formalen Potenzreihenringes
$A[[T]]$ zu. Die generische Faser $A[[T]]_{A^*}$ ist nat"urlich von $Q(A)[[T]]$
verschieden und muss keineswegs faktoriell sein.
Man sagt, dass ein normaler Integrit"atsbereich $A$ eine
{\it diskrete Divisorenklassengruppe} hat, wenn die Abbildung
${\rm DKG}\, A \longrightarrow {\rm DKG}\, (A[[T]])$
bijektiv ist, siehe \cite{fossum}, §19.
Dies ist f"ur eine normale exzellente $k-$Algebra $A$ "uber einem K"orper der
Charakterisik null genau dann der Fall, wenn $A$ 1-{\it rational} ist,
d.h. wenn lokal bei einer Singularit"atenaufl"osung
$X' \longrightarrow {\rm Spek}\, A$ 
gilt $H^1(X',{\cal O}_{X'})=0$, siehe hierzu \cite{binsto2}, 6.1.
1-rationale Singularit"aten besitzen endlich erzeugte Divisorenklassengruppen,
und man kann sie umgekehrt unter gewissen Zusatzbedingungen
dadurch charakterisieren, siehe \cite{binsto2}, 6.3 und 6.5.
Um konkrete Beispiele vor Augen zu haben sei erw"ahnt, dass ein normaler,
lokaler, exzellenter Ring mit Restek"orper der Charakteristik null,
der die Bedingung $S_3$ und $R_2$ erf"ullt,
1-rational ist, \cite{binsto2}, 6.4. Dies wiederum ist beispielsweise
f"ur isolierte Hyperfl"achensingularit"aten der Dimension $\geq 3$ der Fall.
\begin{kor}
Besitzt der normale Integrit"atsbereich $A$
eine diskrete Divisorenklassengruppe, so stimmt die affine
Klassengruppe von $A$ mit der affinen Klassengruppe
des formalen Potenzreihenringes $A[[T]]$ "uberein. \hfill $\Box$
\end{kor}
{\bf Bemerkung}
Es ist offen, ob es bei ${\rm AKG}\, A=0$ Hyperfl"achen in $A[[T]]$ mit
nicht affinem Komplement geben kann.
Es ist auch denkbar, dass es dies sogar bei faktoriellem $A$ geben kann.
Diese Ph"anomene k"onnen zumindest bei der Komplettierung und bei
Grundk"orperwechsel auftreten.
\par\bigskip
\noindent
{\bf Beispiel}
Ein lokaler Ring $A$ kann faktoriell sein, ohne dass die affine Klassengruppe
der Komplettierung verschwindet.
Sei hierzu
$$ A=K[U,V,W][X,Y]/(XY-(U^2-(1+W)V^2))
\mbox{ mit } {\rm char}\, K \neq 2 \, .$$
$A$ ist faktoriell, da in $K[U,V,W] $ das Polynom $F=U^2-(1+W)V^2$
unzerlegbar ist.
"Uber der Komplettierung $K[[U,V,W]]$ 
gilt dagegen 
$$F=(U- \sqrt{1+W}V)(U+ \sqrt{1+W}V) \, .$$
Somit ist in der Komplettierung von $A$
$$ \hat{A}= K[[U,V,W,X,Y]]/(XY-(U- \sqrt{1+W}V)(U+ \sqrt{1+W}V))$$
nicht mehr jede Hyperfl\"ache koaffin,
da darin $( X, U-\sqrt{1+W}V )$ und $(Y,U +\sqrt{1+W}V )$
Primideale der H"ohe eins sind, deren Schnitt durch $(X,Y,U,V)$
beschrieben wird und die Kodimension drei hat.
\par\bigskip
\noindent
{\bf Beispiel}
Auf "ahnliche Weise kann man sich ein Beispiel basteln, bei dem nach
einer Grundk"orpererweiterung aus einem affinen faktoriellen Schema ein
Schema entsteht, bei dem es Hyperfl"achen mit nicht affinem Komplement gibt.
In ${\bf R}[U,V]$ ist $U^2+V^2$ prim und daher ist der
Ring
$$A={\bf R}[U,V,X,Y]/(XY-(U^2+V^2)) $$
faktoriell, es verschwinden also
Divisorenklassengruppe und affine Klassengruppe. Komplexifiziert man diesen
Ring, so erh"alt man
$$A \otimes_{\bf R} {\bf C}={\bf C}[U,V,X,Y]/(XY-(U+iV)(U-iV)) $$
und darin sind $(X,U+iV)$ und $(Y,U-iV)$ Primideale der H"ohe eins mit
einem nicht leeren
Schnitt der H"ohe drei, und die affine Klassengruppe ist $\bf Z$.
Ist der Grundk"orper bereits algebraisch abgeschlossen und von Charakteristik
null, so kann sowas nicht passieren, siehe den Abschnitt 3.5.
\par\smallskip\noindent
Aus diesem Beispiel kann man auch folgendes sehen:
Wenn auf einem affinen Schema $X$ eine Gruppe operiert und im Quotient
alle Hyperfl"achen koaffin sind, so muss das nicht schon in $X$ selbst gelten.
Das Beispiel l"asst sich ${\bf Z}/2{\bf Z}$-graduiert verstehen, und
die Graduierung entspricht einer Operation von
${\rm Spek}\, {\bf C}[{\bf Z}/2]$ auf ${\rm Spek}\, A\otimes_{\bf R}{\bf C}$.
\par\smallskip\noindent
Bei bestimmten Gruppenoperationen aber gelten solche Aussagen, so kann man
die Aussage mit den Vektorb"undeln 3.3.5 so verstehen, dass
unter einer lokal-trivialen Operation der additiven Gruppe $G_a^r$
die affine Klassengruppe des Quotienten mit der affinen Klassengruppe
des Totalraums "ubereinstimmt.
Eine Beziehung gibt es auch bei Operationen
der multiplikativen Gruppe $G_m$, die den $\bf Z$-Graduierungen entsprechen,
siehe hierzu weiter unten im Abschnitt 3.8.
\par
\bigskip
\noindent
\subsection{Folgerungen aus trivialer affiner Klassengruppe}
\par
\bigskip
\noindent
In diesem Abschnitt wollen wir einige Folgerungen aus der Eigenschaft
eines Schemas, dass die affine Klassengruppe verschwindet, zusammenstellen,
insbesondere die Auswirkungen auf rationale Abbildungen und Morphismen
mit einem solchen
Schema als Bildraum besprechen. Ein wichtiger Satz ist dabei der
Reinheitssatz 
f"ur birationale Abbildungen, der unter der st"arkeren Voraussetzung, dass
der Bildraum regul"ar ist, auf van der Waerden zur"uckgeht.
Unter einer {\it rationalen Abbildung} verstehen wir einen auf einer offenen,
nicht leeren Teilmenge eines integren Schemas definierten Schemamorphismus
in ein Bildschema. Ist dieses separiert, so gibt es immer einen gr"o"sten
offenen Definitionsbereich der rationalen Abbildung.
Trivial, aber grundlegend, ist das folgende Lemma.
\begin{lem}
Sei $X$ ein noethersches, separiertes, integres, normales Schema,
$Y={\rm Spek}\, B$ ein affines
Schema und $f:X \cdots \longrightarrow Y$ eine rationale
Abbildung, die in der Kodimension eins definiert sei.
Dann ist $f$ schon global definiert.
\end{lem}
{\it Beweis}.
Es sei $U$ der maximale Definitionsbereich der Abbildung.
Diese Abbildung entspricht einem Ringhomomorphismus
$B \longrightarrow \Gamma(U,{\cal O}_X)$, und da $X$ normal ist und
$U$ die Kodimension eins umfasst, ist
$\Gamma(U,{\cal O}_X) = \Gamma(X,{\cal O}_X)$ und daher ist die
Abbildung global definiert. \hfill $\Box$
\begin{satz}
Seien $X$ und $Y$ noethersche, separierte, integre, normale Schemata und
$f:X \supseteq U \longrightarrow Y$ eine rationale Abbildung mit dem
maximalen Definitionsbereich $U$, die in der
Kodimensions eins definiert sei.
Es sei $Z \subseteq Y$ eine Hyperfl"ache mit affinem Komplement und
$x \not\in U$. Dann ist $x \in \overline{f^{-1}(Z)}$.
Ist $f$ nicht global definiert, so ist $f^{-1}(Z) \neq \emptyset $.
\end{satz}
{\it Beweis}.
Betrachte die Einschr"ankung
$$ X- \overline{f^{-1}(Z)} \supseteq \, U \cap (X-\overline{f^{-1}(Z)})
\longrightarrow Y-Z={\rm Spek}\, B \,  .$$
$U \cap (X-\overline{f^{-1}(Z)})$ umfasst die Punkte der H"ohe eins
von $ X- \overline{f^{-1}(Z)} $ und daher ist nach dem Lemma diese
Abbildung auf $X- \overline{f^{-1}(Z)}$ definiert. Also muss
$x \in \overline{f^{-1}(Z)}$ sein. \hfill $\Box$
\par\bigskip
\noindent
{\bf Beispiel}
Ein einfaches Beispiel, das zeigt,
dass man auf die Voraussetzung $Y-Z$ affin
nicht verzichten kann, ist die offene Einbettung einer punktierten glatten
Fl"a\-che in die an diesem Punkt aufgeblasene Fl"a\-che.
Das Urbild der exzeptionellen Faser und damit auch der Abschluss
davon ist leer.
\par\bigskip
\noindent
{\bf Bemerkung}
Unter den Voraussetzungen des Satzes ist
$X-\overline{f^{-1}(Z)}=U-f^{-1}(Z)$, d.h. man kann die Affinit"at
von dieser offenen Menge allein auf dem Definitionsbereich entscheiden.
\par\smallskip\noindent
Ist $Y$ projektiv "uber einem K"orper,
so ist die Voraussetzung, dass eine rationale Abbildung
eines normalen Schemas $ X \cdots \longrightarrow Y$
in der Kodimension eins definiert ist,
automatisch erf"ullt. Im n"achsten Abschnitt wenden wir den Satz
auf $Y=C \times C$ mit einer projektiven Kurve $C$ und darin liegenden
Kurven mit affinem Komplement an.
\begin{kor}
Seien $X$ und $Y$ noethersche, separierte, integre, normale
Schemata mit ${\rm AKG}\, Y=0$.
Ist dann $f:X \supseteq U \longrightarrow Y$ eine rationale Abbildung,
die "uberall in der Kodimensions eins definiert ist, so liegt
der Undefiniertheitsort im Abschluss des Urbildes einer
jeden Hyperfl"ache von $Y$.
Ist $f$ nicht global definiert, so wird jede Hyperfl"ache erreicht.
\end{kor}
{\it Beweis}.
Die Aussage folgt direkt aus 3.4.2. \hfill $\Box $
\begin{kor}
Sind $X$ und $Y$ noethersche, separierte, integre, normale Variet"aten
"uber einem K"orper $K$ und $Y$ projektiv
mit trivialer affiner Klassengruppe,
so liegt zu einer rationalen Abbildung
$f:X \cdots \longrightarrow Y$ der Undefiniertheitsort im Abschluss
des Urbildes einer jeden Hyperfl"ache. \hfill $\Box$
\end{kor}
\begin{kor}
Ist $C$ eine projektive glatte Kurve und $X \cdots \longrightarrow C$ eine
rationale Abbildung, so liegt der Undefiniertheitsort im Abschluss aller
Fasern.
\end{kor}
{\it Beweis}. Auf einer projektiven Kurve ist das Komplement eines Punktes
affin. 
(F"ur $C={\bf P}^1$ wurde diess Aussage schon in 3.1.2 gezeigt.) \hfill $\Box$
\par\bigskip
\noindent
Eine genauere Aussage zu den Urbildern von Hyperfl"achen folgt
bereits aus 1.12.1.
\begin{kor}
Sei $f:X \longrightarrow Y$ eine Abbildung zwischen noe\-therschen Schemata
mit $Y$ separiert, integer, normal und ${\rm AKG}\, Y=0$.
Dann ist das Urbild einer jeden irreduziblen Hyperfl"ache $D$ entweder leer
oder besitzt die reine Kodimension eins.
\end{kor}
{\it Beweis}.
$U=Y-D$ ist affin, daher ist $U \hookrightarrow Y$ ein affiner Morphismus,
was sich auf $f^{-1}(U) \hookrightarrow X$ "ubertr"agt, so dass 1.12.1 die
Behauptung liefert.
\par\smallskip\noindent
Ein anderer Beweis f"ur $X$ normal geht so:
Man betrachtet die rationale Abbildung $X \cdots \longrightarrow Y-D$.
Ist diese in der Kodimension eins definiert, so ist sie nach dem
Lemma 3.4.1 schon global
definiert, und es ist $f^{-1}(D)= \emptyset$.
Andernfalls gibt es einen Punkt $x$ der H"ohe eins
mit $f(x) \in D$. \hfill $\Box$
\par\bigskip
\noindent
{\bf Beispiel}
Man kann in der Situation des Korollars nicht erwarten, dass der generische
Punkt einer Hyperfl"ache getroffen wird, sobald ihr Urbild nicht leer ist.
Betrachte etwa die
birationale Abbildung
$$ f:{\bf A}^2_K \longrightarrow {\bf A}^2_K \mbox{ mit }
(x,y) \longmapsto (x,xy) \, .$$
Unter dieser Abbildung werden alle Kurven getroffen, das 
Urbild von $V(x)$ ist $V(x)$, das Bild davon ist aber der
Nullpunkt.
\par\bigskip
\noindent
Der folgende Satz ist der Reinheitssatz von van der Waerden.
\begin{satz}
Seien $X$ und $Y$ noethersche, separierte, integre, Schemata und
$f: X \longrightarrow Y$ ein birationaler Morphismus von endlichem
Typ. $Y$ sei normal und in jedem Punkt $y \in Y$ sei
${\rm AKG}\, {\cal O}_y=0$.
Dann besitzt das Komplement des
offenen Ortes $V \subseteq X$,
wo ein lokaler Isomorphismus vorliegt, die reine Kodimension eins.
Ist $f$ in jedem Punkt der H"ohe eins von $X$ ein Isomorphismus, so ist
$V=X$ und $f$ eine offene Einbettung.
\end{satz}
{\it Beweis}.
Da die Abbildung von endlichem Typ ist, ist $V$ offen und nicht leer.
Wir k"onnen $Y$ als affin annehmen mit trivialer affiner Klassengruppe.
Durch Herausnahme von Hyperfl"achen von $X$ k"onnen wir dann annehmen,
dass $V$ alle Punkte der H"ohe eins umfasst, dass also in der Kodimension
eins eine lokale Isomorphie vorliegt, und haben dann zu zeigen,
dass eine offene Einbettung vorliegt.
Die Aussage ist lokal in $X$, deshalb k"onnen wir auch $X$ als affin
annehmen.
\par\smallskip\noindent
Sei weiter $D \subseteq Y$ eine irreduzible Hyperfl"ache. Aufgrund des
Korollars 3.4.6 ist dann $f^{-1}(D)$ leer oder aber es gibt Punkte der
H"ohe eins im Urbild.
Ein solcher Punkt wird aber nach Voraussetzung isomorph auf seinen
Bildpunkt abgebildet, so dass in diesem Fall der generische Punkt der
Hyperfl"ache erreicht wird.
Wir k"onnen daher die endlich vielen Hyperfl"achen, deren generische Punkte
nicht erreicht werden, herausnehmen (das Restschema bleibt affin),
und annehmen, dass alle Punkte der H"ohe eins getroffen werden.
Damit entsprechen sich alle Lokalisierungen der H"ohe eins von $X$ und $Y$,
und insbesondere ist $X$ regul"ar in der Kodimension eins.
Damit stimmt die Normalisierung von $X$ mit $Y$ "uberein und damit
ist bereits $X=Y$. 
\par\smallskip\noindent
Ist $f$ in jedem Punkt der H"ohe eins von $X$ ein Isomorphismus, so liegt
global auf $X$ ein lokaler Isomorphismus vor. Da $X$ separiert ist, muss
es sich um eine offene Einbettung handeln. \hfill $\Box$
\par\bigskip
\noindent
{\bf Bemerkung}
Der vorstehende Satz geht f"ur den Fall, dass $Y$ regul"ar ist, auf van der
Waerden zur"uck und wird in \cite{EGAIV}, 21.12.12, unter etwas schw"acheren
Voraussetzungen als in 3.4.7 gezeigt.\footnote{In \cite{mumcomplex}, 3.20,
findet sich die
Aussage f"ur regul"aren Bildraum unter der Bezeichnung
Zariski´s Main Theorem (smooth case).}
Die dort lokal an das Schema $Y$ gestellte Voraussetzung lautet, f"ur einen
integren lokalen Ring $A$ formuliert, vergl. \cite{EGAIV}, 21.12.8, mit
$Y={\rm Spek}\, A$ und $Y^\times=Y -\{ {\bf m} \}$,
\par\medskip
\noindent
(W) Ist $U \subseteq Y^\times\, , U \neq Y^\times$,
und ist die Einbettung $U \hookrightarrow Y^\times$ ein affiner
Morphismus, so ist $U$ selbst schon affin. 
\par\medskip
\noindent
Diese Bedingung besagt, dass man die Affinit"at ausserhalb des abgeschlossenen
Punktes testen kann, d.h. wenn der affine Ort zu $U=D({\bf a})$
die punktierte Menge $D({\bf m})$ umfasst,
er bereits schon gleich ganz $Y$ ist und daher $U$ affin ist.
Die Bedingung ist nat"urlich erf"ullt, wenn jede Hyperfl"ache
affines Komplement hat. Ist $A$ ein lokaler, normaler, exzellenter
Integrit"atsbereich der Dimension drei, so folgt umgekehrt aus
der Eigenschaft (W), dass die affine Klassengruppe verschwindet.
Ist n"amlich $H$ eine nicht-leere Hyperfl"ache, $U=Y-H$ und $f \in A$ eine
Nichteinheit, so ist $A_f$ zweidimensional, normal und exzellent.
Somit ist $U \cap D(f) \subseteq {\rm Spek}\, A_f$ affin und das Pr"azedens
der Bedingung $(W)$ ist erf"ullt, also muss $U$ affin sein.
\par\smallskip\noindent
Die Eigenschaft (W) reicht auch aus, um folgende Aussage zu erhalten,
vergl. \cite{EGAIV}, 21.12.10.
\begin{kor}
Sei $X$ ein noethersches Schema und $Y$ ein noe\-thersches, separiertes,
integres, normales
Schema mit lokal trivialer affiner Klassengruppe.
$P$ sei ein abgeschlossener Punkt von $Y$, $Y^\times=Y-P$ und
$f:X \longrightarrow Y$ so, dass $f^{-1}(Y^\times) \longrightarrow Y^\times $
eine offene Einbettung ist mit $Z=f^{-1}(P)$ irreduzibel.
Dann hat $Z$ eine Kodimension $\leq 1$ oder ist ein abgeschlossener Punkt.
\end{kor}
{\it Beweis}.
Wir k"onnen annehmen, dass $Y={\rm Spek}\, A$ affin und lokal ist mit trivialer
affiner Klassengruppe, und ebenso, dass $X={\rm Spek}\, B$ affin ist.
$Z$ habe eine Kodimension $\geq 2$. Nach Voraussetzung ist
$X-Z \hookrightarrow Y^\times$ eine offene Einbettung.
Zu jeder irreduziblen Hyperfl"ache $H'$ von $X$ geh"ort der generische Punkt
zu $X-Z$ und damit zu $Y^\times $.
Sei $H$ die entsprechende Hyperfl"ache in $Y$, auf dieser liegt der Punkt $P$.
Dann ist $f^{-1}(H)=H' \cup f^{-1}(P)= H' \cup Z$.
Da $H$ affines Komplement besitzt, ist $Z \subseteq H'$.
$Z$ liegt also auf jeder Hyperfl"ache von $X$, d.h. jede Nichteinheit $f \in B$
liegt im Primideal zu $Z$ und damit ist $B$ ein lokaler Ring mit
maximalem Ideal $Z$. \hfill $\Box$
\par\bigskip
\noindent
{\bf Beispiel}
Typische Beispiele f"ur birationale Morphismen, wo der Reinheitssatz
verletzt ist, sind die sogenannten {\it kleinen Kontraktionen}
$f: Y' \longrightarrow Y$.
Dabei wird eine abgeschlossene Teilmenge $C \subseteq Y'$
der Kodimension $\geq 2$ und der
Dimension $\geq 1$ auf einen abgeschlossenen Punkt kontrahiert, und
au"serhalb liegt ein Isomorphismus vor. Man denke etwa
an die Kontraktion einer Kurve in einer dreidimensionalen Variet"at.
Dabei entsprechen sich dann direkt die Punkte der H"ohe eins,
und damit entsprechen sich die Hyperfl"achen, Divisoren, Divisorenklassengruppe.
Eine solche Kontraktion liefert eine Vielzahl von Hyperfl"achen mit
nicht affinem Komplement im Bild: eine Hyperfl"ache $H$ im Bildraum
durch den Bildpunkt $P=f(C)$ der Kurve kann nur dann affines
Komplement besitzen, wenn die Kurve $C$ ganz auf der entsprechenden
Hyperfl"ache $H'$ liegt. Zu allen Hyperfl"achen $H'$ in $Y'$, die die Kurve
nulldimensional schneiden, ist das Komplement von $H$ nicht affin.
Es gilt ja bei $P \in H$ die Gleichheit $Y'- (H' \cup C) =Y-H$, und bei
$C \not \subseteq H'$ hat das Komplement in $Y'$ eine Komponente
der Kodimension $\geq 2$.
\par\smallskip\noindent
Ist $Y={\rm Spek}\, A$ ein affines, noethersches, integres, normales
Schema mit ${\rm AKG}\, A \neq 0$, und ist $H \subseteq Y$ eine
Hyperfl"ache mit nicht affinem Komplement $U=Y-H$, so ist die Abbildung
${\rm Spek}\, \Gamma(U, {\cal O}_Y)=U^{\rm aff} \longrightarrow Y$ ein Kandidat,
wo der Satz verletzt ist
und eine kleine Kontraktion vorliegen kann.
Diese Abbildung ist birational und auf $U \subseteq U^{\rm aff}$ eine
offene Einbettung, wobei $U$ die Kodimension eins umfasst.
Ist zus"atzlich $U^{\rm aff}$ von endlichem Typ "uber $Y$,
so verletzt $U^{\rm aff} \longrightarrow Y$ die
Aussage des Reinheitssatzes, da das Bild in $U \uplus {\rm Naf}\, (U,Y)$ liegt
und nicht offen ist.
\par\smallskip\noindent
Ist speziell $Y={\rm Spek}\, A$ dreidimensional und exzellent, $U=D({\bf a})$
nicht affin und $U^{\rm aff}$ von endlichem Typ "uber $A$, so liegt sogar
eine kleine Kontraktion vor.
Das Erweiterungsideal ${\bf a}\Gamma(U,{\cal O}_Y)$ hat dann
eine H"ohe von zwei (die H"ohe kann wegen 4.1.4 nicht drei sein) und
beschreibt damit eine Kurve.
Das Bild davon muss im nicht affinen Ort ${\rm Naf}\, (U,Y)=\{ P\}$ liegen.
Diese Kurve wird also auf einen Punkt kontrahiert und au"serhalb liegt eine
offene Einbettung vor.
\par\smallskip\noindent
Als Beispiel betrachten wir in ${\rm Spek}\, K[u,v,x,y]/(ux-vy)=X$ die nicht
affine Teilmenge $U=D(x,y)$ mit der Affinisierung
$$U^{\rm aff}={\rm Spek}\, K[x,y,q=v/x=u/y]={\bf A}_K^3$$
und der Restriktionsabbildung
$${\bf A}^3_K \ni (x,y,q) \longmapsto (qy,qx,x,y)
\in V(ux-vy) \subseteq {\bf A}_K^4 \, . $$
Dabei liegt auf $D(x,y)$ eine Isomorphie vor und die Gerade $x=y=0$
wird auf den Nullpunkt kontrahiert.
\par\bigskip
\noindent
{\bf Bemerkung}
In \cite{EGAIV}, Remarques 21.12.14 (iv) wird gefragt,
ob man die Eigenschaft (W) eines lokalen normalen noetherschen
Ringes $A$ dadurch erhalten kann,
dass f"ur jede Abbildung $X \longrightarrow Y$, f"ur die
$f^{-1}(Y') \longrightarrow Y'$ eine offene Einbettung ist,
das Urbild $f^{-1}(y)$ des abgeschlossenen Punktes $y$
die Kodimension eins hat oder aber ein Punkt ist.
Die Aussagen in diesem Abschnitt gelten aber auch dann noch, wenn man das
Verschwinden der affinen Klassengruppe dadurch ersetzt, dass jedes Ideal
der reinen H"ohe eins auch die Superh"ohe eins besitzt, f"ur diesen
Begriff siehe Kapitel 4.
Zu suchen ist demnach ein solcher (dreidimensionaler) Ring, in dem 
trotzdem nicht alle Hyperfl"achen affines Komplement besitzen.
Am ehesten d"urfte man einen solchen Ring als homogenen Koordinatenring
einer projektiven Fl"a\-che erhalten k"onnen,
auf der sich je zwei (verschiedene) Kurven positiv schneiden, aber
trotzdem nicht jede Kurve ger"aumig ist.
\par
\bigskip
\noindent
\subsection{Die affine Klassengruppe von Produkten}
\par
\bigskip
\noindent
Wir interessieren uns in Verallgemeinerung der Situation
$X \times {\bf A} \longrightarrow X$ f"ur das Verhalten der affinen
Klassengruppe von einer affinen normalen Variet"at $X={\rm Spek}\, A$
"uber einem algebraisch
abgeschlossenen K"orper $K$ zu der affinen Klassengruppe eines Produktes mit
einer weiteren affinen normalen Variet"at, im einfachsten Fall mit einer
affinen glatten Kurve $C$.
Ohne die Voraussetzung, dass der Grundk"orper algebraisch abgeschlossen ist,
gibt es, wie das Beispiel am Schluss von 3.3 zeigt,
sogar nulldimensionale Basiswechsel, unter
denen sich die affine Klassengruppe "andert, weshalb nur unter dieser
Voraussetzung hier etwas zu erwarten ist.
Sind $X$ und $Y$ beide regul"ar, so ist auch das Produkt regul"ar, und
die affine Klassengruppe des Produktes verschwindet.
F"ur die Projektion $X \times Y \longrightarrow X$ gilt,
dass die Urbilder irreduzibler Teilmengen wieder irreduzibel sind.
Ist ${\rm AKG}\, X=0$, so besitzt jede Hyperfl"ache, die in der generischen
Faser $Y_{K(X)}$ durch eine Funktion beschrieben wird, affines Komplement,
diese Faser ist aber in aller Regel nicht faktoriell.
\par\bigskip
\noindent
Wir geben zun"achst einige weitere Korollare zu Satz 3.3.3 an,
bei denen die Hauptschwierigkeit darin besteht, zu zeigen,
dass die Abbildung der Divisorenklassengruppen surjektiv ist.
\begin{kor}
Sei $A$ eine normale $k-$Algebra von endlichem Typ "uber einem algebraisch
abgeschlossenen K"orper $k$ der Charakteristik null,
$k \subseteq K$ sei eine K"orpererweiterung.
Besitzt dann $A$ eine endlich erzeugte Divisorenklassengruppe,
so stimmen die affinen Klassengruppen von $A$ und $A \otimes_kK$ "uberein.
\end{kor}
{\it Beweis}.
Unter der Voraussetzung stimmen die Divisorenklassengruppen "uber\-ein, siehe
\cite{binsto2}, 15.7. \hfill $\Box$
\par\bigskip
\noindent
\begin{kor}
Sei $A$ eine faktorielle $k-$Algebra von endlichem Typ "uber einem
algebraisch abgeschlossenen K"orper $k$ der Charakteristik null und $B$ eine
noethersche, integre, normale $k-$Algebra.
Dann stimmen die affinen Klassengruppen von $B$ und $B \otimes_kA$
"uberein.
\end{kor}
{\it Beweis}. Zur Gleichheit der Divisorenklassengruppen
siehe \cite{binsto2}, 15.10. Man kann auch mit dem Lemma 3.3.6 und dem obigen
Korollar 3.5.1 schlie"sen, da mit $K=Q(B)$ dann zun"achst die generische Faser
${\rm Spek}\, A \times K$ faktoriell ist und damit f"ur
${\rm Spek}\, A \times {\rm Spek}\, B \longrightarrow {\rm Spek}\, B$ die
Voraussetzungen des Lemmas erf"ullt sind.  \hfill $\Box$
\par\bigskip
\noindent
{\bf Beispiel}
Beispiele f"ur die letzte Situation sind etwa $A=k[T_1,...,T_n]$ oder
$A=k[X,Y,Z]/ (X^2+Y^3+Z^5)$ sowie generische homogene Hyperfl"achen
$k[X,Y,U,V]/(F)$ mit $F$ vom Grad $\ge 4$.
In diesen Situationen gilt also f"ur eine
normale $k-$Algebra $B$ die Gleichung
${\rm AKG}\,(B \otimes_kA)={\rm AKG}\, B$.
Beispielsweise ist f"ur jeden normalen exzellenten Ring $B$
der Dimension $\leq 2$ "uber $k$ in ${\rm Spek}\, B[X,Y,Z]/(X^2+Y^3+Z^5)$ jede
Hyperfl"ache koaffin. 
\par\bigskip
\noindent
Sei nun $X$ eine affine normale Variet"at und $C$ eine affine glatte Kurve.
Die einfachsten Hyperfl"achen in $X \times C$ sind die Graphen zu Abbildungen
$X \longrightarrow C$ und allgemeiner die Abschl"usse von Graphen
von rationalen Funktionen in die Kurve. Es ist hier sinnvoll, in der
Fragestellung auch projektive glatte Kurven zuzulassen.
\begin{satz}
Sei $X={\rm Spek}\, A$ eine affine normale Variet"at "uber einem algebraisch
abgeschlossenen K"orper $k$ und $C$ eine affine glatte Kurve oder die
projektive Gerade.
Es sei $f:X \cdots \longrightarrow C$ eine rationale Abbildung derart,
dass $\overline{f}:X \supseteq U \longrightarrow \overline{C}$ ein affiner
Morphismus ist, wobei $U$ der Definitionsort von $\overline{f}$
sei und $\overline{C}$ der projektive glatte Abschluss der Kurve.
Dann hat der Abschluss des Graphen der rationalen Abbildung
affines Komplement.
\end{satz}
{\it Beweis}.
Es sei $\overline{C}$ der Abschluss von $C$ und es seien $P_1,...,P_n$ die
hinzugenommenen Punkte.
Es sei $\overline{\Gamma_f} \subseteq X \times C$ der Abschluss des Graphen von
$f:X \cdots \longrightarrow C$. Es ist dann
$\overline{\Gamma_f}= \overline{\Gamma_{\overline{f}}} \, \cap \, X \times C$,
wie eine einfache topologische "Uberlegung zeigt.
Damit ist die in Frage stehende Menge gleich
\begin{eqnarray*}
X \times C - \overline{\Gamma_f} &=&
X \times C - \overline{\Gamma_{\overline{f}}}\, \cap \, X \times C \\
& = & X \times \overline{C}
- \overline{\Gamma_{\overline{f}}} - X \times \{ P_1,...,P_n \} \\
& =& X \times \overline{C} -
\overline{\Gamma_{\overline{f}} \cup U \times \{ P_1,...,P_n \} } \,  .
\end{eqnarray*}
Die auf $U \times \overline{C}$ (als maximalem Definitionsort)
definierte Abbildung
$$ \overline{f} \times id : U \times \overline{C} \longrightarrow
\overline{C} \times \overline{C} \, $$
ist affin, da nach Voraussetzung
$\overline{f}:U \longrightarrow \overline{C}$ affin ist.
Das  Urbild von $\overline{f} \times id$ von
$\triangle \, \cup \, \overline{C} \times \{P_1,...,P_n \}$
ist dabei gleich $\Gamma_{\overline{f}} \cup U \times \{ P_1,...,P_n \}$.
Ist $C={\bf P}$, so hat in ${\bf P} \times {\bf P}$ die Diagonale affines
Komplement. Im andern Fall ist $n \geq 1$ und dann hat in
$\overline{C}  \times \overline{C}$ die Kurve
$\triangle \cup \overline{C} \times \{ P_1,...,P_n \}$ nach Satz 2.4.3
affines Komplement.
Damit ist dann
$U \times \overline{C} - (\Gamma_{\overline{f}} \cup \{ P_1,...,P_n \})$
affin.
Nach Satz 3.4.2 ist diese Menge aber gleich
$ X \times \overline{C} -
\overline{\Gamma_{\overline{f}} \cup U \times \{ P_1,...,P_n \} } $
und damit gleich $X \times C - \overline{\Gamma_f} $. \hfill $\Box$
\par\bigskip
\noindent
{\bf Bemerkung}
Ist $X$ eine affine Variet"at, $C$ eine beliebige glatte Kurve
und $f:X \longrightarrow C$ ein Morphismus, so ist das
Komplement des Graphen in $X \times C$ immer affin.
Dazu gen"ugt es, die Affinit"at der Abbildung
$$X \times C -\Gamma_f \hookrightarrow X \times C \longrightarrow X$$
nachzuweisen, so dass man zu kleineren affinen Teilmengen von $X$
"ubergehen kann.
Damit kann man annehmen, dass der Punkt $P \in C$ nicht getroffen wird.
Damit ist
$\Gamma_f=(f \times id)^{-1}(\triangle)
=(f \times id)^{-1} (\triangle \cup X \times \{P \}) $
und das Komplement ist affin aufgrund von 2.4.3.
\begin{kor}
Seien $X$ und $C$ wie im Satz und $X$ habe triviale affine Klassengruppe.
Dann hat der Abschluss des Graphen einer rationalen Abbildung
$f:X \cdots \longrightarrow C$ affines Komplement.
\end{kor}
{\it Beweis}. Es ist nur zu zeigen, dass die Ausdehnung $\overline{f}$
ein affiner Morphismus ist. Ist $f$ konstant, so ist das klar, im andern
Fall ist das Urbild
eines Punktes eine Hyperfl"ache, die nach Voraussetzung affines
Komplement hat. \hfill $\Box$
\begin{kor}
Sei $X$ ein affine normale Variet"at und $q \in K(X)$ eine rationale
Funktion mit ${\rm Def}\, (q)$ und ${\rm Def}\, (q^{-1})$ affin.
Dann ist das Komplement des Abschlusses des Graphen der meromorphen
Funktion $q: X \cdots \longrightarrow {\bf P}_K^1$ affin. \hfill $\Box$
\end{kor} 
{\bf Beispiel}
Das Komplement des Abschlusses eines Graphen einer rationalen Abbildung
nach einer projektiven glatten Kurve $C$ von h"oherem Geschlecht muss
nicht affines Komplement haben, da dann die Diagonale nicht ger"au\-mig
ist. Sei $C={\rm Proj}\, S$ eine solche Kurve
und betrachte die Kegelabbildung $f: {\rm Spek}\, S \cdots \longrightarrow C$,
die auf $U=D(S_+)$ definiert ist. $S$ sei als normal vorausgesetzt und
$C$ habe das Geschlecht eins. $D$ sei eine
Gerade auf $C \times C$, die zur Diagonalen disjunkt sei, etwa vom
Typ $(x,x+p)$ mit $p \neq 0$.
Das Urbild der Diagonalen in $U \times C$ hat nicht affines
Komplement. In diesem Komplement liegt n"amlich als
abgeschlossenes Unterschema das Urbild von $D$,
und das ist gleich $\{ (y,z) : z=f(y)+p \}$ und das ist
ebenfalls der Graph einer Abbildung und damit isomorph zum nicht affinen
$U$ (Die Graphen sind zwei disjunkte "ubereinanderliegende
Wendeltreppen mit konstantem Abstand; im Abschluss haben sie die
Schraubungsachse der Wendeltreppe gemeinsam).
\par\bigskip
\noindent
{\bf Bemerkung}
Mit diesem Beispiel ist gezeigt, dass es in projektiven, integren, glatten
Kurven "uber einem affinen normalen Basisschema dominante Hyperfl"achen geben
kann, deren Komplemente nicht affin sind. Es ist offen, ob dies auch sein
kann, wenn die Basis sogar regul"ar ist. Ist die Basis einpunktig,
so ist das Komplement eines Punktes affin, und es geht um die
Frage, inwiefern man diese Aussage auf die relative Kurvensituation
"ubertragen kann. Ist die Basis eine affine regul"are
Kurve (vom endlichen Typ), so zeigt 2.4.3,
dass jede dominante Kurve affines Komplement
hat. Aber schon "uber einem diskreten Bewertungsring ist es nicht klar. 
\par\bigskip
\noindent
{\bf Beispiel}
Besitzen zwei affine normale Variet"aten $X$ und $Y$
triviale affine Klassengruppen,
so muss das nicht f"ur das Produkt $X \times Y$ gelten.
Hierzu betrachten wir die zweidimensionale Variet"at $X={\rm Spek}\, S$
des vorangehenden Beispiels.
Das dortige Schema $X \times C$ kann man
auffassen als ${\rm Proj}\, (S \otimes S)$, wobei man die erste Graduierung
vergisst und die zweite beh"alt, da generell f"ur einen beliebigen Ring $R$
gilt ${\rm Proj}\, (R \otimes S) ={\rm Spek}\, R \times {\rm Proj}\, S$.
Da es in $X \times C$ Hyperfl"achen mit nicht affinem Komplement gibt,
muss es dies auch im affinen Kegel $X \times X$ geben.
\par\bigskip
\noindent
In dem vorangehenden Beispiel l"asst sich die Divisorenklassengruppe
des Produktes nicht auf die Divisorenklassengruppe der Faktoren
zur"uckf"uhren, und in diesen Situationen, wenn man keine
Kontrolle "uber die Divisorenklasengruppe hat, ist es schwierig,
Aussagen "uber die affine Klassengruppe zu machen, bzw. die naheliegenden
Aussagen sind wie im Beispiel schlichtweg falsch.
Auch im Produkt $X \times C$ mit einer affinen glatten Kurve treten
wesentlich neue Divisoren auf, zumindest im zweidimensionalen Fall
aber kann man die affine Klassengruppe trotzdem als trivial erweisen.
\begin{satz}
Sei $X={\rm Spek}\, A$ eine zweidimensionale, affine, normale Variet"at
und $C$ eine affine glatte Kurve
"uber einem algebraisch abgeschlossenen K"orper $k$.
Dann verschwindet die affine Klassengruppe des Produktes $X \times C$,
in $X \times C$ ist also das Komplement jeder Hyperfl"ache affin.
\end{satz}
{\it Beweis}.
Wir f"uhren die Aussage durch Betrachtungen von geeigneten
endlichen normalen Erweiterungen $X' \longrightarrow X$ auf den Satz "uber die
Graphenkomplemente 3.5.3 zur"uck.
Da $X'$ wieder zweidimensional und normal ist, ist stets
${\rm AKG}\, X'=0$, so dass in
$X' \times C$ alle Graphen affines Komplement haben.
Es sei $K=K(X)=Q(A)$.
\par\smallskip\noindent
Sei eine irreduzible Hyperfl"ache $H=V({\bf p}) \subseteq X \times C$
gegeben. Dominiert diese Hyperfl"ache nicht die Basis $X$, so ist sie
das Urbild einer (irreduziblen) Kurve auf $X$, und damit ist ihr
Komplement affin. Sei also weiterhin vorausgesetzt, dass die Hyperfl"ache $X$
dominiert. Dann liegt in den generischen Punkten eine
K"orpererweiterung $K \longrightarrow \kappa({\bf p})=:L$ vor.
Diese K"orpererweiterung ist endlich, da die Gesamtabbildung vom endlichen
Typ und ${\rm Spek}\, L$ ein abgeschlossener Punkt der generischen Faser ist.
Es seien $x_1,...,x_n \in L$ ein $K-$Algebra-Erzeugendensystem von $L$
mit zugeh"origen Ganzheitsgleichungen $F_1,...,F_n \in K[T]$.
Nach dem Satz von Kronecker, siehe \cite{scheja}, 54.10, gibt es einen endlichen
Erweiterungsk"orper $K \subseteq K'$, indem alle diese Polynome $F_i$
in Linearfaktoren zerfallen.
Dann ist jeder Restek"orper von
$K' \otimes_KL=K'[X_1,...,X_n]/(F_1,...,F_n + \mbox{weitere Relationen} )$
isomorph zu $K'$, da in den Restek"orpern jeweils ein Linearfaktor der $F_i$
null werden muss und daher $X_i$ gleich einem Element aus $K'$ sein muss.
\par\smallskip\noindent
Sei nun $A'$ der normale Abschluss von $A$ in $K'$ und $X'={\rm Spek}\, A'$.
Da $A$ exzellent ist, ist $A'$ eine endliche Erweiterung von $A$.
Aufgrund des
Satzes von Chevalley 1.14.1 gen"ugt es, die Affinit"at des
Komplementes des Urbildes von $H$ nachzuweisen.
F"ur $X' \longrightarrow X$ und $X' \times C \longrightarrow X \times C$ gilt
das going down Theorem, siehe \cite{nagloc}, I.10.13, und daher ist das
Urbild von $V=V({\bf r})$ gleich dem Abschluss der Urbildpunkte von $\bf r$,
oder anders gesagt, jede Komponente des Urbildes dominiert das Bild.
Die Faser "uber $\bf p$ ist ${\rm Spek}\, K' \otimes_KL$, dessen Punkte
alle $K'-$Punkte sind. Der Abschluss davon ist die Vereinigung der
Abschl"usse der einzelnen Punkte. Es bleibt also lediglich zu zeigen,
dass in der Situation $X' \times C \longrightarrow X'$ eine Hyperfl"ache $H'$,
deren generischer Punkt gleich $K(X')$ ist, affines Komplement besitzt.
In diesem Fall ist aber $H' \longrightarrow X'$ generisch isomorph, also
eine birationale Abbildung, und es gibt eine rationale Abbildung
$X' \cdots \longrightarrow H' \hookrightarrow X' \times C \longrightarrow C$,
und dabei ist $H'$ der Graph dieser rationalen Abbildung.
Nach Satz 3.5.3 ist dann das Komplement affin. \hfill $\Box$
\par\bigskip
\noindent
{\bf Bemerkung} Das im Satz verwendete Beweisverfahren ist nicht anwendbar
auf den allgemeinen Fall mit ${\rm AKG}\, X=0$ als Voraussetzung,
da sich dies nicht auf $X'$ "ubertr"agt. Vielleicht kann man aber dennoch
mit den Basiswechseln $X' \longrightarrow X$ arbeiten, wenn man
die rationalen Abbildungen $X' \longrightarrow C$ als
affin nachweist.
\par\smallskip\noindent
Die st"arkste Vermutung, die man in diesem Zusammenhang aufstellen kann, ist
folgende Aussage: F"ur eine affine normale
Variet"at $X$ und eine affine glatte Variet"at $Y$ gilt
${\rm AKG}\, (X \times Y)={\rm AKG}\, X$.
\par
\bigskip
\noindent
\subsection{Beispiele}
\par
\bigskip
\noindent
In diesem Abschnitt berechnen wir explizit die affine Klassengruppe von
drei wichtigen Ringklassen, n"amlich von Hyperbeln "uber
faktoriellen Ringen, von Monoidringen und von Determinantenringen.
In all diesen F"allen
ergibt sich, dass die affine Klassengruppe einfach die Divisorenklassengruppe
modulo Torsion ist.
\par
\bigskip
\medskip
\newpage
{\it Hyperbeln}
\par
\bigskip
\noindent
Sei $A$ ein noetherscher faktorieller Ring,
$U_1,...,U_r$ verschiedene, nicht assozierte Primelemente von $A$.
Wir betrachten die (normale) Hyperbel
$$B=A[X,Y]/(XY-f) \mbox{ mit }
f=U_1^{d_1} \cdot ... \cdot U_r^{d_r},\, d_i > 0 \, . $$
${\bf p}_i=(U_i,X)$ und ${\bf q}_i =(U_i,Y)$ sind Primideale der H"ohe eins.
\par\smallskip\noindent
Es ist $B_X=A[X,X^{-1}]$ faktoriell, alle Divisorenklassen werden
also repr"asentiert durch
$n_1{\bf p}_1+...+n_r{\bf p}_r$, wobei man die Divisoren der Einheitengruppe
von $B_X$ rausdividieren muss, also ergibt sich als Divisorenklassengruppe
wegen $(X)=d_1{\bf p}_1+...+d_r{\bf p}_r$
einfach ${\rm DKG}\, B={\bf Z}^r/(d_1,...,d_r) $.
Es ist $(U_i)=(X,U_i)+(Y,U_i)={\bf p}_i+{\bf q}_i$ und insbesondere ist
${\bf p}_i=-{\bf q}_i$ in ${\rm DKG}\, B$.
\par\smallskip\noindent
Das verschwinden der affinen Klassengruppe l"asst sich in dieser Situation
einfach charakterisieren.
\begin{satz}
F"ur $B$ wie oben sind folgende Aussagen \"aquivalent.
\par\smallskip\noindent
{\rm (1)}
Es ist ${\rm AKG}\, B=0$, es sind also alle Komplemente von Hyperfl"achen
affin.
\par\smallskip\noindent
{\rm (2)}
$D({\bf p}_i)$ ist affin f\"ur alle $i=1,...,r$.
\par\smallskip\noindent
{\rm (3)}
Die Primelemente $U_i$ erzeugen in $A$ paarweise das Einheitsideal.
\par\smallskip\noindent
{\rm (4)}
$B$ ist lokal fastfaktoriell.
\end{satz}
{\it Beweis}.
$(1) \Longrightarrow (2)$ ist trivial. $(2) \Longrightarrow (3)$. Wenn
$U_i$ und $U_j,\, i \neq j,$ in $A$ nicht das Einheitsideal erzeugen,
so auch nicht in $B$, es gibt also ein minimales Primoberideal $\bf r$ von
$(X,Y,U_i,U_j)$. Dann gibt es aber die Primidealkette
$(0) \subset (Y,U_j) \subset (X,Y,U_j) \subset {\bf r}$, das letzte
Primideal hat also zumindest die H\"ohe drei.
Damit ist $D(X,U_i)$ nicht affin.
$(3) \Longrightarrow (4)$. Sei $T_i=U_1 \cdot... \check{U_i} ...\cdot U_r $,
wobei $\check{U_i}$ bedeutet, dass dieses Element ausgelassen wird.
Aus der Voraussetzung folgt, dass auch die
$T_1,...,T_r$ zusammen das Einheitsideal erzeugen.
Wir betrachten 
$B_{T_i} =A_{T_i}[X,Y]/XY-U_i^{d_i} \cdot e $ mit $e \in A_{T_i}$ eine
Einheit.
$B_{T_i}$ hat dann die Divisorenklassengruppe ${\bf Z}/d_i $,
ist also eine Torsionsgruppe. Da ${\rm Spek}\, B$ von $D(T_i)$ "uberdeckt
wird, ist $B$ lokal fastfaktoriell.
Von (4) nach (1) gilt immer.     \hfill $\Box$
\par\bigskip
\noindent
Sei nun $A$ lokal. Dann erzeugen in $A$ die $U_i,\, i=1,...,r,$ nicht das
Einheitsideal und daher schneiden sich je zwei Hyperfl"achen $V({\bf p}_i)$
und $V({\bf q}_j)$ f"ur $i \neq j$ in einer Menge der Kodimension $3$. Daher
kann die affine Klassengruppe bei $r \geq 2$ nicht trivial sein.
Die folgenden "Uberlegungen werden zeigen, dass im lokalen Fall die affine
Klassengruppe gleich der Divisorenklassengruppe modulo Torsion ist.
\begin{satz}
Sei $A$ ein lokaler, noetherscher, faktorieller Ring, $U_1,...,U_r$
verschiedene Primelemente von $A$, $f=U_1^{d_1}\cdot ... \cdot U_r^{d_r}$ und
$B=A[X,Y]/(XY-f)$.
Dann ist ein Divisor $D$ genau dann koaffin, wenn er ein Hauptdivisor ist
oder "aquivalent zu einem Divisor $n_1{\bf p}_1+...+n_r{\bf p}_r$ mit
$0 < n_i < d_i$ f"ur $i=1,...,r$.
\end{satz}
{\it Beweis}.
Jeder Divisor ist "aquivalent zu $n_1{\bf p}_1+...+n_r{\bf p}_r$. Durch
Addition von Vielfachen von $d$ betrachten wir den gr"o"stm"oglichen Vertreter
$\leq d$.
\par\smallskip\noindent
Sei zun"achst ein Eintrag nicht positiv.
Es gibt auch positive Eintr"age, sagen wir die ersten $k$ Stellen seien
positiv, die n"achsten von $k+1$ bis $m$ negativ und die restlichen null.
Wir ersetzen die negativ auftretenden ${\bf p}_i$ durch ${\bf q}_i$ und haben
dann einen "aquivalenten effektiven Divisor, wo
der Tr"ager aus ${\bf p}_1,...,{\bf p}_k,{\bf q}_{k+1},...,{\bf q}_m$ besteht
mit $1 \leq k \leq m$.
Ist $m=k$, ist also kein Eintrag negativ, so besitzt
${\bf p}_1 \cap ... \cap {\bf p}_k + {\bf q}_j, \, j > k,$ die H"ohe $\geq 3$
und der Divisor kann nicht koaffin sein.
Sei also nun zus"atzlich ein Eintrag
negativ, und $ m > k$.
Betrachte das Ideal
${\bf r}=(Y+U_1^{d_1}\cdot ... \cdot U_k^{d_k},
X+U_{k+1}^{d_{k+1}} \cdot ... \cdot U_r^{d_r})$. Dieses ist prim in $A[X,Y]$
und wegen
$$XY-U_1^{d_1} \cdot ... \cdot U_r^{d_r}=
X(Y+U_1^{d_1}\cdot ... \cdot U_k^{d_k})
-U_1^{d_1} \cdot ... \cdot U_k^{d_k}
(X+U_{k+1}^{d_{k+1}} \cdot ... \cdot U_r^{d_r})$$
hat dieses Primideal in $B$ die Kodimension eins. Sei
$${\bf g} \supseteq {\bf r}
+{\bf p}_1 \cap ...\cap {\bf p}_k \cap {\bf q}_{k+1} \cap ...\cap {\bf q}_m $$
ein Primoberideal. Dann ist ${\bf p}_i=(U_i,X) \subseteq {\bf g}$ f"ur ein
$i=1,...,k$ oder ${\bf q}_j=(U_j,Y) \subseteq {\bf g}$ f"ur ein
$j=k+1,...,m$. Im ersten Fall ist
\begin{eqnarray*}
{\bf g} \supseteq (U_i,X)+{\bf r}& =& 
(U_i,X,Y+U_1^{d_1}\cdot ... \cdot U_k^{d_k},
X+U_{k+1}^{d_{k+1}} \cdot ... \cdot U_r^{d_r} )\\
&=& (U_i,X,Y,U_{k+1}^{d_{k+1}}\cdot ... \cdot U_r^{d_r})
\end{eqnarray*}
und daher ist ${\bf g} \supseteq (U_i,X,Y,U_j)$
mit einem $k+1 \leq j \leq r$, insbesondere $i \neq j$, also hat
$\bf g$ eine Kodimension von zumindest drei.
Dies ergibt sich auch im anderen Fall, so dass solche Divisoren nicht koaffin
sind.
\par\smallskip\noindent
Seien jetzt alle $n_i$ positiv.
Ist $(n_1,...,n_r)=(d_1,...,d_r)$, so liegt ein koaffiner Hauptdivisor vor.
Sei also $(n_1,...,n_r)$ zumindest an einer Stelle kleiner als das Gradtupel.
Ist zugleich $n_i=d_i$ f"ur ein $i$, so ist $n-d$ nirgendwo positiv,
an einer Stelle null und an einer Stelle negativ,
durch "Ubergang zu ${\bf q}$ erh"alt man dann einen "aquivalenten
Divisor, und man zeigt wie oben, dass dessen Komplement nicht affin sein kann.
\par\smallskip\noindent
Sei also von nun an vorausgesetzt, dass $0 < n < d$ ist, wobei
die Ungleichungen in jeder Komponente echt seien.
Wir wollen zeigen, dass solche Divisoren koaffin sind.
Es ist
$D=n_1{\bf p}_1+...+n_r{\bf p}_r \cong -n_1{\bf q}_1-...-n_r{\bf q}_r=E$.
Wir betrachten dazu das reflexive Ideal
$$ {\bf b}= \{ q \in Q(B) :\, (q) + E\geq 0 \}
={\bf q}_1^{(n_1)} \cap ... \cap {\bf q}_r^{(n_r)}\, . $$
Zur Berechnung des Ideals betrachten wir die (graduierte) Abbildung
$$ A[X,U^dX^{-1}]= B \longrightarrow B_X=A[X,X^{-1}] \, ,$$
wobei im hinteren Ring
das in Frage stehende Ideal einfach das von $U^n$ erzeugte ist, also
${\bf b}= B \cap (U^n)B_X$.
Den Schnitt k"onnen wir stufenweise berechnen. Nat"urlich geh"ort
$U^n \in {\bf b}$ und das erzeugt das Ideal in allen positiven Stufen.
F"ur die negative Stufe ${\bf b}_{-t} \, (t \geq 1)$ ist
${\bf b}_{-t}=A(U^n)X^{-t} \cap B = A(U^dX^{-1})^t$.
Wegen $dt \geq d \geq n$ ist n"amlich
$aU^{dt}X^{-t} = a' U^nX^{-t} \in {\bf b}_{-t}$.
Die negativen Stufen werden also von $Y=U^dX^{-1}$ erzeugt,
und insgesamt ist
$${\bf b}={\bf q}_1^{(n_1)} \cap ... \cap {\bf q}_r^{(n_r)} = (U^n,Y) \, . $$
Zu diesen Erzeugenden geh"oren die effektiven Divisoren
$(U^n)+E= n_1{\bf p}_1+...+n_r{\bf p}_r$ und
$(Y)+E=d_1{\bf q}_1+...+d_r{\bf q}_r-(n_1{\bf q}_1+...+n_r{\bf q}_r)
=(d_1 - n_1){\bf q}_1 +...+ (d_r-n_r) {\bf q}_r$.
Da $n_i$ und $d_i-n_i >0$ sind f"ur $i=1,...,r$ sind die beiden
Komplemente affin.
Nach 3.2.2 ist dann der Divisor $D$ koaffin. \hfill $\Box$
\begin{satz}
Sei $B$ ein Hyperbelring "uber einem lokalen faktoriellen Ring $A$ und $D$
ein Divisor. Dann sind folgende Aussagen "aquivalent.
\par\smallskip\noindent
{\rm (1)} Ein Vielfaches von $D$ ist ein Hauptdivisor.
\par\smallskip\noindent
{\rm (2)} Jedes Vielfache von $D$ ist koaffin.
\par\smallskip\noindent
{\rm (3)} $D$ ist affin-trivial.
\par\bigskip\noindent
Die affine Klassengruppe von $B$ ist die Divisorenklassengruppe modulo
Torsion.
\end{satz}
{\it Beweis}.
Es ist zu zeigen, dass f"ur einen Divisor $D$ mit der Eigenschaft (K)
bereits ein Vielfaches ein Hauptdivisor ist.
Da $D$ insbesondere koaffin ist, ist $D$ "aquivalent zu
$0 < (n_1,...,n_r) < (d_1,...,d_r)$, wobei die Relation in jeder
Komponente erf"ullt ist.
Wir haben zu zeigen, dass die Quotienten
$d_i/n_i$ alle gleich sind. Sei (nach Umordnung) $d_1/n_1 \leq d_i/n_i$
f"ur alle $1 \leq i \leq r$.
Es ist dann $d_1(n_1,...,n_r) \leq n_1(d_1,...,d_r)$, wobei an der ersten
Stelle Gleichheit gilt, und wir haben
$0 \leq d_1(n_1,...,n_r)-n_1(d_1,...,d_r)$.
Dieser Divisor ist effektiv und zu $d_1D$ "aquivalent, da der erste
Eintrag null ist, kann er nur dann koaffin sein, wenn die anderen
Stellen auch null sind. Dann ist es aber ein Hauptdivisor.
\par\smallskip\noindent
Die Aussage "uber die affine Klassengruppe folgt.  \hfill $\Box$
\par\bigskip
\noindent
{\bf Bemerkung}
An dieser Beispielklasse sieht man, dass ein affin-trivialer
Divisor, also einer, der die Eigenschaft (K)
erh"alt, nicht die Koaffinit"at
selbst erhalten muss.
W"urde man affin-trivial und affine Klassengruppe "uber diese
Eigenschaft festlegen, k"onnte in ihr Torsion auftreten.
\par\smallskip\noindent
Betrachte hierzu $K[U_1,U_2][X,Y]/(XY-U_1^3U_2^3)$.
Dann ist ein Vielfaches des Divisors $(1,1)$ ein Hauptdivisor, also ist
$(1,1)$ affin-trivial.
$(1,2)$ ist koaffin, $(1,1)+(1,2)=(2,3)=(-1,0)$ aber nicht, so dass
$(1,1)$ nicht die Koaffinit"at erh"alt.
\par
\medskip
\bigskip
{\it Monoidringe}
\par
\bigskip
\noindent
Wir wollen die affine Klassengruppe zu einem Monoidring $K[M]$ berechnen,
wobei $K$ ein faktorieller Bereich und $M$ ein normales, endlich erzeugtes,
torsionsfreies Monoid ist.
Es sei $\Gamma=\Gamma(M)={\bf Z}^d$ das zugeh"orige
Quotientengitter. Wir fassen kurz einige Tatsachen "uber Monoidringe
zusammen, vergl. \cite{fulton} und \cite{rentzsch}.
\par\smallskip\noindent
Eine von $M$ verschiedene maximale Seite (das ist ein Untermonoid $S$
mit $s+t \in S \Longrightarrow s \in S$ und $t \in S$)
nennt man eine {\it Facette}.
Es seien $F_1,...,F_r$ die Facetten von $M$.
Zu $F_i$ ist ${\bf p}_i=(m \in M:\, m \not\in F_i)$
ein Primideal der H"ohe eins in $K[M]$.
Die Nenneraufnahme $M_{F_i} \cong {\bf Z}^{d-1} \times {\bf N}$
liefert eine Bewertung
$\nu_i: \Gamma \longrightarrow {\bf Z}$ mit $\nu_i(M) \subseteq {\bf N}$,
die der zum diskreten Bewertungsring $K[M]_{{\bf p}_i}$
geh"orenden Bewertung entspricht.
\par\smallskip\noindent
Die Abbildung $\nu=(\nu_1,...,\nu_r): \Gamma \longrightarrow {\bf Z}^r$
ordnet jedem Monom die Ordnung an jeder Facette zu.
Sie ist injektiv mit $M=\Gamma \cap {\bf N}^r$, und ist daher eine Darstellung
mit der Durchschnittseigenschaft, die wir die Divisorenklassendarstellung
von $M$ nennen.
Dabei ist ${\bf Z}/\Gamma$ die Divisorenklassengruppe von $K[M]$.
\par\smallskip\noindent
Die zugeh"orige Ringabbildung
$K[M] \hookrightarrow K[{\bf N}^r]=K[X_1,...,X_r]$, unter der
${\bf p}_i =K[M] \cap (X_i)$ ist, spielt bei der Bestimmung der
koaffinen Divisoren eine entscheidende Rolle.
\par\smallskip\noindent
Zu einem Element $f \in M$ ist $\nu_i(f) >0$ genau dann, wenn $f \in {\bf p}_i$
ist, und genau dann, wenn $f \not\in F_i$ ist.
Wir nennen ${\rm supp}\, (f)=\{F_j: f \not\in F_j \}$
den Tr"ager von $f$.
\begin{lem}
Sei $K[M]$ wie oben.
Dann besitzt ein effektiver Divisor $D=n_1{\bf p}_1+...+n_r{\bf p}_r$
genau dann affines Komplement, wenn der Tr"ager von $D$
gleich dem Tr"ager eines Monoms $f \in M$ ist
\end{lem}
{\it Beweis}.
Sei $D=n_1F_1+...+n_sF_s$ ein effektiver Divisor,
$0 \leq s \leq r,\, n_i >0$
f"ur $i=1,...,s$ (die Facetten seien so geordnet),
dessen Tr"ager $\{F_1,...,F_s \}$ nicht durch eine Funktion
beschrieben werden kann.
Das hei"st, dass jede Funktion $f \in M$, die an den Facetten $F_1,...,F_s$
einen positiven Wert hat, noch auf einer weiteren Facette einen positiven
Wert hat. Dann ist jede monomiale Funktion aus dem beschreibenden Ideal
${\bf a}={\bf p}_1 \cap ... \cap {\bf p}_s$ noch in einem weiteren
Primideal ${\bf p}_j,\, s<j \leq r,$ enthalten.
\par\smallskip\noindent
Andererseits gibt es auf jeder Facette Monoidelemente im Inneren der Facette,
die also nur auf dieser einzigen Facette liegt.
\par\smallskip\noindent
Wir betrachten die Divisorenklassendarstellung
$K[M] \longrightarrow K[X_1,....,X_r]$.
Das Erweiterungsideal ${\bf a}K[X_1,...,X_r]$ liegt einerseits in
$(X_i)$ f"ur $i=1,...,s$
(es ist ${\bf a}=(X_1 \cdot ... \cdot X_s) \cap K[M]$), andererseits
aber auch in
$(X_{s+1},...,X_r)$: ist n"amlich $g \in {\bf a}$, so ist
$g \in {\bf p}_j $ f"ur ein $s+1 \leq j \leq r$
und damit ist $g \in (X_j)$ in $K[X_1,...,X_r]$.
Das Erweiterungsideal kann aber auch nicht in einer
Hyperebene ${\bf a}K[X_1,...,X_r] \subseteq (X_j) \subseteq (X_{s+1},...,X_r)$
liegen, da es f"ur jedes $1 \leq j \leq r$ ein Monom der
Form $X_1^{\nu _1}...X_r^{\nu _r}$ gibt mit $\nu_j=0$ und $\nu_i >0$ f"ur
$i \neq j$. Ist $j \geq s+1$, so geh"ort dieses Monom zum
Erweiterungsideal, aber nicht zu $(X_j)$.
Das Erweiterungsideal hat also nicht die reine Kodimension eins und daher
kann das Komplement des Tr"agers nicht affin sein. \hfill $\Box$
\begin{satz}
In einem Monoidring $K[M]$ zu einem normalen, torsionsfreien, endlich
erzeugten Monoid $M$ "uber einem faktoriellen Bereich $K$ sind f"ur
einen Divisor $D$ folgende Aussagen "aquivalent.
\par\smallskip\noindent
{\rm (1)}
Ein Vielfaches von $D$ ist ein Hauptdivisor.
\par\smallskip\noindent
{\rm (2)}
Jedes Vielfache von $D$ ist koaffin.
\par\smallskip\noindent
{\rm (3)}
$D$ ist affin-trivial.
\par\bigskip\noindent
Die affine Klassengruppe von $K[M]$ ist die Divisorenklassengruppe modulo
Torsion.
\end{satz}
{\it Beweis}.
Von (1) nach (3) und von (3) nach (2) gelten immer.
Sei also $D=n_1F_1+...+n_rF_r$ so, dass alle Vielfachen davon koaffin sind,
so haben wir
zu zeigen, dass ein Vielfaches davon ein Hauptdivisor ist.
Wir gehen zu einem effektiven Vertreter von $D$ "uber und k"onnen
daher die $n_i$ als nichtnegativ annehmen.
Nach dem Lemma ist dann der Tr"ager durch eine Funktion $f \in M$
beschreibbar.
Dann gibt es nat"urliche Zahlen $k$ und $m$ derart,
dass $kD-m(f)$ effektiv ist und
wo zumindest eine null hinzukommt. Dieser neue Divisor ist zu
$kD$ "aquivalent und daher wiederum koaffin, und so gelangt man induktiv
zu einem $t$ mit $tD-s(g)=0$. \hfill $\Box$
\begin{kor}
F"ur einen Monoidring $K[M]$ zu einem normalen Monoid $M$ sind folgende
Aussagen "aquivalent.
\par\smallskip\noindent
{\rm (1)}
$M$ ist simplizial, d.h. die Anzahl der Facetten ist gleich der Dimension.
\par\smallskip\noindent
{\rm (2)}
$K[M]$ ist fastfaktoriell.
\par\smallskip\noindent
{\rm (3)}
${\rm Spek}\, K[M]$ ist eine Quotientensingularit"at.
\par\smallskip\noindent
{\rm (4)}
Die affine Klassengruppe von $K[M]$ verschwindet. 
\end{kor}
{\it Beweis}.
Die "Aquivalenz von (1) und (2) folgt aus der expliziten Beschreibung
der Divisorenklassengruppe als Kokern der exakten
Sequenz $0 \longrightarrow \Gamma \longrightarrow {\bf Z}^r$,
die einer rationalen monomialen Funktion den Hauptdivisor zuordnet.
Der Kokern ist genau dann eine Torsionsgruppe,
wenn die beiden freien Gruppen gleichen Rang haben.
Aus (2) folgt (3) aufgrund der Divisorenklassendarstellung
$K[M] \longrightarrow K[X_1,...,X_r]$, bei der $K[M]$ die nullte Stufe
unter der durch ${\bf Z}^r \longrightarrow {\rm DKG}\, K[M]$
gegeben Graduierung ist.
Diese Graduierung entspricht der Operation des endlichen Gruppenschemas
${\rm Spek}\, K[{\rm DKG}\,(K[M])]$ auf ${\bf A}_K^r$.
und daher liegt eine Quotientensingularit"at vor.
Aus (3) folgt (4) aufgrund des Satzes von Chevalley 1.14.1.
Aus (4) folgt (2) aufgrund des Satzes. \hfill $\Box$
\par
\bigskip
\medskip
\newpage
{\it Determinantenringe}
\par
\bigskip
\noindent
Als dritte Beispielklasse untersuchen wir Determinantenringe und
orientieren uns dabei an \cite{bruns}, II.7.
Sei $K$ ein K"orper und $0<k\leq {\rm min}\,(m,n)$.
Die $k-$Minoren der $(m \times n)-$Variablenmatrix
$$\left (\matrix{ X_{11} & ... & X_{1n} \cr
... & ... & ... \cr X_{m1} & ... & X_{mn} \cr} \right) $$
definieren ein Ideal $I_k$ im Polynomring
$K[X_{ij}:\, 1 \leq i \leq m,\, 1 \leq j \leq n]$.
Die Restklassenringe $R_k=K[X_{ij}]/I_k$ nennt man Determinatenringe
oder Minorenringe.
Die Ideale $I_k$ sind prim der H"ohe $mn-(m+n-k+1)(k-1)$ und die
Determinantenringe $R_k$ sind normale Cohen-Macauley-Bereiche der
Dimension $(m+n-k+1)(k-1)$, siehe \cite{bruns}, Theorem 7.3.1.
Die Divisorenklassengruppe und die affine Klassengruppe l"asst sich wie folgt
beschreiben.
\begin{satz}
Sei $1< k \leq {\rm min}(m,n)$.
Es sei $\bf p$ das Ideal von $R_k$, das durch die
$(k-1)-$Minoren der ersten $(k-1)$ Zeilen erzeugt wird.
Dann gelten folgende Aussagen.
\par\smallskip\noindent
{\rm (1)}
$\bf p$ ist ein Primideal der H"ohe 1 und
${\rm DKG}\, R_k={\bf Z}\cdot {\bf p}$.
\par\smallskip\noindent
{\rm (2)} Die affine Klassengruppe von $R_k$ ist gleich der
Divisorenklassengruppe, also ${\rm AKG}\, R_k={\bf Z}$.
\end{satz}
{\it Beweis}.
F"ur die erste Aussage siehe \cite{bruns}, Theorem 7.3.5.
F"ur die zweite Aussage gen"ugt es zu zeigen, dass das
Komplement von $V({\bf p})$ nicht affin ist.
Dazu betrachten wir die Abbildung
$$K[X_{ij}:\, 1 \leq i \leq m,\, 1 \leq j \leq n]/I_k
\longrightarrow K[X_{ij}:\, 1 \leq i \leq k-1,\, 1 \leq j \leq n]$$
mit $X_{ij} \longmapsto X_{ij}$ f"ur $i \leq k-1$ und
$X_{ij} \longmapsto 0$ f"ur $i \geq k$.
Da nur die ersten $k-1$ Zeilen nicht annuliert werden, verschwinden alle
$k-$Minoren, so dass die Abbildung wohldefiniert ist.
Das Bild von $\bf p$ ist einfach das Ideal $I_{k-1}$ der
$(k-1)-$Minoren im $(k-1)n-$dimensionalen Polynomring.
Dessen H"ohe ist aber $(k-1)n-(k-1+n-(k-1)+1)(k-2)=nk-n-(n+1)(k-2)=n-k+2$.
Wegen $k \leq n$ ist also die H"ohe des Erweiterungsideals $\geq 2$, und
$D({\bf p})$ kann nicht affin sein. \hfill $\Box$
\par\bigskip
\noindent
{\bf Bemerkung}
Anschlie"send an die Beispiele ist zu fragen, welche lokalen
normalen Ringe die Eigenschaft haben, dass die affine Klassengruppe gleich
der Divisorenklassengruppe modulo Torsion ist. Solche Ringe haben insbesondere
die Eigenschaft, dass man aus dem Verschwinden der affinen Klassengruppe
darauf schlie"sen kann, dass sie fastfaktoriell sind.
Die Beispiele sind ("uber einem K"orper) 
rationale Singularit"aten im Sinne von
\cite{binsto2}, es ist zu fragen, ob diese Beschreibung
der affinen Klassengruppe f"ur rationale Singularit"aten immer gilt.
\par
\bigskip
\noindent
\subsection{Koaffine und ger"aumige Divisoren}
\par
\bigskip
\noindent
In diesem Abschnitt untersuchen wir den Zusammenhang von verschiedenen
Begriffen f"ur Geradenb"undel auf eigentlichen Variet"aten, insbesondere
das Verh"altnis von ger"aumig und koaffin. Ferner betrachten wir das
Problem, ob man zu einer offenen affinen Teilmenge einen Divisor mit dem
Komplement als Tr"ager finden kann, der ger"aumig oder koaffin ist.
Hierzu gibt es "uberzeugende Vorarbeiten von Goodman, auf die wir uns
st"utzen k"onnen.
\begin{satz}
Sei $X$ ein eigentliches Schema "uber einem noetherschen Ring $K$ und
$\cal L$ ein Geradenb"undel auf $X$.
Dann sind folgende Aussagen "aquivalent.
\par\smallskip\noindent
{\rm (1)} ${\cal L}$ ist potentiell basispunktfrei und
f"ur alle $n \geq 1$ ist ${\cal L}^n$ koaffin.
\par\smallskip\noindent
{\rm (2)} Die kanonische Abbildung
$f: X \cdots \longrightarrow {\rm Proj}\, \sum_{n \geq 0}{\cal L}^n$ ist global
definiert und ein affiner Morphismus.
\par\smallskip\noindent
{\rm (3)}
Es gibt eine basispunktfreie positive Potenz von $\cal L$
derart, dass jeder Schnitt
daraus eine Nullstelle mit affinem Komplement besitzt.
\par\smallskip\noindent
{\rm (4)} Es gibt eine global definierte, affine Abbildung
$ \varphi: X \longrightarrow {\bf P}^m_{\Gamma(X,{\cal O}_X)}$ mit
$\varphi^{-1}{\cal O}(1)={\cal L}^n,\, (n \geq 1)$.
\par\smallskip\noindent
{\rm (5)} $\cal L$ ist ger"aumig.
\par\bigskip
\noindent
Ist $X$ ein beliebiges noethersches Schema, so sind {\rm (1)} und {\rm (2)}
"aquivalent und es gilt
{\rm (2)} $\Longrightarrow$ {\rm (3)} $\Longrightarrow $ {\rm (4)}
$\Longleftrightarrow $ {\rm (5)}.
\end{satz}
{\it Beweis}.
Wie beweisen zun"achst die Aussagen im allgemeinen Fall.
Die "Aquivalenz von (1) und (2) ist dabei klar, da zu jedem
$s \in {\cal L}^n(X)$ gilt $X_s=f^{-1}(D_+(s))$.
Von (1) nach (3).
Sei ${\cal L}^n$ basispunktfrei, darin hat dann jeder Schnitt die
Affinit"atseigenschaft.
Von (3) nach (4). Man nimmt sich aus ${\cal L}^n$ endlich viele Schnitte
$s_0,...,s_m$ heraus, so dass $X$ von $X_{s_i}$
"uberdeckt wird. Die dadurch definierte Abbildung ist dann affin.
Die "Aquivalenz von (4) und (5).
Die Mengen $X_{s_i}=\varphi^{-1}(D_+(x_i)),\, i=0,...,m,$
sind affin und "uberdecken $X$.
Daraus folgt aber bereits ger"aumig, siehe \cite{EGAII}, 4.5.2.
Umgekehrt bedeutet ger"aumig, dass es eine "uberdeckende affine Familie
von Mengen $X_{s_i}$ gibt mit $s_i \in {\cal L}^n$. Diese definieren dann
einen affinen Morphismus.
\par\smallskip\noindent
Sei nun $X$ zus"atzlich eigentlich "uber einem noetherschen Ring $K$ und
$\cal L$ ger"aumig. Dann ist eine positive Potenz basispunktfrei, und
die kanonische Abbildung ist eine abgeschlossene Einbettung, also erst recht
affin. \hfill $\Box$
\par\bigskip
\noindent
{\bf Bemerkung}
Ohne die Voraussetzung eigentlich sind die anderen Aussagen nicht "aquivalent.
Auf einem quasiaffinen Schema ist jede invertierbare Garbe ger"aumig,
sobald es aber einen Schnitt mit nicht affinem Komplement gibt,
ist (3) (und (1)) nicht erf"ullt. Insbesondere erf"ullt die Strukturgarbe
die Bedingung (3) genau dann, wenn $X$ affin ist.
\par\smallskip\noindent
Es kann auch eine basispunktfreie koaffine invertierbare Garbe geben, ohne
dass jede Potenz davon koaffin ist. Dann ist (3) erf"ullt, aber nicht
(1).
Betrachte dazu das durch die Gleichung $xy-z^2$ "uber einem K"orper
definierte quasiaffine Schema ohne den abgeschlossenen Punkt.
Darin ist jeder positive Divisor koaffin, 
das Schema ist regul"ar, und der Divisor $(x,z)$ ist ein Torsionselement.
Die zugeh"orige invertierbare Garbe $\cal L$ ist basispunktfrei und da sie
nicht trivial ist, f"uhrt jeder Schnitt daraus zu einer echten (nichtleeren)
Kurve mit affinem Komplement. Dagegen ist ${\cal O}_X={\cal L}^2$ trivial,
und davon hat der Einsschnitt leere Nullstellenmenge mit nicht affinem
Komplement.
\par\bigskip
\noindent
Sei $X$ ein normales Schema und $U$ eine offene affine Teimenge mit dem
Komplement $X-U=Y=Y_1 \cup ... \cup Y_r$.
Es stellt sich die Frage, ob es einen koaffinen Divisor mit $Y_1,...,Y_r$ als
Tr"ager gibt.
Ist $X$ affin, so ist hier nichts zu erwarten. Ist $Y\subseteq X$ eine
irreduzible Hyperfl"ache mit nicht affinem Komplement, und so,
dass $Y \cup Z$ affines Komplement besitzt, wobei $Z$ durch ein
Primelement beschrieben sei. Dann findet man keinen koaffinen
Divisor, der den Tr"ager $Y \cup Z$ besitzt.
Ein Divisor $nY+mZ$ mit $n,m >0$
ist ja linear "aquivalent zu $nY$ und kann daher nicht koaffin sein.
\par\smallskip\noindent
Wird dagegen $X$ als projektive Variet"at "uber einem algebraisch
abgeschlossenen K"orper vorausgesetzt, so ist das eine interessante
Frage, und wurde in der Form, ob es einen
ger"aumigen Divisor mit dem Komplement als Tr"ager gibt, von Goodman
studiert.
Es stellt sich heraus, dass das nicht der Fall ist,
und zwar kann man im Allgemeinen weder einen
basispunktfreien noch einen koaffinen Divisor mit dem Komplement
der affinen Menge als Tr"ager finden.
\par\smallskip\noindent
In Anlehnung an der Charakterisierung (4) des Satzes 3.7.1
lautet die Fragestellung, ob es zu $U \subseteq X$ einen global
definierten, affinen Morphismus nach ${\bf P}_K^m$ gibt, unter dem sich $U$
als Urbild einer Menge der Form $D_+(x_0)$ erhalten l"asst.
Diese Formulierung legt die allgemeinere Fragestellung nahe, inwiefern
eine affine Teilmenge sich als Urbild einer affinen Menge unter einem
affinen Morphismus in ein vergleichsweise einfaches Schema erhalten l"asst.
Als Bildr"aume sind da neben dem projektiven Raum beispielsweise
niedrigerdimensionale Schemata untersuchenswert.
Diese Thematik korrespondiert spiegelbildlich zu dem Konzept, die
Nichtaffinit"at durch die Nichtaffinit"at des Urbildes unter einer
affinen Abbildung von einem einfacheren Schema zu erweisen.
In 5.3 werden wir diesen Gesichtspunkt wieder aufgreifen, hier
gehen wir aber nur auf die Frage nach dem Zusammenhang von affin
und eigentlich ein und referieren die Ergebnisse von Goodman.
\begin{satz}
Sei $X$ eine eigentliche normale Variet"at "uber einem
algebraisch abgeschlossenen K"orper $k$ und $U \subseteq X$ affin.
Dann gibt es einen Divisor mit dem Komplement als Tr"ager ohne fixe
Komponente.
\end{satz}
{\it Beweis}.
Ist $Y=X-U$ irreduzibel, so gibt es wegen der Affinit"at eine rationale
Funktion $q$, die l"angs $Y$ ihren einzigen Pol hat, sagen wir der Ordnung $n$.
Dann ist $nY+(q)$ ein zu $nY$ "aquivalenter Divisor, der $Y$ nicht als
Komponente hat. Im Allgemeinen muss man sich eine Funktion basteln, die
genau an den Komponenten von $Y$ einen Pol hat, siehe \cite{goodman},
Prop. 2. \hfill $\Box$
\par\bigskip
\noindent
Im Fl"achenfall kann man mehr erreichen.
\begin{satz}
Sei $U \subseteq X$ eine offene affine Teilmenge einer eigentlichen
Fl"a\-che, derart, dass jeder Punkt des Komplements von $U$
lokal faktoriell ist, dann gibt es einen effektiven ger"aumigen Divisor
mit dem Komplement als Tr\"ager.
\end{satz}
F"ur den Beweis siehe \cite{goodman}, Theorem 2,
oder \cite{haramp}, Theorem II.4.2.
Wir geben einen Beweis f"ur den einfachsten Fall, dass $X$ nichtsingul"ar
und $U$ das Komplement einer irreduziblen Kurve $Y$ ist.
In diesem Fall gibt es zu einem positiven Vielfachen von
$Y$ einen "aquivalenten effektiven
Divisor $D$, der $Y$ nicht als Komponente hat. 
Wegen der Affinit"at von $U$ gilt $Y.C > 0$ f"ur $C \neq Y $ und
damit auch $kY.kY=kY.D >0$, also ist auch der Selbstschnitt positiv.
Aus dem Schnittkriterium folgt dann die Ger"aumigkeit von $Y$. \hfill $\Box$
\par\bigskip
\noindent
Ist $X$ eine eigentliche, aber nicht projektive Variet"at, so gibt es
darauf "uberhaupt keine ger"aumige invertierbare Garbe, und das
Komplement einer affinen Teilmenge kann somit von vornherein nicht Tr"ager
einer ger"aumigen Garbe sein. Der Satz 3.7.3 beweist somit zugleich,
dass eigentliche Fl"a\-chen mit den angegebenen Eigenschaften,
insbesondere glatte Fl"a\-chen, projektiv sind.
Im dreidimensionalen gibt es auch glatte Beispiele f"ur eingentliche, nicht
projektive Variet"aten, siehe \cite{shafarevich}, VI.\S2.3 oder
\cite{haralg}, Ap.B, example 3.4.1.
\par\smallskip\noindent
Die folgenden beiden Beispiele zeigen, dass es auch im Fall von
projektiven Variet"aten nicht immer einen ger"aumigen Divisor mit dem
Komplement als Tr"ager geben muss.
Im ersten Beispiel ist das irreduzible Komplement zwar
potentiell basispunktfrei, es gibt aber andere nicht-affine Vertreter,
so dass der Divisor nicht koaffin ist, w"ahrend im zweiten Beispiel
alle Vielfachen des Divisors eine Kurve als gemeinsamen Basisort
besitzen, der Divisor aber koaffin ist.
\par\bigskip
\noindent
{\bf Beispiel} (vergl. \cite{goodman}, I)
Eine offene Teilmenge in einer projektiven, regul"aren
Variet"at kann affin sein, ohne dass es einen koaffinen Divisor mit dem
Komplement als Tr"ager gibt.
Dazu sei $X' \longrightarrow X$ eine birationale Abbildung zwischen
dreidimensionalen projektiven Variet"aten derart, dass genau eine projektive
Kurve $C\subseteq X'$ auf einen Punkt $P \in X$ kontrahiert wird, und
au"serhalb davon eine Isomorphie vorliegt.
Es liegt also eine kleine Kontraktion vor, siehe den Abschnitt 3.4.
$X'$ sei glatt.
\par\smallskip\noindent
Sei $\cal L$ eine sehr ger"aumige Garbe auf $X$, und $H$ ein irreduzibler
Hyperebenenschnitt
durch $P$. Es ist $X-H$ affin und damit hat auch das Urbild von $H$, die
totale Transformierte $\bar{H}$ von $H$, affines Komplement.
$\bar{H}$ ist die Vereinigung der Kurve
$C$ mit der eigentlich Transformierten $H'$, und daher muss
$C \subseteq H'$ sein.
Also ist $U=X'-H'=X'- \bar{H}$ affin und hat den irreduziblen Divisor $H'$
als Komplement, es kommen also nur Vielfache von $H'$ in Frage.
\par\smallskip\noindent
$H'$ ist aber kein koaffiner Divisor.
Zu einem Hyperebenenschnitt $H_1$ von $\cal L$ auf $X$,
der nicht durch $P$ geht,
ist die totale Transformierte $H_1'$ disjunkt zu $C$ und
"aquivalent zu $H'$. Mit $C \subseteq X'-H_1'=U_1$ kann $H'$ nicht koaffin sein.
Dagegen ist $H'$ basispunktfrei, und damit auch numerisch effektiv, aber
nat"urlich nicht numerisch positiv, wegen $H'.C=0$. Weiterhin ist
der Selbstschnitt $ (H')^3 >0$, wegen asymptotischen Riemann-Roch, siehe
\cite{kollar}, Theorem VI.2.15.
\par\bigskip
\noindent
{\bf Beispiel}
Eine offene Teilmenge in einer projektiven glatten Variet"at kann affin sein,
ohne dass es einen basispunktfreien Divisor mit dem Komplement als Tr"ager
gibt. Zu Einzelheiten und Beweise zu diesem Beispiel vergleiche
vor allem \cite{goodman}, III., ferner \cite{zariskiriemann}, 2, und
\cite{haramp}, II.5.
\par\smallskip\noindent
Es sei $S={\bf P}_K^2 \subseteq {\bf P}^3_K$ "uber einem algebraisch
abgeschlossenen K"orper $K$ der Charakteristik null,
und $C \subseteq S={\bf P}^2_K$ sei eine glatte Kurve vom Grad drei.
Es sei $Q \in C$ ein Punkt, der durch eine Gerade herausgeschnitten wird,
sagen wir $Q=C \cap V_+(x_0)$.
Wir versehen $C$ mit der durch $Q$ als neutralem Element festgelegten
Gruppenstruktur.
Die invertierbare Garbe ${\cal O}(4)$ auf ${\bf P}^2_K$ induziert die
Garbe ${\cal O}_C(4)$ auf $C$, es sei $R$ ein effektiver Vertreter davon
mit den Tr"agerpunkten $R_1,...,R_k$.
\par\smallskip\noindent
Es seien $P_1,...,P_{12}$ Punkte auf $C$ so,
dass sie zu $R_1,...,R_k$
linear unabh"angig sind.
Damit sind dann insbesondere $n_1P_1+...+n_{12}P_{12}$ und $nR$ verschiedene
Divisorenklassen.
\par\smallskip\noindent
Es sei $L$ eine Gerade auf $S$, die die Punkte $P_i$ nicht trifft.
Man betrachtet das lineare System $V \subseteq {\cal O}_{{\bf P}^3_K}(4)$
der Fl"a\-chen im ${\bf P}_K^3$ vom Grad vier, auf denen sowohl $C$
als auch $L$ liegt.
Mit dem Satz von Bertini kann man sehen, dass die generische Fl"a\-che
aus $V$ irreduzibel und an den $P_i$ glatt und nicht tangential an $S$ ist.
\par\smallskip\noindent
Sei nun $X$ die Aufblasung der zw"olf Punkte von ${\bf P}^3_K$,
$\bar{C}, \bar{S}, \bar{D}$ seien die eigentlichen Transformierten davon,
die den Aufblasungen der zw"olf Punkte auf $C,S$ und $D$ entsprechen.
\par\smallskip\noindent
Die irreduziblen Divisoren $n\bar{D},\, n \geq 1,$ besitzen die Kurve $\bar{C}$
als Basiskurve, so dass sie nicht basispunktfrei und nicht ger"aumig sind.
Nach \cite{goodman} gibt es im Komplement von $\bar{D}$
keine projektiven Kurven.
Da der schematheoretische Basisort $B_n$ konstant ist,
folgt aus \cite{haramp} II, Theorem 5.1, dass das
Komplement von $\bar{D}$ affin ist.
\par\bigskip
\noindent
Nach den negativen Resultaten, dass das Komplement einer offenen affinen Menge
nicht Tr"ager eines ger"aumigen Divisors sein muss, zeigt der folgende Satz
von Goodman, dass dies zumindest in einer geeigneten Aufblasung erreicht
werden kann.
\begin{satz}
Eine offene Teilmenge $U \subseteq X$ in einer eigentlichen Variet\"at
ist genau dann affin, wenn es eine Aufblasung $\tilde{X}$ eines Ideals
mit Nullstellenmenge au"serhalb von $U$ gibt derart,
dass das Komplement von $U$ in
der Aufblasung Tr\"ager eines effektiven ger"aumigen Cartierdivisors $D$ ist.
\end{satz}
Zum {\it Beweis} siehe \cite{haramp}, II.Theorem 6.1. oder \cite{goodman},
Theorem 1. \hfill $\Box$
\par\bigskip
\noindent
Aus diesem Resultat kann eine wichtige geometrische Folgerung gezogen
werden, n"amlich, dass das Komplement einer affinen Teilmenge in einer
eigentlichen Variet"at der Dimension $\geq 2$ zusammenh"angend ist.
Dies kann man n"amlich mit dem Aufblasungssatz darauf zur"uckf"uhren, dass
effektive ger"aumige Divisoren zusammenh"angend sind, siehe \cite{haramp},
Cor. II.6.2.
\par
\bigskip
\noindent
\subsection{Affine Klassengruppe von graduierten Algebren}
\par
\bigskip
\noindent
Sei $K$ ein noetherscher normaler Integrit"atsbereich und
$S$ eine positiv graduierte, normale $K-$Algebra, die von der ersten Stufe
endlich erzeugt sei.
Es sei $X={\rm Spek}\, S, \, X'=X-V(S_+)$ und $ Y={\rm Proj}\, S$
das projektive Schema "uber ${\rm Spek}\, K$.
Auf dem punktierten Spektrum $X'$ ist die Kegelabbildung
$p: X' \longrightarrow X$ ein Hauptfaserb"undel,
das lokal die Gestalt
${\rm Spek}\, B[T,T^{-1}] \longrightarrow {\rm Spek}\, B$ hat.
Insbesondere ist die Kegelabbildung treu\-flach, $Y$ ist ebenfalls normal
und lokal stimmen die Divisorenklassen oben und unten "uberein.
Insbesondere ist jede Divisorklasse auf $X'$
"aquivalent zu einer Divisorklasse, die von unten kommt.
\par\bigskip
\noindent
Eine besondere Rolle spielt dabei die durch ein homogenes
Element $x$ vom Grad eins definierte Divisorklasse $H$ auf $Y$.
Sind ${\bf p}_1,...,{\bf p}_r$ die Primideale
der H"ohe eins, die $x$ enthalten, so sind diese alle homogen, und man
nennt den Divisor $n_1Z_1+...+n_rZ_r$, wobei $Z_i=V_+({\bf p}_i)$ ist und
die $n_i$ die Ordnungen von $x$ in $S_{{\bf p}_i}$ sind,
einen Hyperebenenschnitt.
F"ur zwei Elemente $0 \neq x,y \in S_1$ ist der
Quotient $ x/ y$ ein Element des Funktionenk"orpers von $Y$,
und dabei ist $\nu_Z(x/y)= \nu_{p^{-1}(Z)}(x) - \nu_{p^{-1}(Z)}(y)$,
weshalb die Divisorklasse eines Hyperebenenschnittes eindeutig bestimmt ist.
Die zugeh"orige Klasse $p^*H$ auf $X$ ist nat"urlich trivial.
\par\bigskip
\noindent
In diesem Abschnitt wollen wir das Verh"altnis von Koaffinit"atseigenschaften
einer Divisorklasse $D$ auf $Y$ zu solchen von $p^*D$ auf $X'$ bzw. auf $X$
kl"aren und insbesondere untersuchen, wie sich Koaffinit"at und affine
Trivialit"at eines Divisors auf dem Kegel durch Eigenschaften
auf dem projektiven Schema charakterisieren lassen. Dabei behalten wir die
eingangs erw"ahnten Bezeichnungen bei.
\begin{lem}
Sei $D$ ein Weildivisor auf $Y$ mit zugeh"origer reflexiver Modulgarbe
${\cal L}_D={\cal L}$.
Dann gilt auf dem punktierten Kegel $X'$ die Gleichheit
${\cal L}_{p^*D}(X')=\sum_{k \in {\bf Z}} {\cal L}_{D+kH}(Y)$. Hat $V(S_+)$
eine Kodimension $\geq 2$, so gilt dies auch in $X$ selbst.
\end{lem}
{\it Beweis}.
Da $Y$ normal ist, sind die auftretenden Garben in der Kodimension eins
invertierbar. Es sei $W \subseteq Y$ so gew"ahlt, dass sowohl ${\cal L}$
als auch ${\cal H}={\cal L}_H$
darauf invertierbar sind und so, dass $W$ die Kodimension
eins umfasst.
Dann umfasst $W'=p^{-1}(W)$ auf $X'$
ebenfalls die Kodimension eins, und wir k"onnen in der Aussage $Y$ und $X'$
durch $W$ und $W'$ ersetzen.
Auf $W'$ ist ${\cal L}_{p^*D}  =p^*({\cal L}_D)$, also
ist ${\cal L}_{p^*D}(W')=p^*({\cal L}_D)(W') $.
\par\smallskip\noindent
Auf einer affinen offenen Teilmenge $U \subseteq W$ mit affinem Urbild 
$U'$ ist
$$ \Gamma(U',{\cal O}_X)_k={\cal L}_{kH}(U)
=\{ q \in K(Y):\, (q)+kH \geq 0 \mbox{ auf } U\} \, ,$$
wobei ein Element $f$ vom Grad $k$ auf $f/x^k$ abgebildet wird.
Die Abbildung ist injektiv, sei $kH+(q) \geq 0$ mit $q \in K(Y)$.
Dann ist im Kegel ebenfalls $kH+(q) \geq 0$, also ist
$(x^kq)=k(x)+(q) \geq 0$ auf $U$, und damit ist
$x^kq \in \Gamma(U',{\cal O}_X)$.
Dabei muss $x^kq$ aus der $k-$ten Stufe sein.
\par\smallskip\noindent
Auf affinem $U$ ist $ p^*{\cal L}(U')=
{\cal L}(U) \otimes_{\Gamma(U,{\cal O}_Y)} \Gamma(U',{\cal O}_X)$,
und das ist graduiert durch
\begin{eqnarray*}
(p^*{\cal L}(U'))_k & = &
{\cal L}(U) \otimes_{\Gamma(U,{\cal O}_Y)} \Gamma(U',{\cal O}_X)_k \\
&=&{\cal L}(U) \otimes_{\Gamma(U,{\cal O}_Y)} {\cal L}_{kH}(U) \\
&=&({\cal L} \otimes_{{\cal O}_Y} {\cal H}^k)(U) \, .
\end{eqnarray*}
Da die Restriktionen mit diesen Identifikationen vertr"aglich sind, gilt
die Gleichheit auch f"ur $W$ bzw. $W'$.
\par\smallskip\noindent
Hat die Kegelspitze eine Kodimension $\geq 2$, so "andert sich
beim "Ubergang zu $X$ die globalen Schnitte von ${\cal L}_{p^*D}$
nicht. \hfill $\Box$
\par\bigskip
\noindent
{\bf Bemerkung}
Die Voraussetzung an die Kodimension der Spitze ist wesentlich, man denke
etwa an die affine Gerade und einen Punkt als deren projektives Spektrum.
\par\smallskip\noindent
Ist $D$ auf $Y$ koaffin, so muss das nicht f"ur $p^*D$ auf $X'$ gelten.
Das Problem ist hierbei, dass $D-kH$ Schnitte haben kann.
\begin{satz}
In der besprochenen Situation gelten folgende Aussagen.
\par\smallskip\noindent
{\rm (1)}
$p^*D$ ist koaffin auf $X'$ genau dann, wenn f"ur alle $k \in {\bf Z}$
die Divisoren $D+kH$ koaffin auf $Y$ sind.
\par\smallskip\noindent
{\rm (2)}
$p^*D$ erf"ullt die Eigenschaft {\rm (K)} auf $X'$ genau dann,
wenn f"ur alle $m,n \in {\bf Z}$ die Divisoren
$nD+mH$ koaffin oder trivial sind.
\par\smallskip\noindent
{\rm (3)}
$p^*E$ ist genau dann affin-trivial auf $X'$,
wenn zu jedem Divisor $D$, der
die Eigenschaft hat, dass f"ur alle $n,m \in {\bf Z}$ die Divisoren
$nD+mH$ koaffin oder trivial sind, diese Eigenschaft auch f"ur $E+D$ gilt.
\end{satz}
{\it Beweis}.
Zu (1).
Nach dem Lemma sind die globalen Schnitte zu $p^*D$ graduiert
und die homogenen Stufen dazu bestehen aus den globalen Schnitten
von $D+kH$. Die Koaffinti"at kann nach dem Kriterium 3.2.2
mit einem Erzeugendensystem etabliert werden, in der vorliegenden
Situation mit allen homogenen Schnitten. F"ur diese ist nach dem
Satz 2.3.2 die Affinit"at auf dem punktierten Kegel "aquivalent zu der auf dem
projektiven Schema.
\par\smallskip\noindent
Zu (2).
$p^*D$ erf"ulle die Eigenschaft (K), d.h. f"ur alle
$n \in {\bf Z}$ ist $p^*(nD)$ koaffin oder trivial.
Ist $p^*(nD)$ trivial, so ist $nD$ ein Vielfaches von $H$ und
$nD+mH$ ist trivial oder koaffin. Ist $p^*(nD)$ koaffin, so gilt
dies nach Teil (1).
Sei umgekehrt die Bedingung erf"ullt, dass $nD+mH$ koaffin oder trivial ist
f"ur alle $n,m$.
Ohne Einschr"ankung sei $n=1$. Ist dann f"ur alle $m \in {\bf Z}$ der Divisor
$D+mH$ koaffin, so gilt das nach Teil (1) auch f"ur $p^*D$.
Ist dagegen f"ur ein $m$ die Divisorklasse $D+mH$ trivial, so ist $D=-mH$
und es ist auch $p^*D$ trivial.
\par\smallskip\noindent
(3) folgt aus (2), wenn man ber"ucksichtigt, dass alle Divisorklassen
auf $X'$ von unten herkommen. \hfill $\Box$
\par\bigskip
\noindent
{\bf Bemerkung}
In den Aussagen des Satzes kann man nicht den punktierten Kegel
durch den affinen Kegel $X$ ersetzen, auch nicht unter der
Voraussetzung, dass $V(S_+)$ die Kodimension $\geq 2$ hat.
Es folgt zwar dann, dass man aus der Koaffinit"at von $p^*D$ auf $X'$
auf die von $p^*D$ auf $X$ schlie"sen kann, aber nicht umgekehrt.
Dies zeigt schon der triviale Divisor $D$, der ausgenommen den Fall,
dass $Y$ nulldimensional ist, weder auf $Y$ noch auf $X'$ koaffin ist,
auf $X$ dagegen schon. 
\par\bigskip
\noindent
{\bf Beispiel}
Sei $A$ ein noetherscher normaler Integrit"atsbereich.
Die affine Klassengruppe des projektiven Raumes ist dann
${\rm AKG}\,({\bf P}^n_A)={\rm DKG}\, A$ oder $={\rm DKG}\, A+{\bf Z}$.
Sei der Divisor $D+nH$ affin-trivial, $D \in {\rm Div}\, A$.
Da der Hyperebenenschnitt
die Eigenschaft (K) erf"ullt, muss auch $D$,
das von ${\rm Spek}\, A$ kommt und man als effektiv annehmen kann,
als Divisor in ${\bf P}_A^n$ koaffin oder trivial sein, woraus $D=0$ folgt.
Es kann also nur $rH$ affin-trivial sein.
\par\smallskip\noindent
In dem Beispiel scheitert die Affinit"at immer daran, dass es im
Komplement eines Divisors, der vom Basisschema ${\rm Spek}\, A$ herkommt,
projektive Fasern gibt. Die Gr"o"se der affinen Klassengruppe k"onnte
hier f"alschlich suggerieren, die meisten Hyperfl"achen h"atten nicht affines
Komplement. Dem entgegen ist aber, wenn ${\rm AKG}\, A$ trivial ist,
das Komplement jeder Hyperfl"ache, die ${\rm Spek}\, A$ dominiert, affin.
\par\smallskip\noindent
Das Beispiel zeigt insbesondere, dass es im Allgemeinen keine
gruppentheoretische Einschr"ankung an die affine Klassengruppe gibt, im
Gegensatz zur affinen Klassengruppe von affinen Schemata, wo sie
torsionsfrei ist und die Vermutung besteht, dass sie endlich erzeugt ist.
\par\bigskip
\noindent
{\bf Beispiel}
Sei $R$ ein eindimensionaler, normaler, nicht faktorieller Ring und
$S=R[X,Y]$, wobei $X$ und $Y$ den Grad eins haben.
In ${\rm Spek}\, S$ sind nach Satz 3.2.4 alle Hyperfl"achen koaffin.
Die Spitze darin ist eine zu ${\rm Spek}\, R$ isomorphe Kurve,
das projektive Spektrum ist ${\bf P}^1_R$.
Ist $P \in {\rm Spek}\, R$ ein abgeschlossener Punkt,
der nicht durch eine Funktion beschrieben wird,
so ist die Faser dar"uber in $X'$ eine Hyperfl"ache mit nicht affinem
Komplement, ohne dass sie als Divisor trivial ist.
Im Fall eines Grundk"orpers kann diese Situation nicht auftreten.
\begin{prop}
Ist $K$ ein K"orper, so ist ein homogener Divisor $E$ auf
$X$ genau dann koaffin, wenn er trivial ist oder auf $X'$ koaffin ist.
Damit kann man im Fall eines Grundk"orpers in 3.8.2
"uberall $X'$ durch $X$ ersetzen.
\end{prop}
{\it Beweis}.
Sei $E$ koaffin, aber nicht trivial. Wir haben zu zeigen,
dass die Kegelspitze in der Nullstelle von jedem Schnitt liegt, da sich
dann die Komplemente nicht unterscheiden. Dies wiederum gen"ugt f"ur die
homogenen Schnitte zu zeigen.
Ein homogener Schnitt entspricht aber einem effektiven, linear
"aquivalenten, homogenen Divisor, der nach Voraussetzung nicht trivial ist,
also nicht leeren Tr"ager hat. Damit geh"ort die Kegelspitze zum
Tr"ager. \hfill $\Box$
\par\bigskip
\noindent
Im Folgenden betrachten wir nur noch die Situation, wo der Grundring $K$
ein K"orper ist. 
Aus den bisherigen "Uberlegungen erhalten
wir die zum Satz "uber die Vektorraumb"undel
3.2.4 (bzw. 3.3.5) analoge Aussage.
\begin{satz}
Sei $Y={\rm Proj}\, S$ eine projektive normale Variet"at "uber einem
K"orper $K$ mit einer normalen graduierten $K-$Algebra $S$, die von der
ersten Stufe endlich erzeugt sei.
Dann ist genau dann das Komplement jeder Hyperfl"ache von $Y$
affin, wenn dies im affinen Kegel $X={\rm Spek}\, S$ f"ur
jede Hyperfl"ache gilt.
\end{satz} 
{\it Beweis}.
Ist auf dem Kegel jede Hyperfl"ache koaffin, so gilt dies nach Satz 2.3.2
auch f"ur die projektive Variet"at.
Besitze umgekehrt jede Hyperfl"ache auf der projektiven
Variet"at affines Komplement. Nach der vorausgehenden Proposition
und Satz 3.8.2 ist dann auf dem affinen Kegel jeder Divisor vom Typ
$p^*D$ koaffin, also "uberhaupt jeder. \hfill $\Box$
\par\bigskip
\noindent
Mit diesem Satz r"uckt nat"urlich auf einer projektiven Variet"at die von
einem fixierten sehr ger"aumigen Divisor $H$ abh"angige Eigenschaft (K') ins
Interesse, dass zu $m,n \in {\bf Z}$ die Divisoren
$nD+mH$ koaffin oder trivial sind.
In verschiedener Hinsicht scheint f"ur eine projektive Variet"at diese
Bedingung angemessener zu sein als koaffin oder die Eigenschaft (K), und
zwar aus folgendem Grund: ein Divisor, der keine Schnitte hat, ist nach
Definition koaffin, es ist aber gut m"oglich, dass ein numerisch
"aquivalenter Divisor Schnitte hat, die nicht koaffin sind, so dass
in einem gewissen Sinn die Nichtkoaffinit"at mangels Schnitten
verborgen bleiben kann.
\par\bigskip
\noindent
{\bf Beispiel}
Beispielhaft f"ur diese Situation sind die Regelfl"achen
und speziell schon die Produkte $C \times {\bf P}$ mit $g(C) \geq 1$.
Seien $P,Q \in C$ zwei nicht linear "aquivalente Punkte mit den Fasern
$F_P$ und $F_Q$, wobei jede Funktion aus $C-\{ Q\}$ in $Q$ einen Pol
von zumindest zwei habe.
Der Divisor $F_P-F_Q$ ist numerisch trivial, und
$F_P$ ist nicht koaffin, dagegen hat $F_P +(F_P -F_Q)=2F_P-F_Q$
keine Schnitte wegen der Bedingung an $Q$ und ist daher koaffin.
Damit ist der Divisor $F_P-F_Q$ nicht affin-trivial und die affine
Klassengruppe entzieht sich einer Berechnung. Im n"achsten
Abschnitt werden wir sehen, dass im affinen Kegel dar"uber
dieser Divisor affin trival ist und dass der Kegel die affine
Klassengruppe $\bf Z$ hat.
\par\bigskip
\noindent
Dividiert man bei einer projektiven Variet"at diejenigen Divisoren als
{\it affin-trivial im projektiven Sinn} aus der Divisorenklassengruppe
heraus, die die Eigenschaft (K') erhalten, so erh"alt man
die Restklassengruppe ${\rm AGK}'\, (Y)={\rm AKG}'_H\,(Y)$.
Mit dieser Bezeichnung gilt der folgende Satz.
\begin{satz}
Sei $S$ eine normale, positiv graduierte $K-$Algebra "uber einem K"orper $K$,
die von der ertsen Stufe endlich erzeugt sei,
und $H$ der Hyperebenenschnitt (bzw. die sehr ger"aumige Garbe)
auf ${\rm Proj}\, S$.
Dann ist ${\rm AKG}\, S={\rm AKG}'_H\, ({\rm Proj}\, S)$.
\end{satz}
{\it Beweis}. Die Abbildung
$p^*: {\rm DKG}\, ({\rm Proj}\, S) \longrightarrow {\rm DKG}\, S$
ist surjektiv
mit dem Kern ${\bf Z}H$, und $p^*D$ ist genau dann affin-trivial in
${\rm Spek}\, S$, wenn $D$ in ${\rm Proj}\, S$ die angegebene Eigenschaft hat.
Die Aussage gilt dann unmittelbar nach Definition. \hfill $\Box$
\par\bigskip
\noindent
Abgesehen von den Situationen, wo man jede Hyperfl"ache als koaffin erweisen
kann, kann man die affine Klassengruppe des affinen Kegels "uber einer
projektiven Variet"at nur dann explizit ausrechnen, wenn es gelingt,
die relevanten Affinit"atseigenschaften numerisch zu charakterisieren.
In diesen F"allen ist ein Zusammenhang zwischen affiner Klassengruppe
des Kegels und numerischer Klassengruppe zu erwarten.
In diese Richtung weist auch die folgende Aussage.
\begin{prop}
Sei $Y$ eine projektive glatte Variet"at "uber einem algebraisch
abgeschlossenen K"orper,
$H$ ein fixierter ger"aumiger Divisor auf $Y$ und $D$ ein weiterer Divisor.
Sind dann $D$ und $H$ linear abh"angig in der numerischen Klassengruppe
${\rm Num}\, Y$,
so erf"ullt $D$ bez"uglich $H$ die Eigenschaft {\rm (K')}.
\par\smallskip\noindent
Ist dabei $Y={\rm Proj}\, S$ mit einer normalen graduierten
$K-$Algebra $S$, die von der ersten Stufe endlich erzeugt sei,
und $H$ der Hyperebenenschnitt, so erf"ullt f"ur einen zu $H$
numerisch abh"angigen Divisor $D$ der zur"uckgenommene
Divisor $p^*D$ auf dem affinen Kegel $X={\rm Spek}\, S$
die Eigenschaft {\rm (K)} und ist koaffin.
Dies gilt insbesondere f"ur alle numerisch trivialen Divisoren.
\par\smallskip\noindent
Sind in ${\rm Spek}\, S$ alle Divisoren mit der Eigenschaft {\rm (K)}
bereits affin-trivial, so ist ${\rm AKG}\, S$ eine Restklassengruppe
von ${\rm Num}\, Y$ und damit endlich erzeugt.
\end{prop}
{\it Beweis}.
Sei $D \equiv qH$ in ${\rm Num}\, Y \otimes {\bf Q}$
mit $q \in {\bf Q}$.
Wir haben zu zeigen, dass f"ur $n,m \in {\bf Z}$ der Divisor
$nD+mH$ koaffin oder trivial ist.
Es ist $nD+mH \equiv rH,\, r \in {\bf Q},$ und
damit $p(nD+mH) \equiv kH,\, 0 \neq p \in {\bf N},\, k \in {\bf Z}$.
Ist $k >0$, so ist wegen dem Schnittkriterium das $p-$fache und damit
der Divisor selbst ger"aumig, also koaffin.
Ist $k = 0$, so ist $nD+mH$ selbst numerisch trivial und daher entweder
schnittfrei oder trivial, also wieder koaffin oder trivial.
Ist $k< 0$, so ist das $p-$fache schnittfrei, und damit ist der Divisor selbst
schnittfrei, also koaffin. \hfill $\Box$
\par\bigskip
\noindent
{\bf Bemerkung}
Im Beweis st"utzen wir uns als Schnittkriterium nicht auf das
Nakai-Moishezon-Kriterium, da dieses auch die Schnitte $H^{{\rm dim}\, Z}.Z$
f"ur beliebige integre Untervariet"aten $Z \subseteq Y$ miteinbezieht,
w"ahrend die numerische Klassengruppe nur von Kurvenschnitten
abh"angt. Stattdessen ziehen wir das Seshadri-Kriterium f"ur ger"aumige
Divisoren heran, das allein mit Kurvenschnitten arbeitet und damit die
Ger"aumigkeit als eine Eigenschaft der numerischen Klasse erweist,
siehe \cite{kollar}, VI.2.18.
\par\bigskip
\noindent
Sind $D$ und $H$ linear unabh"angig in der numerischen Klassengruppe,
so hat man Ausicht, die Eigenschaft $K'$ zu widerlegen, indem man ein
$nD+mH$ findet mit effektivem, nicht koaffinen Vertreter, wie etwa im
Fall von Regelfl"achen, siehe 3.9.4.
\par
\bigskip
\noindent
\subsection{Affine Klassengruppe von projektiven Fl"a\-chen}
\par
\bigskip
\noindent
Aufgrund den im letzten Abschnitt beschriebenen Ph"anomenen ist die
affine Klassengruppe von eigentlichen Variet"aten schwierig zu berechnen.
In diesem Abschnitt stellen wir f"ur projektive glatte Fl"a\-chen "uber
einem algebraisch abgeschlossenenen K"orper verschiedene Aussagen "uber
die affine Klassengruppe der Fl"a\-chen und ihrer affinen Kegel zusammen.
F"ur die Fl"a\-chen selbst wird die Hauptfrage sein,
ob es nicht koaffine Kurven gibt oder nicht. F"ur den
durch eine sehr ger"aumige Garbe gegebenen affinen Kegel "uber der Fl"a\-che
werden wir auch weitergehende Aussagen "uber die affine Klassengruppe
machen k"onnen. Die affine Klassengruppe des affinen Kegels "uber einer
Regelfl"ache wird berechnet.
\par\bigskip
\noindent
Das Verschwinden der affinen Klassengruppe l"asst sich folgenderma"sen
charakterisieren.
\begin{satz}
Sei $S$ eine projektive glatte Fl"a\-che "uber einem algebraisch abgeschlossenen
K"orper. Dann sind folgende Aussagen "aquivalent.
\par\smallskip\noindent
{\rm (1)} ${\rm AKG}\, S=0$.
\par\smallskip\noindent
{\rm (2)} Jede irreduzible Kurve auf $S$ ist ger"aumig.
\par\smallskip\noindent
{\rm (3)} Jeder positive Divisor ist ger"aumig.
\par\smallskip\noindent
{\rm (4)} Jede irreduzible Kurve $C$ besitzt positiven Selbstschnitt.
\end{satz}
{\it Beweis}.
Von (1) nach (2).
Sei $C$ eine irreduzible Kurve mit $U=X-C$ affin.
Dann hat $C$ mit jeder effektiven, irreduziblen Kurve
$D \neq C$ positiven Schnitt. Es gibt ferner auf $U$ eine nicht konstante
Funktion $q$, die l"angs $C$ einen Pol der Ordnung $n$ hat.
Dann ist $nC+(q)$ ein effektiver, zu $nC$ "aquivalenter Vertreter, in dem
nur Kurven vorkommen mit positivem Schnitt, also ist der Selbstschnitt
auch positiv. Nach dem Schnittkriterium folgt, dass $C$ ger"aumig ist.
\par\smallskip\noindent
Von (2) nach (3) ist wegen dem Schnittkriterium klar, ebenso von (3) nach (4)
und von (3) nach (1).
\par\smallskip\noindent
Sei nun (4) erf"ullt, und $C$ eine irreduzible Kurve mit $C^2 > 0$.
Nach Riemann-Roch gilt dann f"ur $m \geq 0$
$$h^0(mC)+h^0(K-mC) \geq
{m^2 \over 2}C^2 - {m \over 2} C.K + \chi ({\cal O}_S) \, .$$
Die rechte Seite w"achst dabei mit $m$ gegen unendlich, und es ist
$h^0(K-mC)$ durch $h^0(K)$ beschr"ankt, da mit $K-mC+(q)$ auch $K+(q)$
effektiv ist. Also wird $h^0(mC)$ beliebig gro"s, folglich gibt es
auf $U=X-C$ nicht konstante Funktionen.
Sei $q$ eine solche Funktion mit der Polordnung $n$ l"angs $C$.
Der Divisor $nC$ hat dann einen weiteren effektiven Vertreter $nC+(q)$, der
$C$ nicht als Komponente besitzt, daher ist der Divisor komponentenfrei
und ein Vielfaches davon, sagen wir $mC$
ist dann nach dem Satz von Zariski, \cite{zariskiriemann}, Theorem 6.1,
sogar basispunktfrei.
\par\smallskip\noindent
Wir betrachten die durch ein vollst"andiges lineares System aus $mC$ definierte
eigentliche Abbildung nach ${\bf P}^r$.
Das Bild ist nat"urlich nicht nulldimensional, da der Divisor sonst trivial
w"are. Das Bild kann aber auch keine Kurve sein, da sonst die Fasern
alle numerisch "aquivalent w"aren und den Selbstschnitt null h"atten,
was nach Voraussetzung ausgeschlossen ist.
Das Bild ist also eine zweidimensionale, projektive, nicht unbedingt glatte
Fl"a\-che und die Abbildung ist generisch endlich.
Nach \cite{artin}, Theorem 2.3 oder \cite{beau}, Cor. VIII.5,
k"onnen nur Kurven mit negativem Selbstschnitt kontrahiert werden,
also muss die Abbildung quasiendlich und damit endlich sein.
Dann ist aber $C$ nach \cite{haramp}, I. Prop. 4.4. ger"aumig.  \hfill $ \Box$
\par\bigskip
\noindent
{\bf Bemerkung} Die Bedingung $C^2 > 0$ f"ur eine einzige
irreduzible Kurve $C$ reicht nicht f"ur den Schluss auf $X-C$ affin.
Beispiele hierf"ur erh"alt man, wenn man $Y-C$ affin hat und dann
einen Punkt au"serhalb von $C$ aufbl"ast, da sich dabei der
Selbstschnitt nicht "andert. "Ubrigens ist f"ur glattes $C$
die Bedingung $C^2> 0$ "aquivalent zur Lefschetzbedingung
${\rm Lef}\, (X,C)$, siehe \cite{haramp}, IV.1, Exc. 1.7.
Diese Bedingungen sind beide lokal in einer offenen Umgebung von $C$, die
Affinit"at nicht.
\begin{prop}
Sei $S$ eine projektive glatte Fl"a\-che "uber einem algebraisch
abgeschlossenen K"orper, mit ${\rm AKG}\, S=0$.
Dann gelten folgende Aussagen.
\par\smallskip\noindent
{\rm (1)} $S$ ist eine minimale Fl"a\-che.
\par\smallskip\noindent
{\rm (2)} Je zwei Kurven schneiden sich.
\par\smallskip\noindent
{\rm (3)} Jeder Morphismus nach einer Kurve ist konstant.
\par\smallskip\noindent
{\rm (4)} Zu jeder rationalen Funktion $q \in K(S)$ schneiden sich Pol- und
Nullstellendivisor.
\par\smallskip\noindent
{\rm (5)} Es gibt keine rationalen dominanten Abbildungen auf Kurven vom
Geschlecht $\geq 1$. Anders gesagt, alle Unterk"orper von $K(S)$ vom
Transzendenzgrad eins sind Funktionenk"orper in einer Variablen.
\par\smallskip\noindent
{\rm (6)} F"ur alle irreduziblen Kurven $C$ gilt $2g(C)-2 > C.K$, wobei $g(C)$
das arithmetische Geschlecht $h^1(C,{\cal O}_C)$ ist, das auch im
singul"aren Fall definiert ist. Ist $K$ numerisch effektiv, so besitzt
jede Kurve auf der Fl"a\-che ein Geschlecht $\ge 2$.
\par\smallskip\noindent
{\rm (7)}
Die Irregularit"at $q=h^1({\cal O}_S)=h^0(\Omega_S)$ ist ungleich eins.
\par\smallskip\noindent
{\rm (8)} Ist die Kodairadimension positiv,
so muss der kanonische Divisor bereits
ger"aumig sein; insbesondere ist Kodairadimension eins nicht m"oglich.
\end{prop}
{\it Beweis}.
Die ersten vier Aussagen sind klar. Sei $f: S \cdots \longrightarrow C$
eine rationale Abbildung auf eine Kurve. Diese ist in der Kodimension
eins definiert, an den endlich vielen Undefinierbarkeitstellen kann man sie
auf einer geeigneten Aufblasung global fortsetzen, siehe
\cite{beau}, Theorem II.7.
Da $S$ glatt ist, werden dabei Punkte durch projektive Geraden
ersetzt, nach dem Satz
von L"uroth (\cite{haralg}, Example IV.2.5.5)
muss aber deren Bild auf einer Kurve vom Geschlecht $\geq 1$ konstant sein,
und folglich ist die Abbildung schon "uberall definiert, also ein Morphismus
und nach (3) konstant.
\par\smallskip\noindent
Zu (6). Die Absch"atzung mit dem Geschlecht ist eine
Umformulierung von $C^2 > 0$ mittels der Adjunktionsformel, siehe
\cite{haralg}, V.1.5.
\par\smallskip\noindent
Zu (7). Die Aussage zur Irregularit"at zeigen wir im komplexen Fall,
die angegebene
Gleichung beruht auf Hodgetheorie. $q$ ist die
Dimension der Albanesevariet"at ${\rm Alb}\, S$,
siehe \cite{beau}, Theorem V.13,
und es gibt einen (nicht konstanten) Morphismus
$S \longrightarrow {\rm Alb}\, S$, das w"are bei $q=1$
eine Abbildung auf eine Kurve.
\par\smallskip\noindent
Zu (8).
Ist die Kodairadimension positiv, so hat ein Vielfaches $nK$ einen
nicht-trivialen, effektiven Vertreter, und dieser muss ger"aumig sein.
\par\bigskip
\noindent
Nachdem wir einige Folgerungen aus ${\rm AKG}\, S=0$ zusammengetragen haben,
formulieren wir einige hinreichende Bedingungen daf"ur.
\begin{prop}
Sei $S$ eine projektive glatte Variet"at.
Ist die Picardzahl von $S$ gleich eins, so ist ${\rm AKG}\, S=0$.
Dies ist insbesondere bei ${\rm Pic}\, S={\bf Z}$ der Fall
und wird in der komplexen Situation auch durch $b_2=1$ gesichert.
\end{prop}
{\it Beweis}.
Unter der Picardzahl versteht man den Rang der Gruppe der
numerischen Klassen ${\rm Num}\, S$, die man aus der Picardgruppe erh"alt,
wenn man die numerisch trivialen Divisoren herausdividiert.
Ist die Picardgruppe vom Rang eins, so nat"urlich auch ${\rm Num}\, S$.
Nach dem Satz von Neron-Severi ist ${\rm Num}\, S$ eine endlich erzeugte
freie Gruppe. Sei ${\rm Num}\, S={\bf Z}$. Da $S$ projektiv ist, gibt es
einen ger"aumigen Divisor $H$, der nicht numerisch trivial sein kann, durch
ihn sei eine H"alfte von $\bf Z$ als positiv ausgezeichnet.
Sei $C$ eine beliebige effektive Kurve. Auch diese kann nicht numerisch trivial
sein. W"urde die Klasse von $C$ auf der $H$ gegen"uber liegenden
Seite liegen, so h"atten sowohl positive als auch negative Vielfache von $H$
Schnitte $\neq 0$, was nicht sein kann. Also ist die numerische
Klasse von $C$ auf der Seite von $H$ und damit selbst ger"aumig.
\par\smallskip\noindent
In der komplexen Situation betrachten wir die Exponentialsequenz.
Unter dem verbindenden Homomorphismus der langen exakten Kohomologieseqenz
${\rm Pic}\, S \longrightarrow H^2(S,{\bf Z})$
wird die Schnittform erhalten. Ist $H^2(S,{\bf Z})={\bf Z}$, so muss
ein effektiver Divisor numerisch positiv und damit ger"aumig sein. \hfill $ \Box$
\par\bigskip
\noindent
Wir gehen die projektiven glatten Fl"a\-chen gem"a"s der Klassifikation durch
und versuchen, Aussagen "uber die affine Klassengruppe zu machen.
Insbesondere wollen wir diejenigen Fl"a\-chen finden,
f"ur die die affine Klassengruppe verschwindet.
Wir beschr"anken uns auf Minimalfl"achen und, wenn dies die "Uberlegung
vereinfacht, auf den komplexen Fall und orientieren
uns dabei an \cite{beau}. 
\par\bigskip
\noindent
$\bullet $
Sei die Kodairadimension $\kappa = - \infty$. 
Das sind genau die Regelfl"achen, also Fl"a\-chen birational
zu einem Produkt ${\bf P} \times C$.
Die minimalen Modelle von Regelfl"achen sind mit der Ausnahme des ${\bf P}^2$
die sogenannten geometrischen Regelfl"achen, die wiederum
den projektiven B"undeln
${\bf P}_C({\cal E})$ zu einem Vektorb"undel $\cal E$
vom Rang zwei auf der Basiskurve $C$ entsprechen, siehe \cite{beau},
Theorem III.10, oder \cite{haralg}, Prop. V.2.2.
Daraus folgt sofort, dass
unter den Fl"a\-chen mit negativer Kodairadimension nur der ${\bf P}^2$
triviale affine Klassengruppe besitzt, da die anderen direkt als eine Faserung
gegeben sind. 
\par\smallskip\noindent
Die Picardgruppe von $S$
ist gegeben durch ${\bf Z} \oplus p^*({\rm Pic}\, C)$,
wobei $\bf Z$ durch einen
Schnitt $C_0$ erzeugt wird. Die numerische Klassengruppe
ist ${\rm Num}\, S= {\bf Z} \oplus {\bf Z}$, die von $C_0$ und einer
Faser $f$ erzeugt wird. Die Schnitttheorie ist gegeben durch $f^2=0$ und
$C_0.f=1$.
Ist $D$ vom numerischen Typ $aC_0+bf$, so ist $D.f=a$ und $D.C_0=aC_0^2+b$.
Ein Divisor, der von der Basiskurve $C$ kommt, hat $a=0$ und die numerische
Klasse $(0,b)$ mit $b={\rm deg}_C\, D$.
Ist $b$ positiv, so hat zumindest ein Vielfaches davon einen effektiven
Vertreter.
\begin{satz}
Sei $X$ eine geometrische Regelfl"ache "uber der Basiskurve $C$.
Be\-z"uglich einem fixierten ger"aumigen Divisor $H$ ist ein Divisor $D$
genau dann affin-trivial, wenn $D$ und $H$ numerisch abh"angig sind.
Ist $A$ ein normaler homogener Koordinatenring einer geometrischen Regelfl"ache,
so ist ${\rm AKG}\, A={\bf Z}$.
\end{satz}
{\it Beweis}.
Sei $H$ ein fixierter ger"aumiger Divisor auf $X$, und $D$ ein beliebiger
Divisor auf $X$.
Sind die numerischen Klassen ${\bf Q}-$linear abh"angig,
so erf"ullt nach 3.8.6
$D$ bzgl. $H$ die Eigenschaft (K').
Sind dagegen die numerischen Klassen von $D$ und $H$ linear unabh"angig,
so ist eine ganzzahlige Linearkombination davon vom
numerischen Typ $(0,k)$ mit $k$ positiv. Ein solcher Divisor stammt dann von
der Basiskurve her, und ein Vielfaches davon besitzt einen effektiven
Vertreter. Dieser hat dann aber als ein Faserdivisor
nicht affines Komplement.
\par\smallskip\noindent
Wir haben also: $D$ hat genau dann die Eigenschaft (K')
bzgl. $H$, wenn $D$ und $H$ $Q-$linear abh"angig sind.
Dies wiederum ist genau dann der Fall, wenn $D$ affin-trivial bzgl. $H$ ist.
Ist n"amlich $E$ linear abh"angig zu $H$ und 
$D$ linear abh"angig zu $H$, so ist auch $E+D$ linear abh"angig
zu $H$.
\par\smallskip\noindent
Damit ist die affine Klassengruppe des affinen Kegels eine Restklassengruppe
der numerischen Klassengruppe, und der Kern besteht aus den
${\bf Q}-$Vielfachen von $H$, daher ist
${\rm AKG}\, A
=({\bf Z}^2 /H \cdot {\bf Z})/{\rm Torsion} \cong {\bf Z}\, .$  \hfill $\Box$
\par\bigskip
\noindent
Die affine Klassengruppe eines Kegels "uber einer Regelfl"ache sieht also
so aus wie im einfachsten Fall, n"amlich bei $S={\bf P} \times {\bf P}$
mit dem affinen
Kegel ${\rm Spek}\, K[U,V,X,Y]/(UX-VY)$, der uns schon h"aufig begegnet ist.
Ist die Basiskurve nicht rational, so ergibt sich eine Vielzahl von
Beispielen mit vielen nicht trivialen, aber
affin-trivialen Divisoren, ohne dass die affine Klassengruppe verschwindet.
So gesehen handelt es sich neben den in Abschnitt 3.6 besprochenen
Hyperbeln, Monoidringen und Determinantenringen einerseits
und zweidimensionalen, nicht fastfaktoriellen lokalen Ringen andererseits
um einen dritten Typ.
\par\bigskip
\noindent
$\bullet$
Sei die Kodairadimension $\kappa =0$.
Hier gibt es Enriques-Fl\"achen, K3-Fl"a\-chen, biellipitische Fl"a\-chen und
Abelsche Variet\"aten.
\par\bigskip
\noindent
Eine Enriques-Fl\"ache ist definiert durch $q=0,\, p_g=0$ und $2K=0$.
Eine solche besitzt eine elliptische Faserung, d.h. eine Abbildung auf eine
Kurve, deren generische Faser eine elliptische Kurve ist.
Siehe hierzu \cite{beau}, Exc. IX.4, oder \cite{griffiths}, IV.5,
Die affine Klassengruppe kann also nicht verschwinden.
\par\bigskip
\noindent
Sei $S$ eine K3-Fl"a\-che, also eine Fl"a\-che mit $K=0$ und $q=0$.
Jede Enriques-Fl\"ache ist der Quotient einer K3-Fl\"ache unter einer
fixpunktfreien Involution, und umgekehrt f"uhrt eine K3-Fl\"achen
mit einer solchen Involution zu einer Enriques-Fl"a\-che,
so dass solche K3-Fl"a\-chen Kurven besitzen, die disjunkt sind.
Solche K3-Fl"a\-chen k"onnen also nicht triviale affine Klassengruppe
besitzen.
\par\smallskip\noindent
Andererseits gilt nach \cite{griffiths}, IV.5,
dass auf der generischen K3-Fl\"ache jeder Divisor homolog
zu einem Vielfachen der Hyperebenenklasse ist, dass also
${\rm Num}\, S={\bf Z}$
ist und damit auch ${\rm Pic}={\bf Z}$, da $H^1(S,{\cal O}_S)=0$ ist.
Insbesondere gilt dies f"ur Hyperfl"achen vom Grad $4$ im ${\bf P}^3$.
Die Aussage, dass der Restklassenring von $K[x,y,z,u]$ nach einem
generischen homogenen Polynom von Grad vier (und gr"o"ser)
faktoriell ist, nennt man auch den Satz von Noether.
\par\smallskip\noindent
Es gibt nat"urlich auch viele spezielle Gleichungen, wo das nicht der Fall,
etwa die Brieskorn- oder Fermatsingularit"aten.
F"ur $x^4+y^4+u^4+v^4=0$ ist beispielsweise die Existenz von
rationalen Kurven darauf typisch, w"ahrend
auf einer K3-Fl"a\-che mit trivialer affiner Klassengruppe alle Kurven
das Geschlecht $\geq 2$ haben.
Ist $S$ eine generische Hyperfl"ache im ${\bf P}^3$ vom Grad vier, so ist
jede Kurve darauf durch eine Gleichung herausschneidbar, also ein
vollst"andiger Durchschnitt. Nach der Geschlechtsformel daf"ur
erh"alt man sogar $g={1 \over 2}de(d+e-4)+1= 2d^2+1$, also minimal
das Geschlecht drei.
\par\bigskip
\noindent
Eine bielliptische Fl"a\-che hat $q=1$ und daher gibt es eine nicht konstante
Abbildung auf eine elliptische Kurve (die Albanesevariet"at).
Eine solche Fl"a\-che $S$ hat die Gestalt $S=E \times F /G$, wobei
$E$ und $F$ elliptische Kurven sind und $G$ eine endliche Gruppe von
Translationen auf $E$, die auf $F$ operiert mit $F/G={\bf P}^1$, siehe
\cite{beau}, Definition VI.19.
Daraus sieht man ebenfalls, dass es Abbildungen auf Kurven gibt.
Die affine Klassengruppe kann nicht verschwinden.
\par\bigskip
\noindent
Sei schlie"slich $X$ eine abelsche Fl"a\-che und $C$ eine irreduzible Kurve
darauf. Dann ist wieder $C^2=2g(C)-2$, wobei $g(C)$ das arithmetische
Geschlecht bezeichnet, und genau dann ist der Selbstschnitt jeder Kurve
positiv, wenn f"ur jede Kurve das Geschlecht gr"o"ser als eins ist.
Da es auf einer abelschen Variet"at keine
rationalen Kurven gibt, siehe etwa \cite{kollar}, II.3.13,
haben genau die glatten elliptischen Kurven den Selbstschnitt null.
\par\smallskip\noindent
Die Existenz von glatten elliptischen Kurven auf abelschen
Fl"a\-chen kann man noch anders charakterisieren.
Die glatte elliptische Kurve ist dann eine abelsche
Untervariet"at und daher gibt es nach dem Reduzierbarkeitstheorem von
Poincar\'{e}, siehe \cite{mumab}, Kap. 19,
eine weitere abelsche Untervariet"at $Z$ mit endlichem Durchschnitt und mit
$C+Z=S$. Dann ist $S$ isogen zu einem Produkt von elliptischen Kurven.
Ist umgekehrt $S$ isogen zu einem solchen Produkt, d.h. ist
$S/G \cong C \times C'$ mit einer endlichen Untergruppe $G$, so k"onnen
nat"urlich in $S$ nicht alle Kurven affines Komplement haben.
Wir fassen zusammen.
\begin{satz}
Sei $S$ eine abelsche Fl"a\-che. Dann sind folgende Aussagen "aquivalent.
\par\smallskip\noindent
{\rm (1)} Die affine Klassengruppe von $S$ verschwindet.
\par\smallskip\noindent
{\rm (2)} Es gibt auf $S$ keine glatten elliptischen Kurven.
\par\smallskip\noindent
{\rm (3)} $S$ ist einfach, d.h. nicht isogen zu einem Produkt elliptischer
Kurven. \hfill $\Box$
\end{satz}
{\bf Bemerkung}
Bei abelschen Variet"aten gilt allgemein f"ur einen effektiven Divisor
die "Aquivalenz von $D$ ger"aumig und $D^n > 0$, wobei $n$ die
Dimension ist, siehe \cite{iitaka}, Prop. 10.5
\par\bigskip
\noindent
$\bullet$
Nach 3.9.2 kann bei
$\kappa =1$ die affine Klassengruppe nicht
verschwinden, und in der Tat liegt eine elliptische Faserung vor.
\par\bigskip
\noindent
$\bullet$
Bei $\kappa =2$ schlie"slich, den Fl"a\-chen vom allgemeinen Typ,
ist die Situation v"ollig un"ubersichtlich.
Nach der Adjunktionsformel ist $2g(C)-2=C^2+C.K$ und es muss bei
${\rm AKG}\, S=0$ f"ur jede
Kurve dann $2g(C)-2 > C.K$ sein, wobei der Ausdruck rechts dann immer
positiv ist.
Die generische Hyperfl"ache in ${\bf P}^3$ vom Grad $\geq 5$
hat $\kappa =2$ und besitzt wie die K3-Fl"a\-chen Picardgruppe vom Rang eins.

%% file: suphoka.tex
\section{Affinit"at und Superh"ohe eins}
\par
\bigskip
\noindent
Eine Hyperfl"ache mit affinem Komplement in einem affinen Schema hat die
Eigenschaft, dass unter einem Ringwechsel das Urbild wieder eine
Hyperfl"ache ist.
In vielen Beispielen wurde die Nichtaffinit"at eines Hyperfl"achenkomplementes
dadurch erwiesen, dass ein Ringwechsel angegeben wurde, unter dem
diese geometrische Bedingung verletzt wurde. Beispielsweise wurde f"ur eine
Hyperfl"ache eine Fl"ache angegeben, mit der sie einen
Punktschnitt besitzt.
In diesem Kapitel wird untersucht, inwiefern die
Affinit"at eines Hyperfl"achenkomplementes durch die
Eigenschaft, dass das Ideal in jeder Ring\-erweiterung
die Kodimension eins besitzt, gesichert ist.
Insbesondere wird gefragt, wann es f"ur die Affinit"at ein geometrisches
Kriterium gibt.
\par
\bigskip
\noindent
\subsection{Die Superh"ohe eines Ideals}
\par
\bigskip
\noindent
Sei ${\bf a} \subseteq A$ ein Ideal in einem kommutativen Ring und
$A \longrightarrow A'$ ein Ringwechsel. Dann ist das Erweiterungsideal
${\bf a}A'$ definiert, und dieses Ideal beschreibt das Urbild
der offenen Menge $D({\bf a})$ unter der Abbildung
${\rm Spek}\, A' \longrightarrow {\rm Spek}\, A$.
Dabei kann sich die H"ohe von $\bf a$ "andern, und
insbesondere gr"o"ser werden. Die folgenden Begriffe beschreiben die
maximale H"ohe, die bei Ringwechseln aus bestimmten Ringklassen
erreicht werden kann.
Die H"ohe des Einheitsideales definieren wir als null, wenn es sich um den
Nullring handelt, andernfalls als eins.
\par\bigskip
\noindent
{\bf Definition}
Sei ${\bf a} \subseteq A$ ein Ideal in einem kommutativen Ring.
Dann hei"st
\par\medskip\noindent
${\rm supht}\, {\bf a}=
{\rm sup}\, \{   {\rm ht}\, {\bf a}A':\,
A' \,\, {\rm noethersch} \}$
die {\it Superh"ohe} von $\bf a$, bzw. in Abgrenzung zu
den folgenden Begriffen die Superh"ohe von $\bf a$ bez"uglich allen
noetherschen Ringen.
\par\medskip\noindent
${\rm supht}^{\rm end}\, {\bf a}
={\rm sup}\, \{   {\rm ht}\, {\bf a}A': \,
A'\,\, {\rm von \, endlichem \, Typ} \}$ die {\it endliche Superh"ohe}
von $\bf a$.
\par\medskip\noindent
${\rm supht}^{\rm krull}\, {\bf a}
={\rm sup}\, \{   {\rm ht}\, {\bf a}A': \,
A'\,\, {\rm Krullbereich} \}$
die Superh"ohe von $\bf a$ bez"ug\-lich Krullbereichen.
\par\bigskip
\noindent
{\bf Bemerkung}
F"ur ein vom Einheitsideal verschiedenes Ideal muss man zur Bestimmung
der Superh"ohe nur Ringerweiterungen betrachten, wo das Ideal nicht
zum Einheitsideal wird. Die Festlegungen zur H"ohe des Einheitsideals
finden in 4.2 und 4.4 eine Rechtfertigung. Manchmal muss man aber bei
Superh"ohenaussagen das Einheitsideal ausschlie"sen.
\par\smallskip\noindent
Geometrisch betrachtet geht es um eine abgeschlossene Teilmenge
$Y \subseteq X$ in einem affinen Schema $X$, und man fragt sich,
welche Kodimension das Urbild $Y'$ von $Y$ unter einem affinen Morphismus
(mit weiteren Eigenschaften) $X' \longrightarrow X$ bekommen kann.
Dabei "andert sich der Begriff nicht, wenn man beliebige
Schemamorphismen $X' \longrightarrow X$ zul"asst, da die Kodimension
des Urbildes $Y'$ auf einer affinen Teilmenge $W \subseteq X'$
angenommen wird.
\par\smallskip\noindent
Im Prinzip liefert jede Ringklasse einen eigenen Superh"ohenbegriff,
die beiden ersten sind aber die wichtigsten, die Superh"ohe
bez"uglich allen Krullbereichen ist f"ur uns von Bedeutung, da mit ihr
im n"achsten Abschnitt die Affinit"at charakterisiert werden kann.
\par\smallskip\noindent
Der Begriff und die Bezeichnung Superh"ohe gehen auf Hochster zur"uck
und wurden in Zusammenhang mit der direkten Summand-Vermutung
entwickelt, \cite{hochsterbig}. Gem"a"s dieser unbewiesenen Vermutung ist ein
lokaler, regul"arer Ring $A$
in einer endlichen Erweiterung $A \subseteq B$ ein direkter Summand.
Enth"alt $A$ einen K"orper, so ist diese Vermutung
best"atigt, im Allgemeinen ist sie "aquivalent zu einer Superh"ohenaussage
zu einem bestimmten Ideal in einem Hyperfl"achenring,
siehe \cite{hochsterdim}, \cite{hochstercan} (Siehe auch 4.2. und 5.1).
In diesem Kontext wurden von Koh (\cite{koh1}, \cite{koh2})
einige Resultate erzielt, auf die wir uns
hier st"utzen k"onnen, insbesondere 4.1.2 und 4.1.5
\par\smallskip\noindent
Die Sache selbst wurde aber schon fr"uher betrachtet, insbesondere unter
der Fragestellung, wann die Kodimension des Schnittes zweier irreduzibler
Teilmengen die Summe der beiden Kodimensionen ist, und wann sie gr"o"ser
sein kann. In diese Richtung geht das affine
Dimensionstheorem, \cite{haralg}, Prop. I.7.1, und die
allgemeinere Aussage von Serre f"ur regul"are Ringe, siehe den Abschnitt
5.1.
\par\smallskip\noindent
Dagegen ist die Beziehung zu Affinit"atsfragen und zur kohomologischen
Dimension noch nicht eingehend untersucht worden, und soll daher Thema dieses
und des n"achsten Kapitels sein.
\par\bigskip
\noindent
Zu einem Ideal bezeichne ${\rm ara}\, {\bf a}$ die
minimale Anzahl von Funktionen, die das Ideal geometrisch
beschreiben, also die minimale Anzahl von Funktionen, die man braucht,
um das Radikal zu $\bf a$ als Radikal zu beschreiben.
Die folgende Proposition fasst einige einfache Eigenschaften der Superh"ohe
zusammen.
\par\smallskip\noindent
\begin{prop}
Sei $\bf a$ ein Ideal in einem noetherschen Ring $A$.
Dann gelten folgende Aussagen, wobei in {\rm (4)}, {\rm (5)} und {\rm (6)}
${\bf a } \neq A$ vorausgesetzt sei.
\par\smallskip\noindent
{\rm (1)}
Die Superh"ohe h"angt nur vom Radikal ab.
\par\smallskip\noindent
{\rm (2)}
Ist $A \longrightarrow A'$ ein Ringhomomorphismus mit einem
noetherschen Ring $A'$, so
ist ${\rm supht}\, {\bf a} \geq {\rm supht}\, {\bf a}A'$
\par\smallskip\noindent
{\rm (3)}
Es gelten die Absch"atzungen
$$ {\rm ht}\, {\bf a} \leq {\rm alt}\, {\bf a}
\leq {\rm supht}^{\rm end}\, {\bf a}
\leq {\rm supht}\, {\bf a}
\leq {\rm ara}\, {\bf a} \, .$$
Ist $\bf a$ das maximale Ideal in einem lokalen Ring, so gilt "uberall die
Gleichheit.
\par\smallskip\noindent
{\rm (4)}
Die Superh"ohe und die endliche Superh"ohe k"onnen in integren
Ringen so realisiert werden, dass das Erweiterungsideal irreduzibel ist.
\par\smallskip\noindent
{\rm (5)}
Es ist
\begin{eqnarray*}
{\rm supht}\, {\bf a} &= & {\rm sup}\, \{ {\rm dim}\, A':\, 
A' \mbox{ ist noethersch, lokal, komplett, integer,} \\
& & \ \ \ \ \ \ \ \ \ \ \ \ \ \ \ \ \ \ \ \ \ \ \ \ \ \ \ \ \ \ \ \ \ \ 
\mbox{ normal mit }\, V({\bf a}A')=V({\bf m}) \}
\end{eqnarray*}
Die noethersche Superh"ohe wird also als Dimension eines lokalen Ringes
mit den angegebenen Eigenschaften angenommen,
so dass das Urbild von $V({\bf a})$ genau der abgeschlossene Punkt ist.
\par\smallskip\noindent
{\rm (6)}
Die endliche Superh"ohe l"asst sich in einem integren Ring als
H"ohe eines maximalen Ideals realisieren.
\par\smallskip\noindent
{\rm (7)}
Die Superh"ohe "andert sich nicht beim \"Ubergang zu
${\bf a}A_{1+{\bf a}}, A_{1+{\bf a}}$.
Die Superh"ohe eines maximalen Ideals ist gleich der H"ohe.
\par\smallskip\noindent
{\rm (8)}
Die Superh\"ohe lokalisiert in folgender Weise.
Zu ${\bf a}$ und einem Primideal
${\bf p} \in X={\rm Spek}\, A$ sei
${\bf a}_{\bf p}={\bf a}A_{\bf p} \subseteq A_{\bf p}$
die lokalisierte Version. Dann ist
$$ {\rm supht}\, {\bf a}
= {\rm {\rm max}}_{x \in X}\, ({\rm supht}_{A_x}\, {\bf a}_x) \, .$$
\par\smallskip\noindent
{\rm (9)}
Es ist ${\rm supht}\, {\bf a} \leq {\rm dim}\, A$.
\end{prop}
{\it Beweis}.
Zu (1).
Zwei radikalgleiche Ideale sind in jeder Ringerweiterung
radikalgleich und haben daher die gleichen minimalen Primoberideale.
\par\smallskip\noindent
(2) ist klar.
\par\smallskip\noindent
Zu (3).
Die erste Absch"atzung ist nach Definition klar,
die zweite wird unter (4) mitbewiesen.
Die dritte ist klar, die letzte folgt aus dem
allgemeinen Krullschen Hauptidealsatz,
siehe \cite{eisenbud}, Theorem 10.2.\footnote{Im allgemeinen
ist ${\rm ara}\, {\bf a}$
nach dem Satz von Forster $\leq d+1$,
wobei $d$ die Ringdimension ist, siehe \cite{lyub}, Theorem 1.1.
F"ur die Superh"ohe gilt die bessere Absch"atzung (9).}
Der Zusatz ergibt sich, da in einem lokalen Ring f"ur das maximale
Ideal ${\rm ara}\, {\bf m}={\rm ht}\, {\bf a}$ ist.
\par\smallskip\noindent
Zu (4).
Sei ${\bf a}A' \subseteq {\bf p}$ ein minimales Primoberideal in $A'$,
und ${\bf p}_1,...,{\bf p}_r$ die anderen.
Nach dem Satz "uber die Vermeidung von Primidealen gibt es dann ein
$f \not\in {\bf p}$ und $f \in {\bf p}_i$ f"ur $ i=1,...,n$.
Nach der Nenneraufnahme $A' \longrightarrow A'_f$ ist dann ${\bf p}$
das einzige Primoberideal zu ${\bf a}A'_f$, dessen H"ohe nun die H"ohe
von $\bf p$ ist.
Wird die H"ohe von $\bf p$ "uber dem Primideal $\bf q$ der H"ohe null
angenommen, so hat man modulo $\bf q$ den gew"unschten integren Ring.
Beide Prozesse sind von endlichem Typ.
\par\smallskip\noindent
Zu (5). Im Fall der noetherschen Superh"ohe kann man durch
Lokalisierung an $\bf p$ erreichen, dass ${\bf a}A'_{\bf p}$ das maximale
Ideal beschreibt.
Bei der Komplettierung erh"alt sich die Dimension, durch "Ubergang
zur gr"o"sten Komponente erhalten wir einen integren kompletten Ring $R$.
Ein solcher Ring ist exzellent, insbesondere sind die formellen Fasern
geometrisch normal, und damit ist nach \cite{EGAIV}, 7.6.2,
die Normalisierung davon lokal (die Punkte entsprechen den Komponenten
der Komplettierung) und wegen
$(R^{\rm nor})^{\rm kom}=(R^{\rm kom})^{\rm nor}=R^{\rm nor}$ auch komplett.
\par\smallskip\noindent
Zu (6).
Sei ohne Einschr"ankung $A'$ integer und $\bf p$ das einzige
Primoberideal zu ${\bf a}A'$.
Ist $\bf p$ maximal, so ist man fertig, sei also ${\bf q}$ ein direktes
Primoberideal zu $\bf p$.
Sei $x \not\in {\bf p},\, x \in {\bf q}$.
Dann ist ${\bf q}$ modulo $x$ ein minimales Primoberideal zu ${\bf p}A'/x$
(und zu  ${\bf a}A'/x$) und es gelten in $A'/x$ die Absch"atzungen
$$ {\rm supht}^{\rm end}\, {\bf a}
\geq {\rm ht}\, {\bf q}
\geq {\rm dim}\, (A'/x)_{\bf q}
={\rm dim}\, (A'_{\bf q}/x)
\geq {\rm dim}\, (A'_{\bf q}) -1
\geq {\rm ht}\, {\bf p} \, .$$
Dabei gilt die vorletzte Gleichung aufgrund von \cite{eisenbud}, Cor. 10.9.
Da in $\bf p$ schon die endliche Superh"ohe von $\bf a$ angenommen ist,
gilt die Gleichheit, und sukzessive fortfahrend gelangt man zu einem
(beliebig vorgebbaren) maximalen Ideal. 
\par\smallskip\noindent
Zu den letzten drei Aussagen sei $A \longrightarrow R$ ein Ringhomomorphismus
in einen lokalen, noetherschen, integren, kompletten Ring,
wo $\bf a$ (als Radikal) zum maximalen Ideal wird
und wo die Superh"ohe angenommen wird.
\par\smallskip\noindent
Zu (7).
Unter $A \longrightarrow R$ werden alle
Elemente aus dem multiplikativen System $1+{\bf a}$ zu Einheiten, so dass
die Abbildung durch $A_{1+{\bf a}}$ faktorisiert.
F"ur ein maximales Ideal ist $A_{\bf m}=A_{1+{\bf m}}$.
\par\smallskip\noindent
Zu (8).
Die eine Absch\"atzung ist klar, man hat zu zeigen,
dass die Superh\"ohe an einem Punkt angenommen wird.
Sei $A \longrightarrow R$ wie gehabt.
Das maximale Ideal von $R$ bildet auf einen
Punkt $x \in X$ ab, und es liegt eine Faktorisierung
$A \longrightarrow A_x \longrightarrow R$ vor, so dass
${\bf a}_x$ ebenfalls diese Superh"ohe besitzt.
\par\smallskip\noindent
Zu (9).
Ohne Einschr"ankung sei wegen (8) $A$ lokal, $A \longrightarrow R$ wie gehabt.
Die Faser "uber dem maximalen Ideal von $A$ ist nulldimensional
(oder sogar leer) und nach der Dimensionsformel,
\cite{eisenbud}, Theorem 10.10, gilt dann
${\rm supht}\, {\bf a}={\rm dim}\, R \leq {\rm dim}\, A$. \hfill $\Box$
\par\bigskip
\noindent
{\bf Bemerkung} Im Beweis zu (6) wurde allgemein gezeigt, dass es zu einem
Primideal $\bf p$ der H"ohe $n$
in einem noetherschen Ring einen integren Restklassenring
gibt, wo $\bf p$ zu einem maximalen Ideal der H"ohe $n$ wird.
\par\bigskip
\noindent
Eine wichtige Klasse von Ringerwechseln sind die Restklassenabbildungen
$A \longrightarrow A/{\bf b}$, die den abgeschlossenen Einbettungen
$V({\bf b}) \hookrightarrow {\rm Spek}\, A$ entsprechen, und wo es auf die
Kodimension von
$V({\bf b}+{\bf a})=V({\bf b}) \cap V({\bf a})$ in $V({\bf b})$ ankommt.
Diese sind von endlichem Typ, und h"aufig kann man schon mit ihnen
zeigen, dass eine Hyperfl"ache eine Superh"ohe $\geq 2$ hat und somit
nicht affin sein kann, wie etwa bei
$V(x,y) \subseteq {\rm Spek}\, K[u,v,x,y]/(ux-vy)$
mit der Reduktion modulo ${\bf b}=(u,v)$.
Betrachtet man neben abgeschlossenen Unterschemata auch noch deren
Normalisierungen, so erh"alt man sogar schon die endliche Superh"ohe,
wie folgender Satz von Koh zeigt, \cite{koh1}.
\begin{satz}
Sei $A$ ein noetherscher Ring und $\bf a$ ein Ideal.
Dann ist
\begin{eqnarray*}
& & {\rm supht}^{\rm end}\, {\bf a} \\
& & ={\rm sup}\, \{ {\rm alt}\,  {\bf a}A':\,
A' \, \mbox{ist Normalisierung eines integren Restklassenringes} \}.
\end{eqnarray*}
\end{satz}
{\it Beweis}. 
Sei $A \longrightarrow R$ von endlichem Typ derart,
dass $R$ integer und das Erweiterungsideal ${\bf a}R$ geometrisch
ein maximales Ideal $\bf m$ ist, dessen H"ohe die
Superh"ohe von $\bf a$ ist.
Sei $B$ der ganze Abschluss von $A$ in $R$.
F"ur $B$ und $R$ gilt dann der Hauptsatz von Zariski in der Formulierung
des folgenden Lemmas, da $\bf m$ einziger Urbildpunkt von $V({\bf a}B)$ ist
und insbesondere isoliert "uber seinem Bildpunkt liegt.
Somit ist ${\bf m}'={\bf m} \cap B$ minimal "uber ${\bf a}B$
mit gleicher H"ohe wie ${\bf m}$.
\par\smallskip\noindent
Sei ${\bf p}$ der Kern der Abbildung nach $B$ (und $R$) und betrachte die
Faktorisierung
$$ A/{\bf p} \longrightarrow (A/{\bf p})^{\rm nor}
\longrightarrow B \, .$$
Die hintere Abbildung ist injektiv und ganz; $B$ ist der ganze
Abschluss eines normalen Ringes in einer integren Erweiterung, so dass das
going down Theorem gilt, siehe \cite{nagloc}, Theorem 10.13.
Aus dem going down folgt, dass
${\bf m}' \cap (A/{\bf p})^{\rm nor}$ auch minimal "uber
${\bf a}(A/{\bf p})^{\rm nor}$ ist, und zwar mit unver"andeter H"ohe.
Daher ist die Altitude dieses Ideals zumindest so gro"s wie die
endliche Superh"ohe.
\par\smallskip\noindent
F"ur die umgekehrte Absch"atzung gen"ugt es zu zeigen, dass die Altitude
in einer ganzen Erweiterung bereits als Altitude in einer endlichen Erweiterung
angenommen wird.
Sei also $A \longrightarrow B$ eine ganze Erweiterung
und ${\bf a}B \subseteq {\bf q}$ ein minimales Primoberideal.
Durch eine Nenneraufnahme in $B$ und den "Ubergang
zu $A[f,1/f] \longrightarrow B_f$ kann man annehmen,
dass $\bf q$ das einzige Prim\-oberideal ist.
Sei $A \cap {\bf q}={\bf p}=(f_1,...,f_n)$.
Es ist dann $f_1^k,...,f_n^k \in {\bf a}B$. Daher gibt es eine
endlich erzeugte Erweiterung $A'$ von $A$ in $B$, wo ebenfalls
$f_1^k,...,f_n^k \in {\bf a}A'$ gilt. Sei ${\bf q}'={\bf q} \cap A'$
und ${\bf a}A' \subseteq {\bf r} \subseteq {\bf q}'$.
Dann ist
$$(f_1^k,...,f_n^k) \, \subseteq \, {\bf a}A' \cap A
\, \subseteq \, {\bf r} \cap A \, \subseteq \, {\bf p} \, . $$
Diese Ideale sind alle radikalgleich, und da die Fasern der ganzen
Erweiterung nulldimensional sind, folgt ${\bf r}={\bf q}'$
und damit ist ${\bf q}'$ minimal "uber
${\bf a}A'$ mit der gleichen H"ohe wie $\bf q$.  \hfill $\Box$
\par\bigskip
\noindent
Im Beweis wurde der Hauptsatz von Zariski in folgender Version
verwendet, siehe \cite{koh1}.
F"ur den soeben gezeigten Satz braucht man nur die Situation, wo $R$ endlich
erzeugt "uber $B$ ist, f"ur das Beispiel weiter unten die volle
Aussage.
\begin{lem}
Sei $B \subseteq R$ mit $B$ ganz-abgeschlossen in $R$ und $R$ ganz "uber
einer endlich erzeugten $B-$Algebra.
Ist dann ${\bf q} \in {\rm Spek}\, R$ isoliert in der Faser "uber
${\bf p}={\bf q} \cap B$,
so gibt es eine offene Umgebung von $\bf p$,
auf der die Abbildung eine Isomorphie ist.  \hfill $\Box$
\end{lem}
\begin{kor}
Ist $A$ ein noetherscher normaler Integrit"atsbereich der Dimension $d$,
so besitzen alle Ideale mit ${\rm alt}\, {\bf a} < d$ eine endliche
Superh"ohe $< d$.
Insbesondere ist ${\rm supht}^{\rm end}\, {\bf p}=d-1$
f"ur Primideale $\bf p$ der H"ohe $d-1$ {\rm (}also Kurven{\rm )}.
\end{kor}
{\it Beweis}. Nach dem Satz muss man nur schauen, wie sich die H"ohe auf
Normalisierungen von abgeschlossenen Unterschemata verh"alt.
Da $A$ schon normal ist, gilt nach Voraussetzung dort die Absch"atzung, man
muss dann also nur noch Normalisierungen von echten Teilmengen betrachten,
und die haben alle eine Dimension $<d$, wodurch dann auch die Altitude
beschr"ankt ist. \hfill $\Box$
\par\bigskip
\noindent
{\bf Bemerkung}
Der Zusatz besagt im normalen Fl"achenfall, dass Hyperfl"achen (=Kurven)
die endliche Superh"ohe eins besitzen und ist damit eine
schw"achere Aussage als der Satz von Nagata, dass alle
Komplemente von Kurven auf normalen exzellenten Fl"achen affin sind, wobei
allerdings die Voraussetzung exzellent nicht auftaucht.
Eine verwandte kohomologische Aussage ist,
dass f"ur Kurven $C$ in einem kompletten integren Ring die
lokale Kohomologie $H^d_C(X,-)=H^{d-1}(X-C,-)$ verschwindet,
also die kohomologische Dimension von $X-C$ gleich $d-2$ ist, die
kohomologische H"ohe also gleich $d-1$, siehe \cite{peskine},
theorem 3.1 und 5.2.6.
\begin{satz}
Sei $A$ eine endlich erzeugte $K-$Algebra "uber einem K"orper $K$.
Dann ist ${\rm supht}^{\rm end}\, {\bf a}={\rm supht}\, {\bf a}$.
\end{satz}
{\it Beweis}.
Siehe \cite{koh2}, Theorem 1. Dort wird dieser nicht triviale Satz in der Form
${\rm supht}\, ({\bf a}A_{\bf p})={\rm supht}^{\rm end}\, ({\bf a}A_{\bf p})$
f"ur ein beliebiges Primideal ${\bf p} \supseteq {\bf a}$
bewiesen, woraus die Aussage folgt, da es ein geeignetes $\bf p$ gibt
(das Bild des maximalen Ideals in einer Erweiterung) mit
${\rm supht}\, {\bf a}={\rm supht}\, ({\bf a}A_{\bf p})=
{\rm supht}^{\rm end}\, ({\bf a}A_{\bf p})
\leq {\rm supht}^{\rm end}\, ({\bf a}A_{\bf p})$.
Die letzte Absch"atzung gilt dabei, da endlich erzeugte Erweiterungen einer
Lokalisierung von endlich erzeugten Erweiterungen
des Ausgangsringes herkommen. \hfill $\Box$
\par\bigskip
\noindent
{\bf Bemerkung} Der Satz gilt auch, wenn man statt einem K"orper $K$ den
Ring der Witt Vektoren $V_p$ nimmt, siehe \cite{koh2}, in dieser Form ist
die Aussage f"ur die direkte Summand Vermutung von Interesse.
\par\bigskip
\noindent
{\bf Beispiel} (Vergl. \cite{nagloc} und \cite{koh2})
Sei $B=R[X]/(X^2-z)$ der in \cite{nagloc}, A 1, Example 7,
vorgestelle zweidimensionale, lokale,
normale, noetherscher Ring, dessen Komplettierung analytisch reduzibel ist.
Dabei ist $R \subseteq K[[x,y]]$ ein regul"arer Ring, dessen maximales Ideal
durch $(x,y)$ gegeben ist, $w$ eine Potenzreihe aus $xK[[x]]$,
die nicht zu $R$ geh"ort, und $z=(y+w)^2 \in R$.
Die Komplettierung von $B$ ist $K[[x,y]][X]/( (X-(y+w))(X-(y+w)) )$.
Durch $X \longmapsto y+w$ betrachten wir $B \subseteq K[[x,y]]$
in der einen Komponente der Komplettierung.
\par\smallskip\noindent
Es ist $t=(X-y) /x= w/x \in K[[x,y]] \cap Q(B)$, aber $t \not\in B$.
${\rm Spek}\, B[t]$ ist der Abschluss des Graphen der durch $t$
definierten rationalen Funktion, $B[t]$ ist ein
integrer, noetherscher, zweidimensionaler Ring.
Da $t$ im maximalen Ideal von $B$ nicht definiert ist,
liegt "uber dem maximalen Ideal eine affine Gerade.
Das Ideal $(x,y)$ hat also in $B[t]$ die H"ohe eins.
\par\smallskip\noindent
Sei nun $A$ die Normalisierung von $B[t]$.
Dabei ist $A \subseteq K[[x,y]]$, da $ K[[x,y]]$ ganz abgeschlossen ist.
Das Ideal ${\bf a}=(x,y)A$ hat auch in $A$ die reine H"ohe eins.
W"are n"amlich in $V({\bf a})$ der abgeschlossene Punkt $Q$
eine Komponente, so w"are $Q$ unter der Abbildung
${\rm Spek}\, A \longrightarrow {\rm Spek}\, B$ ein isolierter
Punkt "uber dem maximalen Ideal von $B$.
Nach dem Hauptsatz von Zariski w"are dann aber $A=B$, was nicht der
Fall ist.
\par\smallskip\noindent
$V(\bf a)$ ist also eine Hyperfl"ache im zweidimensionalen, normalen,
noetherschen Ring $A$, und damit ist nach dem Korollar 4.1.4
die endliche Superh"ohe gleich eins.
Dagegen wird $\bf a$ unter $A \longrightarrow K[[x,y]]$ zum
maximalen Ideal, so dass die Superh"ohe zwei ist. 
Wir haben also ein Beispiel, wo Superh"ohe und endliche
Superh"ohe auseinander fallen.
Da die Superh"ohe von $\bf a$ zwei ist folgt,
dass $D({\bf a})$ nicht affin ist, und wir haben ein Beispiel, das zeigt,
dass in einem normalen, zweidimensionalen, integren Ring 
nicht jedes Kurvenkomplement affin sein muss.
Im Satz von Nagata kann man also nicht auf die Voraussetzung
exzellent (bzw. analytisch irreduzibel) verzichten.
\par
\bigskip
\noindent
\subsection{Das Erweiterungsideal im Schnittring}
\par
\bigskip
\noindent
Sei ${\bf a} \subseteq A$ ein Ideal in einem kommutativen Ring,
$U=D({\bf a}) \subseteq {\rm Spek}\, A=X$ und
$B=\Gamma(U,{\cal O}_X)$ der globale Schnittring.
In dieser Situation hat man eine offene Einbettung
$U=D({\bf a}B) \hookrightarrow {\rm Spek}\, B$, wobei das Erweiterungsideal
${\bf a}B$ jetzt eine
offene Menge in einem affinen Schema beschreibt, auf der jede Funktion
global ausdehnbar ist. $U$ ist genau dann affin, wenn ${\bf a}B$
das Einheitsideal ist. Ist dies nicht der Fall, so ist unter recht
schwachen Voraussetzungen ${\bf a}B$ ein Ideal der H"ohe $\geq 2$.
\begin{satz}
Sei $A$ ein noetherscher Ring und $\bf a$ ein Ideal, $U=D({\bf a})$.
Dann ist $U$ genau dann affin,
wenn ${\rm supht}^{\rm krull}\, {\bf a} \leq 1$
ist, wenn also unter allen Abbildungen
$A \longrightarrow A'$ mit einem Krullbereich $A'$ das Erweiterungsideal
${\bf a}A'$ die H"ohe $\leq 1$ hat.
\end{satz}
{\it Beweis}.
Ist $U$ affin und $A \longrightarrow A'$ ein solcher Ringwechsel, so ist
auch das Urbild $U'=D({\bf a}A')$ affin und daher muss $V({\bf a}A')$ die
Kodimensionsbedingung 1.11.3 erf"ullen, also eine Altitude $\leq 1$
besitzen.
Sei umgekehrt die angegebene Kodimensionsbedingung erf"ullt. Diese gilt
dann auch f"ur jede Komponente und jede Normalisierung, so dass man
$A$ als einen Krullbereich annehmen darf aufgrund von 1.13.2 und
1.14.2.
Der globale Schnittring $\Gamma(U,{\cal O}_X)$
ist nach 1.11.2 ebenfalls ein Krullbereich,
und das Erweiterungsideal
von $\bf a$ hat nach Voraussetzung eine Altitude $ \leq 1$.
Dies kann aber nur sein, wenn das Erweiterungsideal das Einheitsideal
ist, also ist $U$ affin.  \hfill $\Box$
\par\bigskip
\noindent
{\bf Bemerkung}
Sucht man zu nicht affinem $U$ nach Ringen,
wo die Kodimensionsbedingung verletzt ist, so
ist generell der Schnittring von $U$ ein wichtiger Kandidat, wobei
nat"urlich darauf zu achten ist, zu welcher Ringklasse der Schnittring
geh"ort. Wenn man nicht wei"s, ob $U$ affin ist, so kann man die Frage
so angehen, dass man zu $A$ sukzessive Funktionen aus dem Schnittring
hinzunimmt und schaut, ob das Erweiterungsideal
die Kodimensionsbedingung verletzt oder aber zum Einheitsideal geworden ist.
Ist beides nicht der Fall, muss man weiter machen, allerdings ohne
Garantie, dass man dabei zu einem Ende kommt, insbesondere wenn der
Schnittring nicht vom endlichen Typ ist.
\par\bigskip
\noindent
{\bf Beispiel}
Sei $K$ ein integrer Ring, wir betrachten in
$$A=K[X_1,X_2,Y_1,Y_2]/(X_1^kX_2^k+Y_1X_1^{k+1}+Y_2X_2^{k+1})$$
das Ideal ${\bf a}=(X_1,X_2)$ mit $U=D({\bf a})$. $A$ ist ebenfalls integer.
Auf $U$ ist $Z=-Y_2/ X^k_1=(X_2^k+Y_1X_1)/ X^{k+1}_2$
eine definierte Funktion.
Damit ist $Y_2=-ZX_1^k$ und
$X_1^kX_2^k+Y_1X_1^{k+1}-ZX_1^kX_2^{k+1}=X^k_1(X_2^k+Y_1X_1-ZX_2^{k+1})$.
Adjunktion von $Z$ transformiert damit den Ring zu
$$ K[X_1,X_2,Y_1,Z]/(X_2^k+Y_1X_1-ZX_2^{k+1})\, .$$
Auf $D(X_1,X_2)$ ist jetzt die Funktion
$Y_1 /X_2^k = (ZX_2-1)/ X_1$ definiert, aus der man abliest, dass
eine Darstellung der Eins m"oglich ist.
\par\smallskip\noindent
F"ur $K={\bf Z}$ ist dies der zweidimensionale Spezialfall der
Superh"ohenversion der direkten Summand Vermutung, siehe 5.1 f"ur die
allgemeine Formulierung oder \cite{hochstercan}, Abschnitt 6.
Man kann also den zweidimensionalen Fall der Vermutung leicht
(einfacher als in \cite{hochstercan}) mittels
dem Nachweis der Affinit"at best"atigen.
\par\bigskip
\noindent
{\bf Beispiel}
Wir betrachten jetzt
$$A=K[X_1,X_2,Y_1,Y_2]/(X_1^{k+1}X_2^k+Y_1X_1^{k+1}+Y_2X_2^{k+1})\, .$$
Dann gibt es wieder eine Funktion
$Z=-Y_2/ X^{k+1}_1=(X_2^k+Y_1)/ X^{k+1}_2$.
Adjungiert man $Z$ zu $A$ im globalen Schnittring,
so gelangt man zu
$$ K[X_1,X_2,Y_1,Z]/(X_2^k+Y_1-ZX_2X_2^{k+1})=K[X_1,X_2,Z] \, ,$$
wo die Superh"ohenbedingung verletzt ist.
\begin{satz}
Sei $A$ ein noetherscher integrer Ring und ${\bf a}$ ein Ideal, $U=D({\bf a})$.
Genau dann ist $U$ affin,
wenn der globale Schnittring $\Gamma(U,{\cal O}_X)$ von endlichem Typ
"uber $A$ ist und wenn ${\rm supht}^{\rm end}\, {\bf a} \leq 1$ ist. 
\end{satz}
{\it Beweis}.
Wenn $U$ affin ist, so ist nach 1.4.3
$B=\Gamma(U,{\cal O}_X)$ eine endlich erzeugte
$A-$Alge\-bra.
Sei umgekehrt $U$ nicht affin mit endlich erzeugtem
globalen Schnitt\-ring $B$.
Dann ist $B$ insbesondere noethersch und integer,
das Erweiterungsideal ist nicht das Einheitsideal und besitzt
nach 1.4.4 die Altitude $\geq 2$. \hfill $\Box$
\par\bigskip
\noindent
Wir wollen die Eigenschaft ${\rm supht}\, {\bf a}=1$ im Fall von
integren Algebren von endlichem Typ "uber einem K"orper $K$
genauer untersuchen. Es ergibt sich
eine direkte Beziehung zu der Eigenschaft, affin in der Dimension
zwei zu sein.
Als Korollar ergibt sich, dass man unter der Voraussetzung, dass der
globale Schnittring von endlichem Typ ist, die Affinit"at auf Fl"achen
testen kann. 
\begin{satz}
Sei $A$ eine integre, endlich erzeugte $K-$Algebra,
$D({\bf a})=U \subseteq X={\rm Spek}\,A$
eine offene Teilmenge.
Dann sind folgende Aussagen "aquivalent.
\par\smallskip\noindent
{\rm (1)} ${\rm supht}\, {\bf a} \leq 1$
\par\smallskip\noindent
{\rm (2)} Jede abgeschlossene Untervariet"at der Dimension $\leq 2$ von $U$
ist affin.
\end{satz}
{\it Beweis}.
Sei (1) erf"ullt, f"ur Kurven und Punkte ist die Aussage automatisch erf"ullt,
sei also $S \hookrightarrow U$ eine abgeschlossene integre Fl"ache in $U$
und $S'$ der Abschluss von $S$ in ${\rm Spek}\, A$.
Es ist dann $S' \hookrightarrow X$ ebenfalls eine Fl"ache, da sich im
Fall von Variet"aten die Dimension auf offenen Mengen testen l"asst.
Sei $\tilde{S}$ die Normalisierung von $S'$.
Das Urbild von $Y=V({\bf a})$ unter $\tilde{S} \longrightarrow X$ besitzt nach
Voraussetzung die reine Kodimension eins und hat daher affines Komplement,
was sich auf $S' \cap Y$ "ubertr"agt, also ist $S=S' \cap U =S'-S' \cap Y$
affin.
\par\smallskip\noindent
Sei umgekehrt ${\rm supht}\, {\bf a} \geq 2$ und $A \longrightarrow R$ ein
Ringhomomorphismus von endlichem Typ mit $R$ integer und so,
dass ${\bf a}$ geometrisch zu einem maximalen Ideal $\bf m$ wird, dessen
H"ohe die Superh"ohe von $\bf a$ ist.
Wir k"onnen $R$ durch einen zweidimensionalen integren Restklassenring
ersetzen, worauf $\bf m$ die H"ohe zwei hat.
Eine solche Fl"ache f"uhrt gem"a"s dem Beweis von Satz 4.1.2
zu einer Fl"ache $S' \hookrightarrow X$, in deren Normalisierung das
Ideal die Altitude zwei bekommt.
Damit kann $S' \cap D({\bf a})$ nicht affin sein. \hfill $\Box$
\begin{kor}
Ist $U$ eine quasiaffine Variet"at "uber $K$ mit einer
endlich erzeugten $K-$Al\-ge\-bra als globalem
Schnittring, und sind alle irreduziblen abgeschlossenen Fl"achen
darin affin, so ist $U$ selbst affin.
\end{kor}
{\it Beweis}.
Man kann nach 1.3.6
$U$ als offene Teilmenge in einer affinen Variet"at realisieren,
die Aussage folgt dann direkt aus 4.2.2 und 4.2.3. \hfill $\Box$
\par\bigskip
\noindent
{\bf Bemerkung}
Die Affinit\"at einer offenen Teilmenge kann man (auch unter der
Voraussetzung, dass der globale Schnittring von endlichem Typ ist) nicht
mittels restriktiveren Klassen von Fl\"achen testen, wie die beiden
folgenden Beispiele zeigen.
\par\bigskip
\noindent
{\bf Beispiel}
Wir betrachten ${\rm Spek}\, K[X,Y,U,V]/(UX-VY)$ mit der durch
die Koordinatenfunktion $X$ definierten Abbildung  auf die affine Gerade.
Die Faser \"uber $X=0$ ist gleich ${\rm Spek}\, K[Y,U,V]/(VY)$ und \"uber
einem Punkt $X \neq 0$ gleich ${\rm Spek}\,K[Y,V]$.
Diese Fl"achen haben alle mit $D(X,Y)$ affinen Schnitt, obwohl $D(X,Y)$
selbst nicht affin ist. Es liegt also eine Faserung mit affinen Fl\"achen
als Fasern vor, ohne dass die Menge selbst affin ist.
\par\bigskip
\noindent
{\bf Beispiel}
Sei $S$ die an einem Punkt aufgeblasene projektive Ebene,
$E$ der exzeptionelle Divisor und $C$ eine projektive Gerade,
die nicht durch den Punkt l"auft.
$V=S-(E \cup C)$ ist eine punktierte affine Ebene,
in ihr liegen also keine projektiven Kurven. Sei $A$ ein homogener
Koordinatenring zu $S$, $V=D_+({\bf a})$,
$U=D({\bf a}) \subseteq X={\rm Spek}\, A$.
Sei $V({\bf p})$ eine homogene irreduzible Fl"ache im affinen Kegel.
Die zugeh"orige projektive Kurve schneidet $V_+({\bf a})$ und daher ist
$V({\bf p}) \cap U $ als Urbild von $V_+({\bf p}) \cap V$ affin.
In $U$ ist also jede homogene Fl"ache affin.
Dagegen k"onnen in $U$ nicht alle Fl"achen affin sein.
$U$ ist als Kegel "uber der punktierten Ebene isomorph zu
${\bf A}^\times \times ({\bf A}^2)^{\times}$ und hat in
${\bf A}^{\times} \times {\bf A}^2$ die H"ohe zwei, so dass es eine
Vielzahl nicht affiner Fl"achen darin gibt.
\par\smallskip\noindent
In diesem Beispiel ist $(C+E)^2=0$. Wir werden im Abschnitt 4.5 sehen,
dass man bei $D$ irreduzibel und $D^2=0$ aus dem positiven Schnitt mit allen
Kurven auf Superh"ohe eins schlie"sen kann. 
\par
\bigskip
\noindent
\subsection{Affine und Steinsche offene Teilmengen}
\par
\bigskip
\noindent
Sei $X$ ein affines Schema von endlichem Typ "uber den komplexen Zahlen $\bf C$,
$X$ ist also eine Zariski-abgeschlossene Teilmenge des ${\bf A}^n_{\bf C}$.
Der zugeh"orige komplexe Raum $X^{\rm an}$ ist somit eine abgeschlossene
analytische Teilmenge des ${\bf C}^n$ und damit Steinsch,
siehe \cite{grauert}, V. \S1 Satz 1.
Die Umkehrung gilt dabei nicht, wie das Beispiel weiter unten zeigt.
Aus der Steinschheit einer offenen Menge $U=D({\bf a}) \subseteq X={\rm Spek}\, A$
kann man aber immerhin noch Aussagen "uber die Superh"ohe gewinnen.
\begin{satz}
Sei $A$ eine ${\bf C}-$Algebra vom endlichen Typ und
$U=D({\bf a}) \subseteq {\rm Spek}\, A=X$ eine offene Teilmenge mit
$U^{\rm an}$ Steinsch.
Dann ist ${\rm supht}\, {\bf a} \leq 1$.
\end{satz}
Nach Satz 4.1.5 gen"ugt es, die Aussage f"ur die endliche Superh"ohe
zu zeigen.
Sei also $f:T \longrightarrow X$ ein affiner Morphismus zwischen affinen
${\bf C}-$Variet"aten, mit $f^{-1}(U)=T-\{ P \}$.
Die Abbildung faktorisiert
$T \hookrightarrow X \times {\bf A}^n_{\bf C} \longrightarrow X$
mit einer abgeschlossenen Einbettung und der Projektion,
und dies gilt genauso f"ur die analytische Situation
$T^{\rm an} \hookrightarrow
X^{\rm an} \times {\bf C}^n \longrightarrow X^{\rm an}$
Das Urbild einer Steinschen Menge unter dieser Abbildung ist damit Steinsch,
da sich die Steinschheit auf Produkte und auf abgeschlossene Teilmengen
"ubertr"agt, siehe \cite{grauert}, V. \S1 Satz 1.
Das Komplement einer offenen Steinschen Teilmenge
in einem Steinschen Raum besitzt aber die Kodimension $\leq 1$,
siehe \cite{grauert}, V. \S3, Satz 4,
daher muss ${\rm dim}\, T=1$ sein. \hfill $\Box$
\begin{kor}
Sei $U \subseteq {\rm Spek}\, A $, wobei $A$
eine integre ${\bf C}-$Algebra von endlichem Typ sei.
Der globale Schnittring $\Gamma(U,{\cal O}_X)$ sei ebenfalls
von endlichem Typ "uber $A$. Dann ist $U$ genau dann affin, wenn es
Steinsch ist.
\end{kor}
{\it Beweis}.
Aus der Steinschheit folgt, dass die endliche Superh"ohe eins ist, und aus
der Voraussetzung "uber den Schnittring folgt aus 4.2.2 die
Affinit"at. \hfill $\Box$
\begin{kor}
Sei $U \subseteq {\rm Spek}\, A $, wobei $A$
eine ${\bf C}-$Algebra vom endlichen Typ der Dimension $\leq 2$ sei.
Dann ist $U$ genau dann affin, wenn es Steinsch ist.
\end{kor}
{\it Beweis}.
Im null- und eindimensionalen folgt die Aussage einfach daraus,
dass jede offene Teilmenge affin ist.
Es liege also die zweidimensionale Situation vor.
Die eine Richtung ist klar, sei $U$ Steinsch.
Das "ubertr"agt sich auf die Komponenten und somit k"onnen wir $A$ als integer
annehmen. Sei $\tilde{A}$ die
Normalisierung von $A$ und $\tilde{U}$ die Normalisierung von
$U \subseteq {\rm Spek}\, A$. Mit $U$ ist auch $\tilde{U}$
Steinsch, da die Normalisierung eine endliche Abbildung ist,
\cite{grauert}, V. \S 1, Satz 1, und nach dem Satz von Chevalley gen"ugt es,
$\tilde{U} $ als affin nachzuweisen.
Nach dem Satz hat aber
das Komplement von $\tilde{U}$ in ${\rm Spek}\, \tilde{A}$
die Kodimension eins, also ist nach dem Satz
von Nagata 1.15.2 diese offene Menge affin. \hfill $\Box$
\par\bigskip
\noindent
{\bf Bemerkung}
Die Aussage, dass unter der Bedingung, dass der Schnittring vom
endlichen Typ ist, Steinschheit und Affinit\"at \"aquivalent sind, wird
auch in \cite{bingener}, 5.1, gezeigt.
Dort wird gezeigt, dass unter der Endlichkeitsvoraussetzung
der komplexe Raum zur affinen H\"ulle gleich der Steinschen H\"ulle ist.
Die "Aquivalenzaussage kann man dann mit dem kommutativen Diagramm,
das die algebraische und die analytische Situation verbindet, beweisen.
\par\bigskip
\noindent
Sei $X$ ein komplexer Raum und $Y \subseteq X$ eine abgeschlossene,
analytische Teilmenge.
Die analytische Superh"ohe von $Y$ definieren wir in naheliegender Weise
durch
$$ {\rm supht}^{\rm an}\, (Y,X)  =
{\rm sup}\, \{ {\rm codim}_{x'}\,(f^{-1}(Y),X'):
\, x' \in X',\,   f:X' \longrightarrow X \}\,.$$
{\bf Erl"auterung}
Dabei ist ${\rm codim}_x\, (Y,X)={\rm dim}_x \, X-{\rm dim}_x\, Y$,
und es gilt hierbei
${\rm dim}_x\, X={\rm dim}\, {\cal O}_x={\rm dim}\, \hat{{\cal O}}_x$, siehe
\cite{grauertas}, Kap. II, \S4ff.
Ist $Y_x$ im Punkt $x \in X$ Nullstellenmenge des Ideal $\bf a$,
so ist
${\rm codim}_x\, (Y,X)=
{\rm dim}\, ({\cal O}_{X,x}) -{\rm dim}\, ({\cal O}_{X,x}/{\bf a})$.
Ist $X$ irreduzibel, so ist dies gleich der H"ohe des Ideals $\bf a$,
da die Halmringe katen"ar sind. 
\par\bigskip
\noindent
{\bf Bemerkung}
Der Satz 4.3.1 gilt auch rein analytisch:
ist $X$ ein komplexer Raum, und besitzt $Y \subseteq X$ Steinsches Komplement,
so ist die analytische Superh"ohe $\leq 1$.
Ist n"amlich $f:X' \longrightarrow X $ vorgegeben, so kann man durch
"Ubergang zu einer Steinschen Umgebung von $x'$ annehmen,
dass $X'$ Steinsch ist. 
$f$ faktorisiert durch den abgeschlossenen Graphen\footnote{Alle
komplexen R"aume seien Hausdorffsch.}
und damit ist $f^{-1}(X-Y)$ isomorph zu einer abgeschlossenen
Teilmenge von $X' \times (X-Y)$, also Steinsch, und die
Kodimension von $f^{-1}(Y)$ ist wieder $\leq 1$.
\begin{satz}
Sei $X$ ein affines algebraisches ${\bf C}-$Schema und
$V({\bf a})=Y \subseteq X$
eine algebraische Untervariet"at, so ist
$${\rm supht}\, {\bf a}={\rm supht}^{\rm end}\, {\bf a}
={\rm supht}^{\rm an}\, (Y^{\rm an},X^{\rm an}) \, .$$
\end{satz}
{\it Beweis}.
Die Gleichheit der ersten beiden Superh"ohen beruht auf dem Satz 4.1.5.
Die analytische Superh"ohe kann nat"urlich nicht kleiner sein als die
endliche algebraische Superh"ohe, da die dabei zu
betrachtenden algebraischen Testvariet"aten analytisch interpretierbar sind,
und sich die Dimensionsbegriffe entsprechen.
\par\smallskip\noindent
Sei umgekehrt $f: X' \longrightarrow X^{\rm an}$ ein Morphismus komplexer
R"aume, $x' \in X',\, f(x')=x$.
Wir k"onnen $X'$ als irreduzibel annehmen.
Das Erweiterungsideal ${\bf a}{\cal O}_{X^{\rm an},x}$
besitzt die Nullstelle $Y^{\rm an}$
und $Y'$ wird durch ${\bf a}{\cal O}_{X',x'}$ beschrieben.
Da ${\bf C}-$Stellenalgebren noethersch sind,
\cite{grauertas}, Kap. I, 35.2, Satz 3, und Kap. II, 30, Satz 1, ist
${\rm codim}_{x'}\, (Y',X')={\rm ht}\, ({\bf a}{\cal O}_{X',x'})
\leq {\rm supht}\, {\bf a}$. \hfill $\Box$
\par\bigskip
\noindent
{\bf Bemerkung}
Der Beweis zeigt zugleich, dass zur Berechnung der analytischen
Superh"ohe eines algebraisch definierten komplexen Raumes nur
endliche Morphismen betrachtet werden m"ussen.
Es ist zu fragen, ob dies f"ur beliebige komplexe R"aume gilt.
\begin{kor}
Sei $U=D({\bf a}) \subseteq {\rm Spek}\, A=X$ mit einer integren,
endlich erzeugten ${\bf C}-$Algebra $A$, $Y=X-U$.
Dann sind folgende Aussagen "aquivalent.
\par\smallskip\noindent
{\rm (1)} ${\rm supht}\, {\bf a} \leq 1$.
\par\smallskip\noindent
{\rm (2)} Alle abgeschlossenen Unterschemata der Dimension $\leq 2$ von $U$
sind affin.
\par\smallskip\noindent
{\rm (3)} F"ur jede abgeschlossene analytische Fl"ache
$S \subseteq X^{\rm an}$ ist $S \cap U^{\rm an}$ Steinsch.
\end{kor}
{\it Beweis}.
Die "Aquivalenz von (1) und (2) wurde bereits in 4.2.3 gezeigt.
Sei ${\rm supht}\, {\bf a}=1$ und $S \subseteq X^{\rm an}$ eine abgeschlossene
analytische Fl"ache mit der Normalisierung
$f:\tilde{S} \longrightarrow S \hookrightarrow X^{\rm an} $.
Dann besitzt $f^{-1}(Y)$ auf der normalen Fl"ache $\tilde{S}$
die Kodimension $\leq 1$,
und nach dem Satz von Simha,
siehe \cite{simha},
ist das Komplement Steinsch.
Damit ist die Normalisierung von $U \cap S$ Steinsch, und nach der
analytischen Version des Satzes von Chevalley,
gilt das auch f"ur $U \cap S$ selbst.
\par\smallskip\noindent
Sei umgekehrt (3) erf"ullt, und eine algebraische Fl"ache $S' \subseteq U$
gegeben. Es ist $S'=S \cap U$ mit einer abgeschlossenen Fl"ache $S$ in $X$
und somit ist $U \cap S \subseteq S$ eine offene Teilmenge einer
affin-algebraischen Fl"ache, die nach Voraussetzung Steinsch ist.
Nach dem Korollar 4.3.3 ist dann $U \cap S$ auch affin. \hfill $\Box$
\par\bigskip
\noindent
Wir wollen Beispiele angeben, wo $Y=V({\bf a}) \subseteq X={\rm Spek}\, A$
nicht affines Komplement besitzt, aber trotzdem die Superh"ohe
des Ideals eins ist.
Beispiele hierf"ur ergeben sich nach den bisherigen "Uberlegungen,
wenn $U$ Steinsch, aber nicht affin ist.
Solche Beispiele sind aber bekannt.
Man betrachtet irreduzible Kurven $D \subseteq S$ auf
einer projektiven glatten Fl"ache $S$, die alle
anderen Kurven schneidet, wo $S-D$ Steinsch ist,
der Selbstschnitt aber nicht positiv und daher $S-D$ nicht affin ist.
Interpretiert man diese Ergebnisse im affinen Kegel, so erh"alt man
Beispiele von quasiaffinen Steinschen, nicht affinen Variet"aten.
\par\bigskip
\noindent
{\bf Beispiel}
Sei $C$ eine projektive glatte Kurve "uber $\bf C$, und $\cal E$ eine
lokal freie Garbe auf $C$ vom Rang 2.
Das projektive B"undel $S={\bf P}({\cal E})$ ist eine geometrische
Regelfl"ache, siehe \cite{haralg}, V.2.2.
Es sei eine exakte Garbensequenz
$$0 \longrightarrow {\cal O}_C \longrightarrow {\cal E}
\longrightarrow {\cal O}_C
\longrightarrow 0 $$
gegeben. Die surjektive Abbildung ${\cal E} \longrightarrow {\cal O}_C $
liefert einen Schnitt $s:C \longrightarrow S$,
siehe \cite{haralg}, Prop. II.7.12, und wir setzen $W:= S-s(C)$ und
$D=s(C)$.
\par\smallskip\noindent
${\rm Ext}^1({\cal O}_C,{\cal O}_C) = H^1(C,{\cal O}_C)$ klassifiziert
diese Erweiterungen von ${\cal O}_C$ durch ${\cal O}_C$,
siehe \cite{haralg}, Prop. III.6.3.
Sei $g(C) \neq 0$ und eine nicht triviale Erweiterung durch
ein Element $0 \neq  t \in H^1(C,{\cal O}_C)$ gegeben.
Nach \cite{umemura} ist $D^2=0$ und $D.Z >0$ f"ur jede andere
irreduzible Kurve auf $P$.
Aufgrund von $D^2=0$ kann $W$ nicht affin sein (es gibt nur
konstante Funktionen auf $W$).
\par\smallskip\noindent
Es gibt ein Kriterium von
Umemura, dass in obiger Situation das Komplement von $D$ sogar Steinsch ist.
Die Exponentialsequenz f"uhrt zu
$0 \longrightarrow H^1(C,{\bf Z}) \longrightarrow H^1(C,{\cal O}_C)$,
und dies zu
$H^1(C,{\cal O}_C) \longrightarrow H^1(C,{\cal O}_C)/H^1(C,{\bf Z}) \, .$
Einem Element $t \in H^1(C,{\cal O}_C)$, das eine Erweiterung und
damit einen Schnitt beschreibt, wird also ein Element aus der Jakobischen
Variet"at $J$ zu $C$ zugeordnet. Ist nun
das Bild von ${\bf C}t$ in $J$ kompakt und $t \neq 0$, so ist
$W=S-D$ Steinsch, zum Beweis siehe \cite{umemura}, Th\'{e}or\`{e}me 2.
Die Bedingung ist nat"urlich f"ur eine elliptische Kurve f"ur $t \neq 0$
immer erf"ullt, so dass "uber einer elliptischen Kurve die Komplemente
solcher Schnitte Steinsch, aber nicht affin sind.
F"ur den elliptischen Fall geht dieses Beispiel auf Serre zur"uck und
wird, einschlie"slich dem
Nachweis des Schnittverhaltens und der Steinschheit, auch in
\cite{haramp}, VI.3, besprochen.
\par\smallskip\noindent
Ist nun $A$ ein homogener Koordinatenring von $S$, und ist $W=D_+({\bf a})$
mit einem homogenen Ideal, so ist das Urbild $p^{-1}(W)=D({\bf a})$
unter der Kegelabbildung nach dem Satz 2.3.3 ebenfalls nicht affin,
daf"ur aber als Urbild unter einer Quotientenbildung
nach einer ${\bf C}^\times-$Operation wieder Steinsch.
Damit ist ${\rm supht}\, {\bf a}=1$, aber $D({\bf a})$ nicht affin.
\par
\bigskip
\noindent
\subsection{Die Superh"ohe eines Schemas}
\par
\bigskip
\noindent
Sei ${\bf a}$ ein Ideal im noetherschen Ring $A$. Wir wollen die Superh"ohe
des Ideals $\bf a$ intrinsisch als Invariante zu $U=D({\bf a})$ definieren.
Schon in 4.2.3 haben wir gesehen, dass bei endlich erzeugten $K-$Algebren
die Superh"ohe eins zu der nur von $U$ abh"angigen Aussage "aquivalent ist,
dass alle Fl"achen darin affin sind.
\par\bigskip
\noindent
{\bf Definition}
Sei $X$ ein Schema.
Dann nennen wir die gr"o"ste Zahl $d$ mit der Eigenschaft, dass es ein
noethersches affines Schema $T$
mit einem abgeschlossenen Punkt $P \in T$ der H"ohe $d$
und einen affinen Schemamorphismus $f:T-P \longrightarrow X$ gibt,
die {\it Superh"ohe} von $X$, ${\rm supht}\, X$.
\par\smallskip\noindent
Ist $X$ eine Variet"at "uber einem K"orper $K$, so nennen wir die entsprechende
Zahl unter der Einschr"ankung, dass $T$ eine affine
Variet"at ist, die endliche Superh"ohe von $X$.
\par\bigskip
\noindent
{\bf Bemerkung}
Bei der Definition dieser Superh"ohenzahl kann man sich auf
integre, lokale, komplette, normale Testschemata $T$ beschr"anken.
Im Fall eines lokalen Ringes setzen wir
$T^\times =T -\{ P \}$, wobei $P$ der abgeschlossene Punkt ist.
\par\smallskip\noindent
Ist $X$ separiert, $f: T- \{P \} \longrightarrow X$ ein affiner Morphismus und
${\rm dim}\, T \geq 2$, so ist die Abbildung $f$ nicht auf ganz $T$
ausdehnbar. Gebe es n"amlich eine Ausdehnung $\tilde{f}:T \longrightarrow X$,
so w"are das ebenfalls ein
affiner Morphismus, und f"ur eine affine Umgebung $\tilde{f}(P) \in V$
w"aren $\tilde{f}^{-1}(V)$ und $f^{-1}(V)$ affin, mit
$\tilde{f}^{-1}(V)=f^{-1}(V) \cup \{ P \}$. Ein affines Schema bleibt aber
nach Herausnahme eines abgeschlossenen Punktes der H"ohe $\geq 2$ nicht affin.
\par\smallskip\noindent
Ist $X$ leer, so ist ${\rm supht}\, X = 0$, da in diesem Fall
$T -\{ P\}$ leer sein muss, also ${\rm dim}\, T =0$.
Umgekehrt besitzt ein nicht-leeres Schema eine Superh"ohe $\geq 1$.
Sei ${\rm Spek}\, K \longrightarrow X$ ein Punkt. Dann ist
$$ {\rm Spek}\, K[Y]_{(Y)} \supseteq {\rm Spek}\, K(Y)
\longrightarrow {\rm Spek}\, K \longrightarrow X $$
ein affiner Morphismus mit ${\rm Spek}\, K[Y]_{(Y)}$ eindimensional.
\par\smallskip\noindent
Ist $X$ affin, so ist ${\rm supht}\, X \leq 1$, da in diesem
Fall f"ur einen affinen Morphismus
$T-\{ P \} \longrightarrow X$ bereits $T-\{ P\}$ affin ist, woraus
aus der Kodimensionseigenschaft sofort ${\rm dim}\, T \leq 1$ folgt.
Steinsche, aber nicht affine Beispiele zeigen,
dass die Umkehrung nicht gelten muss.
\begin{prop}
Sei $A$ ein noetherscher Ring.
F"ur eine offene Teilmenge $U=D({\bf a}) \subseteq A$ ist
${\rm supht}\, U={\rm supht}\, {\bf a}$.
\par\smallskip\noindent
Ist $A$ eine endlich erzeugte $K-$Algebra, so gilt das auch f"ur die
endliche Superh"ohe.
\end{prop}
{\it Beweis}.
Sei $A \longrightarrow R$ ein Ringhomomorphismus in einen lokalen, integren,
normalen, noetherschen Ring der Dimension $m={\rm supht}\, {\bf a}$
mit $V({\bf a}R)=V({\bf m}_R)=\{ P \}$.
Dann ist die Abbildung
$f^{-1}(U)={\rm Spek}\, R - \{ P \} \longrightarrow U$ affin und
damit ist ${\rm supht}\, U \geq {\rm supht}\, {\bf a}$.
\par\smallskip\noindent
Sei umgekehrt $f:T- \{ P\} \longrightarrow U$ ein affiner Morphismus mit
$T$ lokal, integer, normal mit $d={\rm dim}\, T = {\rm supht}\, U$.
Ist $d=0$, so ist nichts zu zeigen,
ist $d=1$, so folgt, dass $U$ nicht leer ist.
Das Ideal $\bf a$ ist dann nicht nilpotent und besitzt damit ein
minimales Primoberdieal der H"ohe $\geq 1$, also ist die Superh"ohe des Ideals
$\geq 1$. Sei also $d \geq 2$.
Da $T={\rm Spek}\, R$ normal ist,
haben wir $\Gamma(T^\times,{\cal O}_T)=\Gamma(T,{\cal O}_T)=R$
und zu $f$ geh"ort der Ringhomomorphismus
$A \longrightarrow R$, so dass die Abbildung nach ${\rm Spek}\, A$
fortsetzbar ist.
Dabei muss $\bar{f}(P) \not\in U$ sein, da andernfalls die Abbildung sogar
nach $U$ fortsetzbar w"are, was nach der Bemerkung nicht sein kann.
Unter $A \longrightarrow R$ gilt also f"ur das Erweiterungsideal
$V({\bf a}R) =\{ P\}$ und damit
${\rm supht}\, {\bf a} \geq {\rm supht}\, U$. \hfill $\Box$
\begin{prop}
{\rm (1)}
Ist $X' \longrightarrow X$ ein affiner Morphismus noetherscher Schemata,
so ist ${\rm supht}\, X \geq {\rm supht}\, X'$. Dies gilt insbesondere
f"ur abgeschlossene Unterschemata $X' \hookrightarrow X$.
\par\smallskip\noindent
{\rm (2)}
Die Superh"ohe von $X$ ist gleich dem Maximum der Superh"ohen der Komponenten
von $X$.
\par\smallskip\noindent
{\rm (3)}
Sei $X$ ein noethersches Schema und $Y \subseteq X$ eine abgeschlossene
Teilmenge, $U=X-Y$. Dann gilt die Absch"atzung
$$ {\rm supht}\, X \leq {\rm supht}\, Y + {\rm supht}\, (X-Y) \, .$$
\par\smallskip\noindent
{\rm (4)}
Ist $X=U \cup V$ mit $U,V$ offen, so ist
$$ {\rm supht}\, X \leq {\rm supht}\, V + {\rm supht}\, U \, .$$
{\rm (5)}
F"ur ein noethersches Schema $X$ der Dimension $d$ ist
${\rm supht}\, X \leq d+1$.
\end{prop}
{\it Beweis}.
(1) ist klar.
Zu (2). Bei $T \longrightarrow X$ kann man $T$ als irreduzibel
annehmen, so dass das Bild in einer Komponente von $X$ liegt,
die die Superh"ohe von $X$ besitzt.
\par\smallskip\noindent
Zu (3). Sei $f: T \supseteq T^\times \longrightarrow  X$ ein affiner
Morphismus mit $T={\rm Spek}\, R$ noethersch, lokal, integer, komplett,
wobei die Dimension von $T$ die Superh"ohe von $X$ sei.
Sei $f^{-1}(Y)=V({\bf a})-\{ P\}$ mit einem Ideal ${\bf a} \subseteq R$.
Dann ist einerseits ${\rm dim}\, V({\bf a}) \leq {\rm supht}\, Y$,
wie $V({\bf a})-\{ P \}=f^{-1}(Y) \longrightarrow Y$ zeigt.
Andererseits ist auch $D({\bf a})=f^{-1}(U) \longrightarrow U$ affin und
daher ist ${\rm alt}\,  {\bf a} \leq {\rm supht}\, U$.
Da $R$ komplett ist, ist $R$ katen"ar, siehe \cite{eisenbud}, Cor. 18.10,
und somit gelten mit einem minimalen Primoberideal $\bf p$
von $\bf a$ die Absch"atzungen
\begin{eqnarray*}
{\rm supht}\, X ={\rm dim}\, R &
=& {\rm dim}\, R/{\bf p}+{\rm ht}\, {\bf p} \\
& \leq & {\rm dim}\, R/{\bf a} +{\rm alt}\,  {\bf a} \\
& \leq & {\rm supht}\, Y + {\rm supht}\, (X-Y)
\end{eqnarray*}
Zu (4). Sei $X=U \cup V$. Dann ist $Y=X-U$ eine abgeschlossene Teilmenge
von $V$, nach (1) ist ${\rm supht}\, Y \leq {\rm supht}\, V$
und die Behauptung folgt aus (3).
\par\smallskip\noindent
Zu (5). Wir f"uhren Induktion "uber die Dimension, der Induktionsanfang ist
klar.
Wir k"onnen wegen (2) $X$ als integer annehmen. F"ur eine nichtleere
offene affine Teilmenge $U$ ist dann nach (3)
${\rm supht}\, X \leq {\rm supht}\, U+ {\rm supht}\, (X-U)
\leq 1+d $.  \hfill $\Box$
\par\bigskip
\noindent
{\bf Beispiel}
Sei $Y$ eine projektive Variet"at der Dimension $d$. Dann ist die
Kegelabbildung $X^\times \longrightarrow Y$ eine affine Abbildung, daher
ist die Superh"ohenzahl $\geq d+1$, und wegen (5) muss Gleichheit gelten.
\begin{kor}
Sei $A$ ein noetherscher Ring, $X={\rm Spek}\, A$, ${\bf a} $ ein Ideal.
Dann gelten folgende Aussagen.
\par\smallskip\noindent
{\rm (1)}
Ist ${\bf b}$ ein weiteres Ideal mit ${\bf b} \subseteq {\bf a}$, so ist
${\rm supht}_A\, {\bf a}
\leq {\rm supht}_{A/{\bf b}}\, {\bf a} +{\rm supht}_A \, {\bf b} $.
\par\smallskip\noindent
{\rm (2)}
Ist ${\bf a}={\bf b}+ {\bf c}$ {\rm (}als Radikal{\rm )},
so ist
${\rm supht}\, {\bf a} \leq {\rm supht}\, {\bf b} + {\rm supht}\, {\bf c} $.
\end{kor}
{\it Beweis}. Die Aussagen folgen direkt aus der Proposition. \hfill $\Box$
\par\bigskip
\noindent
{\bf Bemerkung}
Die Eigenschaft ${\rm supht}\, X=1$ f"ur eine Variet"at $X$ 
kann im nicht quasiaffinen Fall nicht dadurch charakterisiert werden,
dass alle abgeschlossenen Fl"achen darin affin sind, 4.2.3 und 4.2.4
lassen sich also nicht verallgemeinern.
Ein Beispiel dazu ist die Steinsche, nicht affine Teilmenge einer Regelfl"ache.
Diese hat als globalen Schnittring die Konstanten, so dass auch
Satz 4.2.2 nicht ohne die quasiaffine Voraussetzung gilt.
Dasgleiche gilt f"ur die Aussagen 4.3.2, 4.3.3 und 4.3.5.
\par\bigskip
\noindent
{\bf Bemerkung}
${\rm supht}\, X \leq 1$ l"asst sich manchmal dadurch erweisen,
dass jede Abbildung $T - \{P \} \longrightarrow X$ mit ${\rm dim}\, T =2$
fortsetzbar ist und daher nicht affin sein kann, siehe Satz 4.5.2.
Man kann aber die Superh"ohe nicht dadurch charakterisieren,
dass f"ur ${\rm dim}\, T > {\rm supht}\, X$ jede au"serhalb
eines Punktes $P$ definierte Abbildung nach $X$ fortsetzbar ist.
\par\smallskip\noindent
Sei $R=K[X_1,...,X_n,Y_1,...,Y_n](X_iY_j-X_jY_i:\, 0 \leq i,j \leq n)$.
Dann ist $q=X_1/Y_1=...=X_n/Y_n$ eine auf $D(X_i,Y_j)$ definierte
meromorphe Abbildung nach ${\bf P}^1_K$, die nicht in den abgeschlossenen
Punkt hinein ausdehnbar ist.
Die Superh"ohe der projektiven Gerade ist aber eins.
\begin{kor}
Sei $X$ ein noethersches Schema mit ${\rm supht} \leq d$.
Dann besitzt das Komplement von $X$ in jeder offenen Einbettung
$X \subseteq X'$ mit $X'$ noethersch und separiert eine
Kodimension $\leq d$.
\end{kor}
F"ur eine affine Teilmenge $U$ von $X'$ ist $U \cap X \hookrightarrow X$
ein affiner Morphismus, so dass sich die Voraussetzung auf $U \cap X$
"ubertr"agt. Da die Behauptung lokal ist, kann man $X'$ als affin
annehmen. Dann ist aber die Aussage klar aufgrund von 4.4.1. \hfill $\Box$
\par\bigskip
\noindent
{\bf Bemerkung}
Die Umkehrung des Korollars gilt nicht.
In 1.12.3 hatten wir gezeigt, dass eine semiaffine Variet"at
in einem separierten noetherschen Schema ein Komplement der Kodimension
$\leq 1$ besitzt. Eine semiaffine Variet"at kann aber projektive Kurven
enthalten und die Superh"ohe kann $\geq 2$ sein.
\par
\bigskip
\noindent
\subsection{Die Superh"ohe in graduierten Ringen}
\par
\bigskip
\noindent
Im Beispiel am Schluss von 4.2 haben wir gesehen, dass in einem
graduierten Ring ein homogenes Ideal die endliche Superh"ohe zwei haben kann,
obwohl es auf allen homogenen Fl"achen die H"ohe eins besitzt.
Wir besprechen in diesem Abschnitt die Frage, inwiefern sich die
Superh"ohe homogen messen l"asst und beschreiben eine
Situation, wo man aus der Affinit"at aller homogenen Fl"achen auf
die Superh"ohe eins schlie"sen kann. Diese Situation ist insbesondere
im Beispiel am Ende von 4.3 erf"ullt, so dass man dort
die Superh"ohe eins rein algebraisch und
unabh"angig von der Steinschheit erh"alt.
\begin{prop}
Sei $A$ ein ${\bf Z}-$graduierter Ring, $\bf a$ ein homogenes
Ideal. Dann wird die noethersche und die endliche Superh"ohe in
homogenen Erweiterungen angenommen.
\end{prop}
{\it Beweis}.
Sei $f: A \longrightarrow R$ ein Ringhomomorphismus,
wo die Superh"ohe von $\bf a$
als H"ohe eines maximalen Ideals $\bf m$ angenommen wird.
Die Abbildung
$A \longrightarrow A[T,T^{-1}],\, A_d \ni a_d \longmapsto a_dX^d$
(die Kooperation) f"uhrt zu einer Faktorisierung
$$A \longrightarrow R[T,T^{-1}] \longrightarrow R \, ,$$
wobei vorne $a_d \longmapsto f(a_d)T^d$ und hinten
$T \longmapsto 1$ geschickt wird.
Sei ${\bf a}=(a_1,...,a_s)$ mit homogenen Elementen.
In $R[T,T^{-1}]$ ist dann
$${\bf a}R[T,T^{-1}]=(f(a_1)T^{d_1},...,f(a_s)T^{d_s})=
(f(a_1),...,f(a_s))={\bf m}R[T,T^{-1}] \, ,$$
da die letzte Gleichung bereits in $R$ gilt (als Radikal).
Es hat aber ${\bf m}R[T,T^{-1}]$ die gleiche H"ohe wie $\bf m$. \hfill $\Box$
\par\bigskip
\noindent
{\bf Bemerkung}
Die Aussage erlaubt nicht darauf zu schlie"sen, dass die Superh"ohe
als H"ohe eines maximalen homogenen Ideales in einem graduierten Ring
angenommen wird.
Insbesondere kann man die endliche Superh"ohe $k$ nicht homogen in
der Dimension $k$ realisieren.
\begin{satz}
Sei $S$ eine projektive glatte Fl"ache "uber einem algebraisch abgeschlossenen
K"orper $K$, $Y \subset S$ ein Primdivisor mit $Y^2 \geq 0$ und
$Y.C >0$ f"ur alle Kurven $C \neq Y$, $U=S-Y$.
Dann ist jede Abbildung $T \supseteq T-\{P\} \longrightarrow U$,
wobei $T$ eine zweidimensionale, normale, integre, affine Variet"at ist,
fortsetzbar.
Insbesondere ist ${\rm supht}^{\rm end}\, U=1$.
Ist $S={\rm Proj}\, A$ mit einer endlich erzeugten graduierten
$K-$Algebra $A$ und mit $U=D_+({\bf a})$, so ist ${\rm supht}\, {\bf a}=1$.
\end{satz}
{\it Beweis}.
Die Zus"atze sind klar: Gibt es einen affinen Morphismus
$f:T^\times=T-\{ P \} \longrightarrow U$ mit ${\rm dim}\, T \geq 2$, so gebe es
auch einen mit ${\rm dim}\, T=2$ (eine abgeschlossene Teilmenge),
was nach der ersten Bemerkung von 4.4. der Fortsetzbarkeit
widerspricht.
Damit besitzt auch der affine Kegel $U'$ die endliche Superh"ohe eins
und damit nach Satz 4.1.5 auch die Superh"ohe eins.
\par\smallskip\noindent
Ist sogar $Y^2 >0$, so ist $Y$ ger"aumig und $U$ affin, woraus
sich sofort die Superh"oheneigenschaft und auch die
Fortsetzbarkeit ergibt. Sei also im Folgenden
$Y^2=0$.
\par\smallskip\noindent
Sei also $f:T^\times \longrightarrow U$ ein Morphismus einer integren,
normalen, affinen Fl"ache $T$, ohne Einschr"ankung sei $T^\times$ regul"ar.
Ist $f(T^\times)$ ein Punkt, so ist $f$ nat"urlich ausdehnbar.
Liegt $f(T^\times) \subseteq U$ in einer irreduziblen Kurve $C$, so ist
$C$ nach Voraussetzung nicht projektiv, also affin.
Dann ist aber $f$ als Abbildung nach $C \subseteq U$ nach 3.4.1
ausdehnbar.
Es sei also das Bild von $f$ zweidimensional und $f$ dominant.
\par\smallskip\noindent
Sei $T \hookrightarrow T'$ eine offene Einbettung in eine projektive Fl"ache
mit dem Komplement $D'$.
Es sei $p:\tilde{T} \longrightarrow T'$ eine Singularit"atenaufl"osung von
$T'$ als auch eine Aufl"osung der Undefinierbarkeitsstellen von
$f:T' \supseteq T^\times \longrightarrow S$, siehe \cite{haralg}, 3.8.1
und Theorem V.5.5.
Es ist also
$\tilde{f}: \tilde{T} \longrightarrow S$ eine Fortsetzung von
$f$ auf $T^\times \cong \tilde{T}-p^{-1}(P)-p^{-1}(D')$.
Es seien $C_1,...,C_n$ die irreduziblen eindimensionalen
Komponenten von $p^{-1}(P)$ und
$D_1,...,D_m$ entsprechend die Komponenten von $p^{-1}(D')$.
\par\smallskip\noindent
Da $\tilde{f}$ surjektiv ist, induziert $\tilde{f}$ eine
Abbildung $\tilde{f}^*$ der Divisoren (=Cartier-Divisoren).
Sei
$\tilde{f}^*(Y)=C+D$ mit $C=k_1C_1+...+k_nC_n$ und $D=l_1D_1+...+l_mD_m$
der Urbilddivisor von $Y$,
andere Komponenten kommen nicht in Frage.
Es ist dann $\tilde{f}_*(C)$ ein Vielfaches von $Y$,
da eine Komponente aus $\tilde{f}^*(Y)$ ihr Bild in der
irreduziblen Kurve $Y$ hat.
Daher ist
$$0=\tilde{f}_*(C).Y=C.\tilde{f}^*(Y)=C.C +C.D=C.C \,\, .$$
Nach \cite{artin}, Theorem 2.3, besitzt aber ein Divisor $\neq 0$,
dessen Komponenten alle (singul"ar) kontrahierbar sind, negativen
Selbstschnitt, und $C$ wird ja unter $p$ auf den Punkt $P$ kontrahiert.
Also ist $C=0$ und das hei"st, dass
in $\tilde{f}^*(Y)$ keine Komponente oberhalb von $P$
vorkommt und dass f"ur diese Komponenten $C_i$ immer
$\tilde{f}(C_i) \not\subseteq Y$ gilt.
\par\smallskip\noindent
Das Urbild von $Y$ unter $\tilde{T}-p^{-1}(D') \longrightarrow S$
besteht also h"ochstens aus einzelnen Punkten (auf den $C_i$).
Da $S-Y$ lokal affin in $S$ ist, muss das Urbild von $Y$ "uberall
die Kodimension eins haben, so dass das nicht sein kann.
\par\smallskip\noindent
Es muss also $\tilde{f}(C_i) \subseteq U=S-Y$ f"ur alle Kurven $C_i$
oberhalb von $P$ sein.
Da es in $U$ keine projektiven Kurven gibt,
werden alle diese Kurven
unter $\tilde{f}$ auf Punkte in $U$ kontrahiert, und zwar auf einen
einzigen, da die $C_i$ zusammenh"angen.
Damit ist aber schon $f$ selbst in $P$ fortsetzbar. \hfill $\Box$
\par\bigskip
\noindent
{\bf Bemerkung}
Die Voraussetzungen des Satzes lassen sich wohl nicht abschw"achen.
Dass man bei $Y$ reduzibel nicht auf Superh"ohe eins
schlie"sen kann, hat schon das Beispiel am Schluss von 4.2 gezeigt.
Ist $Y \cap C = \emptyset $, so liegen in $U$ projektive Kurven.
Deren Kegelabbildungen sind Morphismen auf
zweidimensionalen Fl"achen, die nicht in die Kegelspitze hinein
ausdehnbar sind.
Auch auf die Bedingung $Y^2 \geq 0$ kann man vermutlich nicht verzichten.
Ist $K={\bf C}$ und $Y^2 <0$, so kann man $Y$ nach einem
Satz von Grauert, \cite{grauertmod},
kontrahieren mit einem komplexen Raum als Bildraum.
Dieser komplexe Raum definiert dann aber eine nur im Kontraktionspunkt
nicht definierte Abbildung zur"uck.
\par\bigskip
\noindent
{\bf Beispiel}(Vgl. \cite{goodlandman}, example 6.10)
Sei $D_0 \subseteq {\bf P}^2_K$ eine glatte Kurve vom Grad drei, also eine
elliptische Kurve, und $P_1,...,P_n$ Punkte auf $D_0$.
Es sei $S$ die Fl"ache, die entsteht, wenn man an den $n$ Punkten
aufbl"ast, und $D$ sei die eigentliche Transformierte von $D_0$.
Es ist $D^2=9-n$.
$D$ schneidet nat"urlich alle exzeptionellen Divisoren, der Schnitt mit
den anderen Kurven h"angt von der Auswahl der Punkte ab.
Es sei $Q \in D_0$ ein Punkt, der durch eine lineare Form aus der
Kurve herausgeschnitten wird, sagen wir $D_0 \cap V_+(x)=\{ Q\}$
(also ein inflection point, siehe \cite{haralg}, II. Exc. 6.19.2).
Durch die Isomorphie $P \longmapsto P-Q \in {\rm Div}^0/ {\rm HD} $ wird dann
die Kurve mit einer Gruppenstruktur mit $Q$ als neutralem Element
versehen.
\par\smallskip\noindent
Es seien nun die Punkte so gew"ahlt, dass sie ${\bf Z}-$linear
unabh"angig in dieser Gruppenstruktur sind.
Eine beliebige Kurve $C \subseteq {\bf P}^2_K$ trifft dann
$D_0$ nicht nur in den ausgew"ahlten Punkten.
Ist n"amlich $C=V_+(p)$ und $p$ vom Grad $d$, und w"are
$V_+(p) \subseteq \{ P_1,...,P_n\}$, so w"are
$ (p/x^d )= r_1P_1+...+r_nP_n -3dQ$, und es w"urde in der
Gruppenstruktur eine Relation vorliegen.
Das hei"st f"ur die eigentliche Transformierte $D$, dass sie mit jeder
Kurve einen positiven Schnitt besitzt.
\par\smallskip\noindent
Ist $n=9$, so ist der Selbstschnitt der Kurve $=0$ und es sind die
Voraussetzungen des Satzes erf"ullt. In einem affinen Kegel dar"uber
hat also das beschreibende Ideal die Superh"ohe eins, dagegen ist
die offene Menge nicht affin.
\par\bigskip
\noindent
{\bf Bemerkung}
Umemura fragt in \cite{umemura}, probl\`{e}me A,
ob in der komplexen Fl"a\-chensituation generell f"ur eine irreduzible Kurve $D$
mit $D^2=0$ und $D.C >0$ f"ur alle Kurven $C \neq D$ folgt,
dass $U=S-D$ Steinsch ist.
Es ist mir weder im vorangehenden Beispiel
noch in der Situation des Beispiels am Ende von 4.2, wenn das Kriterium
von Umemura nicht erf"ullt ist, klar, ob $U$ bzw.
das Urbild $W=p^{-1}(U)$ davon im affinen Kegel Steinsch ist oder nicht.
Aus dem Satz ergibt sich, dass $U$ und $W$ die Superh"ohe eins besitzen
und dass daher nach Satz 4.3.5 alle analytischen Fl"achen in $W$
(zumindest diejenigen, die einen analytischen Abschluss im Kegel haben)
Steinsch sind.
\par\smallskip\noindent
Mit dieser Problematik h"angt direkt das sogenannte
Hyperfl"achenproblem der komplexen Analysis zusammen,
vergl. \cite{diederich1}, 3.4.
Es geht dabei um die Frage, ob eine offene Menge $U \subseteq X$ in einem
Steinschen Raum $X$ der Dimension $\geq 3$
selbst Steinsch ist, wenn f"ur jede Hyperfl"ache
$H \subseteq X$ der Durchschnitt $U \cap H$ Steinsch ist.
In den obigen Beispielen ist (im Kegel) also die Voraussetzung des
Hyperfl"achenproblems erf"ullt.
Ein Gegenbeispiel zu diesem Problem findet sich in \cite{diederich2} und
\cite{diederich1}, worauf zu verwandten Fragestellungen aus der komplexen
Analysis hier generell verwiesen sei.
\par\bigskip
\noindent
{\bf Bemerkung}
Ist f"ur eine (nicht notwendigerweise irreduzible) Kurve $D$ auf
einer projektiven Fl"ache $S$ die offene Menge $U=S-D$ semiaffin, so ist $U$
bereits affin, wenn in $U$ keine projektiven Kurven liegen,
vergleiche 2.4.5 bzw. \cite{goodman}, Prop. 3.10.
\par
\bigskip
\noindent
\subsection{"Uberdeckungen mit trivialer affiner Klassengruppe}
\par
\bigskip
\noindent
Wir haben gesehen, dass im Allgemeinen die geometrische
Bedingung, dass die noethersche Superh"ohe eins ist, nicht hinreichend
f"ur die Affinit"at ist. Trotzdem kann man in vielen Situationen
die Affinit"at geometrisch sichern.
In diesem Abschnitt interessieren wir uns exemplarisch f"ur den
bestm"oglichsten Fall,
dass es f"ur ein affines Schema $X= {\rm Spek}\, A$
eine Erweiterung $A \longrightarrow B$ gibt mit ${\rm AKG}\, B=0$ und
so, dass die Affinit"at von $U \subseteq X$ aus der des Urbildes $U'$ folgt.
Das Ziel ist dabei, die Affinit\"atsfrage f\"ur Hyperfl"achen
in ${\rm Spek}\, A$ auf eine geometrische
Kodimensionsfrage in ${\rm Spek}\, B$ zur\"uckzuf\"uhren.
Als Ringhomomorphismen, die den R"uckschluss auf die Affinit"at zulassen,
haben wir in Abschnitt 1.14 endliche Erweiterungen (Satz von Chevalley),
treuflache Erweiterungen
und direkte Summanden kennengelernt.
Als Beispiele werden wir hier die Normalisierung von Fl"achen
und f\"ur Monoidringe die Darstellungen mit der
Durchschnittseigenschaft besprechen.
\par
\bigskip
\medskip
\noindent
{\it Fl"achen}
\begin{satz}
Sei $A$ ein integrer, exzellenter, zweidimensionaler Ring.
Dann besitzt eine Kurve $C \subseteq X={\rm Spek}\, A$ genau dann affines
Komplement, wenn das Urbild $\tilde{C}$ in der Normalisierung $\tilde{X}$
die reine Kodimension eins hat.
\end{satz}
{\it Beweis}.
Die eine Richtung ist klar, habe also das Urbild die Kodimension eins.
Nach dem Satz von Nagata hat dann $\tilde{C}$ affines Komplement, und
nach dem Satz von Chevalley auch $C$ selbst. \hfill $\Box$
\begin{satz}
Sei $A$ ein lokaler, integrer, exzellenter, zweidimensionaler Ring.
Genau dann ist die Normalisierung lokal (das bedeutet, dass $A$ analytisch
irreduzibel ist),
wenn in ${\rm Spek}\, A$ jede Kurve affines Komplement besitzt.
\end{satz}
{\it Beweis}.
Sei zun\"achst die Normalisierung $B$ lokal mit maximalem Ideal
${\bf n}$ und ${\bf p}$ ein Primideal der H\"ohe eins in $A$.
Es ist dann ${\bf p}B \neq {\bf n}$, da dieses Erweiterungsideal das
Urbild von $V({\bf p})$ beschreibt, und damit einen Punkt enthalten muss,
der auf $\bf p$ abbildet, das maximale Ideal $\bf n$ aber auf
$\bf m$ runterschneidet.
Daraus folgt, dass ${\bf p}B$ die reine Kodimension eins besitzt und daher
ist $U$ affin.
\par\smallskip
\noindent
Seien jetzt umgekehrt \"uber dem abgeschlossenen Punkt $P$ zumindest
zwei Punkte $Q$ und $R$. Wir zeigen, dass es dann Kurven mit nicht
affinem Komplement gibt. Jede Bildkurve einer
Kurve, die durch den einen, aber nicht durch den anderen Punkt geht,
ist daf"ur ein Kandidat, doch muss man etwas sorgf"altiger argumentieren.
\par\smallskip\noindent
Wegen der Exzellenz gibt es eine nichtleere, offene, affine Teilmenge $U$,
auf der die Normalisierung eine Isomorphie ist. Es seien $Z_1,...,Z_n$
die Komponenten von $\tilde{X} - U$ und es sei $\tilde{C}$ eine
Kurve $\neq Z_i,\, i=1,...,n$, die durch $R$ aber nicht durch $Q$ geht.
Die Existenz einer solchen Kurve wird so gesichert:
sei $V$ eine affine Umgebung von $Q$, die $R$ nicht enth\"alt, aber die
die generischen Punkte der $Z_i$ enth\"alt.
Das Komplement von $V$ besteht
dann aus Kurven, und $\tilde{C}$ sei eine Komponente davon, auf der
$R$ liegt.
\par\smallskip\noindent
Das Bild von $\tilde{C}$ sei $C$. Da $\tilde{C}$ generisch in $U$ verl"auft,
besitzt das Urbild von $C$ nur $\tilde{C}$ als gro"se Komponente,
w"ahrend $Q$ ein isolierter Punkt im Urbild ist. Also ist das Komplement
von $C$ nicht affin. \hfill $\Box$
\par\bigskip
\noindent
{\bf Bemerkung} Die R\"uckrichtung des obigen Satzes gilt in jeder
Dimension $\geq 2$. Wenn \"uber dem maximalen Ideal mehrere Punkte
liegen, so muss man f"ur die allgemeine Argumentation
irreduzible Hyperfl\"achen statt wie oben Kurven durch den einen Punkt
nehmen, die von den Komponenten des nicht normalen Ortes verschieden sind,
und die den anderen Punkt nicht treffen.
Das Bild dieser Hyperfl"ache hat dann nicht affines Komplement,
da das Urbild eine Komponente besitzt, die nicht die H"ohe eins hat.
\par\bigskip
\noindent
Der folgende Satz liefert eine gro"se Anzahl von Kurven mit nicht affinem
Komplement auf einer Fl\"ache, die "uber einem abgeschlossenen Punkt in
der Normalisierung mehrere Punkte besitzt (also eine analytisch nicht
irreduzible Singularit\"at ist).
\begin{satz}
Sei $X={\rm Spek}\, A$ eine affine, integre, noethersche, exzellente Fl\"a\-che
und $P$ ein abgeschlossener Punkt,
\"uber dem in der Normalisierung zumindest zwei Punkte liegen. Dann hat
eine irreduzible Kurve $C$, die durch $P$ geht, die als Kurve regul"ar ist
{\rm (}oder h"ochstens kuspidal{\rm )} und die generisch im normalen Ort
verl"auft, nicht affines Komplement.
\end{satz}
{\it Beweis}.
Zum Urbild $C'$ von $C$ in der Normalisierung geh"oren alle Punkte oberhalb
von $P$. Da der generische Punkt von $C$ im offenen normalen Ort liegt,
worauf die Normalisierung eine Isomorphie ist, gibt es "uber
diesem Punkt nur einen Punkt der H"ohe eins, $C'$ besitzt also nur eine
Komponente der Dimension eins, die mit $\tilde{C}$ bezeichnet sei.
W"urden alle Urbildpunkte von $P$ auf $\tilde{C}$ liegen, w"urden unter
der eingeschr"ankten birationalen Abbildung
$\tilde{C} \longrightarrow C$ verschiedene Punkte identifiziert, was
der Voraussetzung regul"ar (bzw. kuspidal) widerspricht. \hfill $\Box$
\par\bigskip
\noindent
{\bf Bemerkung}
Betrachtet man die Normalisierung aufgrund von 2.5.7 als eine
Aufblasung mit dem Ideal $I$, so besagt die generische Bedingung, dass der
Aufblasungsort $Z=V(I)$ nicht die Kurve $C$ umfasst, sondern nur einzelne
Punkte mit ihr gemeinsam hat.
Es ist dann $C'$ die totale Transformierte  und
$\tilde{C}$ die eigentiche Transformierte von $C$.
Diese ist die Aufblasung von $C$ an dem Erweiterungsideal $I$ auf $C$ mit
dem Aufblasungsort $Z \cap C$.
Nur dann kann $C$ affines Komplement besitzen, wenn bei dieser
Kurvenaufblasung bereits alle Punkte oberhalb von $P$ dazukommen,
die bei der Fl"achenaufblasung von $X$ oberhalb von $P$ entstehen.
Im gewissen Sinn kann man also sagen, dass die Kurve $C$ nur dann affines
Komplement besitzt, wenn ihre Singularit"at so stark ist wie die
Singularit"at der Fl"ache, auf der sie liegt. In eine "ahnliche
Richtung weist folgender Umstand: Besitzt die Kurve nicht nur affines
Komplement, sondern wird sie sogar
idealtheoretisch durch eine einzige Funktion beschrieben, so ist
die Einbettungsdimension der Kurve nur um eins geringer als die
Einbettungsdimension der Fl"achensingularit"at, die Differenz zwischen
Einbettungsdimension und Dimension also gleich.
\par\bigskip
\noindent
{\bf Beispiel}
Sei $A=K[y,yt,t^2]=K[y,x,z]/(x^2-y^2z) $ mit der
Normalisierung $K[y,t]$.
Der singul"are Ort ist dabei die Gerade $Z=V(x,y)$
(dies ist auch der Aufblasungsort).
"Uber einem Punkt aus dieser Gerade liegen bei $z=0$ ein Punkt und
bei $z \neq 0$ zwei Punkte  (Es sei ${\rm char}\, K \neq 2 $).
Die Fl"ache entsteht aus der affinen Ebene,
indem man auf der $t-$Achse ($y=0$) die gegen"uberliegenden Punkte
identifiziert.
\par\smallskip\noindent
Das Erweiterungsideal von $(x,y)$ in der
Normalisierung ist das Hauptideal $(y)$, damit ist
$D(x,y)$ affin.
Auf die Voraussetzung, dass die Kurve nicht ganz im singul"aren
Ort verl"auft, kann man also nicht verzichten.
\par\smallskip\noindent
Wegen $x^2 -y^2z= (x-y)(x+y)-(z-1)y^2$ liegt auf der Fl"ache die
regul"are Kurve (eine Gerade)
$C=V(x-y,z-1)$. Sie schneidet den singul"aren Ort in
$x=y=0, z=1$, die eigentliche Transformierte ist also einfach wieder
eine affine Gerade, und nach dem Satz kann das Komplement von $C$ nicht affin
sein.
In der Tat ist das Erweiterungsideal gleich $(yt-y,t^2-1)$,
und es ist $V(yt-y,t^2-1)=V(t-1) \cup V(t+1,y)$.
\par\smallskip\noindent
Betrachte auf ${\bf A}^2$ die durch $y=t^2-1$ definierte Kurve, die man durch
$t \longmapsto (t^2-1,t)$ parametrisieren kann.
Die Bildkurve $C$ hat dann die Parametrisierung
$t \longmapsto (t^2-1,(t^2-1)t,t^2)$ und dabei werden die beiden
Punkte $t=1$ und $t=-1$ verklebt. $C$ ist eine nodale Kurve und
hat affines Komplement, da die
eigentliche Transformierte beide Punkte oberhalb von $(0,0,1)$ trifft.
"Ubrigens wird die Kurve durch die Gleichung $z-y+1$ allein beschrieben,
und hat auch deshalb affines Komplement (Mit dieser Ebene hat unsere
Fl"ache den irreduziblen Schnitt $V(x^2-y^2(y-1))$).
\par
\bigskip
\medskip
\newpage
\noindent
{\it Monoidringe}
\par
\bigskip
\noindent
Sei $M$ ein normales, torsionsfreies, endlich erzeugtes Monoid im
Quotienten\-gitter $\Gamma={\bf Z}M \cong {\bf Z}^d$.
Ferner sei $M$ {\it spitz} (oder {\it positiv}),
d.h. das neutrale Element sei die einzige Einheit.
Dann gibt es eine Einbettung $ M \hookrightarrow {\bf Z}^k$ so,
dass $M= \Gamma \cap {\bf N}^k$ ist, siehe \cite{bruns}, exc. 6.1.10
oder die Divisorenklassendarstellung in 3.6.
Man spricht von einer Darstellung des Monoids mit der 
{\it Durchschnittseigenschaft}. Umgekehrt ist ein so als
Durchschnitt gegebenes Monoid normal und spitz.
\par\smallskip\noindent
Die Monoideinbettung $M \hookrightarrow {\bf N}^k$ liefert
eine Ringeinbettung
$$K[M] \hookrightarrow K[{\bf N}^k]=K[T_1,...,T_k] \, . $$
Die Restklassenabbildung ${\bf Z}^k \longrightarrow {\bf Z}^k/\Gamma=:D$
liefert eine $D-$Graduierung auf $K[T_1,...,T_k]$, dessen Unterring der
nullten Stufe gerade $K[M]$ ist. Das hei"st insbesondere, dass man
Monoidringe zu einem normalen spitzen Monoid als direkte Summanden
einer Polynomalgebra erhalten kann.
\par\smallskip\noindent
Mit der Proposition 1.14.4 kann man die Affinit"atsfrage f"ur Monoidringe
auf die Frage nach der Kodimension im Polynomring zur"uckf"uhren.
\begin{satz}
Sei $M$ ein endlich erzeugtes, torsionsfreies Monoid und
$K$ ein normaler noetherscher Ring mit ${\rm AKG}\, K=0$.
Dann besitzt der Monoidring $K[M]$ eine Einbettung
$K[M] \hookrightarrow B$ in eine endlich erzeugte $K-$Algebra von
Laurentpolynomen derart,
dass eine offene Teilmenge $D({\bf a}) \subseteq {\rm Spek}\, K[M]$
genau dann affin ist, wenn das Erweiterungsideal ${\bf a}B$ die reine
Kodimension eins besitzt {\rm (}oder das Nullideal ist{\rm )}.
\end{satz}
{\it Beweis}.
Sei $\tilde{M}$ die Normalisierung von $M$ und
$\tilde{M}={\bf Z}^s \times M'$ mit $M'$ spitz,
siehe \cite{bruns}, Theorem 6.1.4 und Prop. 6.1.3.
$M' \hookrightarrow {\bf N}^k$ sei eine Darstellung mit der
Durchschnittseigenschaft. Dann ist
$$ K[M] \longrightarrow K[\tilde{M}] \cong
K[{\bf Z}^s][M'] \longrightarrow 
K[V_1,...,V_s,V_1^{-1},...,V_s^{-1}][T_1,...,T_k]=B$$
insgesamt eine Abbildung von endlichem Typ, die als Hintereinanderschaltung
einer endlichen Abbildung und eines direkten Summanden erlaubt, auf die
Affinit"at zur"uckzuschlie"sen. Mit $K$ besitzt auch $B$ triviale
affine Klassengruppe, siehe Satz 3.3.5. 
Ist ${\bf a} \subseteq K[M]$ ein Ideal, das in $B$ die Altitude
$\leq 1$ hat, so definiert dort das Erweiterungsideal eine affine offene
Menge, und das muss dann bereits f"ur $D({\bf a})$ gelten. \hfill $\Box$
\par\bigskip
\noindent
{\bf Beispiel}
Sei $K[T_1,...,T_n,S_1,...,S_m]$ gegeben mit der durch
$T_i \longmapsto 1,\, S_j \longmapsto -1$ definierten
${\bf Z}-$Graduierung.
Der Ring der nullten Stufe ist dann
$$A=K[T_iS_j: 1 \leq i \leq n,\, 1 \leq j \leq m ]
\cong K[X_{i,j}]/(X_{i_1,j_1}X_{i_2,j_2}-X_{i_1j_2}X_{i_2,j_1}) \, .$$
Sei ${\bf p}_i=(T_iS_1,...,T_iS_m)$ und ${\bf q}_j=(T_1S_j,...,T_nS_j)$.
Betrachte
$ Y=V({\bf p}_{i_1})\cup ... \cup V({\bf p}_{i_r}) \cup
V({\bf q}_{j_1}) \cup ... \cup V({\bf q}_{j_s})$.
Das Urbild von $V({\bf p}_i)$ ist $V(T_i) \cup V(S_1,..,S_m)$, das
Urbild von $V({\bf q}_j)$ ist entsprechend $V(S_j) \cup V(T_1,..,T_n)$.
Die Superh"ohe von ${\bf p}_i$ ist damit m, die von ${\bf q}_j$ ist n.
Sei $Y \neq \emptyset $, und sei $V({\bf p}_i)$ eine Komponente davon.
Dann geh"ort $V(S_1,..,S_m)$ zum Urbild von $Y$ und das Urbild kann nur
dann die reine Kodimension eins haben, wenn in $Y$ auch ein
$V({\bf q}_j)$ vorkommt, und in diesem Fall liegt dann auch die Kodimension
eins vor.
Genau dann besitzt also $Y$ affines Komplement, wenn
in $Y$ beide Sorten vertreten sind.
\par\smallskip\noindent
Unser Standardbeispiel $A=K[u,v,y,x]/(ux-vy)$ ist von diesem
Typ, und zwar ist $A \cong K[T_1S_1,T_1S_2,T_2S_1,T_2S_2]$.
$D(x,y)$ ist nicht affin, das Erweiterungsideal davon ist
$(T_2S_1,T_2S_2)=(T_2) \cap (S_1,S_2)$, was eine Komponente der H"ohe zwei
hat.
Dagegen hat das Ideal $(x+u)(x,y)=(x^2+ux,xy+uy)$ das Erweiterungsideal
$((T_2S_2+T_1S_1)T_2S_2,\, (T^2S_2+T_1S_1)T_2S_1)$ mit den minimalen
Primoberidealen $(T_2)$ und $(T_2S_2+T_1S_1)$. Es liegt als Kodimension
eins vor und daher ist $D((x+u)(x,y))$ affin.

%% file: divhoka.tex
\section{Superh"ohe und kohomologische Dimension}
\par
\bigskip
\noindent
In diesem abschlie"senden Kapitel untersuchen wir verschiedene
Begriffsbildungen, die geeignet sind, die Eigenschaft einer Hyperfl"ache,
affines Komplement zu besitzen, auf irreduzible Teilmengen h"oherer
Kodimension entsprechend zu "ubertragen. Im Mittelpunkt stehen dabei
verschiedene H"ohenbegriffe zu einem Primideal, Absch"atzungen zwischen
ihnen und die Frage, wann sie gleich der Kodimension sind.
\par
\bigskip
\noindent
\subsection{Die Superh"ohe in regul"aren Ringen}
\par
\bigskip
\noindent
Das klassische Superh"ohenresultat ist das affine Dimensionstheorem.
\begin{satz}
Sei ${\bf A}_K^n$ der affine Raum "uber einem algebraisch
abgeschlossenen K"orper und
$Y$ und $Z$ irreduzible Variet"aten der Dimension $r$ und $s$.
Dann besitzt eine irreduzible Komponente von $Y \cap Z$ zumindest eine
Dimension $r+s-n$.
\end{satz}
Zum {\it Beweis} siehe \cite{haralg}, Prop. I.7.1 oder
 \cite{mumcomplex}, Prop. 3.28. \hfill $\Box$
\par\bigskip
\noindent
{\bf Bemerkung} Das affine Dimensionstheorem ist in folgender Weise als
ein Superh"ohenresultat zu interpretieren.
Ist $Y=V({\bf p})$, $Z=V({\bf q})$, und
${\bf p}+{\bf q} \subseteq {\bf r}$ ein minimales Primoberideal, so besagt
der Satz einfach
${\rm ht}\, {\bf r} \leq {\rm ht}\, {\bf p} + {\rm ht}\, {\bf q}$.
Modulo $\bf q$ ist $\bf r$ ein minimales Primoberideal von $\bf p$, und
dort ist also
${\rm ht}_{A/{\bf q}}\, {\bf p}
\leq {\rm ht}_{A/{\bf q}}\, {\bf r}
\leq {\rm ht}\, {\bf r}-{\rm ht}\, {\bf q} \leq {\rm ht}\, {\bf p}$.
Das hei"st, dass die Superh"ohe von $\bf p$
bezogen auf abgeschlossene Einbettungen gleich der H"ohe von $\bf p$ ist.
Der Satz von Serre verallgemeinert das affine Dimensionstheorem
auf beliebige regul"are Ringe.
\begin{satz}
Sei $A$ ein regul"arer Ring, ${\bf p},\, \bf q$ Primideale in $A$
und $\bf r$ ein minimales Primoberideal von ${\bf p}+{\bf q}$.
Dann gilt die Absch"atzung
$${\rm ht}\, {\bf r} \leq {\rm ht}\, {\bf p} +{\rm ht}\, {\bf q} \, .$$
\end{satz}
F"ur den {\it Beweis} siehe
\cite{serre}, V.B.6, Th\'{e}or\`{e}me 3. \hfill $\Box$
\par\bigskip
\noindent
Dieser Satz gilt f"ur beliebige Ideale $\bf a$ und $\bf b$, wenn man
darin die H"ohe durch die Altitude ersetzt.
\begin{satz}
Sei $A$ ein noetherscher regul"arer Ring. Dann gilt f"ur
ein Ideal $\bf a$ die Gleichheit
${\rm alt}\, {\bf a} ={\rm supht}\, {\bf a}$.
\end{satz}
{\it Beweis}.
Es ist zu zeigen, dass f"ur ein Ideal $\bf a$ die Superh"ohe durch
die Altitude beschr"ankt ist.
Die endliche Superh"ohe von $\bf a$ ist durch das Schnittverhalten von
$V({\bf a})$ mit anderen abgeschlossenen Untervariet"aten
in ${\rm Spek}\, A[T_1,...,T_n]$ bestimmt. Da der Polynomring wieder regul"ar
ist, folgt die Aussage f"ur die endliche Superh"ohe aus 5.1.2.
Nach dem Beweis zu 4.1.2 gilt die Superh"ohenabsch"atzung dann auch f"ur alle
ganzen Ringwechsel.
\par\smallskip\noindent
Wir nehmen nun $A$ als lokal an und betrachten die
Komplettierung $\hat{A}$ mit dem Erweiterungsideal
$\hat{\bf a}={\bf a}\hat{A}$.
Ein minimales Primoberideal von $\hat{\bf a}$ schneidet auf ein minimales
Primoberideal von $\bf a$ herunter.
Da $A$ formal
"aquidimensionell ist und sich diese Eigenschaft
auf "aquidimensionelle Restklassenringe
"ubertr"agt, siehe \cite{EGAIV}, Cor. 7.1.4 und Prop. 7.1.10,
besitzen die minimalen Primoberideale der Erweiterung eines Primideals
alle die gleiche H"ohe wie dieses Ideal.
Damit ist die Altitude von $\hat{\bf a}$ gleich der Altitude von $\bf a$.
Ferner "andert sich nach 5.5.2 die Superh"ohe bei der Komplettierung
nicht.
Wir k"onnen also $A$ als komplett, lokal und regul"ar annehmen.
Nach dem Cohenschen Struktursatz, \cite{nagloc}, 31, ist dann
$A=B[[T_1,...,T_n]]$, wobei $B$ ein K"orper oder ein diskreter
Bewertungsring ist mit Restek"orper $K=B/(p)=A/{\bf m}_A$
(ist $B$ ein K"orper, so sei $p=0$).
\par\smallskip\noindent
Sei $A \longrightarrow R$ ein Ringhomomorphismus
in einen lokalen, noetherschen, kompletten Ring, wo $\bf a$ geometrisch zum
maximalen Ideal wird.
Da die Lokalisierung nach einem Primideal wieder regul"ar ist, k"onnen
wir annehmen, dass ${\bf m}_R$ auf ${\bf m}_A$ herunterschneidet.
Es sei $K \longrightarrow L=R/{\bf m}_R$ die Erweiterung der Restek"orper
und $L \ni x_i,\, i \in I,$ eine Familie algebraisch unabh"angiger
Elemente so, dass $L$ algebraisch "uber $K(x_i)$ ist. Es seien $x'_i \in R$
Urbilder von $x_i$.
\par\smallskip\noindent
Wir betrachten die Erweiterung
$B(X_i:\, i \in I)=B[X_i:\, i \in I]_{p B[X_i]}=B'$,
wobei $B'$ wieder ein K"orper oder ein diskreter Bewertungsring ist.
Ist $F=\sum b_\nu X^\nu \not\in pB[X_i]$, so ist $b_\nu \in B^\times$
f"ur ein $\nu$  und damit ist $\sum \bar{b}_\nu x^\nu \neq 0$ in $L$,
da die $x^\nu$ linear unabh"angig "uber $K$ sind.
Damit geht unter der Abbildung
$B[X_i,\, i \in I] \longrightarrow R,\, X_i \longmapsto x_i'$
ein Polynom $F \not\in pB[X_i]$ auf eine Einheit und
die Abbildung faktorisiert durch $B'$.
\par\smallskip\noindent
Beim "Ubergang $B[[T_1,...,T_n]] \hookrightarrow B'[[T_1,...,T_n]]$
"andert sich die Altitude des Ideals und die Superh"ohe nicht.
Wir k"onnen also zus"atzlich annehmen, dass bei $A \longrightarrow R$
die Erweiterung der Restek"orper algebraisch ist. Es ist dann zu zeigen,
dass diese Abbildung ganz ist.
Sei $x \in R$ und $A'$ der von $A$ und $x$ in $R$ erzeugte Unterring.
Da ${\bf m}_R$ auf ${\bf m}_A$ herunterschneidet, ist
$A'/({\bf m}_A \cap A')=K[x]$ ein endlicher Erweiterungsk"orper.
Da andererseits die Faser "uber ${\bf m}_A$ nur aus ${\bf m}_R \cap A'$
besteht, ist nach \cite{nagloc}, Theorem 30.6, die Komplettierung
von $A'$ endlich.
Damit liegt $x$ in einer endlichen Ringerweiterung von $A$,
und $R$ ist ganz "uber $A$. \hfill $\Box$
\par\bigskip
\noindent
{\bf Bemerkung}
Die Eigenschaft, dass f"ur jedes Primideal H"ohe und Superh"ohe
"ubereinstimmt, gilt beispielsweise auch f"ur zweidimensionale, normale,
exzellente Ringe.
Aus ${\rm AKG}\, A=0$ folgt generell, dass alle Primideale der H"ohe eins
auch die Superh"ohe eins besitzen, und f"ur maximale Ideale stimmt
immer H"ohe mit Superh"ohe "uberein.
\par\smallskip\noindent
F"ur eine dreidimensionale, normale, integre $K-$Algebra
von endlichem Typ sind auch die Primideale der H"ohe zwei unproblematisch, 
ihre Superh"ohe ist nach 4.1.4 und 4.1.5 gleich der H"ohe.
Besitzt eine solche Algebra triviale affine Klassengruppe,
so gilt darin f"ur jedes Primideale die Gleichheit von H"ohe und Superh"ohe.
In einer dreidimensionalen faktoriellen $K-$Algebra von endlichem Typ
gilt damit die Aussage des Satzes genau so.
Im Allgemeinen kann es aber in einem faktoriellen Ring Primideale geben,
deren Superh"ohe gr"o"ser als die H"ohe ist, siehe das n"achste Beispiel.
\par\bigskip
\noindent
{\bf Beispiel}
$A=K[X_1,...,X_6]/(X_1X_2+X_3X_4+X_5X_6)$ ist eine faktorielle,
f"unfdimensionale, isolierte Singularit"at.
Jede Hyperfl"ache hat affines Komplement und daher die Superh"ohe
eins. Das Primideal $(X_1,X_2,X_3)$ hat die H"ohe zwei, aber die Superh"ohe
drei, wie Reduktion modulo $(X_4,X_5,X_6)$ zeigt.
\par\bigskip
\noindent
{\bf Beispiel}
${\bf p}=(X_1,...,X_n)$ ist ein Primideal der H"ohe $n-1$ im
Hyperfl"achenring
$$K[X_1,...,X_n,Y_1,...,Y_n]/(f) \mbox{ mit }
f=X_1^k\cdot...\cdot X_n^k+Y_1X_1^{k+1}+...+Y_nX_n^{k+1}\, .$$
F"ur $K={\bf Z}$ ist die Aussage, dass die Superh"ohe von $\bf p$
ebenfalls $n-1$ ist, "aquivalent zur Monomialvermutung, dass in einem lokalen
noetherschen Ring zu einem Parametersystem $x_1,...,x_n$ das Monom
$x_1^k...x_n^k$ nicht in dem von $x_1^{k+1},...,x_n^{k+1}$
erzeugten Ideal enthalten ist. Im K"orperfall ist dies richtig.
Zur Monomialvermutung und zur "aquivalenten direkten Summand Vermutung
siehe \cite{roberts3}, \cite{hochsterdim}, \cite{hochstercan}
und \cite{bruns}, 9.2.
\par\bigskip
\noindent
{\bf Beispiel}
Sei $K$ ein algebraisch abgeschlossener K"orper mit ${\rm char}\, K \neq 2$.
Betrachte die Abbildung ${\bf A}_K^d \longrightarrow {\bf A}_K^{2d-1}$ mit
$$(r,t_1,...,t_{d-1}) \longmapsto 
(r^2,t_1,...,t_{d-1},rt_1,...,rt_{d-1})
=(x,y_1,...,y_{d-1},z_1,...,z_{d-1})\, .$$
Das Bild davon ist eine integre, nicht normale,
$d-$dimensionale, affine Variet\"at $X$.
Wir betrachten darin die $(d-1)-$dimensionale lineare Variet\"at 
$Y=V(x-1,y_1-z_1,...,y_{d-1}-z_{d-1})$.
Wir bestimmen die Urbildpunkte von $Y$.
Aus $r^2=x=1$ folgt $r=1$ oder $r=-1$.
Ist $r=1$, so kann man die $t_i$ frei w\"ahlen, was zugleich zeigt,
dass $Y$ ganz im Bild liegt und daher eine Hyperfl"ache von $X$ ist.
Ist dagegen $r=-1$, so muss $t_i=0$ sein f"ur alle $i$,
das Urbild von $Y$ hat also einen isolierten Punkt der H\"ohe $d$.
Wir haben damit also eine Hyperfl"ache in einem $d-$dimensionalen Spektrum
mit Superh"ohe $d$.
\par\bigskip
\noindent
{\bf Beispiel}
$A=K[X_1,...,X_n,Y_1,...,Y_n]/(X_iY_j-X_jY_i)$
ist eine normale Determinantensingularit"at der Dimension $d=n+1$.
Die beiden Hyperfl"achen $V(X_1,...,X_n)$ und $V(Y_1,...,Y_n)$ schneiden
sich darin im Nullpunkt und haben damit die Superh"ohe $n=d-1$.
Dies ist die maximale Superh"ohe, die eine Hyperfl"ache im Fall einer
normalen affinen Variet"at haben kann.
\par\bigskip
\noindent
{\bf Beispiel}
Wir wollen in einem noetherschen integren Ring Primideale finden mit
${\bf p} \subset {\bf q}$,
aber so, dass die Superh\"ohe von $\bf p$
gr"o"ser ist als die von $\bf q$.
\par\smallskip\noindent
Sei $K$ ein K"orper, wir betrachten dar"uber
den zwei-Minoren-Ring $A$ zur Matrix
$$\left(\matrix{ X_1 & X_2 & X_3 \cr  Y_1 & Y_2 & Y_3
\cr Z_1 & Z_2 & Z_3 \cr} \right) \, .$$
Sei $B=A[U,V]$. Mit $A$ ist auch $B$ ein integrer, normaler, noetherscher
Ring. Darin ist ${\bf q}=(U,V)$ ein Primideal der H"ohe=Superh"ohe zwei.
Sei $h={X_1 \over X_3}U+V$ und
$h_1 =X_3h=X_1U+X_3V,\,  h_2 =Y_3h=Y_1U+Y_3V,\,  h_3=Z_3h=Z_1U+Z_3V $.
Es gelten dann die Relationen
$Y_ih_1=X_ih_2,\, Z_ih_1=X_ih_3,\, Z_ih_2= Y_ih_3 \mbox{ f"ur } i=1,2,3$.
Auf $B_{X_2}=A_{X_2}[U,V]=K[X_1,X_2,X_2^{-1},X_3,Y_1,Z_1][U,V]$
ist $h_1$ prim, sei ${\bf p}=(h_1) \cap B$.
\par\smallskip\noindent
Es ist $(h_1,h_2,h_3) \subseteq {\bf p}$.
Wir zeigen, dass dies das einzige minimale Prim\-oberideal davon ist.
Sei $\bf r$ ein anderes. 
$\bf p$ wird auf $D(X)$ zu jeder Unbestimmten $X$
der Matrix durch eine der Funktionen
$h_1,h_2,h_3$ beschrieben, wie die Relationen zeigen.
Auf $D(X)$ ist also ${\bf p} \subseteq {\bf r}$, damit die Gleichheit
nicht gilt, muss $\bf r$ jede Unbestimmte enthalten.
Andererseits ist aber
${\bf p} \subseteq (X_1,X_2,X_3,Y_1,Y_2,Y_3,Z_1,Z_2,Z_3)$,
so dass ${\bf p} \subset {\bf r}$ folgt.
Damit ist also $V(h_1,h_2,h_3)=V({\bf p})$ und ${\bf p} \subset {\bf q}$.
\par\smallskip\noindent
Die Substitution
$U \longmapsto 1,\, V \longmapsto 0$ und
$ X_2,X_3,Y_2,Y_3,Z_2,Z_3 \longmapsto 0 $ nach
$K[X_1,Y_1,Z_1]$ sind mit den Minorenrelationen vertr"aglich
und zeigt, dass $(h_1,h_2,h_3)$ die Superh"ohe drei besitzt.
\par
\bigskip
\noindent
\subsection{Die kohomologische Dimension}
\par
\bigskip
\noindent
Die folgende Definition geht auf R. Hartshorne zur"uck, vergleiche
\cite{harcodi}.
\par\bigskip
\noindent
{\bf Definition}
Sei $X$ ein Schema. Dann hei"st die Zahl
$${\rm cd}\, X=
{\rm max}\, \{ n: H^n(X,{\cal F}) \neq 0 \mbox{ f"ur eine quasikoh"arente
Garbe } {\cal F}  \} $$
die {\it kohomologische Dimension} von $X$.
Daneben betrachten wir die Zahl ${\rm kohoht}\, X= {\rm cd}\, X+1$
und nennen sie die {\it kohomologische H"ohe} von $X$, da diese mit der
Superh"ohe besser zu vergleichen ist.
\footnote{In \cite{harspeiser} und \cite{lyub}
findet sich daf"ur der Begriff local
cohomological dimension, {\rm lcd}.}
Die kohomologische Charakterisierung affiner Schemata 1.13.1 besagt, dass ein
noethersches Schema $X$ genau dann affin ist, wenn ${\rm cd}\, X=0$ ist.
\begin{prop}
{\rm (1)}
Die kohomologische Dimension eines noetherschen Schemas kann man
mit koh"arenten Garben allein testen.
\par\smallskip\noindent
{\rm (2)}
Ist $f : X \longrightarrow Y $ ein affiner Morphismus noetherscher
und separierter Schemata, so ist $ {\rm cd}\, X  \leq  {\rm cd}\, Y$.
Ist $A \longrightarrow A'$ ein Ringhomomorphismus und ${\bf a} \subseteq A$
ein Ideal, so ist
${\rm cd}\, D({\bf a}A') \leq {\rm cd}\, D({\bf a})$.
\par\smallskip\noindent
{\rm (3)}
F"ur ein Ideal ${\bf a}$ in einem noetherschen Ring $A$ mit $X={\rm Spek}\, A$
ist
$${\rm cd}\, D({\bf a}) = {\rm max}_{x \in X}\, {\rm cd}\, D({\bf a}_x) \, .$$
\end{prop}
{\it Beweis}.
Zu (1).
Jede quasikoh"arente Garbe ist direkter Limes von koh"arenten
Garben, und auf einem  noetherschen Schema vertauscht die
Kohomologie mit direkten Limiten, siehe \cite{haralg}, Prop. III.2.9.
Ist also $H^r(X,{\cal F})=0 $ f"ur alle koh"arenten Garben, so auch f"ur
alle quasikoh"arenten Garben.
\par\smallskip\noindent
Zu (2).
Unter einem affinen Morphismus zwischen noetherschen separierten Schemata
ist $H^{i}(Y,f_*({\cal G}))=H^{i}(X,{\cal G})$ f"ur eine quasikoh\"arente
Garbe $\cal G$ auf $X$, siehe \cite{haralg}, Exc. III.4.1.
Ist nun $\cal G$ eine solche Garbe mit $H^{i}(X,{\cal G}) \neq 0 $, so
ist auch $H^{i}(Y,f_*({\cal G})) \neq 0 $.
\par\smallskip\noindent
Zu (3). Nach (2) ist die eine Absch\"atzung klar,
man hat zu zeigen,
dass die kohomologische Dimension an einem Punkt angenommen wird.
F"ur eine quasikoh"arente Garbe $\cal F$ auf $U=D({\bf a})$ gilt
$ H^n(U_{\bf p}, {\cal F}_{\bf p})=(H^n(U,{\cal F}))_{\bf p} $, da
Lokalisierung ein flacher Basiswechsel ist.
Ist $H^n(U,{\cal F}) \neq 0$, so muss dies auch in einem Punkt gelten,
da man das Verschwinden eines Moduls lokal testen kann. \hfill $\Box$
\par\bigskip
\noindent
{\bf Bemerkung}
Ein Spezialfall der letzten Aussage ist,
dass die Affinit\"at eine lokale Eigenschaft ist.
\begin{satz}
Ist $X$ ein noetherscher topologischer Raum der Dimension $d$, so ist
f"ur jede Garbe $\cal F$ von abelschen Gruppen $H^{i}(X,{\cal F})=0$
f"ur $i >d$.
Ist $X$ ein noethersches Schema, so ist ${\rm cd}\, X \leq {\rm dim}\, X$.
\end{satz}
F"ur den {\it Beweis} siehe \cite{tohoku}, Theorem 3.6.5 oder
\cite{haralg}, Theorem III.2.7. \hfill $\Box$
\par\bigskip
\noindent
F"ur offene Unterschemata affiner Schemata kann man diese Aussage um eins
verbessern.
\begin{kor}
Ist $A$ ein noetherscher Ring der Dimension $d$,
und $D({\bf a})=U \subseteq {\rm Spek}\, A$ offen,
so ist $H^{i}(U,{\cal F})=0$
f"ur $i \geq d $ und jede
quasikoh"arente Garbe $\cal F$. Es ist also ${\rm cd}\, U < {\rm dim}\, A$.
\end{kor}
{\it Beweis}.
Nach Proposition 5.2.1. kann man sich auf den lokalen Fall
beschr"anken.
F"ur einen lokalen Ring der Dimension $d$ besitzt $D({\bf a})$
bei ${\bf a} \subseteq {\bf m}$ die Dimension $\leq d-1 $, und die Aussage
folgt aus dem Satz. Ist dagegen $\bf a$ das Einheitsideal,
so ist $U=X$ affin mit ${\rm cd}\, U=0$. \hfill $\Box$
\par\bigskip
\noindent
{\bf Bemerkung}
Aufgrund der exakten Sequenz (mit $Y=X-U$) f"ur die lokale Kohomologie
$$ H^{i-1}(U,{\cal F}) \longrightarrow H^{i}_Y(X,{\cal F})
\longrightarrow H^{i}(X,{\cal F})
\longrightarrow H^{i}(U,{\cal F}) $$
verschwinden in einem Ring der Dimension $d$ f"ur $i > d$
die lokalen Kohomologiegruppen $H^{i}_Y(X,{\cal F})$.
F"ur das punktierte Spektrum eines lokalen Ringes ist die $d-$te lokale
Kohomologie nicht null, wie der folgende Satz von Grothendieck zeigt.
\begin{satz}
Sei $A$ ein noetherscher lokaler Ring der Dimension $d$ mit maximalem
Ideal $\bf m$, $U={\rm Spek}\, A - \{ {\bf m} \}$.
Dann ist die kohomologische Dimension des
punktierten Schemas gleich $d-1$, die kohomologische H"ohe gleich der
Dimension.
\end{satz}
{\it Beweis}.
Sei zun"achst $A$ ein noetherscher lokaler Ring und $M$ ein endlich erzeugter
$A-$Modul der Tiefe $t$, $X={\rm Spek}\, A$, $P$ der abgeschlossene
Punkt und ${\cal F}=\tilde{M}$.
Dann ist nach \cite{SGA2}, Exp. III. Exemple 3.4, 
$H^{i}_P(X,{\cal F}) =0$ f"ur $i < t$ und $H_P^t(X,{\cal F}) \neq 0$.
Insbesondere ist $H_P^d(X,{\cal O}_X) \neq 0$
f"ur einen Cohen-Macauley Ring der Dimension $d$.
Da unter einem affinen Morphismus die kohomologische Dimension nicht
gr"o"ser werden kann, k"onnen wir zur Komplettierung und auf eine
Komponente maximaler Dimension "ubergehen.
Sei also $A$ lokal, komplett und integer.
Dann ist $A$ nach dem Normalisierungssatz f"ur komplette Ringe
endlich "uber einem formalen Potenzreihenring, insbesondere
"uber einem regul"aren Ring $R$, dessen abgeschlossener Punkt mit $Q$
bezeichnet sei.
\par\smallskip\noindent
$R$ ist insbesondere ein Cohen-Macauley Bereich der Dimension $d$ und daher ist
$H^d_Q(R) \neq 0$.
Wir zeigen auf ${\rm Spek}\, R$, dass f"ur den endlich erzeugten $R-$Modul $A$
gilt $H^d_Q(A) \neq 0$.
Es ist $A \hookrightarrow Q(A) =Q(R)^r$ und damit ist auch
$aA \hookrightarrow R^r $ f"ur ein geignetes $a \in A$.
Der Restklassenmodul ist ein Torsionsmodul, hat also einen Tr"ager kleinerer
Dimension, und daher ist der letzte Teil der langen Kohomologiesequenz
gleich
$H^d_Q(A) \longrightarrow H^d_Q(R^r) \longrightarrow 0 $, so dass
$H^d_Q(A) \neq 0$ ist. Dieser Modul ist aber gleich $H^d_P(A)$. \hfill $\Box$
\par\bigskip
\noindent
Der n"achste Satz ist das kohomologische Analogon zur Additivit"at
der Superh"ohe 4.4.2.
Diese so naheliegende Aussage findet sich nicht
in \cite{harcodi}, daf"ur kann man sie auf eine Proposition
von dort zur"uckf"uhren.
\begin{satz}
Sei $X$ ein noethersches separiertes Schema, $U \subseteq X$ eine
offene Teilmenge und $Y=X-U$.
Dann ist
${\rm kohoht}\, X \leq {\rm kohoht}\, U + {\rm kohoht}\, Y$.
Ist $X=U \cup V$ mit $U,V$ offen, so ist
${\rm kohoht}\, X \leq {\rm kohoht}\, U + {\rm kohoht}\, V$.
\end{satz}
{\it Beweis}.
Die zweite Aussage folgt direkt aus der ersten wie in 4.4.2.
Sei $r= {\rm kohoht}\, U$ und $s={\rm kohoht}\, Y$.
\par\smallskip\noindent
Wir zeigen, dass f"ur $j >r$ und einen koh"arenten Modul $\cal F$
die Garben $H_Y^{j}({\cal F})$ versschwinden. Diese Garbe ist die Vergarbung
der Pr"agarbe $W \longmapsto H^{i}_{W\cap Y}(W,{\cal F})$,
siehe \cite{SGA2}, I.2.4.
Daher k"onnen wir die Sache auf affinen Teilmengen $W \subseteq X$ betrachten.
Da $X$ separiert ist, ist die Einbettung $W \hookrightarrow X$ ein affiner
Morphismus und daher ist ${\rm kohoht}\, (W \cap U) \leq {\rm kohoht}\, U =r$.
Die lange exakte Sequenz f"ur die lokale Kohomologie auf $W$ liefert
$$ H^{j-1}(W \cap U,{\cal F}) \longrightarrow
H^{j}_{Y \cap W}(W,{\cal F}) \longrightarrow H^{j}(W,{\cal F}) \, .$$
F"ur $j >r$ ist der linke Modul null wegen der kohomologischen H"ohe von
$U$ und $U \cap W$, der rechte Modul ist null, da $W$ affin ist, und damit
ist auch der mittlere null.
\par\smallskip\noindent
Wir haben also insgesamt ${\rm cd}\, Y =l \leq s-1,\, {\rm cd}\, U=r-1$ und
$H^{i}_Y({\cal F})=0 $ f"ur $i >r$.
Damit ist nach \cite{harcodi}, Prop. 1.4,
${\rm cd}\, X \leq {\rm max}\, (l+r,r-1)= l+r \leq s-1+r$,
also ist
${\rm kohoht}\, X \leq {\rm kohoht}\, U+ {\rm kohoht}\, Y$. \hfill $\Box$
\par\bigskip
\noindent
{\bf Bemerkung}
Der letzte Satz umfasst als Spezialfall die Aussage, dass f"ur ein separiertes
noethersches Schema, das durch $n$ affine Teilmengen "uberdeckt wird,
die kohomologische H"ohe h"ochstens gleich $n$ ist. Diese Aussage kann man
nat"urlich auch \v{c}echkohomologisch beweisen.
Aus dem Satz kann man wie in 4.4.2 f"ur ein noethersches Schema
${\rm cd}\, X \leq {\rm dim}\, X$ folgern unabh"angig von 5.2.2.
\begin{satz}
Sei $A$ ein kompletter, lokaler, noetherscher, integrer Ring der Dimension $d$
und ${\bf a}$
ein Ideal in $A$ so, dass jede Komponente von $V({\bf a})$ eine Dimension
von zumindest eins besitzt. Dann ist f"ur einen $A-$Modul $M$
$H^{d-1}_{\bf a}(M)=0$. Die kohomologische Dimension von $D({\bf a})$ ist also
h"ochstens $d-2$.
\par\smallskip\noindent
Ist $A$ lokal und noethersch, und $\bf a$ ein Ideal, so gilt die Aussage
ebenso, wenn $\bf a$ in jeder integren Komponente der Komplettierung
zumindest die Dimension eins hat.
\end{satz}
{\it Beweis}. F"ur einen Beweis siehe \cite{harcodi}, Theorem 3.1,
\cite{peskine}, Theorem 3.7, oder \cite{schenzel}, Theorem 2.14. \hfill $\Box$
\par\bigskip
\noindent
{\bf Bemerkung} Diesen Satz nennt man den Verschwindungssatz von
Lichten\-baum-Hartshorne. Ist speziell $A$ komplett, integer und
zweidimensional, so besagt er, dass das Komplement einer (reinen) Kurve
affin ist und beinhaltet somit einen neuen Beweis f"ur den Satz
von Nagata f"ur den kompletten Fall, auf den man aber auch den normalen
exzellenten Fall zur"uckf"uhren kann.
Der Satz ist analog zu dem (freilich einfacheren) Superh"ohenresultat
4.1.4, dass f"ur Kurven in normalen noetherschen Integrit"atsbereichen
die Superh"ohe gleich der H"ohe ist.
Mit dem Satz erh"alt man die folgende Aussage, die von Lichtenbaum
vermutet wurde und erstmals von Grothendieck
in \cite{grolocal}, Th\'{e}or\`{e}me 6.9, bewiesen wurde.
\begin{kor}
Ein quasiprojektives integres Schema "uber einem K"orper $K$
ist genau dann projektiv, wenn
die kohomologische Dimension gleich der Dimension ist.
\end{kor}
{\it Beweis}.
Ist $X$ projektiv mit Dimension $d$, so f"uhrt ein sehr
ger"aumiges Geradenb"undel zu einem punktierten affinen Kegel, auf dem
nach Satz 5.2.4 die $d-$te Kohomologie nicht verschwindet,
was dann wegen der Affinit"at der Kegelabbildung auch f"ur $X$ gilt.
\par\smallskip\noindent
Ist umgekehrt $X$ nicht projektiv, so gibt es eine offene
Einbettung von $X$ in ein projektives Schema $Y$, sei $Z=Y-X \neq \emptyset$.
Dann sind unter der Kegelabbildung mit einem affinen Kegel zu $Y$ die Fasern
"uber $Z$ eindimensional und das Urbild von $Z$ erf"ullt die
Voraussetzungen des Satzes 5.2.6. Damit ist die kohomologische Dimension von
$p^{-1}(X)$ und $X$ kleiner als $d$. \hfill $\Box$
\par\bigskip
\noindent
Als Korollar ergibt sich, dass die Aussage genauso gilt,
wenn man die kohomologische H"ohe durch die Superh"ohe
ersetzt. Ein projektives Schemas hat auf Grund des affinen Kegels eine
Superh"ohe $=d+1$. Besitzt umgekehrt das quasiprojektive Schema eine
solche Superh"ohe, so ist nach der Absch"atzung in 5.3.1 auch
die kohomologische Dimension $\geq d$ und der Satz von Lichtenbaum liefert
die Projektivit"at.
\par
\bigskip
\noindent
\subsection{Absch"atzungen f"ur Superh"ohe und kohomolo\-gi\-sche Dimension}
\par
\bigskip
\noindent
In diesem Abschnitt besprechen wir Absch"atzungen zwischen Superh"ohe,
kohomologischer H"ohe und weiteren H"ohenbegriffen.
F"ur ein Schema $X$ bezeichne ${\rm afra}\, X$ die
minimale Anzahl von affinen offenen Teilmengen von $X$,
die f\"ur eine affine \"Uberdeckung gebraucht wird. 
\begin{prop}
Sei $X$ ein noethersches separiertes Schema. Dann gelten die Absch"atzungen
$$ {\rm afra}\, X \geq {\rm kohoht}\, X \geq {\rm supht}\, X \, .$$
\end{prop}
{\it Beweis}.
Die erste Absch"atzung
ist \v{C}ech-Kohomologie bzw. der Additionssatz 5.2.5
f"ur die kohomologische H"ohe.
Sei $T$ ein lokales, noethersches, normales, integres Schema und
$T^\times \longrightarrow X$ ein affiner Morphismus.
Dann ist nach dem Satz von Grothendieck 5.2.4 und der Proposition 5.2.1
${\rm dim}\, T = {\rm kohoht}\, T^\times \leq {\rm kohoht}\, X$.
Also ist ${\rm supht}\, X \leq {\rm kohoht}\, X$. \hfill $\Box$
\begin{prop}
Sei ${\bf a}$ ein Ideal in einem noetherschen Ring $A$, $U=D({\bf a})$.
Dann gelten folgende Aussagen.
\par\smallskip\noindent
{\rm (1)} $ {\rm ara}\, {\bf a}  \geq {\rm afra}\, U
\geq {\rm kohoht}\, U \geq {\rm supht}\, U \geq {\rm alt}\, {\bf a} \,.$
\par\smallskip\noindent
{\rm (2)} Ist $A$ fastfaktoriell, so ist ${\rm ara}\, {\bf a}= {\rm afra}\, U$.
\par\smallskip\noindent
{\rm (3)} Ist $A$ lokal und ${\bf a}$ das maximale Ideal,
so sind die Zahlen alle gleich der Ringdimension.
\end{prop}
{\it Beweis}.
Zu (1).
Die erste Absch\"atzung beruht einfach darauf, dass die Basismengen $D(f)$
affin sind, die mittleren wurden in obiger Proposition gezeigt, die letzte
ist klar. 
\par\smallskip\noindent
(2) Ist $A$ fastfaktoriell, so wird jedes Komplement einer affinen Teilmenge
durch eine Funktion beschrieben.
\par\smallskip\noindent
(3) wurde schon in 4.1.1 gezeigt.  \hfill $\Box$
\par\bigskip
\noindent
{\bf Beispiel} 
In eindimensionalen, nicht fastfaktoriellen Ringen gibt es abgeschlossene
Punkte, deren Komplemente affin sind, die aber nicht geometrisch durch
eine Funktion beschrieben werden k"onnen.
Im Allgemeinen ist also die erste Absch"atzung echt.
\par\smallskip\noindent
Ein Beispiel f"ur eine solche affine Kurve ist etwa eine punktierte
elliptische Kurve. Ein Punkt $P$ auf
$C= \bar{C} - \{ Q\} \subseteq \bar{C}$
ist genau dann durch eine Funktion $f$ beschreibbar, wenn der
Hauptdivisor $(f)= nP-nQ$ ist. Dies ist genau dann der Fall, wenn
auf der durch $Q$ als neutralem Element definierten Gruppenstruktur
der Punkt $P$ ein Torsionselement ist, was aber nicht f"ur jeden Punkt
sein muss.
\par\smallskip\noindent
Das gleiche gilt f"ur Kurven in lokalen, zweidimensionalen, normalen,
exzellenten, nicht fastfaktoriellen Ringen, etwa den Kegelringen "uber
projektiven Kurven von Geschlecht $\geq 1$.
\par\bigskip
\noindent
{\bf Beispiel}
Die Absch"atzung der Superh"ohe durch die kohomologische H"ohe
enth"alt als Spezialfall die Aussage, dass die Superh"ohe eines
affinen Schemas $\leq 1$ ist.
Die im Abschnitt 4.3 besprochenen Beispiele von Steinschen,
nicht affinen Teilmengen zeigen, dass diese
Absch"atzung echt sein kann.
\par\bigskip
\noindent
{\bf Beispiel}
Im vierdimensionalen affinen Raum hat
die Vereinigung $Y$ zweier affiner
Ebenen, die sich in einem Punkt schneiden, die kohomologische H\"ohe
drei, aber die Superh\"ohe zwei.
Nach \cite{haralg}, III. Exc. 4.9., ist n"amlich einerseits
$H^2({\bf A}^4- Y,{\cal O}_X) \neq 0$.
Andererseits ist nach 5.1.3 die Superh"ohe von $Y$ zwei.
Insbesondere gilt bei ${\rm dim}\, A \geq 4$
5.1.3 nicht mit der kohomologischen H"ohe statt der Superh"ohe.
\par\bigskip
\noindent
{\bf Beispiel}
Sei $A={\bf Z}[X_1,X_2,X_3,Y_1,Y_2,Y_3]/(f)$
mit $f=X_1^2X_2^2X_3^2+Y_1X_1^3+Y_2X_2^3+Y_3X_3^3$
und ${\bf a}=(X_1,X_2,X_3)$.
Gem"a"s der Monomialvermutung m"usste das Ideal die
Superh"ohe 2 besitzen.
In \cite{roberts4} wird gezeigt,
dass $H^3_{(X_1,X_2,X_3)}(A)=H^2(D({\bf a}), {\cal O}_X)$
nicht verschwindet. Damit ist die kohomologische H"ohe des Ideals 3
und man kann diese Vermutung nicht kohomologisch best"atigen.
Vergleiche hierzu 4.2 und \cite{hochstercan}, insbesondere Abschnitt 6.
\par\bigskip
\noindent
{\bf Bemerkung}
Die Echtheit der Absch"atzung ${\rm afra}\, U \geq {\rm kohoht}\, U$
macht einige Probleme.
Ist ${\rm kohoht}\, U=1$, so ist $U$ affin und nat"urlich ${\rm afra}\, U=1$.
Sei ${\rm kohoht}\, U=2$.
Die einfachsten Beispiele daf"ur sind Komplemente $U$ von einem abgeschlossenen
Punkt auf einer (nicht lokalen) Fl"ache $X$.
${\rm afra}\, U \geq 3$
w"urde bedeuten, dass das Komplement nicht von zwei affinen Teilmengen
"uberdeckt werden kann.
Nehmen wir, um afra besser in den Griff zu bekommen, zus"atzlich an, dass
${\rm AKG}\, X=0$ ist, dann ist afra einfach die Anzahl der ben"otigten
Kurven, um $P$ mengentheoretisch herauszuschneiden.
Ist $X={\rm Spek}\, A$ affin und faktoriell,
so ist ${\rm afra}\, U={\rm ara}\, U$.
Zu untersuchen sind hier etwa glatte und
faktorielle Hyperfl"achen vom Grad $d \geq 4$ im ${\bf A}_K^3$.
\par\smallskip\noindent
Eine Dimension h"oher geht es um Kurven in einer
affinen, faktoriellen, dreidimensionalen Variet"at. Nach 5.2.6
besitzt das Komplement der Kurve die kohomologische H"ohe zwei,
dagegen gibt es keinen Grund, warum sie durch zwei Fl"achen
herausschneidbar sein sollte.
Schon f"ur den dreidimensionalen affinen Raum ist das ein offenes
Problem, siehe etwa \cite{lyub}.
Macht man in $K[X,Y,Z]$ alle $X-a,\, a \in K,$ zu Einheiten
($K$ algebraisch abgeschlossen), so ergibt sich eine affine Fl"ache,
die Kurve wird zu einem Punkt, und es ist dann
zu fragen, ob dieser Punkt durch Kurven herausschneidbar ist.
\par\bigskip
\noindent
Die Superh"ohe eines Ideals ${\bf a} \subseteq A$ gibt an,
welche H"ohe das Erweiterungsideal in einem noetherschen Ring h"ochstens
haben kann und ist durch die Ringdimension von $A$ beschr"ankt.
Die umgekehrte Frage, ob $\bf a$ bereits in einem niederdimensionalen
Ring definiert werden kann, f"uhrt zu einigen weiteren H"ohenbegriffen.
\par\bigskip
\noindent
{\bf Definition}
(1) Sei $U$ ein noethersches separiertes Schema. Dann nennt man die
minimale Dimension eines noetherschen separierten Schemas $U'$ mit
der Eigenschaft, dass es einen affinen Schemamorphismus
$f:U \longrightarrow U'$ gibt, die {\it Definitionsdimension} von $U$,
${\rm defdim}\, U$.
Die Definitionsdimension plus $1$ nennen wir die {\it Definitionsh"ohe},
${\rm defht}\, U$
\par\smallskip\noindent
(2) Ist $U \subseteq X$ eine offene Teilmenge eines noetherschen separierten
Schemas, so nennt man die minimale Dimension von $X'$ derart,
dass es einen affinen Morphismus $f:X \longrightarrow X'$ und
eine offene Teilmenge $U' \subseteq X'$ mit $U=f^{-1}(U')$ gibt,
die {\it eingebettete Definitionsdimension}, ${\rm defdim}\, (U,X)$.
\par\smallskip\noindent
(3) Ist $A$ ein noetherscher Ring und $\bf a$ ein Ideal, so nennt man die
minimale Ringdimension eines noetherschen Ringes $B$ derart, dass es
ein Ideal ${\bf b} \subseteq B$ und einen
Ringhomomorphismus $B \longrightarrow A$ gibt mit $V({\bf a})=V({\bf b}A)$ ,
die {\it Ringdefinitionsdimension} des Ideals,
${\rm defdim}^{\rm ring}\, {\bf a}$.
\begin{prop}
{\rm (1)}
Ist $U$ ein noethersches separiertes Schema, so ist
$$ {\rm dim}\, U \geq {\rm defdim}\, U \geq {\rm cd}\, U \, .$$
{\rm (2)}
Ist $U \subseteq X$ mit $X$ separiert und noethersch, so ist
$$ {\rm defdim}\, (U,X) \geq {\rm defdim}\, U \, .$$
{\rm (3)}
Ist $A$ ein noetherscher Ring und $\bf a$ ein Ideal, so ist
$$ {\rm defdim}^{\rm ring}\, {\bf a}
\geq {\rm defdim}\, (D({\bf a}),{\rm Spek}\, A) \, .$$
Ist $A$ lokal und $\bf a$ nicht das Einheitsideal, so ist sogar
$$ {\rm defdim}^{\rm ring}\, {\bf a} \geq {\rm defht}\, D({\bf a})
\geq {\rm supht}\, {\bf a}\, .$$
\par\smallskip\noindent
{\rm (4)}
Ist $K$ ein K"orper,
$A$ eine noethersche $K-$Algebra und $\bf a$ ein Ideal, so
ist ${\rm ara}\, {\bf a} \geq {\rm defdim}^{\rm ring}\, {\bf a}$.
Ist $A$ zus"atzlich lokal, so gelten die
Absch"atzungen 
$$ {\rm ara}\, {\bf a} \geq {\rm defdim}^{\rm ring}\, {\bf a}
\geq {\rm defht}\, D({\bf a}) \geq {\rm kohoht}\, {\bf a}
\geq {\rm supht}\, {\bf a} \, . $$
\end{prop}
{\it Beweis}.
Zu (1).
Ist $f:U \longrightarrow U'$ affin, so ist
${\rm dim}\, U' \geq {\rm cd}\, U' \geq {\rm cd}\, U$,
also ist ${\rm defdim}\, U \geq {\rm cd}\, U$.
\par\smallskip\noindent
Zu (2). Es ist ${\rm defdim}\, (U,X) =
{\rm dim}\, X' \geq {\rm dim}\, U' \geq {\rm defdim}\, U$.
\par\smallskip\noindent
Zu (3). Die erste Absch"atzung ist klar,
sei $A$ lokal und $B \longrightarrow A$ ein Ringhomomorphismus mit einem
noetherschen Ring $B$ und einem definierenden Ideal $\bf b$ mit
${\rm defdim}^{\rm ring}\, {\bf a}={\rm dim} \, B$.
Man kann $B$ ebenfalls als lokal annehmen (mit ${\bf b} \neq B)$.
Dann ist $D({\bf a}) \longrightarrow D({\bf b})$ affin und
${\rm defdim}\, D({\bf a}) \leq {\rm dim}\, D({\bf b}) < {\rm dim}\, B$.
\par\smallskip\noindent
Zu (4).
Ist $V({\bf a})=V(f_1, ...,f_r)$, so definieren diese Funktionen
eine $K-$Ab\-bil\-dung ${\rm Spek}\, A \longrightarrow {\bf A}_K^r $, und dabei
ist $V({\bf a})$ das Urbild des Nullpunktes des r-dimensionalen affinen Raumes.
Das Ideal l"asst sich also im $r-$dimensionalen Polynomring definieren.
Ist $A$ lokal, so folgen die Absch"atzungen alle aus den
bisherigen "Uberlegungen. \hfill  $\Box$
\par\bigskip
\noindent
{\bf Beispiel}
Ein affines Schema kann die Definitionsdimension $-1,0$ oder $1$ haben,
je nachdem, ob es leer ist, ob es einen Morphismus
auf das Spektrum eines K"orpers gibt oder nicht.
Entscheidend ist dabei der kanonische Morphismus nach ${\rm Spek}\, {\bf Z}$.
Ist beispielsweise $R$ ein lokaler Ring mit gemischter Charakteristik,
so ist die Definitionsdimension 1 (und die Definitionsh"ohe 2).
Umgekehrt muss ein Schema, das affin "uber einer Kurve ist, nat"urlich nicht
selbst affin sein.
\par\bigskip
\noindent
{\bf Beispiel}
Besitzt $U \subseteq X$ die eingebettete Definitionsdimension 1, so
gibt es einen affinen Morphismus
$f:X \longrightarrow C$ auf ein eindimensionales
noethersches Schema mit einer offenen Menge $V \subseteq C$ und
mit $U=f^{-1}(V)$.
Da $V \subseteq C$ lokal geometrisch durch eine Funktion beschrieben
werden kann, gilt dies auch f"ur $Y=X-U$.
Wird umgekehrt lokal $Y \subseteq X$ durch eine Funktion beschrieben
und ist $X$ ein $K-$Schema,
so liegt nat"urlich lokal die eingebettete Definitionsdimension eins
vor.
\par\smallskip\noindent
Ist $X$ ein affines integres $K-$Schema und $q$ eine rationale Funktion
mit ${\rm Def}\, (q) \cup {\rm Def}\, (q^{-1})=X$, so ist
$q: X \longrightarrow {\bf P}^1_K$ ein affiner Morphismus,
vergleiche 1.4.4 und 3.2.
\par\bigskip
\noindent
{\bf Beispiel} Sei $A$ eine noethersche $K-$Algebra und
$V({\bf a})=V(f_1,...,f_n)$. Genau dann besteht zwischen
den $f_1,...,f_n$ eine algebraische Relation, wenn
die durch $(f_1,...,f_n)$ definierte Abbildung auf den $n-$dimensionalen
affinen Raum nicht dominant ist.
In diesem Fall ist die Ringdefinitionsdimension echt kleiner als ${\rm ara}$.
Beispielsweise ist das Komplement von $V(f,g)$ in ${\rm Spek}\, A$ affin,
wenn zwischen $f$ und $g$ eine algebraische Relation besteht.
\par\bigskip
\noindent
{\bf Bemerkung}
Das Problem, ob eine irreduzible Kurve im dreidimensionalen affinen Raum
Durchschnitt zweier Fl"achen ist, kann man in zwei Probleme unterteilen,
n"amlich, ob sie die Ringdefinitionsdimension zwei besitzt,
also eine Faser "uber einem
Punkt einer Abbildung auf eine affine Fl\"ache $S$ ist, und zweitens, ob jeder
Punkt einer affinen Fl"ache durch zwei Funktionen beschreibbar ist.
\par
\bigskip
\noindent
\subsection{Die projektive Dimension}
\par
\bigskip
\noindent
Als eine weitere H"ohenzahl zu einem Ideal $\bf a$ besprechen wir die
projektive Dimension und ihr Verh"altnis zu Superh"ohe und kohomologischer
Dimension.
F"ur den Zusammenhang mit der Superh"ohe k"onnen wir uns auf \cite{bruns}
st"utzen.
\par\bigskip
\noindent
{\bf Definition}
Sei $A$ ein kommutativer Ring und $M$ ein $A-$Modul.
Die {\it projektive Dimension} von $M$, ${\rm pd}\, M$, ist die minimale
L\"ange $r$ einer projektiven Aufl\"osung 
$$ 0 \longrightarrow F_r \longrightarrow ... \longrightarrow F_1
\longrightarrow F_0 \longrightarrow M \, .$$
\par\bigskip
\noindent
{\bf Bemerkung}
Eine endliche projektive Aufl"osung muss nicht existieren.
Uns interessieren projektive Aufl"o\-sung von einem Ideal $\bf a$ bzw.
dem Restklassenring $A/{\bf a}$. Dabei ist
$ {\rm pd}\, (A/{\bf a}) = {\rm pd}\, {\bf a} +1 $.
Diese Zahl ist von uns von Interesse und ist mit Superh"ohe und
kohomologischer H"ohe zu vergleichen, und wir werden sie die
{\it projektive H"ohe} des Ideals nennen.
Sie ist gleich eins genau dann, wenn $\bf a$ projektiv ist,
und dies bedeutet, dass das Ideal lokal durch einen Nichtnullteiler
erzeugt wird, also invertierbar ist. Ein solches Ideal hat damit
trivialerweise affines Komplement, so dass auf dieser Ebene
(dem Niveau eins) keine interessante Beziehung zu Affinit"atsfragen
besteht und die projektive Dimension daher in den bisherigen Kapiteln
unber"ucksichtigt blieb.  
\par\smallskip\noindent
Der projektiven Dimension fehlen im Gegensatz zur kohomologischen H"ohe
und zur Superh"ohe einige strukturelle Eigenschaften.
Sie h"angt nicht nur vom Radikal ab, sondern wesentlich vom Ideal selbst,
sie muss nicht endlich sein und verh"alt sich unter Ringhomomorphismen
unkontrollierbar. Um so erstaunlicher ist die Absch"atzung nach der
Superh"ohe.
\par\bigskip
\noindent
Das Haupthilfsmittel bei der Absch"atzung zwischen Superh"ohe und
projektiver Dimension ist folgender Satz, der den Namen
new intersection theorem tr"agt, und vollst"andig erstmals von
Roberts 1987 in \cite{roberts} bewiesen wurde. F"ur den Beweis dieses
Satzes sowie auf verwandte Resultate sei hier generell auf \cite{bruns},
9.4 und \cite{roberts} verwiesen, f"ur den Zusammenhang zu anderen
homologischen Vermutungen auf \cite{roberts2} und \cite{roberts3}
\begin{satz}
Sei $A$ ein noetherscher lokaler Ring,
$F_\bullet: 0 \longrightarrow F_s \longrightarrow ... \longrightarrow F_0$
ein Komplex endlicher freier $A-$Moduln.
Der Komplex sei exakt f"ur alle ${\bf p} \neq {\bf m}$.
Ist dann $s < {\rm dim}\, A$, so ist $F_\bullet$ exakt.
\end{satz}
Zum {\it Beweis} siehe \cite{bruns}, Cor. 9.4.3. \hfill $\Box$
\par\bigskip
\noindent
Aus diesem Satz folgt nun leicht folgendes Superh"ohenresultat, das
vormals die homological height conjecture war.
\begin{satz}
Sei $A$ ein noetherscher Ring, $M \neq 0$ ein endlicher $A-$Modul.
Dann ist ${\rm pd}\, M \geq  {\rm supht}\,({\rm Ann}\, M )$.
Insbesondere ist f\"ur ein Ideal
${\rm pd}\,(A/{\bf a}) \geq {\rm supht}\, {\bf a} .$
\end{satz}
{\it Beweis}.
Sei die projektive Dimension von $M$ gleich $r$ und sei
$0 \longrightarrow F_r \longrightarrow ... \longrightarrow F_0
\longrightarrow M$ eine projektive Aufl"o\-sung.
Da die Superh"ohe in einem Punkt angenommen wird, k"onnen wir $A$ als lokal
und die Moduln $F_i$ als frei annehmen.
Sei
$A \longrightarrow B$ ein Ringhomomorphismus in einen lokalen integren Ring
der Dimension $d$ derart, dass $V({\rm Ann}\, M) B= V({\bf m}_B)$ ist.
Dabei ist einerseits $({\rm Ann}\, M)B \subseteq {\rm Ann}\, (M \otimes B)$,
es muss aber als Radikal die Gleichheit gelten, da andernfalls $M \otimes B$
der Nullmodul w"are. Dies kann aber nicht sein,
da ${\bf p}:= \varphi({\bf m}_B) \in V({\rm Ann}\, M)$ ist
und daher $M_{\bf p} \neq 0$ ($\varphi$ sei die Spektrumsabbildung).
\par\smallskip\noindent
Wir betrachten den tensorierten Komplex
$$ F_\bullet \otimes_AB : 0 \longrightarrow
F_r \otimes B \longrightarrow ...\longrightarrow
F_1 \otimes B \longrightarrow F_0 \otimes B \longrightarrow 0 \,  .$$
Sei ${\bf m} \neq {\bf q} \in {\rm Spek}\, B$ und $\varphi ({\bf q})={\bf r}$.
Es ist dann $M \otimes _A A_{\bf r}=0$ und daher ist
$F_\bullet \otimes_A A_{\bf r}$ eine exakte freie Aufl"osung des
Nullmoduls "uber einem lokalen Ring.
Nach \cite{eisenbud}, Lemma 20.1, ist dann dieser Komplex trivial und
spaltet.
Damit ist auch
$$ (F_\bullet \otimes_A B)_{\bf q}
=(F_\bullet \otimes_A A_{\bf r}) \otimes_{A_{\bf r}} B_{\bf q}$$
exakt, d.h. $F_\bullet \otimes_AB$ ist in jedem vom maximalen Ideal verschiedenen
Primideal exakt.
\par\smallskip\noindent
Der Komplex selbst ist aber nicht exakt, da wegen der Rechtsexaktheit
des Tensorproduktes, siehe \cite{scheja}, Satz 82.9,
die Sequenz $F_1 \otimes B \longrightarrow F_0 \otimes B \longrightarrow M \otimes B
\longrightarrow 0$ exakt ist und $M \otimes B \neq 0$ ist.
Damit liegt die Situation des obigen Satzes vor,
und es ergibt sich $ r \geq {\rm dim}\, B$. \hfill $\Box$
\par\bigskip
\noindent
{\bf Bemerkung}
Besitzt das Ideal ${\bf a}$ in $A$ selbst keine endliche projektive
Aufl"o\-sung, ist es aber (geometrisch) das
Erweiterungsideal unter einem Ringhomomorphisms
$B \longrightarrow A$ eines Ideal ${\bf b} \subseteq B$ mit einer
projektiven Aufl"osung der L"ange $r$, so "ubertr"agt sich die
Superh"ohenabsch"atzung auf $\bf a$.
Auch in diesem Kontext kann es sich also lohnen, wie bei der Definitionsh"ohe
zu schauen, wo ein Ideal herkommt.
\par\bigskip
\noindent
Aus dem Satz folgt insbesondere, dass f"ur ein Ideal, das eine
endliche projektive H"ohe gleich der H"ohe besitzt, auch die Superh"ohe
mit der H"ohe "ubereinstimmt. Dies ist unter recht allgemeinen
Voraussetzungen der Fall.
\begin{satz}
Sei $A$ ein lokaler noetherscher Cohen-Macauley Ring der Dimension $d$,
${\bf a}$ sei ein Ideal mit endlicher projektiver Dimension und so, dass
der Restklassenring $A/{\bf a}$ ebenfalls Cohen-Macauley ist.
Dann ist die H"ohe gleich der projektiven H"ohe und gleich der Superh"ohe.
\end{satz}
{\it Beweis}.
Nach der Formel von Auslander\--Buchsbaum,
siehe \cite{eisenbud}, Theorem 19.9, ist
${\rm pd}\, (A/{\bf a})=d-{\rm tf}_A (A/{\bf a})
=d-{\rm dim}\, (A/{\bf a}) = {\rm ht}\, {\bf a} $.
Damit stimmt diese Zahl auch mit der Superh"ohe "uberein.
\par\bigskip
\noindent
{\bf Beispiel}
Ohne die Voraussetzung, dass der Restklassenring
Cohen-Macauley ist, gilt der letzte Satz nicht.
Sei $A=K[x,y]$ und ${\bf a}=(x^2,xy)$.
Das Ideal wird geometrisch durch $x$ beschrieben, daher ist die Superh"ohe
und die kohomologische H"ohe des Ideals eins.
Die projektive H"ohe ist dagegen zwei.
In $A/{\bf a}$ ist das maximale Ideal gleich dem Annulator von $x$
und enth\"alt demnach keine Nichtnullteiler, daher ist die Tiefe null.
Als projektive Dimension von $A/{\bf a}$ (= projektive H"ohe von $\bf a$)
kommt nach der Auslander\--Buchsbaum-Formel
nur zwei in Frage.
\par\bigskip
\noindent
{\bf Beispiel}
Auch f"ur Primideale in einem regul"aren Ring muss die
Aussage des Satzes nicht gelten.
Wir betrachten den Monoidring $A=K[x^2,x^3,xy,y]=K[M]$.
Dieser Ring ist regul\"ar in der
Kodimension eins, aber nicht normal und nicht Cohen-Macauley.
Sei $\bf p$ der Kern der Restklassendarstellung
$R=K[u,v,w,y] \longrightarrow A$ mit
$u \longmapsto x^2, v \longmapsto x^3, w \longmapsto xy $.
Die Tiefe von $A=K[u,v,w,y]/{\bf p}$ ist 1, also muss die projektive
H"ohe 3 sein.
\par\smallskip\noindent
Die Superh"ohe des Primideals ist zwei, ebenso ist die kohomologische
H"ohe zwei, letzteres folgt aus Theorem 2.1 in \cite{huneke}.
\par\bigskip
\noindent
{\bf Bemerkung}
Wir stellen die Frage, ob in einem lokalen noetherschen Ring $A$
ein Primideal der H"ohe eins
mit endlicher projektiver Dimension affines Komplement besitzen muss oder
zumindest die Superh"ohe eins haben muss.
Dies steht in Zusammenhang mit einem Multiplizit"atenproblem,
ob f"ur einen Modul $M$ mit endlicher projektiver Dimension und einen
weiteren Modul $N$ derart, dass $M \otimes_A N$ von endlicher L"ange ist, die
Absch"atzung
$${\rm dim}\, M+ {\rm dim}\, N \leq {\rm dim}\, A$$
gelten muss, siehe \cite{roberts3}.
Ist diese Dimensionsaussage f"ur $M=A/{\bf p}$ mit ${\rm ht}\, {\bf p}=1$
und endlicher projektiver Dimension richtig,
so folgt f"ur eine Fl"ache $V({\bf q}) \subseteq {\rm Spek}\, A$,
dass $V({\bf q}) \cap V({\bf p}) \neq {\bf m}_A$ ist, da andernfalls
$A/{\bf p} \otimes_A A/{\bf q} = A/ ({\bf p}+{\bf q})$ endliche L"ange h"atte.
\par\smallskip\noindent
Ist ${\rm pd}\, A/{\bf p}=1$, so ist $\bf p$
invertierbar und $D({\bf p})$ affin.
Ist ${\rm pd}\, A/{\bf p}=2$, so ist nach
dem Satz von Hilbert-Burch, siehe \cite{eisenbud}, Theorem 20.15,
${\bf p}=aJ$ mit einem Nichtnullteiler $a \in A$ und einem
(Fitting-)Ideal $J$ der Tiefe zwei.
Dann kann aber $\bf p$ nicht die H"ohe eins haben.
\par\bigskip
\noindent
Wir besprechen nun den Zusammenhang von projektiver und kohomologischer
H"ohe eines Ideals.
\par\smallskip\noindent
Auf der einen Seite kann man ${\rm Ext}^r(V,W)$ mittels einer projektiven
Aufl"osung $... \longrightarrow F_1 \longrightarrow F_0$ von $V$ berrechnen.
Besitzt $V$ die projektive Dimension $r$, so ist ${\rm Ext}^s(V,W)=0$
f"ur $s > r$ und jeden Modul $W$.
Dies gilt insbesondere f"ur den Fall $V=A/{\bf a}$, wenn $\bf a$ die
projektive H"ohe $r$ besitzt.
\par\smallskip\noindent
Auf der anderen Seite geh"ort zu zwei Idealen ${\bf b} \subseteq {\bf c}$
die surjektive Abbildung $A/ {\bf b} \longrightarrow A/{\bf c}$
und dazu f"ur einen $A-$Modul $M$ die Abbildungen der Ext-Moduln
${\rm Ext}^{i} (A/{\bf c}, M) \longrightarrow {\rm Ext}^i(A/{\bf b},M) $.
\par\smallskip\noindent
Insbesondere liefern alle zu einem Ideal $\bf a$ radikalgleichen Ideale ein
gerichtetes Sytem von Ext-Moduln, und es gilt mit $Y=V({\bf a})$ und
${\cal F}=\tilde{M}$, vergleiche \cite{grolocal}, Theorem 2.8
oder \cite{SGA2}, Exp. II, Th\'{e}or\`{e}me 6,
$$\lim_{\bf b}\, {\rm Ext}^{i}(A/{\bf b},M)=
\lim_n \,{\rm Ext}^{i}(A/{\bf a}^n,M)=
H^{i} _Y(X,{\cal F}) \cong H^{i-1}(D({\bf a}),\cal{F}) \,  .$$
Besitzt $Y$ ein beschreibendes Ideal mit endlicher projektiver Dimension,
so kann man daraus allein noch keine Schl"usse f"ur die kohomologische
Dimension ziehen, allerdings unter folgender st"arkeren Bedingung.
\begin{satz}
L"asst sich die abgeschlossene Teilmenge $Y\subseteq X={\rm Spek}\, A$
geometrisch durch ein Ideal $\bf b$ mit der Eigenschaft beschreiben,
dass f\"ur alle {\rm (}es gen\"ugt f\"ur alle hinreichend gro"sen{\rm )}
Potenzen von $\bf b$ die projektive Dimension ${\rm pd}\,(A/{\bf b}^n) \leq r$
ist, so ist die kohomologische H\"ohe von $\bf a$ kleiner gleich $r$.
\end{satz}
{\it Beweis}.
F"ur $s >r$ verschwinden die einzelnen Glieder und damit auch der Limes,
also ist $H^{i}(D({\bf a}),{\cal F})=0$ f"ur
$i \geq r$, und die kohomologische H"ohe ist kleiner gleich $r$. \hfill $\Box$
\par\bigskip
\noindent
Sei im folgenden $A$ ein Ring von endlicher Charakteristik $p$, d.h.
$A$ ist eine ${\bf Z}/p-$Algebra.
In dieser Situation ist die Abbildung $F: x \longmapsto x^p$ ein
Ring\-endomorphismus von $A$, der der {\it Frobeniusmorphismus}
genannt wird.
Der Frobeniusmorphismus f"uhrt zu einer $A-$Algebra-Struktur auf $A$,
dieser Ring sei mit $A^F$ bezeichnet.
Die Zuordnung ${\cal F}: M \longmapsto A^F \otimes_A M$ legt einen Funktor,
den sogenannten {\it Frobeniusfunktor} von der Kategorie der $A-$Moduln
in die Kategorie der $A-$Moduln fest, siehe \cite{bruns}, 8.2.
Dabei ist ${\cal F} (A/{\bf a})= A/ ({\bf a}^{[p]})$, wobei
${\bf a}^{[p]}$ das von allen $p-$ten Potenzen von $\bf a$ erzeugte Ideal
bezeichnet, siehe \cite{bruns}, Prop. 8.2.2.
\begin{satz}
Sei $A$ ein lokaler noetherscher Ring der Charakteristik $p >0$.
Besitzt dann das Ideal $\bf a$ die endliche projektive H"ohe $r$,
so ist $r \geq {\rm kohoht}\, {\bf a}$.
\end{satz}
{\it Beweis}.
Ist $0 \longrightarrow M_r \longrightarrow ... \longrightarrow M_1
\longrightarrow M_0 \longrightarrow A/{\bf a} \longrightarrow 0$
eine projektive Aufl"o\-sung mit endlichen freien Moduln $M_i$, so ist
auch diese Aufl"osung nach Anwendung des Frobeniusfunktors exakt.
Dies ist f"ur die kleinen Stufen $i \leq 0$ klar aus der Rechtsexaktheit
des Tensorproduktes und folgt f"ur die anderen Stufen aus dem
Azyklit"atskriterium \cite{bruns}, Theorem 1.4.12, siehe \cite{bruns},
Theorem 8.2.7. 
\par\smallskip\noindent
Dies liefert dann eine freie Aufl"osung gleicher
L"ange von ${\cal F} (A/{\bf a})= A/({\bf a}^{[p]})$.
Damit sieht man iterativ, dass f"ur alle $n \in {\bf N}$
die Ideale ${\bf a}^{[p^n]}$ die projektive H"ohe $r$ haben
(siehe auch \cite{bruns}, Exc. 8.2.8).
F"ur jedes zu $\bf a$ radikalgleiche Ideal ist
${\bf a}^{[p^n]} \subseteq {\bf b}$ f"ur ein geeignetes $n$
und daher bilden diese Ideale
bzw. ihre Restklassenringe ein hinreichend feines Subsystem, und man kann
mit ihnen den Limes, also die lokale Kohomologie berechnen.
F"ur einen $A-$Modul $V$ ist dann
$H_{{\bf a}}^{i}(V)= H^{i-1}(D({\bf a}), \tilde{V})=0$
f"ur $i >r$, also ist die kohomologische H"ohe $\geq r$. \hfill $\Box$
\par\bigskip
\noindent
{\bf Beispiel}
Ist die Charakterisitk des Grundk"orpers null, so kann die kohomologische
H"ohe gr"o"ser als die projektive H"ohe sein, wie das folgende Beispiel zeigt.
\par\smallskip\noindent
Sei $A =K[X,Y,Z,U,V,W]$ und $\bf a$ das Ideal, der
$2 \times 2$ Unterdeterminanten der Matrix
$$ \left(\matrix{ X & Y & Z \cr U & V & W \cr } \right) \, .$$
Es ist also ${\bf a}=(XV-YU,XW-ZU,YW-ZV)$.
$\bf a$ ist ein Primideal der H"ohe zwei, der Restklassenring ist
als Minorenring normal und Cohen-Macauley der Dimension 4.
Daraus folgt sofort, dass die projektive Dimension dieses Restklassenringes
zwei ist. Die projektive H"ohe ist also zwei, nach 5.4.2 ist das auch die
Superh"ohe, und besitzt $K$ die Charakteristik $p$, so ist auch die
kohomologische H"ohe zwei.
In diesem Fall ist also f"ur jede quasikoh"arente Garbe
$H^2(D({\bf a}),{\cal F})=0$.
Wir werden zeigen, dass dies "uber den komplexen Zahlen nicht der Fall ist.
\par\smallskip
\noindent
Wir erinnern an die Segre-Einbettung, vergl. \cite{haralg},
I. Exc. 2.14 und II. Exc. 5.11.
Die Segre-Einbettung ist die abgschlossene Einbettung
$${\bf P}_K^r \times {\bf P}_K^s \longrightarrow {\bf P}^n \,
\mbox{ mit }\, (a_0,...,a_r;b_0,...,b_s) \longmapsto ( ...,a_ib_j,...) \,\, .$$
Das Bild wird dabei durch die Gleichungen
$T_{i,j} \cdot T_{k,l}- T_{i,l} \cdot T_{j,k}$ beschrieben. Bei $r=1$
und $s=2$ wird $(a_0,a_1;b_0,b_1,b_2)$ auf
$(a_0b_0, a_0b_1, a_0b_2, a_1b_0, a_1b_1, a_1b_2) =(X,Y,Z,U,V,W)$ ge\-schickt
und die Relationen sind dabei die von oben.
Projektiv betrachtet liegt also die Situation
$V_+({\bf a})={\bf P}_K^1 \times {\bf P}_K^2 \subseteq {\bf P}_K^5$ vor.
\par\smallskip\noindent
Da sich die kohomologische Dimension unter einer Kegelabbildung zu
einer positiv graduierten $K-$Algebra nicht "andert,
haben wir zu zeigen, dass f"ur $D_+({\bf a}) \subseteq {\bf P}_K^5$
die zweite Kohomologie im Fall $K={\bf C}$ nicht verschwinden kann.
\par\smallskip\noindent
Es gilt folgender Satz "uber $K={\bf C}$, der die kohomologische
Dimension "uber die algebraische de Rahm-Kohomologie mit
der singul"aren Kohomologie in Verbindung bringt.
\begin{satz}
Sei $X$ ein eigentliches Schema von endlichem Typ \"uber
${\bf C}$ der Dimension $n$. $Y$ sei ein abgeschlossenes Unterschema
und $X-Y$ sei regul\"ar. Ist dann ${\rm cd}\,(X-Y) < r$, so sind die
nat\"urlichen Abbildungen der singul\"aren Kohomologie
$$H^i(X^{\rm an},{\bf C}) \longrightarrow H^i(Y^{\rm an} ,{\bf C}) $$
f\"ur $i < n-r$ bijektiv und f\"ur $i =n-r$ injektiv.
\end{satz}
F"ur den Beweis siehe \cite{haramp}, III. Theorem 7.4. \hfill $\Box$
\par\bigskip
\noindent
In unserem Beispiel ist $X={\bf P}^5$, also $n=5$, und
$Y \cong {\bf P}^1 \times {\bf P}^2$. 
W\"are die kohomologische Dimension gleich eins, so erh\"alt man aus dem
Satz, wenn man $r=2$ setzt, f\"ur $i=0,1,2$ Bijektionen
und f\"ur $i=3$ noch eine Injektion.
Insbesondere h"atte man dann f\"ur $i =2 $ die Isomorphie
$$H^2({\bf P}^5,{\bf C}) \longrightarrow 
H^2( {\bf P}^1 \times {\bf P}^2,{\bf C}) \,  ,$$
und die zweiten Betti-Zahlen von ${\bf P}^5$ und von
${\bf P}^1 \times {\bf P}^2$ w"aren gleich.
F\"ur die komplexen projektiven R\"aume
ist aber $H_r({\bf P}_{\bf C}^d)= {\bf Z}$ f\"ur jede gerade Zahl $r$
unterhalb der reellen Dimension, siehe \cite{stocker}, Satz 9.9.10.
Nach der K"unneth-Formel f"ur Produkte $X \times Y$,
siehe \cite{stocker},  12.4.3,
gibt es einen injektiven Homomorphismus
$$H_0(X) \otimes H_n(Y) \oplus ... \oplus H_n(X) \otimes H_0(Y)
\longrightarrow H_n (X \times Y ) \, .$$
Damit ist die zweite Bettizahl von ${\bf P}_1 \times {\bf P}_2$
zumindest zwei (sie ist auch genau zwei, da in unserem Fall
der Komplement\"arsummand verschwindet).
\par\smallskip\noindent
Wir fassen zusammen: "uber $K={\bf C}$ ist
${\rm ara}\, {\bf a}={\rm afra}\, {\bf a}={\rm kohoht}\, {\bf a}=3$,
aber 
$ {\rm proht}\, {\bf a}={\rm supht}\, {\bf a}={\rm ht}\, {\bf a} =2 $.
\par
\bigskip
\noindent
\subsection{Verhalten bei Ringhomomorphismen}
\par
\bigskip
\noindent
Analog zu der Eigenschaft eines affinen Morphismus
${\rm Spek}\, B \longrightarrow {\rm Spek}\, A$, dass Urbilder affiner
Mengen wieder affin sind, steht die Eigenschaft, dass die besprochenen
H"ohenzahlen unter einem beliebigen
Ringhomomorphismus nicht gr"o"ser werden k"onnen.
Viele Ringhomomorphismen  gestatten es aber auch umgekehrt, von der Affinit"at
von $p^{-1}(U)$ auf die Affinit"at von $U$ selbst zu schlie"sen, wie
endliche Erweiterungen, treuflache Erweiterungen oder direkte
Summanden, siehe 1.14..
In diesem Abschnitt besprechen wir f"ur welche Ringhomomorphismen
kohomologische
Dimension und Superh"ohe nicht kleiner werden k"onnen.
\begin{prop}
Sei $A \longrightarrow A'$ eine treuflache Abbildung zwischen noe\-ther\-schen
Ringen. Dann gilt f"ur ein noethersches separiertes Schema $X$,
das von endlichem Typ "uber ${\rm Spek}\, A$ ist, die Gleichheit
${\rm cd}\, X={\rm cd}\, X'$ mit $X'=X \times_A {\rm Spek}\, A'$.
Speziell gilt f"ur ein Ideal $\bf a$ in $A$ die Gleichheit
${\rm cd}\, D({\bf a})={\rm cd}\, D({\bf a}A')$.
\end{prop} 
{\it Beweis}.
Da $X' \longrightarrow X$ affin ist, kann die kohomologische Dimension beim
"Ubergang zu $X'$ nicht gr"o"ser werden.
Schon im Beweis der entsprechenden Aussage f"ur die Affinit"at haben wir
$ H^{i}(X',{\cal F}') =H^{i}(X,{\cal F}) \otimes_A A'$ verwendet.
Steht links null, so muss auch $H^{i}(X,{\cal F})$ schon null sein,
da die Abbildung treu ist. Also kann die kohomologische
Dimension nicht kleiner werden. \hfill $\Box$
\begin{prop}
Sei $A$ ein noetherscher Ring, ${\bf a} \subseteq A$ ein Ideal und
$A \longrightarrow B$ eine
treuflache Erweiterung mit einem noetherschen Ring $B$. Dann gelten
folgende Aussagen.
\par\smallskip\noindent
{\rm (1)} Ist $B$ endlich erzeugt "uber $A$, so gilt
${\rm supht}\, {\bf a}={\rm supht}\, {\bf a}B$.
\par\smallskip\noindent
{\rm (2)} Ist $A$ lokal und $B$ die Komplettierung von $A$, so ist
${\rm supht}\, {\bf a}={\rm supht}\, {\bf a}B $.
\end{prop}
{\it Beweis}.
Sei $A \longrightarrow R$ ein Ringwechsel mit einem noetherschen,
lokalen, kompletten Ring $R$,
wo das Ideal seine Superh"ohe annimmt als H"ohe
des maximalen Ideals $\bf m$.
In (1) betrachtet man $R'=R \otimes_AB$. $R'$ ist von endlichem Typ "uber
$R$ und daher noethersch,
und $R \longrightarrow R'$ ist auch treuflach.
Daher gilt das going down Theorem und die H"ohe von ${\bf m}R'$
ist zumindest so gro"s wie die H"ohe  von $\bf m$.
In (2) nimmt man die Komplettierung $\hat{R}$ von $R$.
Dabei "andert sich die Dimension nicht, und der komplettierte
Ringhomomorphismus $\hat{A}=B \longrightarrow \hat{R}$ zeigt,
dass die Superh"ohe nicht kleiner wird.  \hfill $\Box$
\par\bigskip
\noindent
{\bf Bemerkung}
Es ist unklar, ob bei einer beliebigen
treuflachen Erweiterung $A \longrightarrow B$ die Superh"ohe kleiner
werden kann. Die Schwierigkeit liegt darin, dass zu einem
Testmorphismus $A \longrightarrow R$ das Tensorprodukt
$R \otimes_AB$ nicht noethersch sein muss.
\par\bigskip
\noindent
Wir wenden uns nun den endlichen surjektiven
Schemamorphismen $X' \longrightarrow X$ zu,
also der Situation, wo es der Satz von Chevalley erlaubt, von der
Affinit"at von $X'$ auf die Affinit"at von $X$ zu schlie"sen.
In dieser Situation bleiben kohomologische Dimension und Superh"ohe konstant.
\begin{prop}
Sei $f: X' \longrightarrow X$ ein endlicher surjektiver Morphismus
von noetherschen separieten Schemata. Dann ist
${\rm cd}\, X'={\rm cd}\, X$.
\end{prop}
{\it Beweis}. Siehe \cite{harcodi}, Prop. 1.1. \hfill $\Box$
\begin{prop}
Sei $X' \longrightarrow X$ eine endliche surjektive Abbildung
zwischen noetherschen separierten Schemata.
Dann ist die Superh"ohe von $X'$ gleich der Superh"ohe von $X$.
\end{prop}
{\it Beweis}.
Da die Abbildung affin ist, wird die Superh"ohe nicht gr"o"ser.
Wir zeigen die Umkehrung. Sei
${\rm Spek}\, R=T \supset T^\times \longrightarrow X$
ein affiner Morphismus mit $R$ lokal, integer, noethersch, normal und
${\rm dim}\, R=n={\rm supht}\, X$.
Es sei ${\rm Spek}\, L$ ein abgeschlossener Punkt von
${\rm Spek}\, Q(R) \times_X X'$, so dass $Q(R) \subseteq L$
eine endliche K"orpererweiterung ist.
Sei $R'$ der normale Abschluss von $R$ in $L$.
Damm gibt es auf $D({\bf m}_R R')$ eine affine kommutierende
Abbildung nach $X'$, wie Betrachtung der Situation auf affinen
St"ucken von $X$ zeigt.
Da $R \longrightarrow R'$ injektiv und ganz ist, besitzt
${\bf m}_R R'$ ebenfalls die H"ohe $n$. \hfill $\Box$
\par\bigskip
\noindent
{\bf Bemerkung} Liegt eine endliche Ringerweiterung $A \longrightarrow A'$
mit surjektiver Spektrumsabbildung vor, so kann man die Konstanz
der Superh"ohe f"ur eine offene Teilmenge $U \subseteq {\rm Spek}\, A$
direkt durch Tensorierung eines Testhomomorphismus
$A \longrightarrow R$ beweisen.
\par\bigskip
\noindent
Als dritte Klasse von Ringwechseln, die den R"uckschluss auf die
Affinit"at gestattet, haben wir direkte Summanden kennengelernt, die
wir jetzt in Hinblick auf das Verhalten der H"ohenzahlen untersuchen.
\begin{prop}
Seien $A$ und $B$ noethersche Ringe mit $B=A+V$, wobei $V$ ein $A-$Modul
sei. Dann gilt f"ur eine offene Menge
$U \subseteq {\rm Spek}\, A$ mit dem Urbild $U'$ die Gleichheit
${\rm cd}\, U = {\rm cd}\, U'$. 
\end{prop}
{\it Beweis}.
Ist $M$ ein $A-$Modul, so gelten f"ur
$M_B=M \otimes_A B= M \oplus (M\otimes_AV)$
die Gleichungen
\begin{eqnarray*}
H^k(U',\tilde{M_B}) & = & H^k(U,p_*(\tilde{M_B}))\\
& =& H^k(U, \tilde{M} \oplus \widetilde{(M\otimes_AV)}) \\
&=& H^k(U,\tilde{M}) \oplus H^k(U,\widetilde{M\otimes_AV}) \, .
\end{eqnarray*}
Ist $H^k(U, \tilde{M}) \neq 0$, so ist auch $H^k(U',\tilde{M_B}) \neq 0$, und
die kohomologische Dimension kann nicht kleiner werden. \hfill $\Box$
\par\bigskip
\noindent
Besitzt die Abbildung $A \longrightarrow B$ sogar einen Schnitt, ist also 
der andere Summand ein Ideal,
so kann nat"urlich die Superh"ohe nicht kleiner werden.
Dies gilt trivialerweise f"ur jede H"ohenzahl, die unter einem affinen
Morphismus nicht gr"o"ser werden kann.
Das folgende Beispiel zeigt, dass bei einem beliebigen
direkten Summanden sich die Superh"ohe "andern kann.
\par\bigskip
\noindent
{\bf Beispiel}
Sei $B=K[S_1,S_2,T_1,T_2]$ graduiert, wobei $S_1$ und $S_2$ den Grad eins
und $T_1$ und $T_2$ den Grad minus eins haben.
Der Ring der nullten Stufe ist der Ring
$A=K[S_1T_1,S_1T_2,S_2T_1,S_2T_2]=K[X,Y,V,U]$ mit der Relation
$UX =VY$. Als Ring der nullten Stufe
einer ${\bf Z}-$Graduierung ist $A$ ein direkter Summand.
\par\smallskip\noindent
Der Punkt $(U,X,V,Y)$ hat auf $ UX-VY$ als abgeschlossener
Punkt auf einer dreidimensionalen affinen Variet"at
die H\"ohe, Superh\"ohe und kohomologische H\"ohe drei.
Als Radikal l\"asst sich der Punkt durch die drei Funktionen
$X+U,V,Y$ beschreiben und daher ist ara ebenfalls drei.
\par\smallskip\noindent
Das Urbild dieses Punktes ist $(S_1T_1,S_1T_2,S_2T_1,S_2T_2)
=(S_1,S_2) \cap (T_1,T_2) ={\bf a}$.
Nach 5.5.5 hat dann ${\bf a}$ die kohomologische
H\"ohe drei, vergleiche auch das dritte Beispiel in 5.3.
Damit ist auch ${\rm ara}\, {\bf a}=3$.
Dagegen ist die H"ohe von $\bf a$ zwei und nach 5.1.3
auch die Superh\"ohe.

%% file: problem.tex
\section*{Probleme und offene Fragen}
Wir versammeln zum Schluss einige offene Fragen und Probleme, die uns schon im
Verlauf der Arbeit begegnet sind.
Es ist m"oglich, dass einzelne von ihnen eine einfache Antwort besitzen,
wir verzichten auf eine weitere Kommentierung.
\par\bigskip
\noindent
1)
Ist ein noethersches, separiertes, integres, eindimensionales Schema
mit eindimensionalem globalen Schnittring affin?
\par\bigskip
\noindent
2)
Seien $X$ und $Y$ Variet"aten "uber einem K"orper und
$f:X \longrightarrow Y$ eine Hom"oomorphie.
Folgt aus $X$ affin, dass $Y$ affin ist?
\par\bigskip
\noindent
3)
Man gebe ein handliches Kriterium f"ur affin-triviale Divisoren an.
\par\bigskip
\noindent
4)
Ist $p:{\rm Spek}\, A' \longrightarrow {\rm Spek}\, A$ ein
(flacher) Ringhomomorphismus zwischen noetherschen, normalen, integren
Ringen.
Ist zu einem affin-trivialen Divisor $D$ auf ${\rm Spek}\, A$
der zur"uckgenommene
Divisor $p^*D$ ebenfalls affin-trivial?
Gibt es einen kanonischen Gruppenhomomorphismus
${\rm AKG}\, A \longrightarrow {\rm AGK}\, A'$?
\par\bigskip
\noindent
5)
Sei $A$ ein noetherscher, integrer, normaler, lokaler Ring (mit zus"atzlichen
Eigenschaften: exzellent, im wesentlichen von endlichem Typ, komplett).
Ist ${\rm AKG}\, A$ endlich erzeugt?
\par\bigskip
\noindent
6)
Sei $A$ ein integrer, normaler, noetherscher Ring mit ${\rm AKG}\, A=0$.
"Ubertr"agt sich dies auf den formalen Potenzreihenring $A[[T]]$?
Gilt zumindest bei $A$ faktoriell ${\rm AKG}\, A[[T]]=0$?
\par\bigskip
\noindent
7)
Ist $X$ eine normale affine Variet"at und $Y$ eine glatte affine Variet"at
"uber einem algebraisch abgeschlossenen K"orper.
Ist dann ${\rm AKG}\, (X \times Y)={\rm AKG}\, X$?
Gilt dies, wenn $X$ eine Fl"ache ist?
Gilt dies, wenn $Y$ eine Kurve ist?
\par\bigskip
\noindent
8)
Sei $X$ ein normaler, komplexer, Steinscher Raum, auf dem jedes
Hyperfl"achenkomplement Steinsch ist. Gilt dies dann auch
f"ur $X \times {\bf C}$?
Gilt dies zumindest, wenn $X$ eine normale Steinsche Fl"ache ist,
so dass die Voraussetzung automatisch erf"ullt ist.
Gilt dies sogar f"ur $X \times Y$, wenn $Y$ eine beliebige Steinsche
Mannigfaltigkeit ist?
\par\bigskip
\noindent
9)
Gilt f"ur eine lokale, (1-) rationale Singularit"at ${\rm Spek}\, A$,
dass die affine Klasengruppe gleich der Divisorenklassengruppe modulo
Torsion ist?
Wie sieht das f"ur eine isolierte Hyperfl"achensingularit"at im ${\bf A}^4$
aus.
Man berechne die affine Klassengruppe von Brieskorn- und insbesondere
Fermatsingularit"aten.
\par\bigskip
\noindent
10)
Ist $A$ der Kegelring "uber einer (glatten) projektiven Variet"at $Y$ mit
der Kegelabbildung $p:{\rm Spek}\, A - \{ (A_+) \} \longrightarrow Y$.
Ist das Urbild eines numerisch trivialen Divisors $D$ auf $Y$
affin-trivial?
Man kl"are das Verh"altnis von numerischer und affiner Klassengruppe.
\par\bigskip
\noindent
11)
"Ubertr"agt sich die Eigenschaft ${\rm AKG}\, A=0$ auf (generische)
Hyperfl"achenschnitte (unter gewissen Dimensions- und Gradbedingungen)?
\par\bigskip
\noindent
12)
Verschwindet in einem dreidimensionalen, normalen, noetherschen, integren Ring,
indem jedes Primideal der H"ohe eins die Superh"ohe eins besitzt
(damit gilt diese Gleichheit "uberhaupt f"ur jedes Primideal),
die affine Klassengruppe?
Gibt es projektive glatte Fl"achen, auf denen sich je zwei verschiedene
irreduzible Kurven positiv schneiden, der Selbstschnitt jeder Kurve
nicht negativ ist, aber nicht jedes Kurvenkomplement affin ist
(etwa alle Kurvenkomplemente Steinsch, aber nicht alle affin)?
\par\bigskip
\noindent
13)
Man kl"are zu einem Ideal das Verh"altnis von Superh"ohe und analytic spread.
\par\bigskip
\noindent
14) Kann man unter den noetherschen, separierten Schemata (ohne die
Voraussetzung quasiaffin) die affinen durch die Superh"ohe eins bzgl. allen
Krullbereichen charakterisieren. (F"ur quasiprojektive Schemata ist das richtig)
\par\bigskip
\noindent
15)
Definiert ein Primideal der H"ohe eins mit endlicher projektiver Dimension
eine affine Teilmenge?
\par\bigskip
\noindent
16)
Wie kann man die Eigenschaft eines quasiaffinen Schemas $U$, dass es eine
Einbettung $U=D(I) \subseteq {\rm Spek}\, A$ gibt mit einem definierenden
Ideal $I$ endlicher projektiver Dimension intrinsisch in $U$ beschreiben,
etwa in der Auswirkung auf die Kohomologiemoduln?
L"asst sich diese Eigenschaft (wie die Superh"ohe)
f"ur beliebige Schemata sinnvoll definieren
und gilt die entsprechende Absch"atzung?
\par\bigskip
\noindent
17)
Kann sich unter einem treuflachen Ringwechsel die Superh"ohe
"andern?
Kann bei der Komplettierung eines exzellenten Ringes die endliche
Superh"ohe gr"o"ser werden?
Dies ist die Hauptschwierigkeit bei der Frage von Koh, ob
f"ur exzellente Ringe die Gleichheit von Superh"ohe und
endlicher Superh"ohe gilt.
\par\bigskip
\noindent
18)
Man charakterisiere die Ringwechsel, die den R"uckschluss auf die Affinit"at
erlauben bzw. unter denen sich die kohomologische Dimension
einer offenen Teilmenge nicht "andert.

%% file: lebenslauf.tex
\par
\bigskip
\bigskip
\bigskip
\par
\bigskip
\bigskip
\noindent
LEBENSLAUF VON HOLGER BRENNER
\par\bigskip
\noindent
Ich wurde am 3.11.1966 in Calw geboren als Sohn von Willy Brenner
und Irmgard Brenner, geborene Pross.
\par\smallskip\noindent
Ab 1973 besuchte ich die Grundschule Wildberg, ab 1977
das Otto-Hahn-Gymnasium Nagold, wo ich 1986 das Abitur machte.
\par\smallskip\noindent
1986-88 leistete ich in Essen meinen Zivildienst ab.
\par\smallskip\noindent
Von 1988 an studierte ich an der Ruhr-Universit"at Bochum Mathematik mit
Nebenfach Philosophie und erlangte 1994 das Diplom der Fakult"at f"ur
Mathematik.
\par\smallskip\noindent
Seit 1995 bin ich Doktorand bei Herrn Prof. Dr. U. Storch in Bochum.
Dabei hatte ich von Jan 1995-Sept 1995 eine Stelle als wissenschaftliche
Hilfskraft, von Okt 1995- April 1997 war ich Stipendiat im Graduiertenkolleg
Geometrie  und mathematische Physik, seit April 1997 bin ich wissenschaftlicher
Mitarbeiter an der Ruhr-Universit"at Bochum.